\documentclass[reqno,12pt,a4paper]{article}
\usepackage{sectsty}
\usepackage{amsmath}
\usepackage{amsthm}        
\usepackage{amssymb}       
\usepackage{amsfonts}
\usepackage{mathtools}
\usepackage{enumitem}
\usepackage{float}%for including figures at specific point
\usepackage{adjustbox}
\usepackage[svgnames]{xcolor}
\usepackage{titletoc}
\usepackage[titletoc]{appendix} % package for appendices in toc
\usepackage{dsfont}
\usepackage{caption}
\usepackage{subcaption}
\usepackage{anyfontsize}
\usepackage{tikz-cd}
\usepackage[hyperfootnotes=true,draft=false,pdftex, colorlinks=true, linkcolor=Black, citecolor=Black, urlcolor=DarkBlue, filecolor=Black, hypertexnames=false, % linkbordercolor=red
	]{hyperref}
\sectionfont{\fontsize{15}{15}\selectfont}
\numberwithin{equation}{section}%label equation by sections

\theoremstyle{plain} %%% Plain Theorem Styles.
\newtheorem{theorem}{Theorem}[section]
\newtheorem{lemma}[theorem]{Lemma}
\newtheorem{corollary}[theorem]{Corollary}   
\newtheorem{prop}[theorem]{Proposition}      

\theoremstyle{definition} %%%% Definition-like Commands  
\newtheorem{definition}[theorem]{Definition} 

\theoremstyle{remark}  %%%% Remark-like Commands

\newtheorem{remark}[theorem]{Remark}
\newtheorem{example}[theorem]{Example}
\newtheorem{examples}[theorem]{Examples}

\newcommand{\nid}{\noindent}
\newcommand{\ra}{\rightarrow}

\makeatletter
\newcommand{\markthis}[3]{% #1 = marker, #2 = label, #3 = relation
  \overset{
    \textup{\makebox[0pt]{\normalfont\footnotesize(#1)}}
    \def\@currentlabel{#1}%
    \ltx@label{#2}%
  }{
    #3
  }
}
\makeatother %equation number above equality sign.

\newcommand\eqdef{\mathrel{\overset{\makebox[0pt]{\mbox{\normalfont\tiny\sffamily def}}}{=}}} %equation with text "def" above

\newcommand\eqabove[1]{\mathrel{\overset{\makebox[0pt]{\mbox{\normalfont\footnotesize #1}}}{=}}} %equation with text "#1" above

\newcommand{\parll}{{\mathpalette\rotamp\relax}} %par operator
\newcommand{\rotamp}[2]{\rotatebox[origin=c]{180}{$#1\&$}}
\DeclareMathOperator{\parLL}{\parll}

\newcommand{\distl}{\delta^l} %left distributor
\newcommand{\distr}{\delta^r} %right distributor

\def\cA{\mathcal A}\def\cB{\mathcal B}\def\cC{\mathcal C}\def\cD{\mathcal D}
\def\cE{\mathcal E}\def\cF{\mathcal F}\def\cG{\mathcal G}\def\cH{\mathcal H}
\def\cI{\mathcal I}

\def\cV{\mathcal V}\def\cW{\mathcal W}

\author{Max Demirdilek and
  Christoph Schweigert\\ \\ \multicolumn{1}{p{.7\textwidth}}{\centering \emph{Fachbereich Mathematik\\ Universit\"at Hamburg\\Bereich Algebra und Zahlentheorie\\
Bundesstra{\ss}e 55\\ D-20146 Hamburg, Germany}}}
\title{Surface Diagrams for\\ Frobenius Algebras and Frobenius-Schur Indicators\\in Grothendieck-Verdier Categories}
\date{}

\begin{document}
\emergencystretch 3em	

\thispagestyle{empty}

\maketitle

\begin{abstract}\nid Grothendieck-Verdier categories (also known as $\ast$-autonomous categories) generalize rigid monoidal categories, with notable representation-theoretic examples including categories of bimodules, modules over Hopf algebroids, and modules over vertex operator algebras.

In this paper, we develop a surface-diagrammatic calculus for Grothendieck-Verdier categories, extending the string-diagrammatic calculus of Joyal and Street for rigid monoidal categories into a third dimension. This extension naturally arises from the non-invertibility of coherence data in Grothendieck-Verdier categories.
 
We show that key properties of Frobenius algebras in rigid monoidal categories carry over to the Grothendieck-Verdier setting. Moreover, we introduce higher Frobenius-Schur indicators for suitably finite $k$-linear pivotal Grothendieck-Verdier categories and prove their invariance under pivotal Frobenius linearly distributive equivalences. 

The proofs are carried out using the surface-diagrammatic calculus. To facilitate the verification of some of our results, we provide auxiliary files for the graphical proof assistant homotopy.io.
\end{abstract}

%%%%%%%%%%%%%%%%%%%%%%%%%%%%%%%%%%%%%%%%%%%%%%%

\pagebreak
\setcounter{tocdepth}{2}
\tableofcontents

%%%%%%%%%%%%%%%%%%%%%%%%%%%%%%%%%%%%%%%%%%%%%%%

\section{Introduction}

\nid\textbf{Background.} Contragredient representations are a common phenomenon, often acting as rigid dual objects in monoidal categories. Rigidity is, however, too restrictive for many interesting representation-theoretic examples.

Consider, for instance, the monoidal category of finite-dimensional bimodules over a finite-dimensional algebra $A$ over a field $k$. Although only some objects in this monoidal category have rigid duals, the $k$-linear dual of any $A$-bimodule is still an $A$-bimodule. The $k$-linear dual of the regular bimodule $A$ is not necessarily isomorphic to $A$ itself. It is a \emph{dualizing object}, and thus endows the monoidal category of finite-dimensional $A$-bimodules with the structure of a \emph{Grothendieck-Verdier category}, see \cite{fuchs2024grothendieckverdierdualitycategoriesbimodules}. More generally, finite-dimensional modules over a Hopf algebroid with an antipode carry Grothendieck-Verdier structures \cite{allen2024hopfalgebroidsgrothendieckverdierduality}.

Further examples arise in conformal field theory. A large class of representation categories of vertex operator algebras (VOAs), such as the $\cW_{2,3}$ triplet model, are monoidal categories in which not all objects possess rigid duals \cite{HLZ}. However, any choice of conformal structure on a VOA $\cV$ still allows for the definition of a grade-wise dual of a $\cV$-module. The grade-wise dual of $\cV$ itself is a dualizing object, and the monoidal category of $\cV$-modules admits the structure of a Grothendieck-Verdier category \cite[Thm. 2.12]{allen2021duality}. 

By definition, the dualizing object of a Grothendieck-Verdier category $\cC$ induces an antiequivalence $D$ on $\cC$. Trivializations of the double duality functor $D^2\colon \cC\ra \cC$ yield \emph{pivotal structures} on $\cC$. Both the category of finite-dimensional bimodules and representation categories of VOAs naturally carry such pivotal structures.

\medskip

Inspired by these examples, we address the following two questions in this paper.

\smallskip

\nid\textbf{Question 1.} \textit{What is the theory of Frobenius algebras in Grothendieck-Verdier categories?}

\smallskip

\nid We provide motivations for focusing on Frobenius algebras later.

\medskip

\nid\textbf{Question 2.} \textit{What invariants can be defined for objects in Grothendieck-Verdier categories?}

\smallskip

\nid We focus on higher Frobenius-Schur indicators for \emph{pivotal} Grothendieck-Verdier categories.

\medskip

The following notion plays a central role in our investigation.

\nid\textbf{Linearly distributive categories.} 
The duality functor of a Grothendieck-Verdier category is generally not opmonoidal, leading to the introduction of a second monoidal product on the category. Specifically, as recalled in Theorem \ref{thm:GV correspond LD neg}, Grothendieck-Verdier categories correspond to linearly distributive categories with negation. A linearly distributive category \cite{WDC}, see also Definition \ref{def:LDcat}, is a category endowed with two monoidal structures $(\otimes, 1)$ and $(\parLL, K)$, referred to as otimes and par, respectively. It also comes with two coherent natural transformations \begin{align*}
    \distl_{X,Y,Z}:\: X\otimes{(Y \;{\parLL}\; Z)} \,\longrightarrow \, (X \otimes Y) \;{\parLL}\; Z,\\
    \distr_{X,Y,Z}:\: (X\;{\parLL}\;Y) \otimes Z \, \longrightarrow \, X \;{\parLL}\;{(Y \otimes Z)},
\end{align*} called \emph{distributors}. The distributors are not required to be invertible and, in fact, are non-invertible in relevant examples. A linearly distributive category has a notion of duals and is called a linearly distributive category with negation if every object possesses two-sided linearly distributive duals; see Definition \ref{def:LDduals}.

\medskip

\nid\textbf{Frobenius algebras.} Linearly distributive categories provide a natural framework for defining Frobenius algebras. An LD-Frobenius algebra in a linearly distributive category $\cC$, as recalled in Definition \ref{def:FrobAlg}, is a quintuple $(A,\mu,\eta,\Delta,\epsilon)$, where the triple $(A,\mu,\eta)$ forms a unital algebra in the monoidal category $(\cC,\otimes,1)$, and the triple $(A,\Delta,\epsilon)$ forms a counital coalgebra in the other monoidal category $(\cC,\parLL,K)$. The algebra and coalgebra structures are thus defined with respect to different monoidal products. Additionally, we impose the following relations:
\begin{equation*}
    (\mu \parLL A)\circ \distl_{A,A,A}\circ (A \otimes  \Delta)\,\;\markthis{F1}{eq:F1}{=}\;\,(A\parLL \mu)\circ \distr_{A,A,A}\circ (\Delta \otimes A) \,\;\markthis{F2}{eq:F2}{=}\;\, \Delta \circ \mu.
\end{equation*}
These LD-Frobenius relations involve both monoidal products and are obtained by replacing the associators in the standard monoidal Frobenius relations with the appropriate distributors, adjusting the monoidal products accordingly.

LD-Frobenius algebras arise in disguise in various contexts, as we will see in Examples \ref{example ld-frobenius algebras}. For instance, they appear in the study of quantales \cite{EggerFrobenius,FrobeniusStructures}, the classification of K3 surfaces up to isogeny \cite{MotivicGlobalTor}, and the investigation of Chow motives of smooth cubic fourfolds \cite{CubicFourfolds}. Our primary interest in LD-Frobenius algebras stems from their possible role in logarithmic conformal field theory. In the better-understood setting of rational conformal field theories, Frobenius algebras in rigid monoidal categories play a significant role, e.g. \cite{FuchsRunkelSchBoundariesDefects}, \cite{TCFTRFFS}, \cite[\S 3]{stigner2012hopffrobeniusalgebrasconformal}, \cite{fuchs2023algebraicstructurestwodimensionalconformal}. It has been proposed in \cite[\S 3]{LogBulkBound} that in a general logarithmic conformal field theory, boundary conditions yield a structure that we identify as an LD-Frobenius algebra in a Grothendieck-Verdier category.

\medskip

\nid\textbf{Frobenius-Schur indicators.} Classical Frobenius-Schur indicators are invariants of irreducible complex representations of a finite group that determine whether a representation can be realized over the real numbers. These indicators can, for instance, be used to prove classical results in group theory, such as the Brauer-Fowler Theorem. 

Frobenius-Schur indicators have been generalized in various directions. In one direction, higher-order generalizations of Frobenius-Schur indicators for pivotal rigid monoidal categories exist \cite{KaSoZh,SchauNgHigherFrob}. In another direction, Shimizu has studied indicators for categories with a very general notion of duality structure \cite{ShimizuFrobSchur}. Frobenius-Schur indicators find applications in modular fusion categories relevant to conformal field theory \cite{Bantay, BantayCurrentExt, BantayKleinBottle, MasonFrob}, as well as in the study of modules over semisimple quasi-Hopf algebras \cite{MasonNg,CentralInvariantsAndHigher}. Given these applications, it is natural to explore generalizations of higher Frobenius-Schur indicators to pivotal Grothendieck-Verdier categories.
\medskip

\nid\textbf{Problem.}
The many coherence axioms for distributors in linearly distributive categories and Grothendieck-Verdier categories make computations involved and tedious. In the special case of monoidal categories, all coherence morphisms are typically neglected, which is justified by Mac Lane's strictification theorem. 

No analogous strictification result holds for linearly distributive categories. Specifically, any linearly distributive category with invertible distributors is Frobenius LD-equivalent to a shift monoidal category; see Remark \ref{strictifiability of LD-cats}. In a shift monoidal category, the second monoidal product $\parLL$ is, by definition, given by an $S$-shifted monoidal product ${X \parLL Y=Y\otimes (S \otimes X)}$, for $S$ an $\otimes$-invertible object. Shift-monoidal categories do not suffice in practice.

In fact, linearly distributive categories with negation are not even fully coherent in the sense that there exist formal diagrams that do not commute in general; see Example \ref{ex:LD not fully coherent}.

\medskip

\nid\textbf{Surface diagrams.} Given the complexity of the linearly distributive coherence axioms, a tool to aid in calculations is desirable. For strict monoidal categories, this is provided by Joyal and Street's string-diagrammatic calculus. In linearly distributive categories, the presence of non-invertible coherence data and the lack of a full coherence theorem suggests a graphical language that explicitly represents both coherence morphisms and their axioms.

We develop such a diagrammatic language by ‘moving up a dimension’. Surface diagrams replace string diagrams, with the key insight being that the monoidal $2$-category of categories, functors, and natural transformations provides an appropriate setting for this language. This insight allows us to present a three-dimensional diagrammatic language for linearly distributive categories using the standard graphical calculus for strict monoidal $2$-categories; see Section \ref{sec:surfdiagforGray} for a discussion of the calculus. Variants of the calculus have been applied in various contexts \cite{TrimbleSurface,FunctorialCalcMonBi, WillertonHopfM, douglas2018fusion, barrett2024gray}. Similar ideas have appeared in \cite{StreetSurfaceProc,ProfSemantics, DunnVicarySurface, DunnVicaryFrob}.

In our diagrammatic language, the coherence morphisms of a linearly distributive category become stratified $2$-manifolds embedded into framed three-space, see, e.g., Figures \ref{fig:alpha} and \ref{fig: Distributors} for the geometry relevant for associators and distributors. Intuitively, these surfaces consist of sheets ($2$-strata) glued together in various configurations. The coherence axioms amount graphically to the invariance under certain geometric moves. These can be described as moving gluing boundaries of a sheet along a surface; see Appendix \ref{coherenceLD}. 

The surface-diagrammatic calculus also allows us to represent functors-with-structure between linearly distributive categories. Their surface diagrams arise from drawing string diagrams onto the stratified surfaces representing the coherence morphisms, see Figure \ref{fig: Relations FrobLD-fun}. 

\medskip

Surface diagrams are used in the proofs of our main results.

\smallskip

\nid\textbf{Main results.} First, we present three novel results for LD-Frobenius algebras, which extend well-known results for Frobenius algebras defined with respect to a single monoidal structure. Our results are non-trivial since they involve the non-invertible distributors $\distl$ and $\distr$.
\vspace{-0.17cm}
\begin{itemize}
    \item Theorem \ref{FrobRelationsNotIndep} shows that the LD-Frobenius relations are not independent, meaning that the LD-Frobenius relation \eqref{eq:F1} implies the LD-Frobenius relation \eqref{eq:F2}.
    \vspace{-0.17cm}
    \item Theorem \ref{Frobenius via ideals} characterizes LD-Frobenius algebras in abelian Grothendieck-Verdier categories as algebras endowed with a form (with values in the dualizing object) whose kernel has no non-trivial one-sided ideals. Our result generalizes a theorem of Walton and Yadav \cite[Thm. 5.3]{walton2022filtered} from the rigid to the Grothendieck-Verdier setting.
    \vspace{-0.17cm}
    \item Theorem \ref{thm:Frob Alg modules isomorphic to comodules} provides an isomorphism between the category of modules and the category of comodules of an LD-Frobenius algebra.
\end{itemize}
\nid Second, we provide an invariant of $k$-linear pivotal Grothendieck-Verdier categories with finite-dimensional hom-spaces:
    \begin{itemize}
    %\vspace{-0.2cm}
    \item Section \ref{sec:FrobSchurIndic} generalizes higher Frobenius-Schur indicators from $k$-linear pivotal rigid monoidal categories to $k$-linear pivotal Grothendieck-Verdier categories. Theorem \ref{thm: FrobSchur categorical invariants} establishes graphically the invariance of our generalized Frobenius-Schur indicators under equivalences of $k$-linear pivotal Grothendieck-Verdier categories.
\end{itemize}

\nid\textbf{Proof techniques.}
Many of our surface-diagrammatic proofs are higher-dimensional analogues of string-diagrammatic proofs in the strict monoidal setting, which can be recovered by viewing surface diagrams from the the `front face', see, e.g. Figure \ref{fig:string diag surf diag}.

This observation provides a strategy for extending string-diagrammatic statements and proofs from (rigid) monoidal categories to linearly distributive categories (with negation). First, draw each string diagram in the monoidal setting. Next, convert it into a surface diagram for linearly distributive categories by adding three-dimensional building blocks that represent the coherence data. Finally, identify the correct moves between the diagrams, corresponding to the coherence axioms. This strategy is described in detail in Remark \ref{proof strategy}. For a concrete example illustrating this approach, we refer to the figures in the proof of Theorem \ref{FrobRelationsNotIndep}.

\medskip

\nid\textbf{Overview.} The paper is organized as follows: In Section \ref{sec:Categorical prelims}, we provide the categorical ingredients required for the main results of this paper. We begin with a review of the graphical calculus for strict monoidal $2$-categories, applying the calculus to a strictification of the monoidal $2$-category of (enriched) categories, functors, and natural transformations. From there, we derive a three-dimensional diagrammatic language for linearly distributive categories. 

Using this language, we recall Frobenius LD-functors between linearly distributive categories, their morphisms, and the notion of dual objects in a linearly distributive category. We conclude Section \ref{sec:Categorical prelims} by reviewing Grothendieck-Verdier categories and their correspondence to linearly distributive categories with negation, proving parts of this correspondence graphically. Finally, we show that the duality functor of a Grothendieck-Verdier category carries the structure of a strong Frobenius LD-functor.

In Section \ref{sec:LD-FrobAlg}, we study LD-Frobenius algebras in linearly distributive categories, proving several key properties graphically and providing examples. We then investigate Frobenius forms and establish an isomorphism between the category of modules and the category of comodules of an LD-Frobenius algebra.

Section \ref{sec:FrobSchurIndic} begins in Theorem \ref{equivalent def of pivotal structure} with a new characterization of pivotal Grothendieck-Verdier categories, which we show to be equivalent to Boyarchenko and Drinfeld’s definition \cite[Def 5.1]{BoDrinfeld}, though in a form more suited for the surface-diagrammatic calculus. We then offer a graphical definition of pivotal Frobenius LD-functors between pivotal Grothendieck-Verdier categories. Finally, we introduce higher Frobenius-Schur indicators for suitably finite $k$-linear pivotal Grothendieck-Verdier categories and discuss their key properties. 

Throughout the paper, we provide STL files for 3D models of relevant surface diagrams. Additionally, we supply various HOM files, the file format used by the web-based graphical proof assistant homotopy.io for finitely-presented globular $n$-categories. These resources are available on the first \href{https://maxdemirdilek.github.io/Research/SurfaceDiagrams}{author's website}.

\newpage

\section{Categorical preliminaries}\label{sec:Categorical prelims}
In this paper, $n$-category and $n$-functor mean \emph{strict} $n$-category and \emph{strict} $n$-functor. Following \cite[Def. 2.6]{GPSCoTri}, a one-object tricategory which is locally a $2$-category is called a \emph{monoidal $2$-category}. One-object $3$-categories are referred to as \emph{strict monoidal $2$-categories}.
\subsection{Surface diagrams for monoidal bicategories}\label{sec:surfdiagforGray}
Our presentation of the three-dimensional graphical calculus for strict monoidal $2$-categories is informal. More detailed accounts of the calculus for the more general notion of Gray categories (also known as \emph{semistrict tricategories}) appear in \cite{HummonPhD} and \cite{barrett2024gray}. %Manifold diagrams, which generalize surface diagrams to higher dimensions, are discussed in \cite{dorn2023associative, Dorn:2022vwk}.

\subsubsection{The macrocosm $\cV{\textnormal-}\mathsf{Cat}$}\label{main example:V-Cat}
We aim to apply the surface-diagrammatic calculus to ordinary categories, preadditive categories, and $k$-linear categories. Formally, this can be set up as follows:

Let $(\cV,\otimes,J)$ be a symmetric monoidal category. Consider the monoidal $2$-category $\cV{\textnormal-}\mathsf{Cat}$ of small $\cV$-enriched categories, $\cV$-enriched functors, and $\cV$-enriched natural transformations. The monoidal product of $\cV{\textnormal-}\mathsf{Cat}$ is the Kelly tensor product of $\cV$-enriched categories \cite[\S 1.4]{KellyEnCat}. Let $\cI$ denote the \emph{unit $\cV$-category}, which consists of the single object $\ast$ with $\operatorname{Hom}_{\cI}(\ast,\ast)=J$. 

By standard arguments \cite[Cor. 8.4]{GPSCoTri}, we may assume, without loss of generality, that $\cV{\textnormal-}\mathsf{Cat}$ is a strict monoidal $2$-category. We will maintain this strictness assumption from here on, as it is necessary for the surface-diagrammatic calculus.

Unless otherwise stated, we let $(\cV,\otimes,J)$ be the cartesian monoidal category of sets $(\mathsf{Set},\times,\cI)$. In other words, we will mainly work with $\mathsf{Cat}$, i.e. a strictification of the monoidal $2$-category of small categories, functors, and natural transformations. We do this to avoid clumsy terminology and notation, even though our graphical proofs go through to the enriched setting. We denote by $\times$ the cartesian product of categories. The terminal category is the monoidal unit $\cI$ of the strict monoidal $2$-category $\mathsf{Cat}$. Occasionally, we will let $(\cV,\otimes,J)$ be the symmetric monoidal category of abelian groups or (finite-dimensional) vector spaces.

\subsubsection{Categories}\label{sec:cats}
We shall develop a graphical calculus involving categories, functors, and natural
transformations. The calculus is drawn in framed three-space, which we build up gradually. 

\begin{figure}[H]
     \centering
     \begin{subfigure}[b]{0.32\textwidth}
         \centering
         \includegraphics[width=0.65\textwidth]{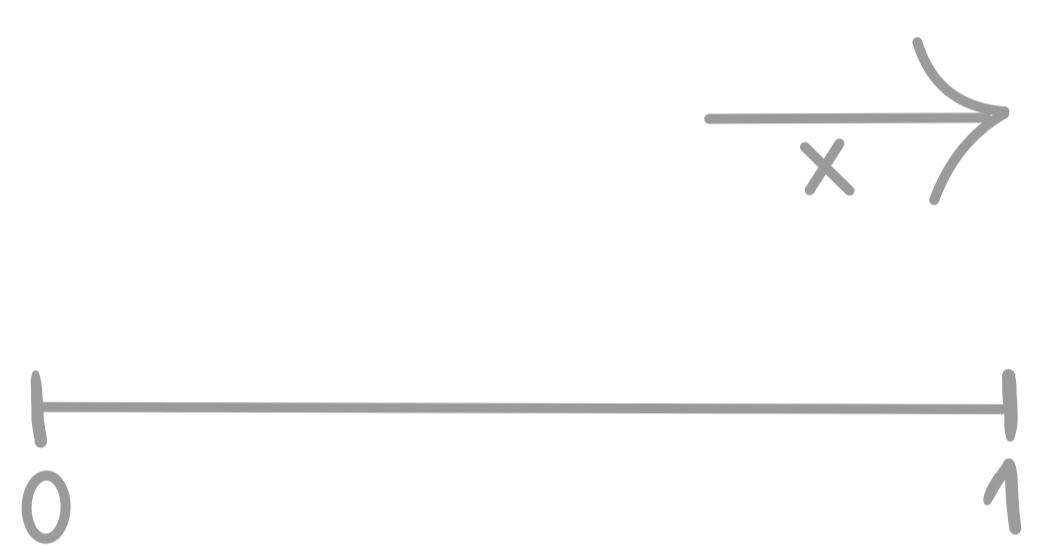}
         \caption{The canvas.}
         \label{fig:canvas1}
     \end{subfigure}
     \hfill
     \begin{subfigure}[b]{0.32\textwidth}
         \centering
         \includegraphics[width=0.65\textwidth]{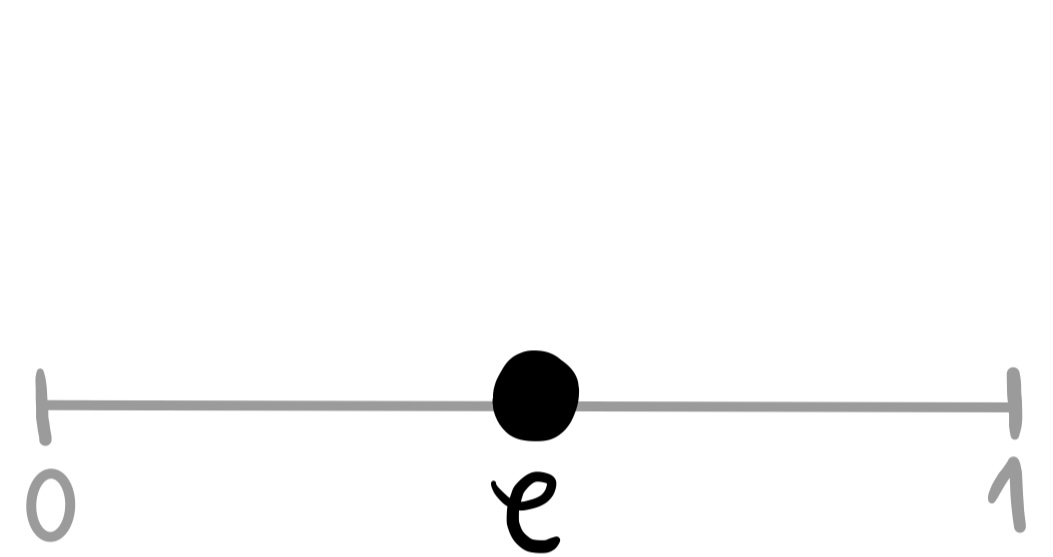}
         \caption{A category $\cC$.}
         \label{fig:a cat C}
     \end{subfigure}
          \hfill
     \begin{subfigure}[b]{0.32\textwidth}
         \centering
         \includegraphics[width=0.65\textwidth]{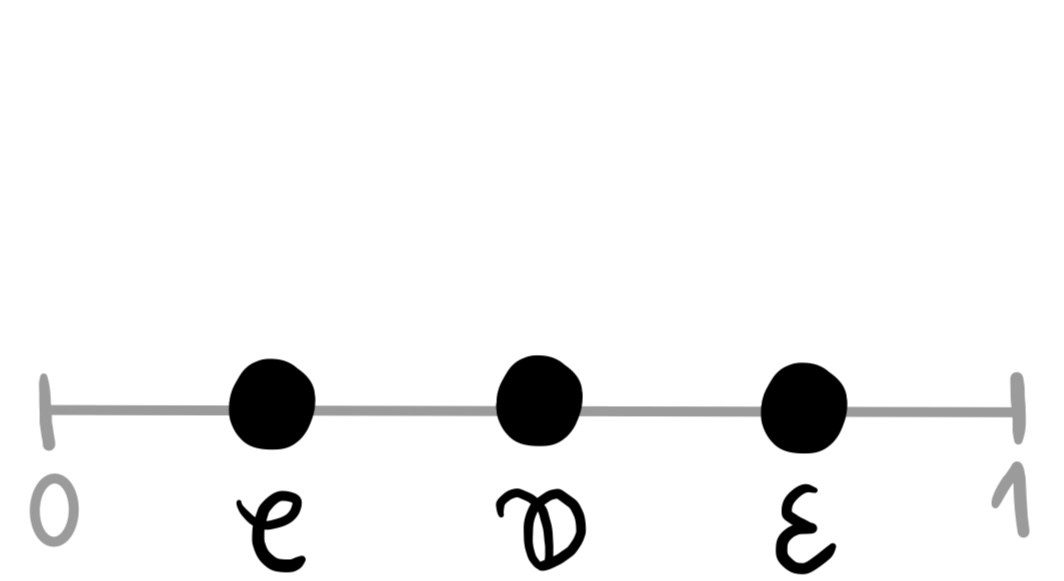}
         \caption{$\cC\times \cD \times \cE$.}
         \label{fig:cart prod cat}
     \end{subfigure}
        \caption{Depicting categories.}
        \label{fig:categories}
\end{figure}

We start with categories only. Ultimately, categories will correspond to $2$-strata of stratified surfaces. For the moment, we have only one dimension, a framed closed unit interval $I:=[0,1]\subset \mathbb{R}$, as in Figure \ref{fig:canvas1}. This interval has finitely many marked points that are labeled by categories, see, e.g., Figure \ref{fig:a cat C}. Categories can be composed using the cartesian product (the Kelly tensor product, more generally). Geometrically, this composition is depicted by juxtaposition, as indicated in Figure \ref{fig:cart prod cat}. The chosen $1$-framing of the interior $(0,1)$, shown in the top right corner of Figure \ref{fig:canvas1}, indicates the direction of composition. The unit category $\cI$ is transparent. Note that we have already used our assumption that our monoidal $2$-category is strict.

\subsubsection{Functors}\label{sec:funs}

\begin{figure}[H]
     \centering
     \begin{subfigure}[b]{0.27\textwidth}
         \centering
         \includegraphics[width=0.999\textwidth]{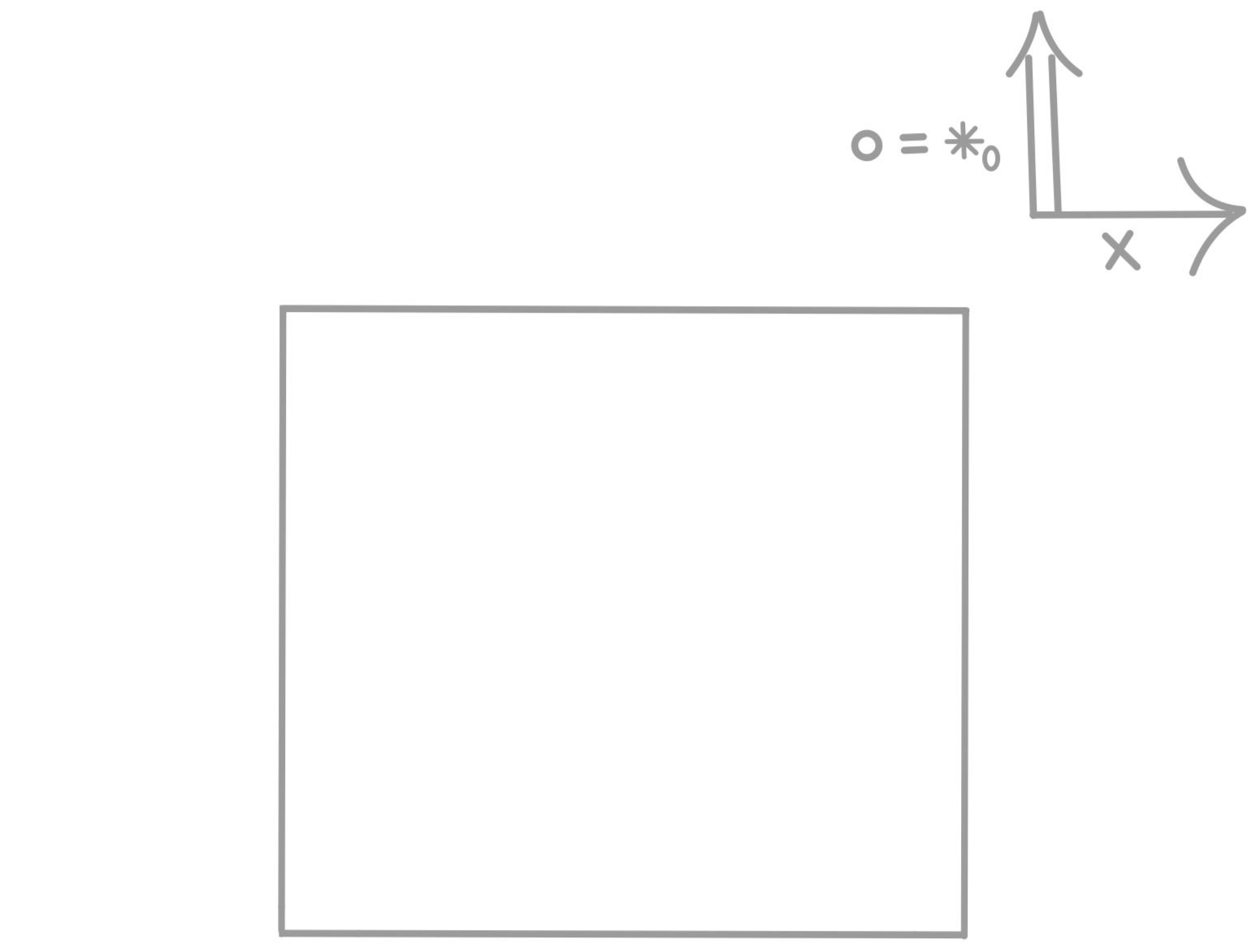}
         \caption{The canvas.}
         \label{fig:canvas2}
     \end{subfigure}
     \hfill
     \begin{subfigure}[b]{0.23\textwidth}
         \centering
         \includegraphics[width=0.9\textwidth]{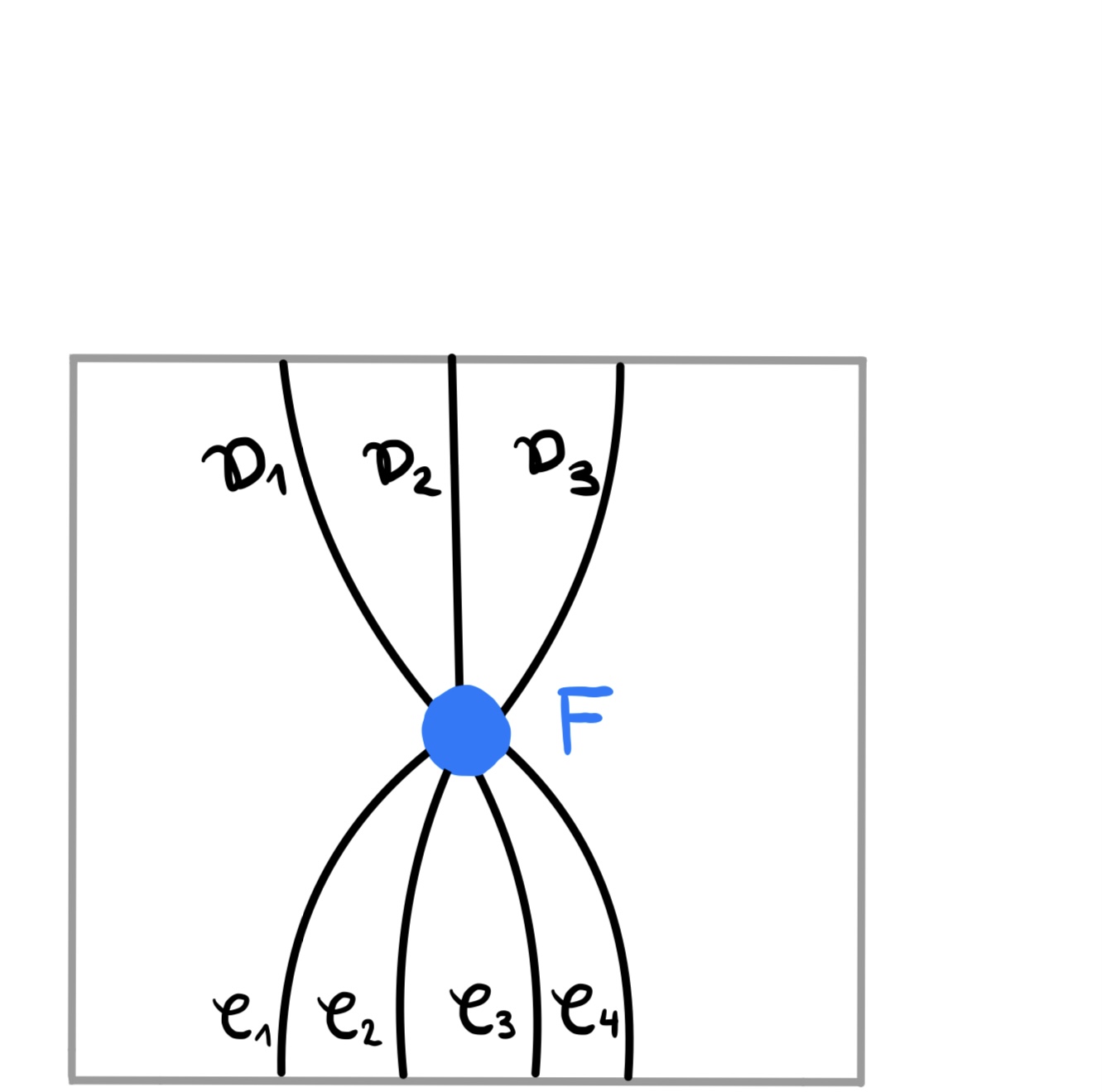}
         \caption{The functor $F$.}
         \label{fig:fun F}
     \end{subfigure}
          \hfill
     \begin{subfigure}[b]{0.23\textwidth}
         \centering
         \includegraphics[width=0.94\textwidth]{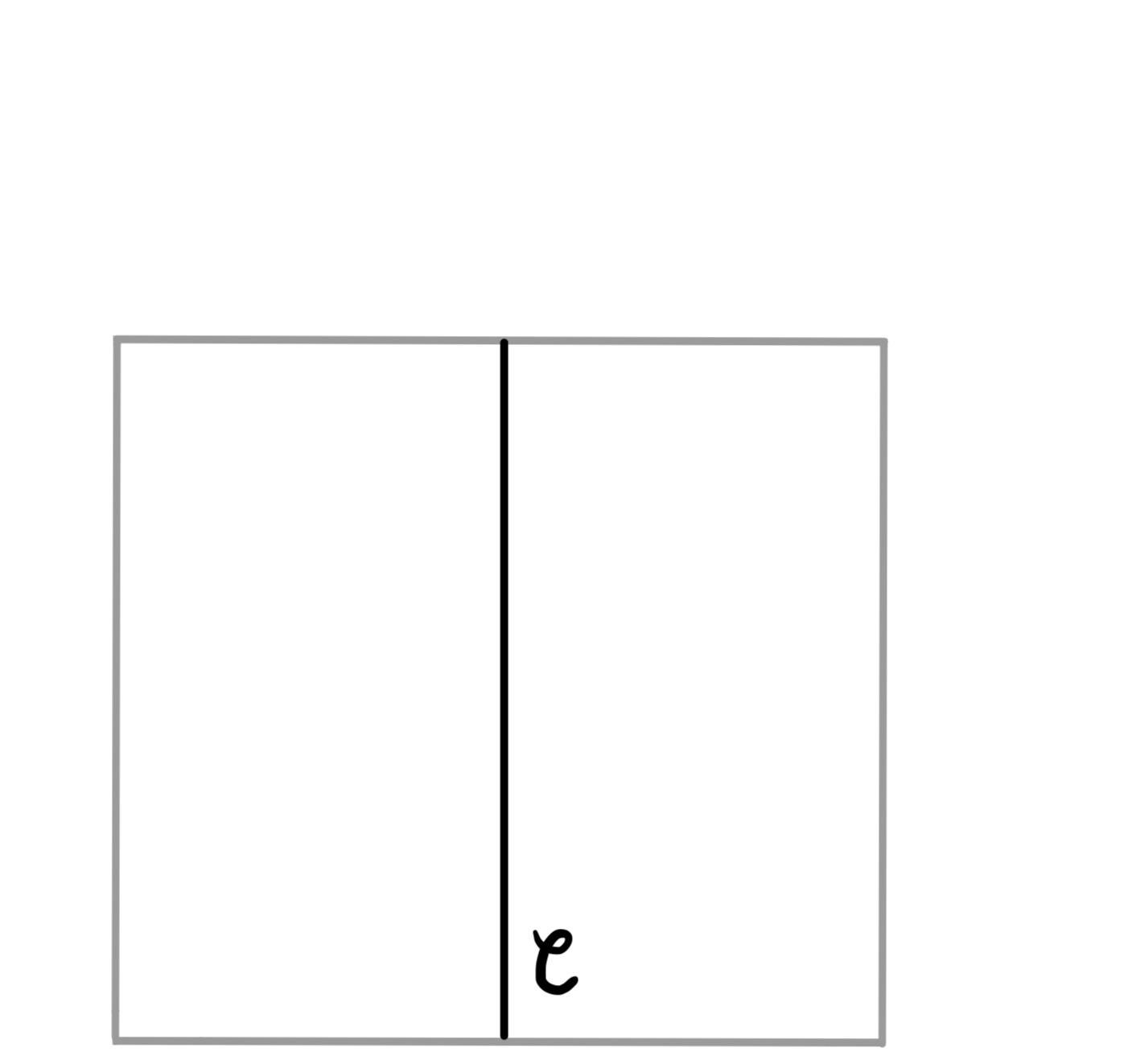}
         \caption{$\operatorname{id}_{\cC}$.}
         \label{fig:idC}
     \end{subfigure}
               \hfill
     \begin{subfigure}[b]{0.23\textwidth}
         \centering
         \includegraphics[width=0.96\textwidth]{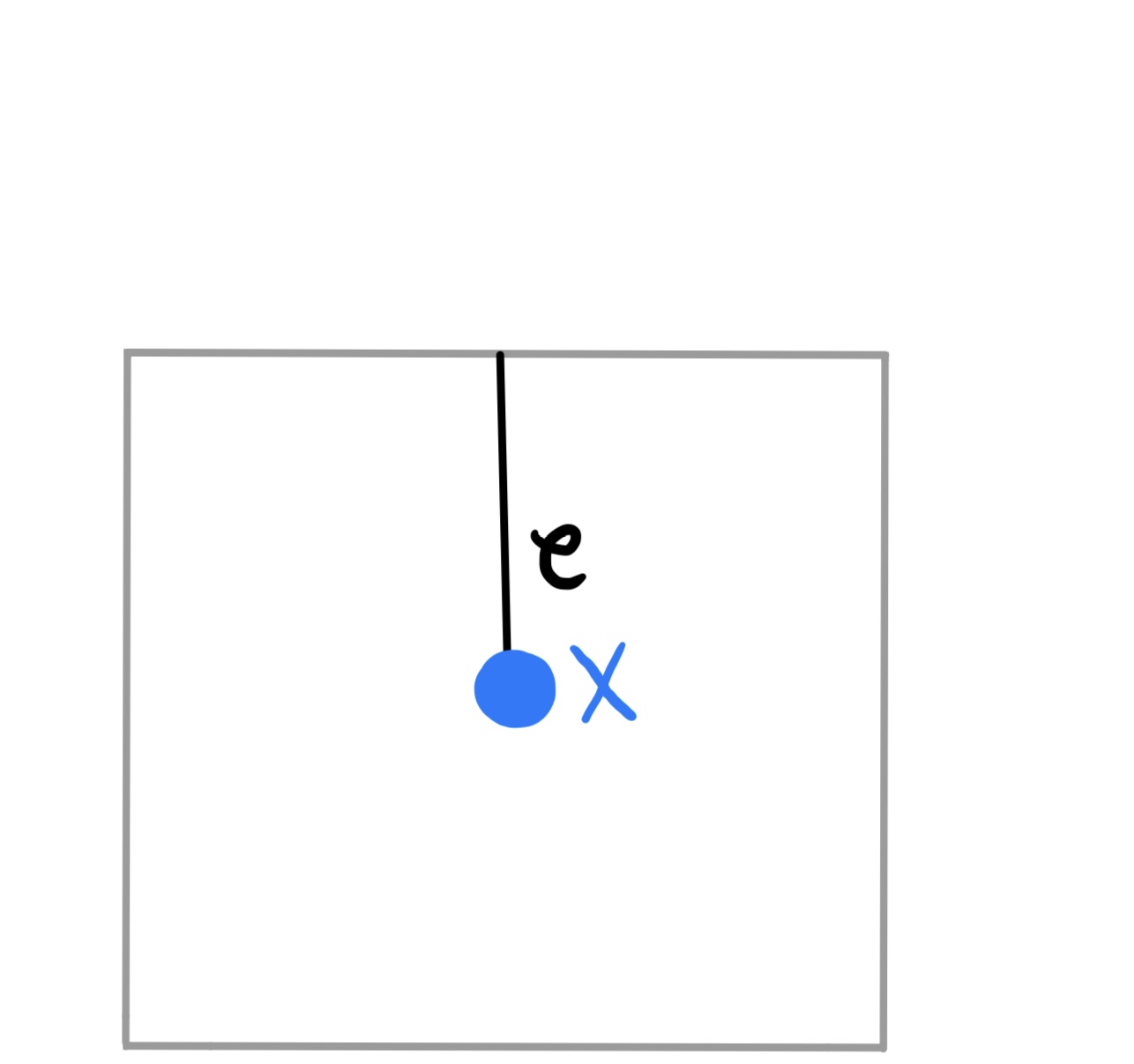}
         \caption{An object $X\in \cC$.}
         \label{fig:obX}
     \end{subfigure}
        \caption{Depicting functors.}
        \label{fig:functors}
\end{figure}

\paragraph{The canvas.}
Let us also include functors, and hence move up one dimension. Consider the closed unit square $I^2$. We equip its interior $(0,1)^2$ with the $2$-framing indicated in the upper right corner of Figure \ref{fig:canvas2}. 

\paragraph{Graphical depiction.}
Given categories $\cC_i$ and $\cD_j$, Figure \ref{fig:fun F} shows a functor 
    \begin{equation*}
    F\colon\,\bigtimes_{i=1}^4\cC_i\,\longrightarrow \,\bigtimes_{j=1}^3\cD_j.
    \end{equation*} 
Figure \ref{fig:fun F} is obtained as follows: In the way explained in Section \ref{sec:cats}, we first draw the source and target of $F$ onto the bottom edge $I\times \{0\}$ and the top edge $I\times \{1\}$, respectively. We then connect the source and target to the blue \emph{singular vertex} $F$, a process that we might call \emph{conification}. The identity functor on a category $\cC$ is transparent; see Figure \ref{fig:idC}. A (generalized) object $X$ in a category $\cC$, i.e. a functor $\cI\ra \cC$ from the unit category, is depicted in Figure \ref{fig:obX}. 

\paragraph{Compositions.}

\begin{figure}[H]
    \centering
\includegraphics[scale=0.13]{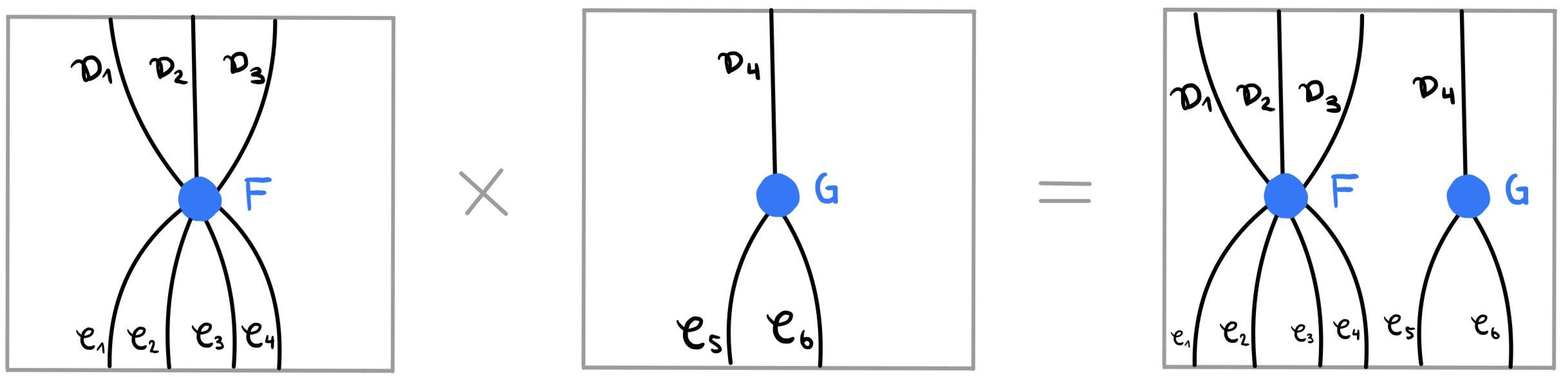}
 \caption{Cartesian product $\times$ of functors.}
  \label{fig:cart prod fun}
\end{figure}

\begin{figure}[H]
    \centering
\includegraphics[scale=0.13]{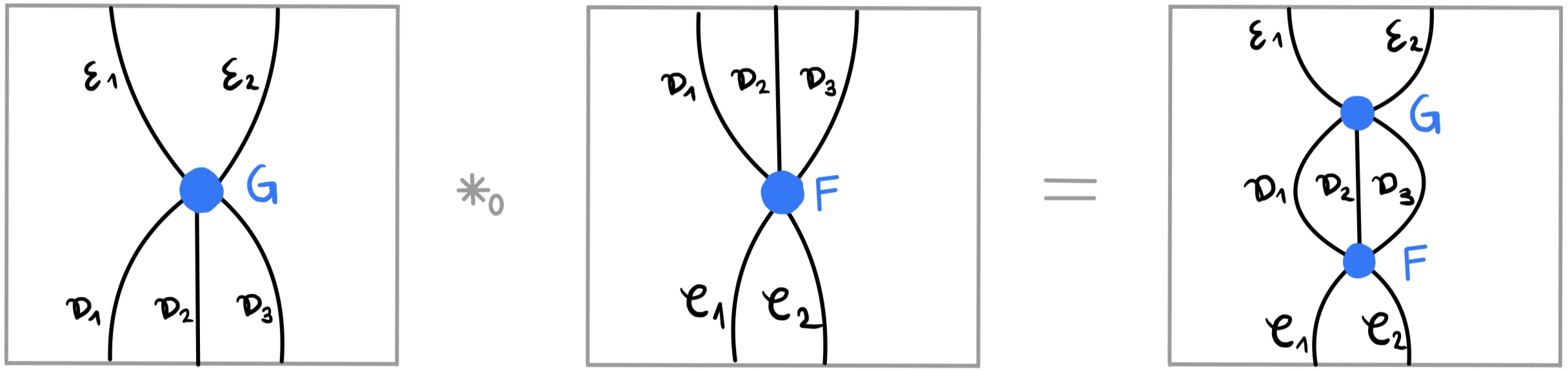}
 \caption{Composition $\ast_0$ of functors.}
  \label{fig:comp of functors}
\end{figure}

Functors can be composed in two ways. Their cartesian product is given by horizontal juxtaposition as shown in Figure \ref{fig:cart prod fun}. The vertical merging of functor wires describes the composition $\circ=\ast_0$ of functors along (cartesian products of) categories; see Figure \ref{fig:comp of functors}.

\subsubsection{Natural transformations.}
\paragraph{The canvas.}\label{subsec:The canvas}
We finally want to include natural transformations. This forces us to move up one more dimension. The ultimate canvas, which will allow us to also treat natural transformations, is the standard closed unit cube $I^{3}$ in three-space $\mathbb{R}^{3}$. The sets $\{1\}\times I \times I$, $I \times \{0\} \times I$ and $I\times I \times \{0\}$ are called \emph{front face}, \emph{side face} and \emph{bottom face}, respectively. 

\begin{figure}[H]
    \centering
\includegraphics[scale=0.1]{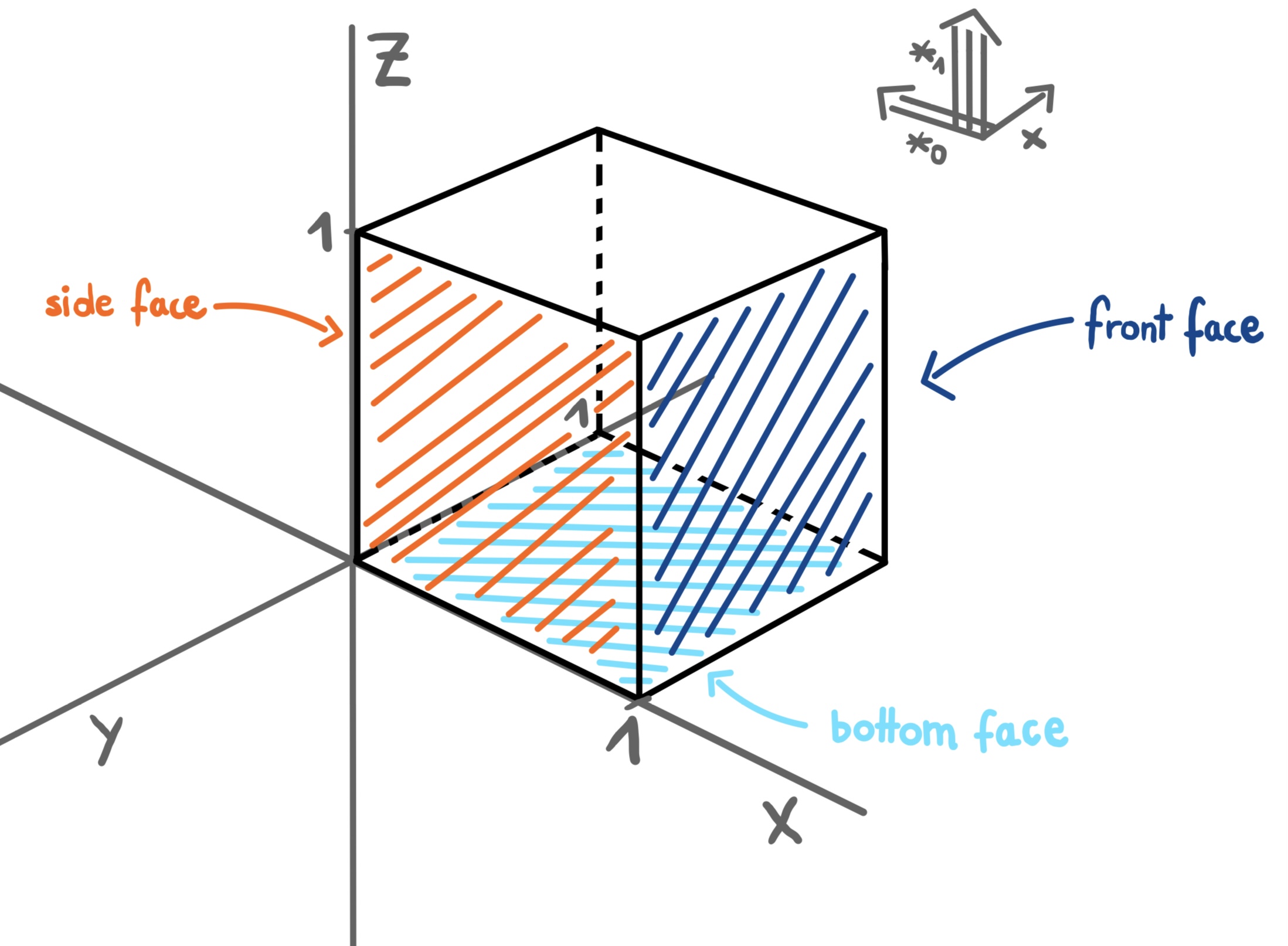}
 \caption{The $3$-framed canvas.}
  \label{fig:canvas}
\end{figure}

We endow the three-manifold $(0,1)^{3}\subset \mathbb{R}^3$ with a $3$-framing. That is, for any point ${x\in (0,1)^{3}}$, we pick an ordered basis of the tangent space at $x$ as indicated in the upper right corner of Figure \ref{fig:canvas}.

The $3$-framing on our canvas specifies the directions in which we read our three-dimensional depictions of natural transformations: The three ordered basis vectors of each tangent space correspond to the three possible compositions of natural transformations in the strict monoidal $2$-category $\mathsf{Cat}$. More precisely, the $Y$-direction (single arrow), negative $X$-direction (double arrow) and the $Z$-direction (triple arrow) are used for the cartesian product $\times$, the (horizontal) composition $\ast_0$ along a category, and the (vertical) composition $\ast_1$ along a functor, respectively. 

\begin{figure}[H]
     \centering
     \begin{subfigure}[b]{0.48\textwidth}
         \centering
         \includegraphics[width=0.58\textwidth]{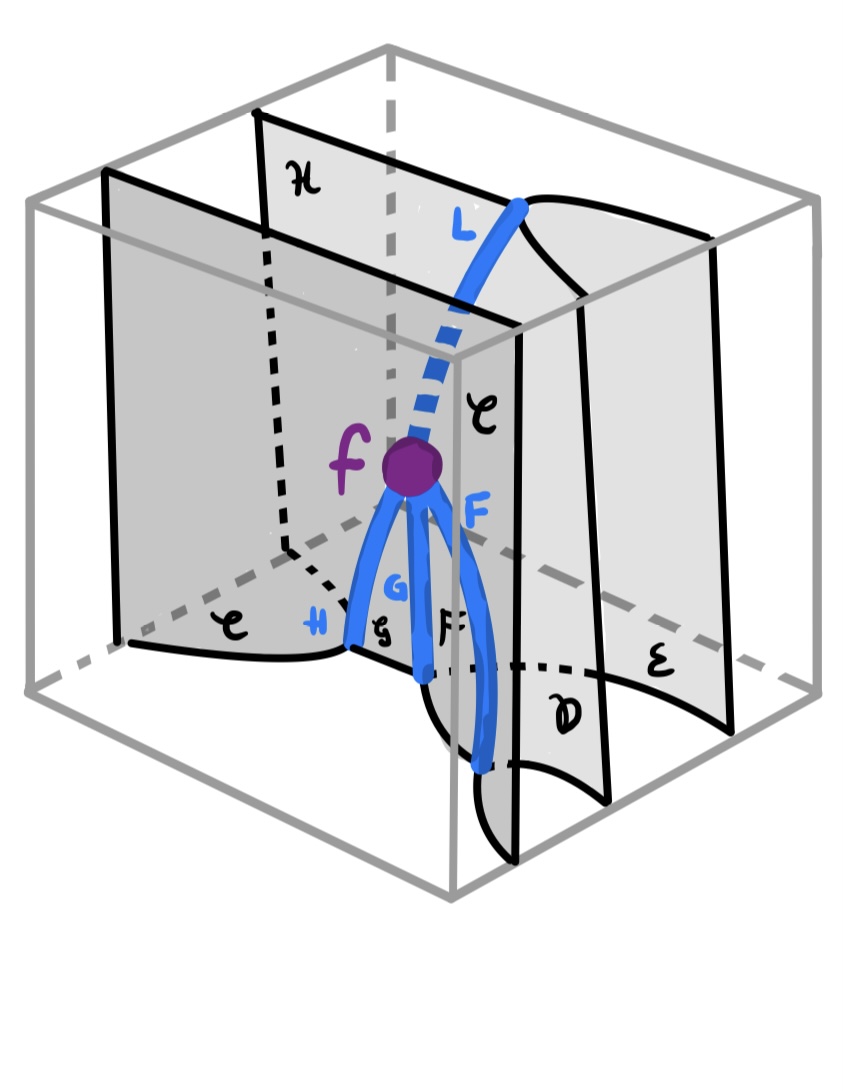}
         \caption{Surface diagram.}
         \label{fig:natTransfof}
     \end{subfigure}
     \hfill
     \begin{subfigure}[b]{0.48\textwidth}
         \centering
         \includegraphics[width=0.38\textwidth]{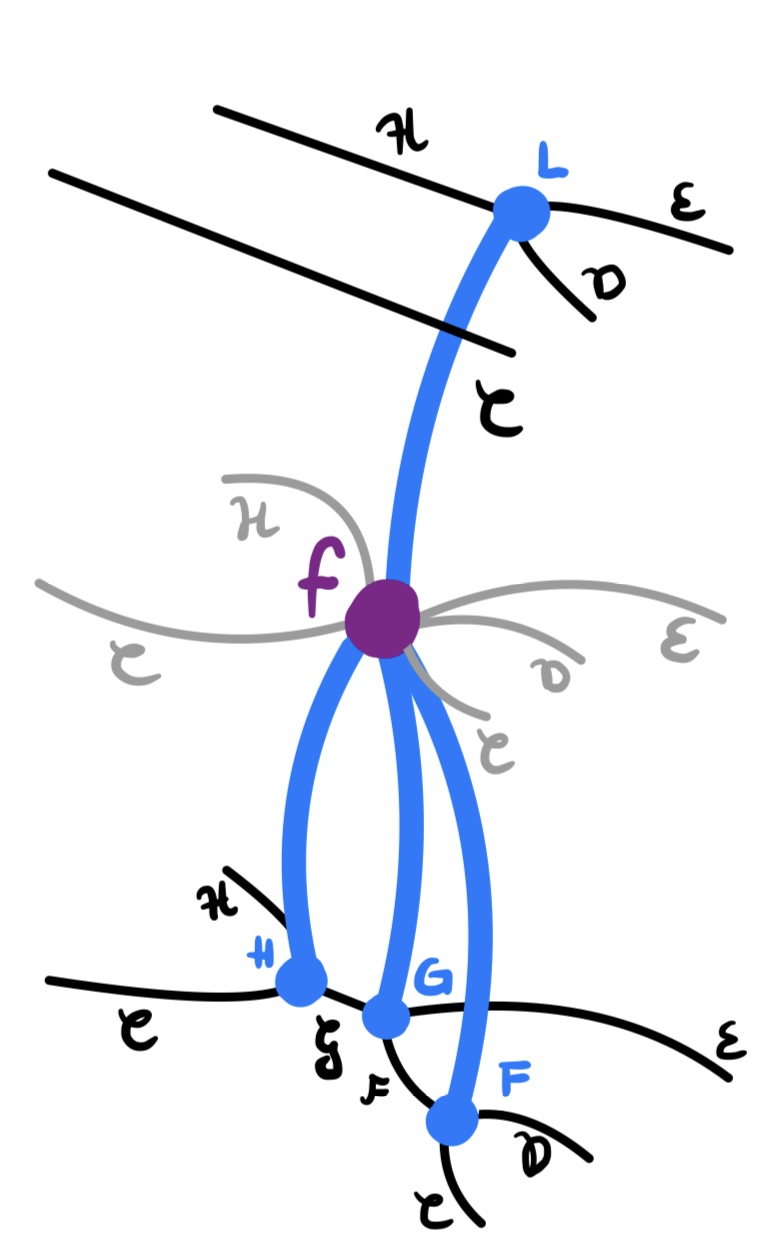}
         \caption{Constructing the surface diagram.}
         \label{fig:natTransfofIllust}
     \end{subfigure}
     \caption{The natural transformation $f$.}
\end{figure}

\paragraph{Graphical depiction.} Natural transformations are represented by labeled stratified surfaces. In these \emph{surface diagrams}, closed $2$-cells, $1$-cells and $0$-cells are labeled by categories, functors, and natural transformation, respectively. We illustrate this with an example. Figure \ref{fig:natTransfof} shows the surface diagram for a natural transformation
    \begin{equation*}
        f\colon H\circ G\circ(F\times \operatorname{id}_{\cE})\ra {\operatorname{id}_{\cC}}\times {L},
    \end{equation*} 
given functors $F\colon \cC \times \cD \ra \cF$, $G\colon \cF\times \cE\ra \cG$, $H\colon \cG\ra \cC\times \cH$, $L\colon \cD\times \cE \ra \cH.$ To obtain Figure \ref{fig:natTransfof}, we draw the source functor $H\circ G\circ(F\times \operatorname{id}_{\cE})$ onto the bottom face. Similarly, we draw the target functor onto the top face. As shown in Figure \ref{fig:natTransfofIllust}, we then connect the blue vertices to \emph{singular vertex} $f$, a process that can be again called \emph{conification}. Finally, we extend the $1$-cells on both the bottom and top faces (labeled by categories) to $2$-cells.

\begin{figure}[H]
\begin{subfigure}[b]{0.24\textwidth}
         \centering
         \includegraphics[width=0.74\textwidth]{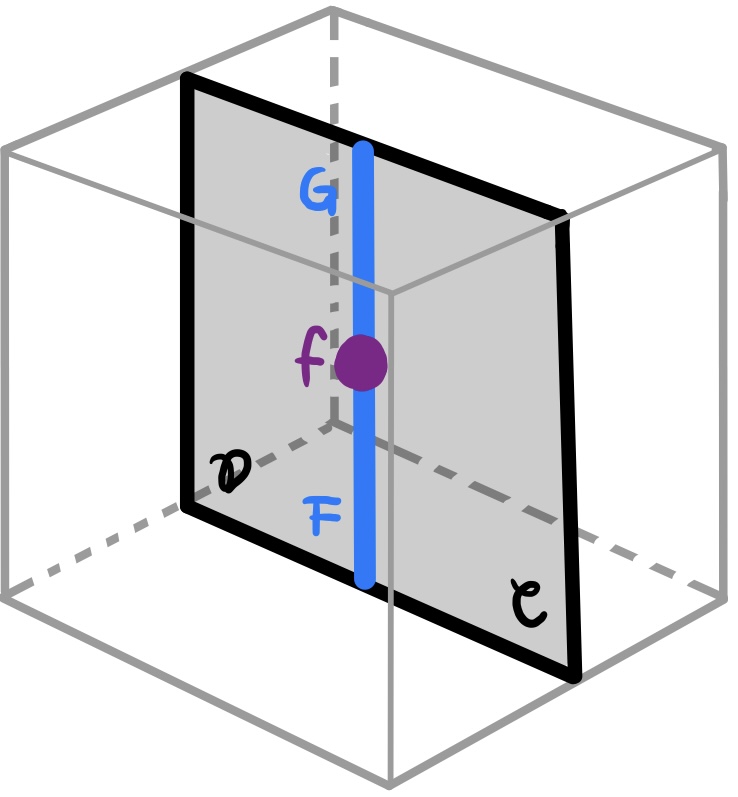}
         \caption{$f\colon F\ra G$.}
        \label{fig:f}
\end{subfigure}
\hfill
\begin{subfigure}[b]{0.24\textwidth}
         \centering
         \includegraphics[width=0.74\textwidth]{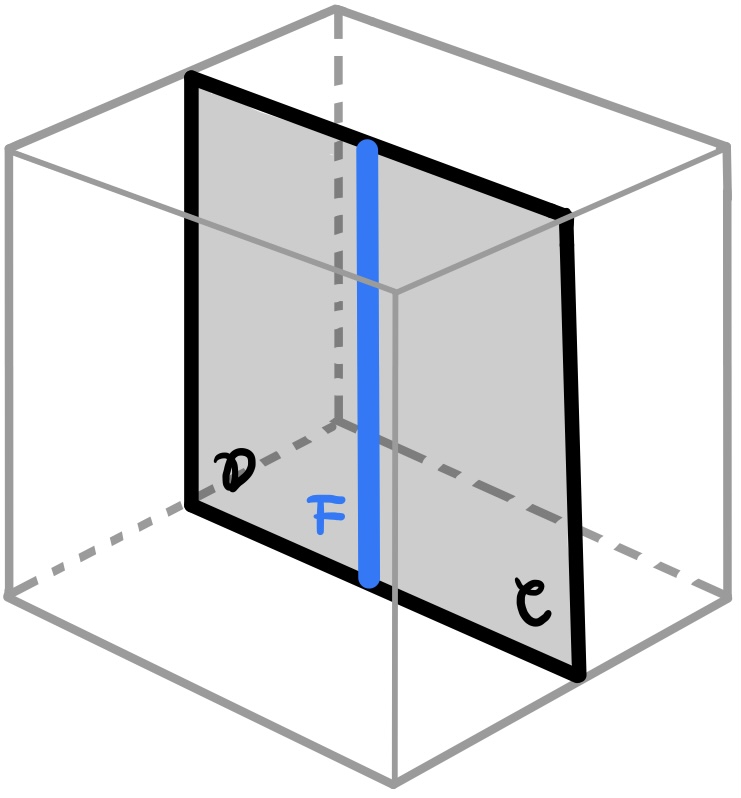}
         \caption{$\operatorname{id}_{F}$.}
         \label{fig:idF}
     \end{subfigure}
        \hfill
        \hfill
\begin{subfigure}[b]{0.24\textwidth}
         \centering
         \includegraphics[width=0.74\textwidth]{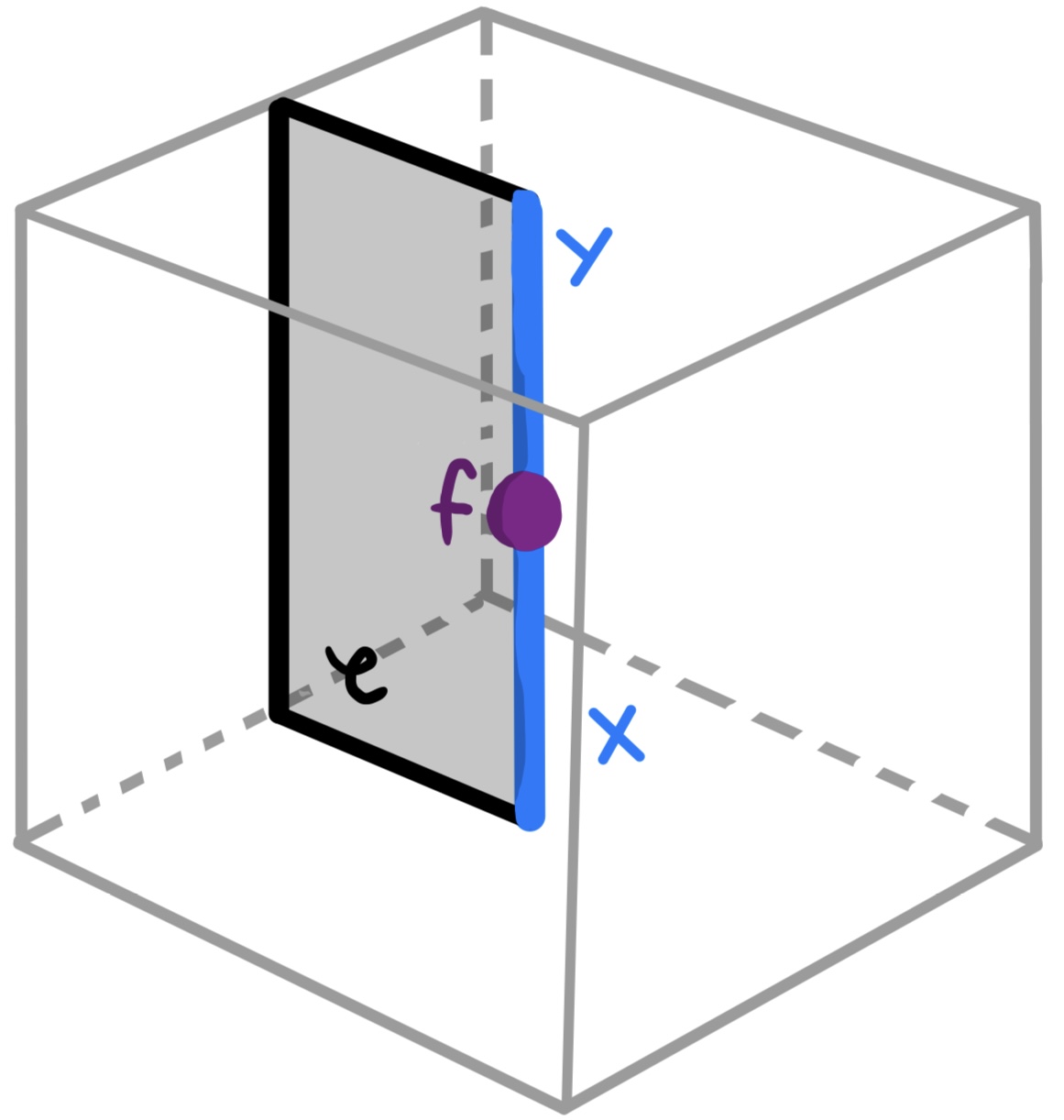}
\caption{$f\in \operatorname{Hom}_{\cC}(X,Y)$.}
\label{fig:morphf}
     \end{subfigure}
\begin{subfigure}[b]{0.24\textwidth}
         \centering
         \includegraphics[width=0.74\textwidth]{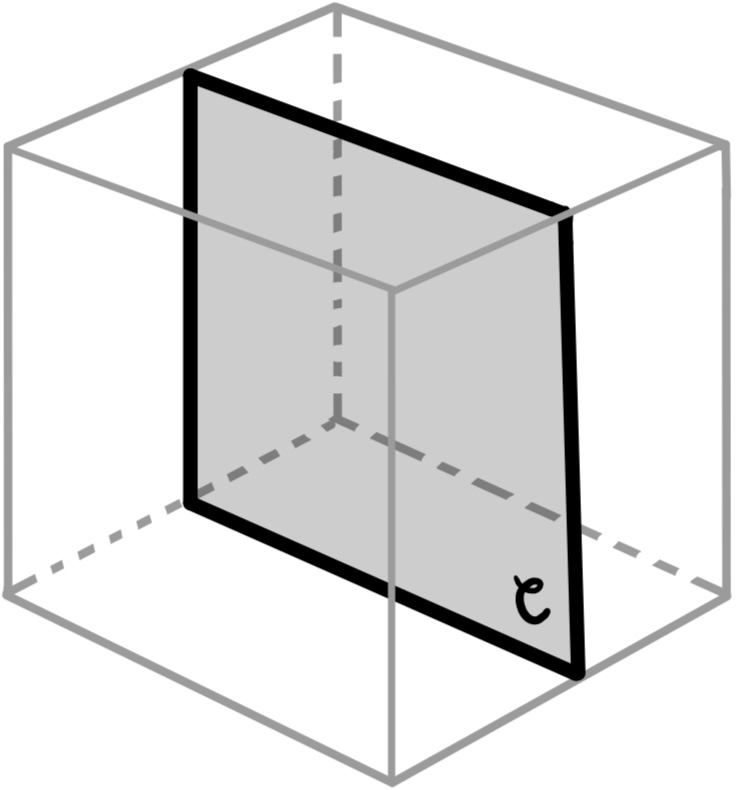}
\caption{$\operatorname{id}_{\operatorname{id}_{\cC}}$.}
\label{fig:ididF}
     \end{subfigure}
        \caption{Examples of natural transformations.}
        \label{fig:Natural transfos}
\end{figure}

\begin{remark}
 The bottom face of the canvas in Figure \ref{fig:natTransfof} shows the source functor of the natural transformation $f$. The source category $\cC\times\cD\times\cE$ of the source functor ${H\circ G\circ(F\times \operatorname{id}_{\cE})}$ is read off similarly by looking at the intersection of the front and bottom faces. Analogously, by considering other (subsets of) faces of our cubical canvas, the target operation in $\mathsf{Cat}$ can be understood graphically. The globular identities relating the source and the target operation in $\mathsf{Cat}$ imply that the intersection of the front face (and its opposite face) with the underlying surface of any surface diagram consists of a finite number of straight horizontal lines. This number can be zero if the unit category $\cI$ appears, e.g., see Figure \ref{fig:morphf}.
\end{remark} 

\begin{remark}
The surface diagram of a natural transformation is the Poincaré dual of its pasting diagram (as used, for example, in \cite{GPSCoTri}) within the one-object $3$-category $\mathsf{Cat}$.
\end{remark}

\paragraph{Examples.}
Figure \ref{fig:f} shows a natural transformation $f\colon F\ra G$. The identity natural transformation on a functor $F$ is, by convention, transparent; see Figure \ref{fig:idF}. A (generalized) morphism ${f\in \operatorname{Hom}_{\cC}(X,Y)}$ in a category $\cC$, i.e. a natural transformation from the functor $X\colon \cI\ra \cC$ to the functor $Y\colon \cI\ra \cC$, is depicted in Figure \ref{fig:morphf}. Finally, Figure \ref{fig:ididF} shows the identity natural transformation on the identity functor on a category $\cC$.

\begin{figure}[H]
    \centering
\includegraphics[scale=0.125]{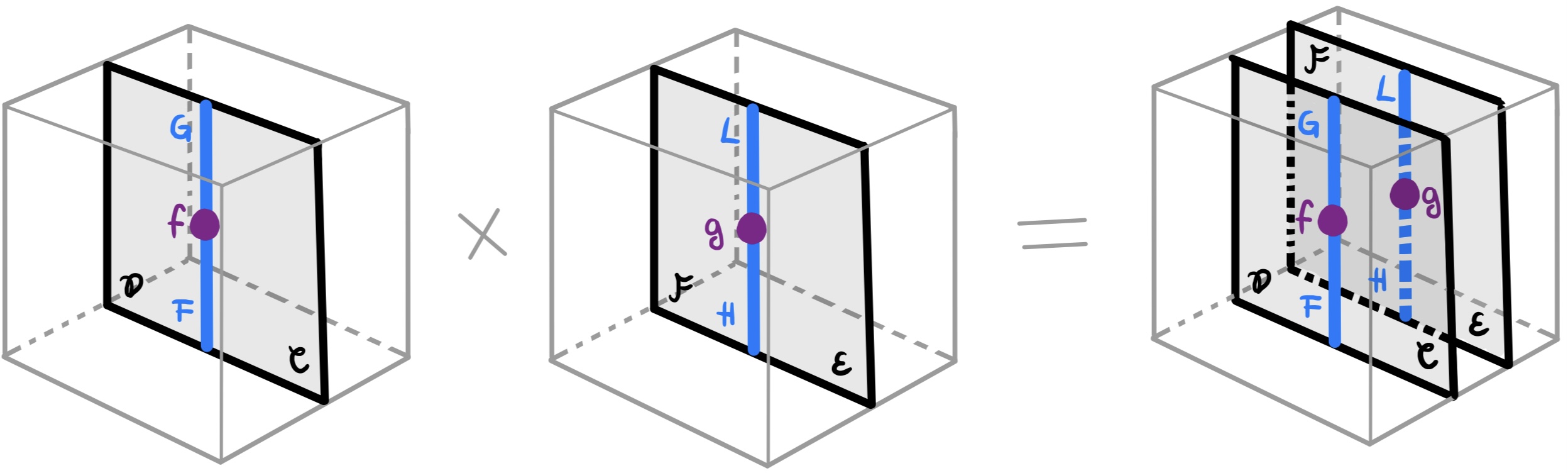}
 \caption{Composition $\times$ of natural transformations by stacking sheets.}
  \label{fig:cart prod nat}
\end{figure}

\begin{figure}[H]
    \centering
\includegraphics[scale=0.125]{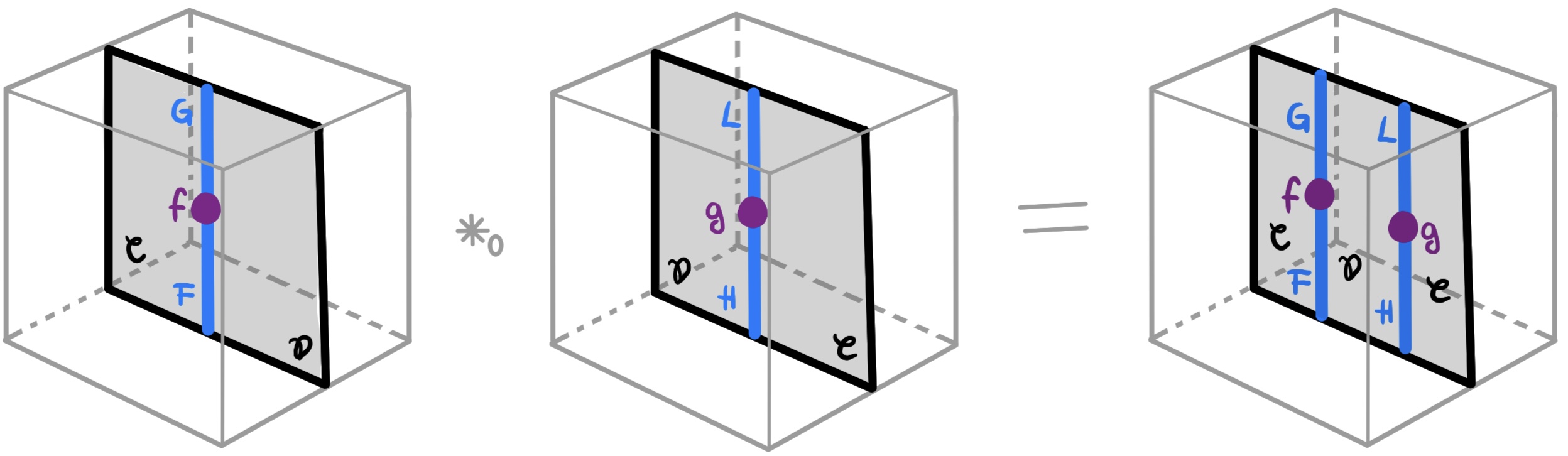}
 \caption{Horizontal composition $\ast_0$ by juxtaposing along closed $2$-cells.}
  \label{fig:hor comp of nats}
\end{figure}

\begin{figure}[H]
    \centering
\includegraphics[scale=0.125]{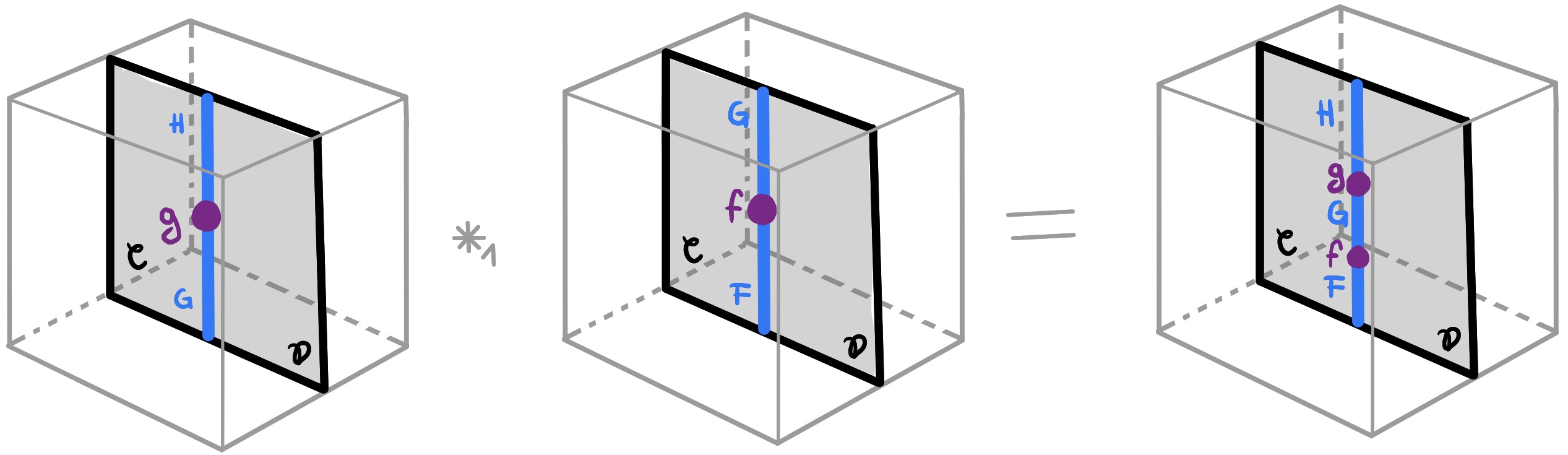}
 \caption{Vertical composition $\ast_1$ by juxtaposing along closed $1$-cells.}
  \label{fig:vert comp of nats}
\end{figure}

\paragraph{Compositions.} Natural transformations can be composed in three different ways. The juxtaposition in the Y-direction represents their cartesian product $\times$; see Figure \ref{fig:cart prod nat}. The horizontal composition $\ast_0$ along categories is shown in Figure \ref{fig:hor comp of nats}. Finally, Figure \ref{fig:vert comp of nats} describes the vertical composition $\ast_1$ along a functor; see also \cite[Fig. 9]{barrett2024gray}.

\begin{remark}
We can recover the well-known string diagrammatic calculus for the underlying $2$-category of our strict monoidal $2$-category. We achieve this by projecting the surface diagrams of the underlying $2$-category onto the side face.
\end{remark}

At this point, it is important to emphasize that we work with strict monoidal $2$-categories, rather than just semistrict ones:
\begin{remark}
    As explained in \cite[Notation 2.1.2]{douglas2018fusion}, for semistrict monoidal $2$-categories, the monoidal product of two $1$-morphisms or two $2$-morphisms is not uniquely defined. However, the monoidal product of an object with a $1$-morphism or a $2$-morphism is uniquely defined. Following \cite[Def. 2.1.1]{douglas2018fusion}, we denote, for an object $V$ in a semistrict monoidal $2$-category, its left and right monoidal product $2$-functors by $V\,\square\,-$ and $-\,\square\,V$, respectively. For $1$-morphisms $F\colon X\ra Y$, $G\colon U\ra V$, the convention is to define the symbol $F\,\square\,G$ as $(F\,\square\,V)\circ (X\,\square\,G)$, or alternatively, $(Y\,\square\,G)\circ (F\,\square\,U)$. One similarly defines the monoidal product of two $2$-morphisms. Graphically, for semistrict monoidal $2$-categories, this means that the two parallel $1$-cells and $0$-cells on the right-hand side of Figure \ref{fig:cart prod nat} do not lie on the same $X$-coordinate but are slightly nudged. In contrast, for our macrocosm $\mathsf{Cat}$, and more generally for any strict monoidal $2$-category, nudging is unnecessary, as the product of two functors (and also of two natural transformations) is uniquely defined. 
\end{remark}

\paragraph{Conventions and terminology.} Hereafter, we will not draw the cubical canvas. That is, we only depict the embedded surfaces diagram. Additionally, we will often not color 2-cells. If clear from context, we also leave out labels. For a given surface diagram, we refer to a closed $1$-cell representing a functor as a \emph{functor line}. We will often depict the functor line of a functor $\cI\ra \cC$, i.e. an object $X\in \cC$, by a squiggly black line. 

\subsubsection{Application: Monoidal categories and their functors}\label{application:moncats and functors}
The graphical calculus is best understood by an example. This example will reappear later.

\begin{figure}[H]
     \centering
     \begin{subfigure}[b]{0.26\textwidth}
         \centering
         \includegraphics[width=0.68\textwidth]{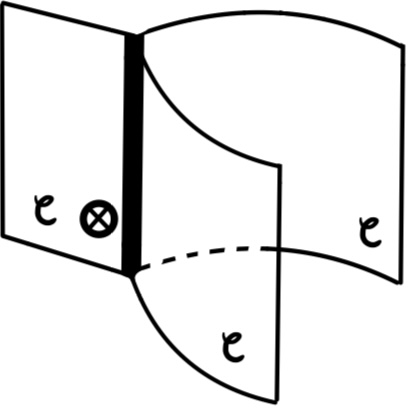}
         \caption{Monoidal tuning fork.}
         \label{fig:tunfork}
     \end{subfigure}
     \hfill
     \begin{subfigure}[b]{0.23\textwidth}
         \centering
         \includegraphics[width=0.58\textwidth]{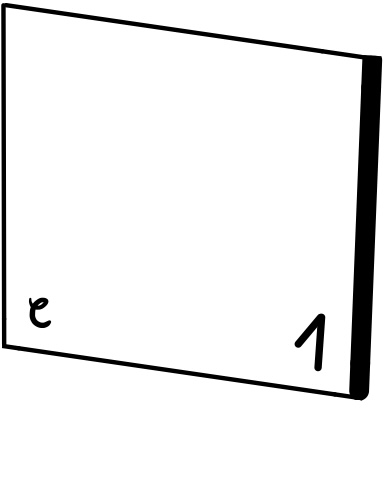}
         \caption{Monoidal unit $1$.}
         \label{fig:monUnit1}
     \end{subfigure}
          \hfill
     \begin{subfigure}[b]{0.23\textwidth}
         \centering
         \includegraphics[width=0.63\textwidth]{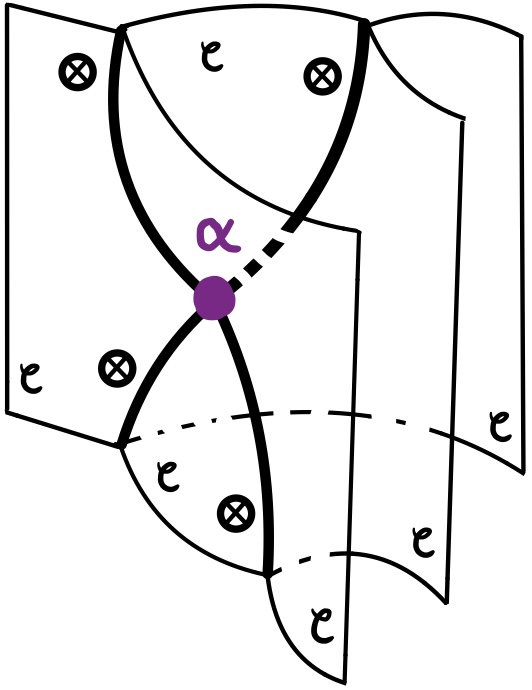}
         \caption{Associator $\alpha$.}
         \label{fig:alpha}
     \end{subfigure}
               \hfill
     \begin{subfigure}[b]{0.23\textwidth}
         \centering
         \includegraphics[width=0.68\textwidth]{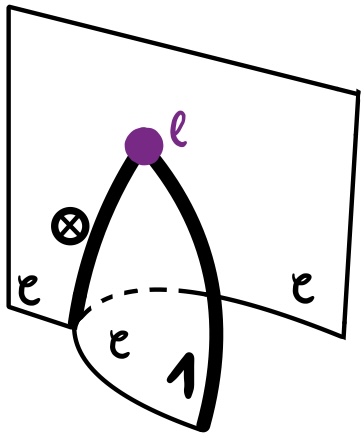}
         \caption{Left unitor $l$.}
         \label{fig:Left unitor}
     \end{subfigure}
        \caption{Data of a monoidal category $\cC$.}
        \label{fig:Monoidal categories.}
\end{figure}

\paragraph{Monoidal categories.}
Let $\cC$ be a monoidal category with monoidal product $\otimes \colon\cC\times \cC\ra \cC$ and monoidal unit $1\colon \cI\ra \cC$. The functor lines of the monoidal product $\otimes$ and the monoidal unit $1$ will both be colored black. Thus, we draw the identity natural transformation on the monoidal product $\otimes$ as shown in Figure \ref{fig:tunfork}. The underlying surface of the surface diagram in Figure \ref{fig:tunfork} is called \emph{monoidal tuning fork}. The identity morphism on the monoidal unit $1$ is shown in Figure \ref{fig:monUnit1}. Figures \ref{fig:alpha} and \ref{fig:Left unitor} depict the associator ${\alpha\colon 
    {\otimes}\; {\circ}\; (\otimes \times \operatorname{id_\cC})\ra {\otimes}\;{\circ}\; (\operatorname{id_\cC} \times \;\otimes)}$ and the left unitor $l\colon {\otimes}\;{\circ}\; (1 \times \operatorname{id_\cC}) \ra \operatorname{id_\cC}$, respectively.         
We obtain the surface diagram for the right unitor $r\colon {\otimes}\;{\circ}\; (\operatorname{id_\cC} \times \;1) \ra \operatorname{id_\cC}$ by reflecting the underlying surface of the surface diagram of Figure \ref{fig:Left unitor} along the side face. 

The reader is invited to check that the underlying associator surface of Figure \ref{fig:alpha} is made up of six sheets, all meeting at the singular vertex labeled $\alpha$. These sheets correspond to the six indices appearing in the $6j$-symbols; see \cite[Ex. 4.9.3]{EGNO}, \cite{baez}, \cite[Appendix F]{TuraevVir}. The underlying surface in Figure \ref{fig:alpha} is the Poincaré dual of a tetrahedron; see \cite{baez}.

\begin{remark}
STL files for 3D models of the surface diagrams in Figure \ref{fig:Monoidal categories.} can be viewed \href{https://maxdemirdilek.github.io/Research/SurfaceDiagrams}{here}.
\end{remark}

\begin{figure}[H]
     \centering
     \begin{subfigure}[b]{0.48\textwidth}
         \centering
         \includegraphics[width=0.9\textwidth]{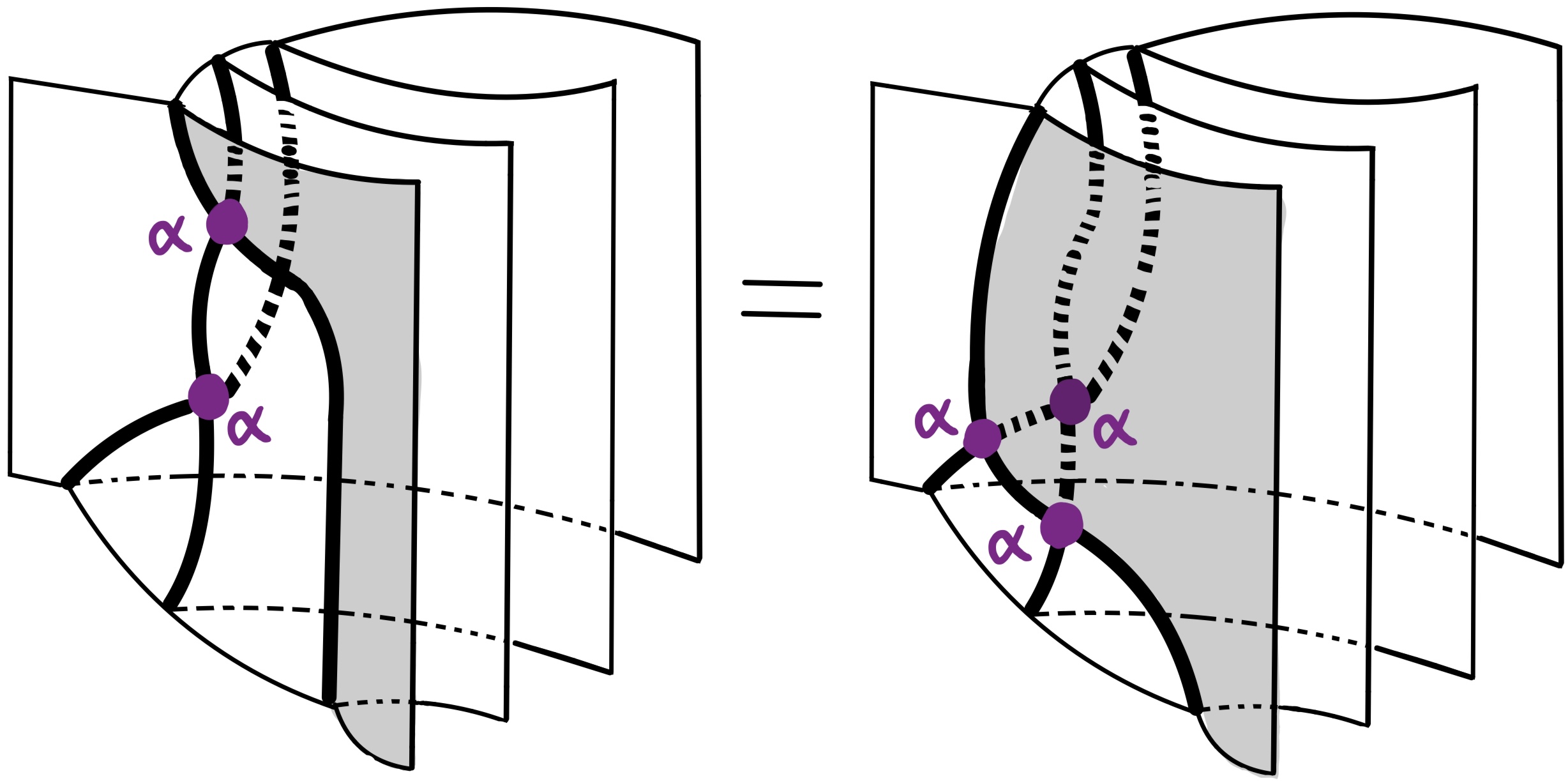}
         \caption{Mac Lane's pentagon.}
         \label{fig:MacPentagon}
     \end{subfigure}
     \hfill
     \begin{subfigure}[b]{0.48\textwidth}
         \centering
         \includegraphics[width=0.8\textwidth]{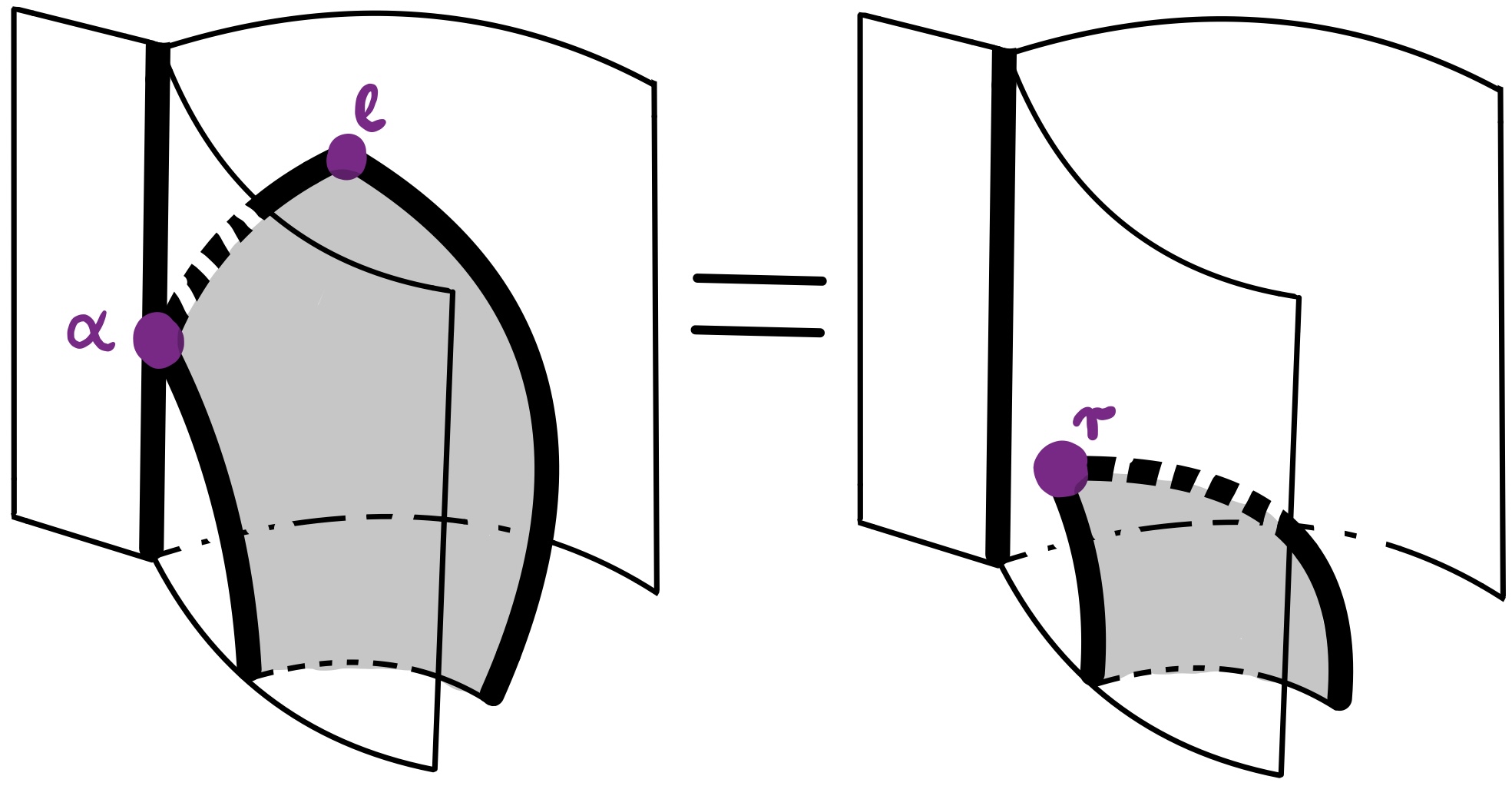}
         \caption{Mac Lane's triangle.}
         \label{fig:MacTriangle}
     \end{subfigure}
        \caption{Coherence axioms.}
        \label{fig:Coherence axioms.}
\end{figure}

\paragraph{Coherence axioms.}
Figure \ref{fig:Coherence axioms.} shows the coherence axioms for the monoidal category $\cC$. Graphically, they amount to moving a gluing boundary of a sheet along a surface. We have shaded a sheet being moved in light grey, purely to enhance readability.

The intersections of the two surface diagrams in Figure \ref{fig:MacPentagon} with the bottom face are identical and represent the same functor. The same applies to the intersection with the top face. This is an instance of globularity and applies to all our diagrams.

\medskip

\begin{figure}[H]
    \centering
\includegraphics[scale=0.09]{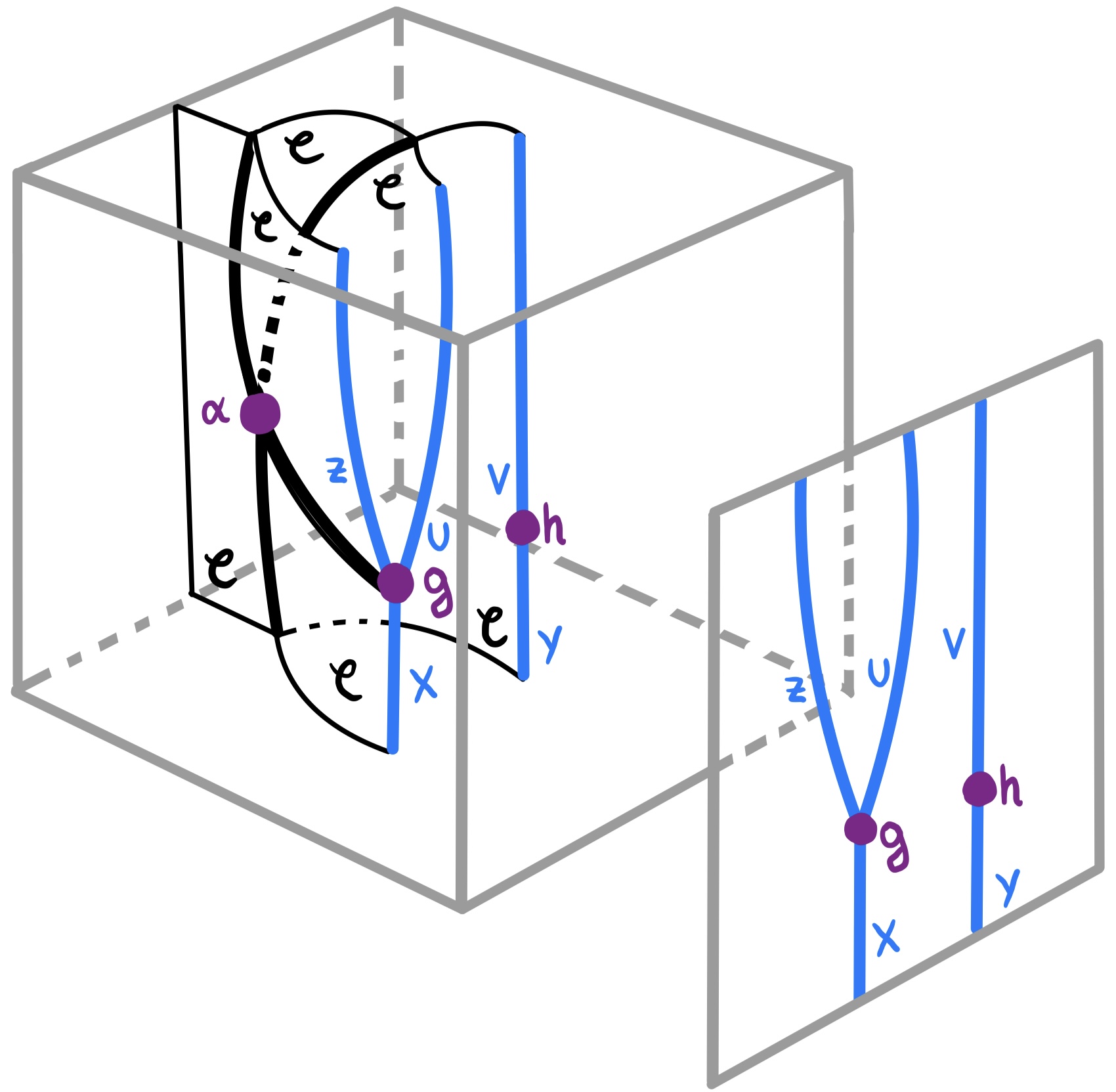}
 \caption{The front face of $f$ from equation (\ref{morph f}) yields a string diagram.}
  \label{fig:string diag surf diag}
\end{figure}

\begin{remark}
   A HOM file for the signature of monoidal categories is available \href{https://maxdemirdilek.github.io/Research/SurfaceDiagrams}{here.} In the spirit of homotopy.io, the coherence axioms are implemented not as equalities but as invertible cells between natural transformations.        
\end{remark}
\paragraph{Extracting string diagrams.}
The front face of surface diagrams yields string diagrams:
\begin{remark}\label{Mac Lane's strictification theorem}
Let $f$ be a formal morphism in a monoidal category $\cC$. That is, $f$ is made up of composites and monoidal products of morphisms in $\cC$, as well as the coherence data of $\cC$, e.g.
    \begin{equation}\label{morph f}
        f:=\alpha_{Z,U,V}\circ(g\otimes h)\in \operatorname{Hom}_{\cC}(X\otimes Y,Z\otimes (U\otimes V)),
    \end{equation}
for $g\in\operatorname{Hom}_{\cC}(X, Z\otimes U)$ and $h\in\operatorname{Hom}_{\cC}(Y,V)$. As illustrated in Figure \ref{fig:string diag surf diag}, when viewing the surface diagram of $f$ from the front face, we see the (light blue) string diagram for the morphism $f$ considered in a strictification of $\cC$. Informally: Surface diagrams for morphisms in a monoidal category are string diagrams with coherence morphisms in the back. From this viewpoint, Mac Lane's coherence theorem states that, if the string diagrams in the front faces of two surface diagrams (for two formal morphisms in the monoidal category $\cC$) agree up to planar isotopy and the source and target functor of these surface diagrams are equal, then we can rearrange sheets in the background so that the surface diagrams agree.
\end{remark}

An important feature of the three-dimensional calculus for the strict monoidal $2$-category $\mathsf{Cat}$, which was central to \cite{WillertonHopfM}, is its capacity to handle functors-with-structure. For future reference, we now discuss such functors-with-structure between monoidal categories.

\paragraph{Lax and oplax monoidal functors.}
Let $(\cC,\otimes,1)$ and $(\cD,\otimes,1)$ be monoidal categories. Let $F\colon \cC\ra \cD$ be a lax monoidal functor with not necessarily invertible natural transformations ${\varphi^2\colon\, {\otimes} \circ {(F\times F)} \ra F\circ \otimes}$ and $\varphi^0\colon \,1\ra F(1)$, called \emph{multiplication (morphism)} and \emph{unit (morphism)}, respectively. 

\begin{figure}[H]
     \centering
     \begin{subfigure}[b]{0.48\textwidth}
         \centering
         \includegraphics[width=0.6\textwidth]{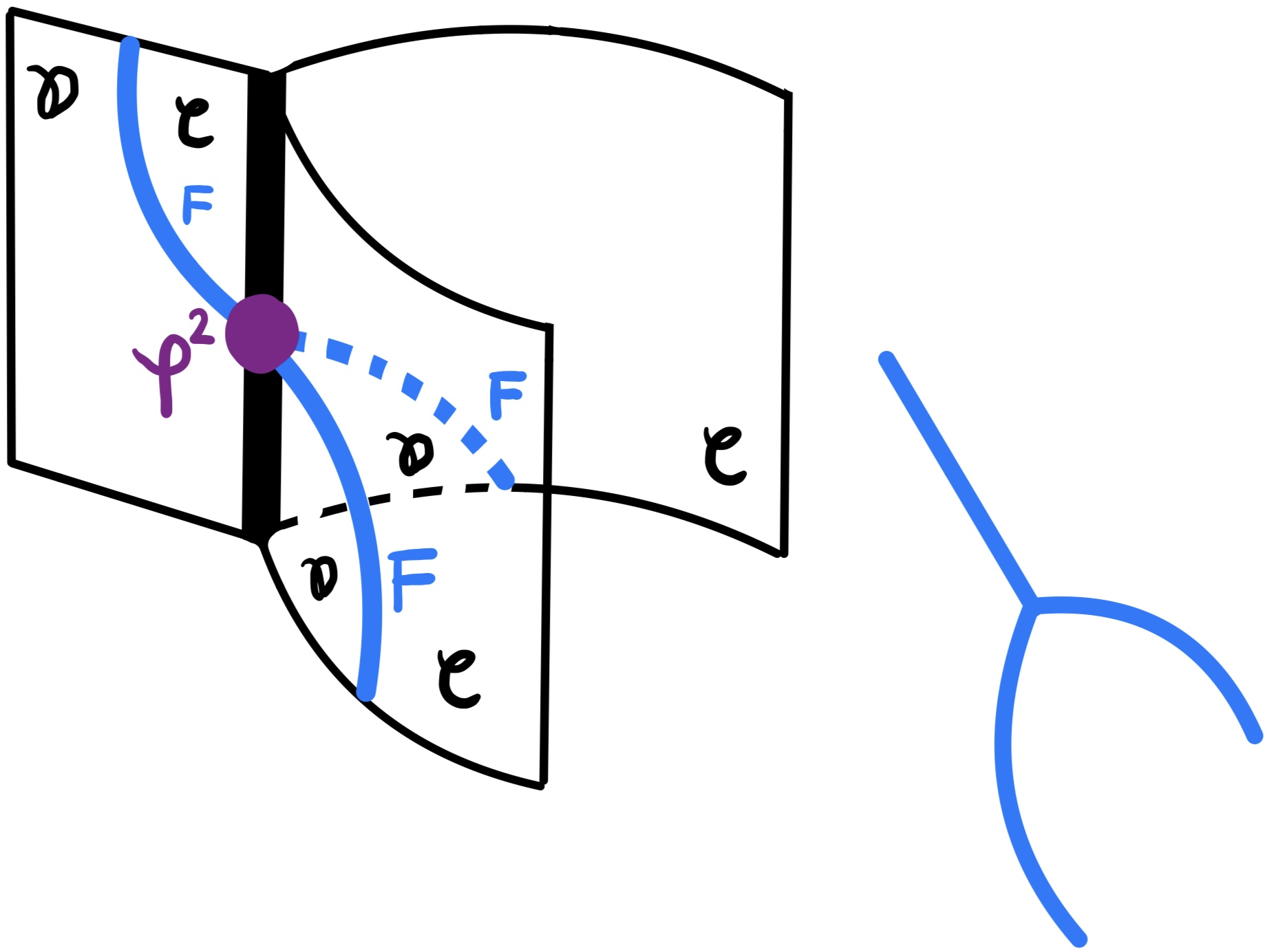}
         \caption{Multiplication $\varphi^2$.}
         \label{fig:varphi2}
     \end{subfigure}
     \hfill
     \begin{subfigure}[b]{0.48\textwidth}
         \centering
         \includegraphics[width=0.5\textwidth]{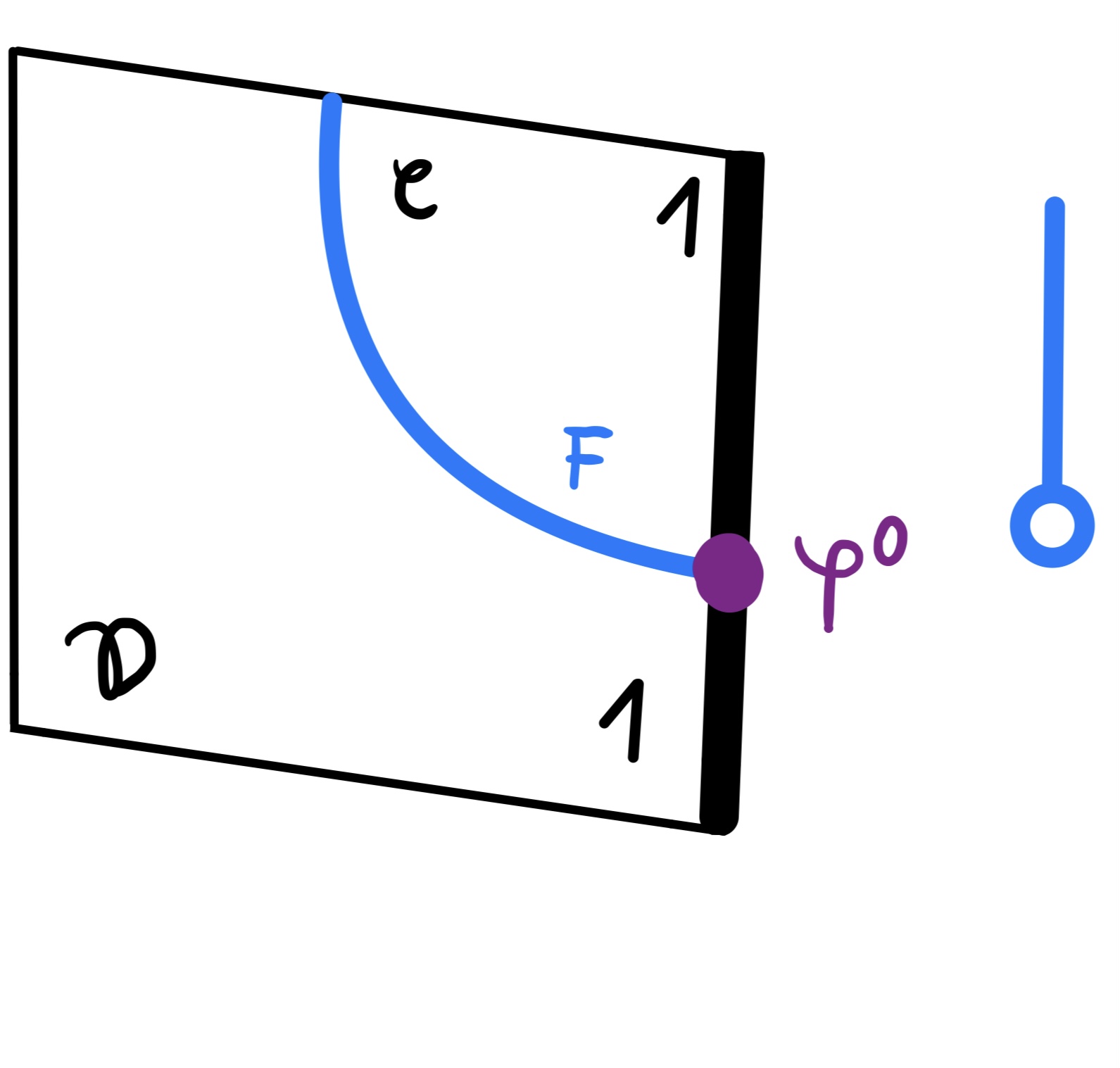}
         \caption{Unit $\varphi^0$.}
         \label{fig:varphi0}
     \end{subfigure}
        \caption{Data of a lax monoidal functor $(F,\varphi^2,\varphi^0)$ and their extracted string diagrams.}
        \label{fig:LaxMonFunctF.}
\end{figure}

The left-hand side of Figure \ref{fig:varphi2} shows the surface diagram of the multiplication $\varphi^2$. It is obtained by drawing the (blue) bident on the right-hand side of Figure \ref{fig:varphi2} onto the monoidal tuning fork. The singular vertex of the bident becomes the purple intersection point of the blue functor line of $F$ with the black functor line of $\otimes$. Figure \ref{fig:varphi0} shows the unit $\varphi^0$. The coherence axioms for the multiplication $\varphi^2$ and unit $\varphi^0$ in the language of surface diagrams can be found in \cite[\S 1.2.2.]{WillertonHopfM}. Informally, they arise from drawing the string diagrams encoding the associativity and unitality of an algebra in a strict monoidal category onto suitable stratified surfaces. We obtain the surface diagrams for an oplax monoidal functor by reflecting those for a lax monoidal functor along the bottom face.

\begin{remark}
The STL and HOM files for lax monoidal functors are available \href{https://maxdemirdilek.github.io/Research/SurfaceDiagrams}{here}. In the spirit of homotopy.io, the coherence axioms are implemented not as equalities, but as invertible cells between natural transformations.
\end{remark}

\subsection[Linearly distributive categories (LD-categories) and their functors]{Linearly distributive categories and their functors}\label{LD-cats and functors}
As mentioned in the introduction, we will see in Theorem \ref{thm:GV correspond LD neg} that Grothendieck-Verdier categories correspond to so-called linearly distributive categories with negation.

\begin{definition}\label{def:LDcat}(\cite[Def. 2.1]{WDC}). A (non-symmetric) \emph{linearly distributive category}, or \emph{LD\nobreakdash-category} for short, consists of a \mbox{category $\cC$} together with the following additional data:
\begin{itemize}
\item Two monoidal structures $(\otimes, 1)$ and $(\parLL,K)$ on $\cC$. 
\item Two natural transformations between functors 
$\cC\times \cC\times \cC \ra \cC$
    \begin{equation*}
        \distl \colon\; {\otimes}\;{\circ}\; (\operatorname{id_\cC} \times \parLL) \,\longrightarrow \, {\parLL}\; {\circ}\; (\otimes \times \operatorname{id_\cC}),
    \end{equation*}
    \begin{equation*}
        \;\distr \colon\; {\otimes}\;{\circ}\;(\parLL \times\; {\operatorname{id_\cC}}) \, \longrightarrow \, {\parLL}\; {\circ}\;(\operatorname{id_\cC}\times \; {\otimes),}\;
    \end{equation*}
called the \emph{left distributor} and the \emph{right distributor}.\footnote{This terminology originates from Cockett and Seely \cite{WDC}. Distributors can be viewed as resource-sensitive one-directional distributive laws. That is, unlike the distributive laws that arise in ring theory, in the definition of distributors, each variable occurs only once in source and target.} They are required to satisfy the coherence axioms in Appendix \ref{coherenceLD}. The axioms consist of four triangle and six pentagon diagrams. These encode the compatibility of the distributors with the unitors and associators of both monoidal structures. Distributors do not need to be invertible.
\end{itemize}
\end{definition}

Throughout the paper, we invoke the following example. For additional examples (which carry even more structure), see Examples \ref{pivotality of bimodule example}.
\begin{example}\label{bimodule example}
Consider the category $A\operatorname{-bimod}$ of bimodules over a finite-dimensional algebra $A$ over a field $k$. The usual right exact tensor product $\otimes_A$ of $A$-bimodules endows $A\operatorname{-bimod}$ with a monoidal structure with unit $A$. The $k$-linear dual $DM:=\operatorname{Hom}_{k}(A,k)$ of an $A$-bimodule $M$ carries an $A$-bimodule structure given by 
    \begin{equation*}
        (x.f.y)(m)\,:=\,f(y. m . x),
    \end{equation*}
for $x,y\in A$, $m\in M$ and $f\in DM$. Any $A$-bimodule becomes a bicomodule over the dual coalgebra $DA$. The left exact tensor product of $DA$-bicomodules, also known as a \emph{cotensor product}, endows $A\operatorname{-bimod}$ with a second monoidal structure $\otimes^A$ and monoidal unit $DA$. Together, these two monoidal structures make $A\operatorname{-bimod}$ into an LD-category. For further details and the definition of distributors, see \cite{fuchs2024grothendieckverdierdualitycategoriesbimodules}. (Although the paper \cite{fuchs2024grothendieckverdierdualitycategoriesbimodules} restricts attention to the full subcategory $A\operatorname{-bimod}^{\operatorname{f.d.}}$ of finite-dimensional $A$-bimodules, the LD-structure extends directly to the entire category $A\operatorname{-bimod}$.)
\end{example} 

Let us discuss surface diagrams for LD-categories: 
\begin{remark}
The $\otimes$-monoidal structure is depicted as explained in Section \ref{application:moncats and functors}. Specifically, the functor lines of the monoidal product $\otimes$ and the monoidal unit $1$ are colored black, while the functor lines of the other monoidal product $\parLL$ and its monoidal unit $K$ are colored red. As a result, the surface diagrams for the distributors feature the same configuration of sheets as those for the associators, but with a different coloring of the functor lines:

\begin{figure}[H]
    \centering
    \begin{subfigure}[b]{0.48\textwidth}
         \centering
         \includegraphics[width=0.43\textwidth]{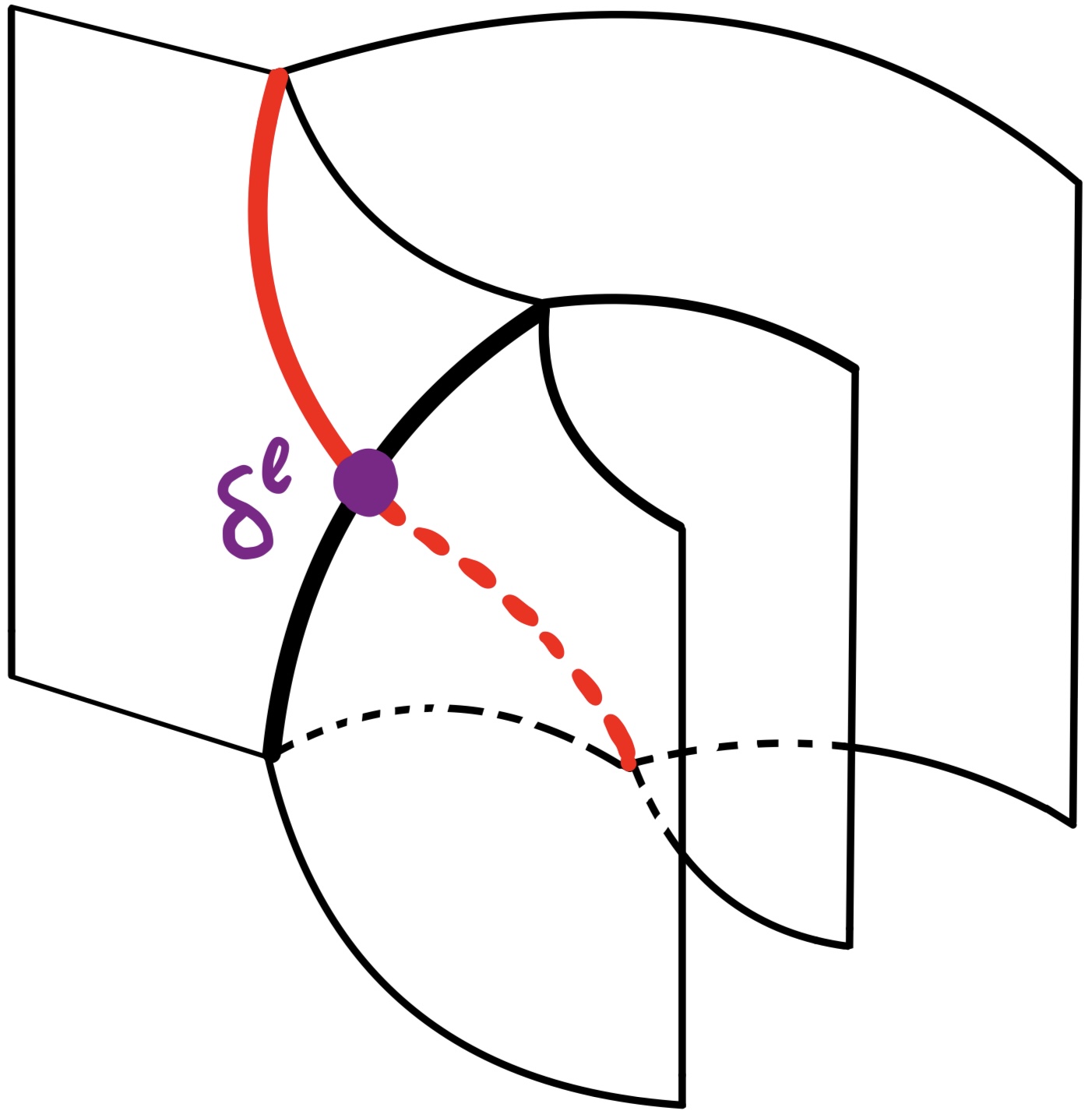}
         \caption{Left distributor $\distl$.}
    \end{subfigure}
    \begin{subfigure}[b]{0.45\textwidth}
         \centering
         \includegraphics[width=0.41\textwidth]{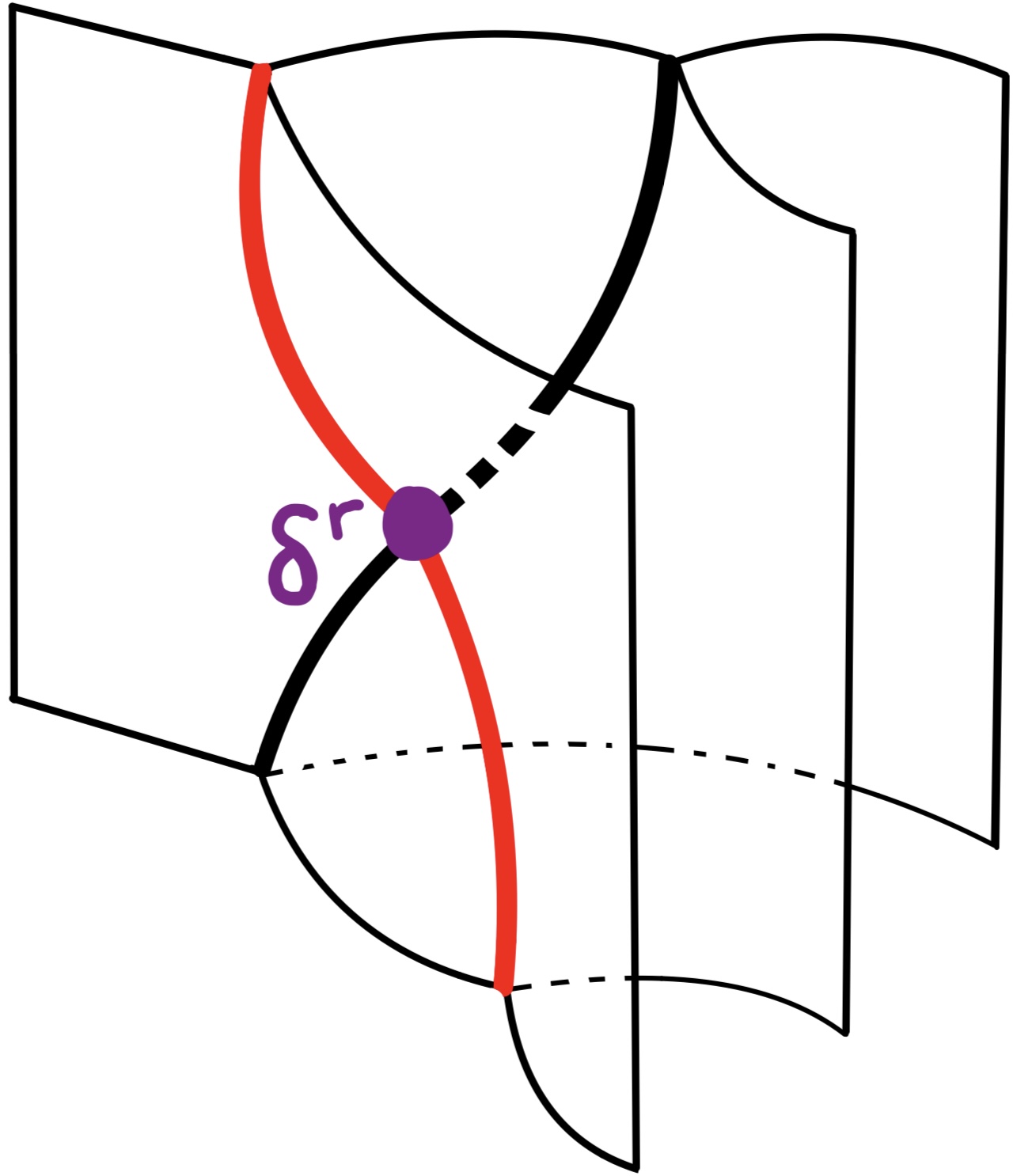}
         \caption{Right distributor $\distr$.}
    \end{subfigure}
    \caption{Surface diagrams for the distributors.}\label{fig: Distributors}
\end{figure}

The surface diagrams in Figure \ref{fig: Distributors} are obtained from each other by reflecting along the side face. Three-dimensional STL models of the (monochromatic variants of the) surfaces in Figure \ref{fig: Distributors} can be viewed \href{https://maxdemirdilek.github.io/Research/SurfaceDiagrams}{here}. The website also provides a HOM file for the signature of LD-categories.
\end{remark}

As was first observed by Cockett and Seely \cite[Prop. 5.4]{WDC}, the following type of LD-category is the only example where the distributors are invertible: 

\begin{example}\label{degnerateLDC}
    Given a monoidal category $(\cC,\otimes,1)$, choose an $\otimes$-invertible object $S$ in $ \cC$ with inverse $S^{-1}$. Then the bifunctor $(-\otimes S)\otimes -$, together with appropriate coherence data, defines a second monoidal structure on $\cC$ with monoidal unit $S^{-1}$. It endows $\cC$ with a linearly distributive structure \cite[\S 5.2]{WDC}, making it a \emph{shift-monoidal category}. As a special case, taking the monoidal unit $S=1$, one can view any monoidal category as an LD-category with coinciding monoidal products.
\end{example}

There exist various notions of functors-with-structure between LD-categories. The following functors-with-structure constitute a subclass of the \emph{linear functors} defined in \cite{LD-functors}.
\begin{definition}\label{def:FrobLin}(\cite[Def. 3.1.]{DaggerLinLogic}).
Let $\cC=(\cC,\otimes,1, \parLL,K)$ and $\cD=(\cD,\otimes,1, \parLL,K)$ be LD-categories. A \emph{Frobenius linearly distributive functor} from $\cC$ to $\cD$, or \emph{Frobenius LD-functor} for short, consists of a functor $F\colon\cC\ra\cD$ between the underlying categories together with the following additional data:
    \begin{itemize}
        \item A lax $\otimes$-monoidal structure $(\varphi^2,\varphi^0)$.
        \item An oplax $\parLL$-monoidal structure $(\upsilon^2,\upsilon^0)$. 
    \end{itemize}
 
\nid These structures are required to satisfy the following \emph{Frobenius relations}, for all $X,Y,Z\in \cC$:
    \begin{align}
     \tag{F1} 
        \upsilon^2_{X\otimes Y,Z}\circ F(\distl_{X,Y,Z})\circ \varphi^2_{X,Y\parLL Z}\:=\:\big(\varphi^2_{X,Y}\parLL F(Z)\big)\circ \distl_{F(X),F(Y),F(Z)}\circ \big(F(X)\otimes \upsilon^2_{Y,Z}\big),\label{eq:F1 Frob LD}\\[1.2ex]
    \tag{F2}
        \upsilon^2_{X,Y\otimes Z}\circ F(\distr_{X,Y,Z})\circ \varphi^2_{X \parLL Y,Z}\:=\:\big(F(X)\parLL{\varphi^2_{Y,Z}}\big)\circ \distr_{F(X),F(Y),F(Z)}\circ \big(\upsilon^2_{X,Y}\otimes F(Z)\big).\label{eq:F2 Frob LD}
    \end{align}

\nid The surface diagrams for the Frobenius relations look as follows:

    \begin{figure}[H]
        \centering
        \begin{subfigure}[b]{0.46\textwidth}
            \centering
            \includegraphics[width=0.99\textwidth]{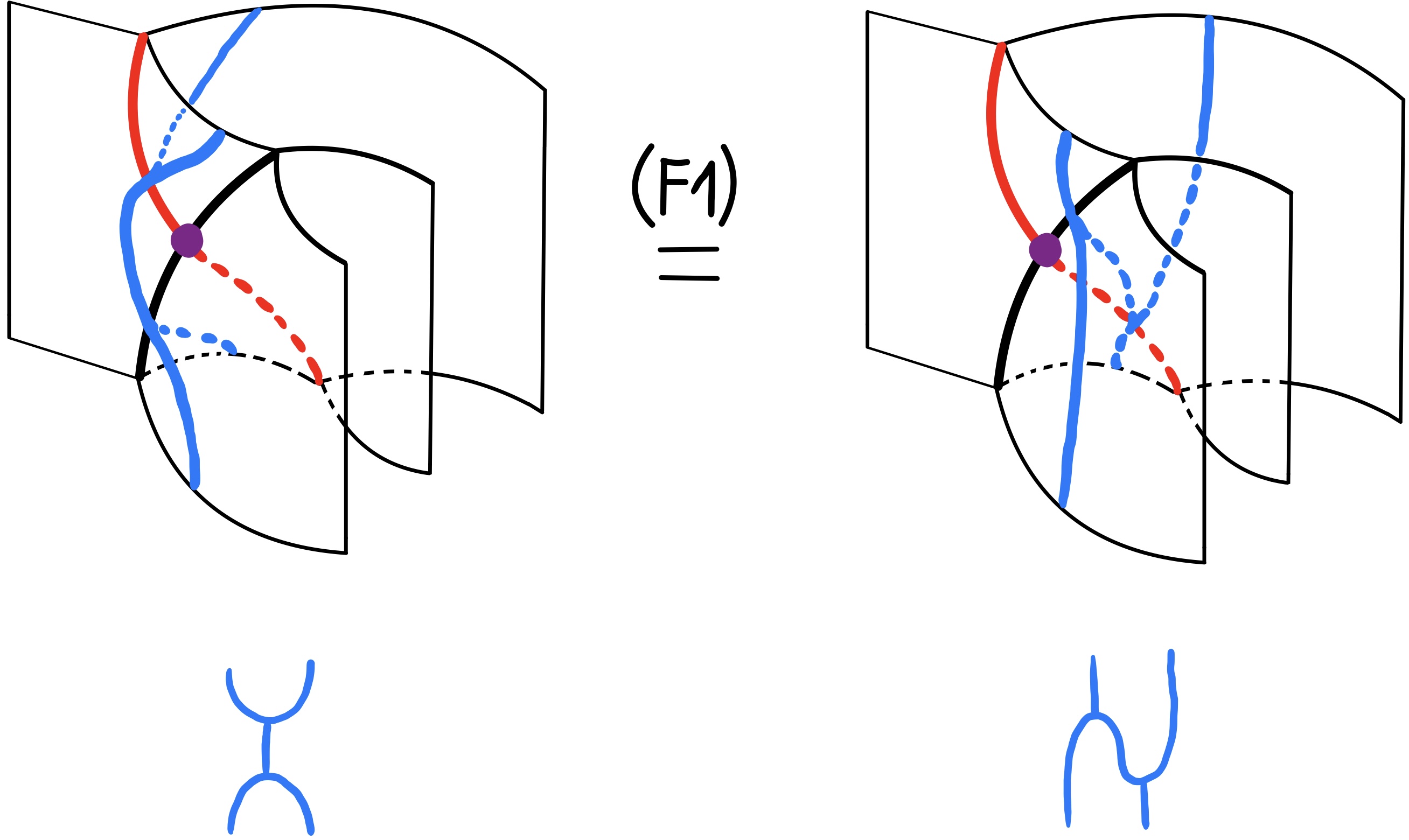}
            \caption{Relation \eqref{eq:F1 Frob LD}.}\label{rel F1}
        \end{subfigure}
     \hfill
        \begin{subfigure}[b]{0.46\textwidth}
         \centering
         \includegraphics[width=0.99\textwidth]{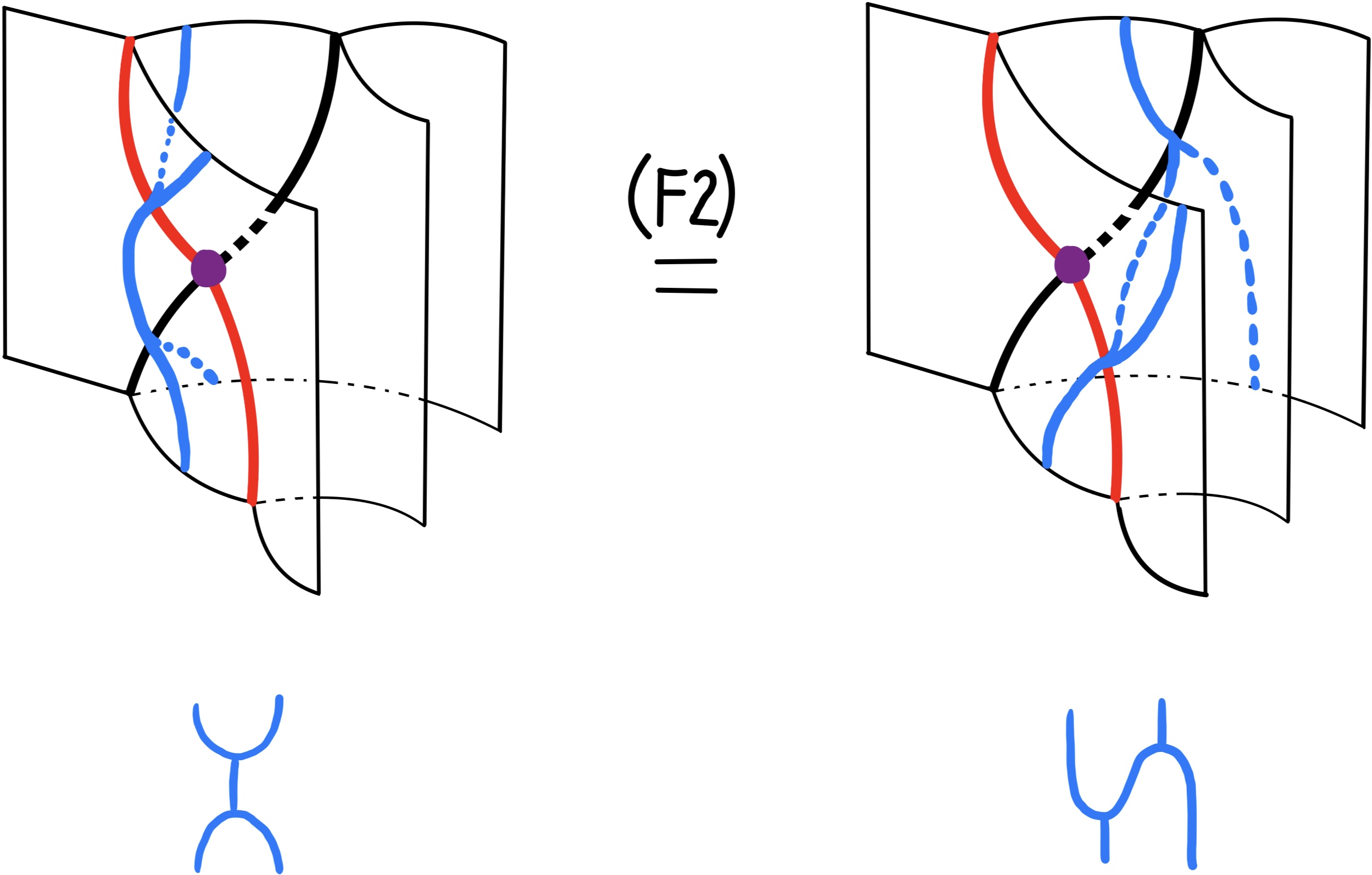}
         \caption{Relation \eqref{eq:F2 Frob LD}.}\label{rel F2}
        \end{subfigure}
     \caption{The Frobenius relations. The extracted string diagrams are shown below.}\label{fig: Relations FrobLD-fun}
    \end{figure}
\nid Here, the line showing the functor $F$ is depicted in blue. Below the surface diagrams, we have drawn the blue lines of the functor $F$ when viewed from the front face.
\end{definition}

\begin{remark}
The surface diagrams for the above Frobenius relations are obtained by drawing onto the distributor surface diagrams from Figure \ref{fig: Distributors} the string diagrams encoding the monoidal Frobenius relations (as depicted in blue at the bottom of Figure \ref{fig: Relations FrobLD-fun}).
Note that the blue functor lines for the relations \eqref{eq:F1 Frob LD} and \eqref{eq:F2 Frob LD} are drawn onto different surfaces. In fact, their surface diagrams are obtained from each other by reflection along the side face.
\end{remark}

\begin{examples}\label{ex:restriction of scalars}
    \begin{enumerate}[label=(\roman*)]
        \item  Frobenius LD-functors between monoidal categories are often called \emph{Frobenius monoidal functors}; see e.g. \cite{NoteOnFrobMonFunctors, Runkel-Kong, McSt, FrobeniusMonDijkgraaf, fuchs2023stringnetmodelspivotalbicategories, Yadav}.
        \item Restriction along a morphism $\varphi\colon A \rightarrow B$ of finite-dimensional $k$-algebras defines a functor
        \begin{equation*}
        \operatorname{Res}_{\varphi}\colon \, B\operatorname{-bimod}\, \longrightarrow \, A\operatorname{-bimod}.
        \end{equation*}
         It is a Frobenius-LD functor between the LD-categories from Example \ref{bimodule example}. Let ${M,N\in B\operatorname{-bimod}}$. The coequalizer surjection ${M\otimes_k N\twoheadrightarrow M\otimes_B N}$ induces an $A$-bimodule morphism $ {\operatorname{Res}_{\varphi}(M)\otimes_A \operatorname{Res}_{\varphi}(N) \rightarrow \operatorname{Res}_{\varphi}(M\otimes_B N)}.$
        Together with ${\varphi\colon A \rightarrow \operatorname{Res}_{\varphi}(B)}$, this defines a lax $\otimes_A$-monoidal structure on $\operatorname{Res}_{\varphi}$. Dually, the equalizer inclusion ${M\otimes^B N\hookrightarrow M\otimes_k N}$ yields an $A$-bimodule morphism 
        ${\operatorname{Res}_{\varphi}(M\otimes^B N)\rightarrow \operatorname{Res}_{\varphi}(M)\otimes^A \operatorname{Res}_{\varphi}(N)}.$ Together with the $k$-dual ${\operatorname{Res}_{\varphi}(D(B))\rightarrow DA}$ of $\varphi$, this defines an oplax $\otimes^A$-monoidal structure on $\operatorname{Res}_{\varphi}$. The Frobenius relations follow from the definition of distributors in \cite[\S 6]{fuchs2024grothendieckverdierdualitycategoriesbimodules}.
        \item Other examples of Frobenius LD-functors between genuine LD-categories (i.e., LD-categories where the two monoidal categories do not necessarily coincide) are discussed in this paper; see Proposition \ref{duality functor is FrobLD} or Examples \ref{example ld-frobenius algebras}, for instance.
    \end{enumerate}   
\end{examples}

Day and Pastro \cite[Prop. 3]{NoteOnFrobMonFunctors} prove the following proposition using commutative diagrams. We provide a surface-diagrammatic proof to get accustomed to surface diagrams:

\begin{prop}\label{strong mon implies frob mon}
Let $\cC$ and $\cD$ be monoidal categories. Any strong monoidal functor ${F\colon\cC\ra \cD}$ is a Frobenius monoidal functor.    
\end{prop}

\begin{proof}
We show that the strong monoidal functor $F\colon \cC\ra \cD$ satisfies Frobenius relation \eqref{eq:F1 Frob LD}: 
    \begin{figure}[H]
    \centering
\includegraphics[width=0.96\textwidth]{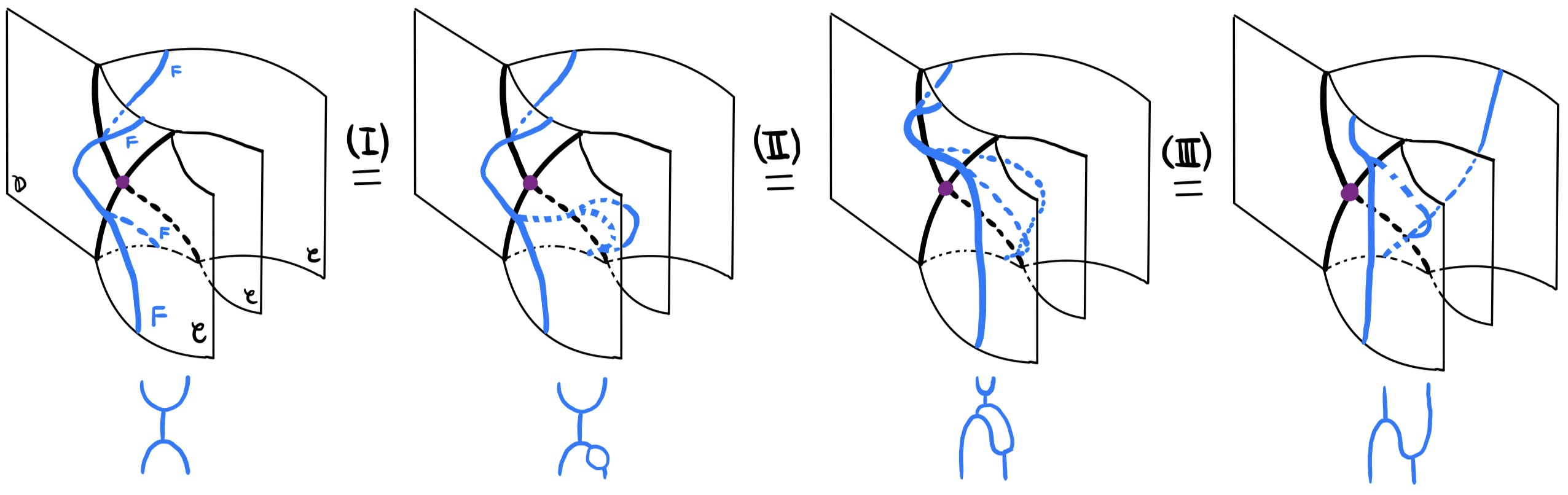}.
\end{figure}
%\vspace{-0.3cm}
Equations (I) and (III) hold since the monoidal structure on $F$ is strong, while equation (II) follows from the associativity axiom of the (lax) monoidal structure.  
The proof of the Frobenius relation \eqref{eq:F2 Frob LD} follows by reflecting each surface diagram along the side face.
\end{proof} 

The converse of Proposition \ref{strong mon implies frob mon} does not hold; e.g. \cite[Prop. 2.23]{Runkel-Kong}.

\begin{definition}\label{def:strong Frobenius LD-functor}
    A Frobenius LD-functor is called \emph{strong} (resp. \emph{strict}) if the lax $\otimes$-monoidal structure and the oplax $\parLL$-monoidal structure are strong (resp. strict). A \emph{Frobenius LD-equivalence} is a strong Frobenius LD-functor whose underlying functor is an equivalence of categories.
\end{definition}

Using the previous definition, we can now make precise the claim that LD-categories with invertible distributors are shift-monoidal categories in the sense of Example \ref{degnerateLDC}:

\begin{remark}\label{strictifiability of LD-cats}
Any LD-category with invertible distributors is Frobenius LD-equivalent to a shift-monoidal category; see \cite{MyBScThesis} for a detailed proof.
\end{remark}

Shift-monoidal categories do not suffice in practice:

\begin{example}
Recall from Example \ref{bimodule example} the LD-category $A\operatorname{-bimod}$ of bimodules over a finite-dimensional $k$-algebra $A$. Its right exact tensor product $\otimes_A$ is generally, e.g. for $A=k[X]/ \langle X^2\rangle$ the $k$-algebra of dual numbers, not left exact. If the right exact tensor product $\otimes_A$ is not left exact, a tensor product obtained by shifting the $\otimes_A$-monoidal structure by an $\otimes_A$-invertible object is not left exact either. However, the $\parLL$-tensor product $\otimes^A$ on $A\operatorname{-bimod}$ is left-exact. Thus, by Remark \ref{strictifiability of LD-cats}, the distributors of the LD-category $A\operatorname{-bimod}$ are generally non-invertible; see \cite[Lemma 5.5.]{fuchs2024grothendieckverdierdualitycategoriesbimodules} for more information.
\end{example}

In fact, LD-categories do not even satisfy a full coherence theorem, i.e. there exist (small) LD-categories in which not every formal diagram\footnote{To define formal diagrams, consider the forgetful functor from the category of (small) LD-categories and strict Frobenius LD-functors to the category of sets, which sends each category to its underlying set of objects. This functor admits a left adjoint, called \emph{free} LD-category functor. Adapting the terminology from \cite[\S 1]{MalPo}, we call a diagram in a (small) LD-category $\cC$ \emph{formal} if it lifts to the free LD-category on the underlying set of objects of $\cC$.} commutes:
\begin{example}\label{ex:LD not fully coherent}
Let $\cC$ be an LD-category. The following endomorphism in $\cC$
\begin{equation}\label{tripledualmorph}
{(({1} \parLL {1})\otimes K)}\parLL {1}\xrightarrow{\distr\parLL 1} {({1} \parLL {(1\otimes K)})}\parLL {1} \cong ({1} \parLL {1}) \otimes ({K}\parLL {1})\xrightarrow{\distl} {(({1} \parLL {1})\otimes K)}\parLL {1},
\end{equation}
where the isomorphism in the middle is obtained by applying multiple unitors, is in general not equal to the identity morphism on ${(({1} \parLL {1})\otimes K)}\parLL {1}$.\footnote{We thank Aaron David Fairbanks \cite{DaFa} for pointing out this example.} For instance, in the LD-category $A\operatorname{-bimod}$ of bimodules over the three-dimensional $k$-algebra $A:=k[x,y]/\langle x^2,y^2,xy\rangle$, the first map $\distr_{A,A,DA}\otimes^A A$ (and thus the map \eqref{tripledualmorph}) is non-injective; see \cite[\S 6]{fuchs2024grothendieckverdierdualitycategoriesbimodules}.
\end{example}
 
Having defined Frobenius LD-functors between LD-categories, we propose the following definition for morphisms between them:

\begin{definition}\label{def:MorphFrob} A \emph{morphism of Frobenius LD-functors} $F,G\colon \cC \ra \cD$ is a natural transformation $f\colon F\ra G$ that has the property of being $\otimes$-monoidal and $\parLL$-opmonoidal.
\end{definition}

\begin{remark}
    Up to recoloring, the surface diagrams for a $\parLL$-opmonoidal natural transformation can be found in \cite[\S 1.2.3]{WillertonHopfM}. We obtain the surface diagrams for $\otimes$-monoidal natural transformations from these by recoloring the functor lines and reflecting along the bottom face.
\end{remark}

With surface diagrams at our disposal, it is immediate that the composite of Frobenius LD-functors carries a natural Frobenius LD-structure. The graphical proof consists of drawing a second functor line parallel to the blue functor lines in Figure \ref{fig: Relations FrobLD-fun} and then using the Frobenius relations for each of these functor lines separately. It is a similarly simple graphical verification to show that morphisms of Frobenius LD-functors are compatible with both horizontal and vertical composition. We thus find:

\begin{prop}\label{LDC2Cat}
    LD-categories, together with Frobenius LD-functors, and their morphisms form a $2$-category $\mathsf{Frob}$.
\end{prop}

By the following proposition, equivalences in the $2$-category $\mathsf{Frob}$ are precisely Frobenius LD-equivalences. A proof of this proposition is outlined in Appendix \ref{sec:appendixproof}.

\begin{prop}\label{quasi-inverse gets Frobenius LD-structure}
Let $F\colon \cC \ra \cD$ be a Frobenius LD-equivalence between LD-categories. The quasi-inverse $G$ of the underlying functor of $F$ carries a unique structure of a Frobenius LD-functor such that the unit $\eta\colon \operatorname{id}_{\cC}\xrightarrow{\simeq}GF$ and the counit $\epsilon\colon FG\xrightarrow{\simeq}\operatorname{id}_{\cD}$ are isomorphisms of Frobenius LD-functors. This Frobenius LD-structure on $G$ is strong.
\end{prop}

Given a monoidal category $(\cC,\otimes,1)$, we define the reversed monoidal product \begin{equation*}
X \otimes^{\text{rev}}Y\,:=\,Y \otimes X,\end{equation*} for $X,Y\in \cC$. With this notion in hand, we introduce for future use:

\begin{definition}\label{three2-functors}
Let $\cC=(\cC,\otimes,1,\parLL,K)$ be an LD-category.
    \begin{enumerate}[label=(\roman*)]
        \item The category $\cC$ carries another linearly distributive structure
            \begin{equation*}
            \cC^{\text{rev}}\,:=\,(\cC,\otimes^{\text{rev}},1,\parLL^{\text{rev}},K).
            \end{equation*} 
        \item The opposite category $\cC^{\text{op}}$ carries a linearly distributive structure
            \begin{equation*}
            \cC^{\text{cop}}\,:=\,(\cC^{\text{op}}, \parLL,K,\otimes,1).
            \end{equation*}
        \item The opposite category $\cC^{\text{op}}$ also carries a second linearly distributive structure
            \begin{equation*}
            \cC^{\text{lop}}\,:=\,(\cC^{\text{op}},\parLL^{\text{rev}},K,\otimes^{\text{rev}},1).
            \end{equation*} 
        We refer to $\cC^{\operatorname{lop}}$ as the \emph{linearly distributive opposite category} of the LD-category $\cC$. 
    \end{enumerate}
\end{definition}

Each of the assignments $\cC\mapsto \cC^{\text{rev}}$, $\cC\mapsto \cC^{\text{cop}}$ and $\cC\mapsto \cC^{\text{lop}}$ extends to a $2$-endofunctor on the $2$-category $\mathsf{Frob}$. Each of these $2$-functors is a $2$-involution, meaning that it is an involution on objects, $1$-morphisms and $2$-morphisms. Graphically, these $2$-involutions are constructed from reflections of the cubical canvas and specific recolorings within the cube:

\begin{remark}
    \begin{enumerate}[label=(\roman*)]
        \item The $2$-functor $(-)^{\text{rev}}$ reflects all surface diagrams along the side face.
        \item The $2$-functor $(-)^{\text{cop}}$ reflects all surface diagrams along the bottom face and swaps the red and black coloring of the functor lines. Put differently, applying the $2$-endofunctor $(-)^{\text{cop}}$ corresponds to changing the $3$-framing of the canvas so that the triple arrow in Figure \ref{fig:canvas} points, after an application of $(-)^{\text{cop}}$, in the negative Z-direction. 
        \item The $2$-functor $(-)^{\text{lop}}$ rotates all occurring surface diagrams by $180^\circ$ and swaps the red and black coloring of the functor lines. Of course, the rotation by $180^\circ$ can be expressed as a composite of two reflections.
    \end{enumerate}
\end{remark}

\subsection{Dual objects in LD-categories}\label{duality in LD-categories}
Hereafter, we fix an LD-category $\cC$. An LD-category comes with a natural notion of duality. We present this notion in Definition \ref{def:LDduals}. In doing so, we will employ the following concepts:

\begin{definition}\label{pairing/copairing,side-inverse}
    Let $X,Y$ be objects in $\cC$. 
    \begin{enumerate}[label=(\roman*)]
        \item An \emph{LD-pairing} of $X$ and $Y$ is a morphism $\kappa\in \operatorname{Hom}_{\cC}(X\otimes Y,K)$; an \emph{LD-copairing} of $X$ and $Y$ is a morphism $\overline{\kappa}\in \operatorname{Hom}_{\cC}(1, Y\parLL X)$.
        \item An LD-pairing $\kappa\colon X\otimes Y\ra K$ and an LD-copairing $\overline{\kappa}\colon 1\ra Y\parLL X$ are called \emph{side-inverse} to each other if they satisfy the following two \emph{snake equations} involving the distributors:
    \begin{align}\label{firstzigzag}
     \tag{S1}({\kappa} \parLL {X})\circ\distl_{X,Y,X}\circ(X\otimes {\overline{\kappa}})\;=\;{(l^{\parLL}_{X})^{-1}}\circ {r_X^{\otimes}},\\[1.5ex] \label{secondzigzag}\tag{S2}
    ({Y} \parLL {\kappa})\circ\distr_{Y,X,Y}\circ({\overline{\kappa}}\otimes {Y})\;=\;{(r^{\parLL}_Y)^{-1}}\circ {l^{\otimes}_Y}.
    \end{align}
    \nid Here, the symbols $l^{\otimes}$, $r^{\otimes}$, $l^{\parLL}$ and $r^{\parLL}$ denote the left and right unitors for the $\otimes$-monoidal and $\parLL$-monoidal structure.
    \end{enumerate}
\end{definition}

The snake equations are depicted as follows:

\begin{figure}[H]
    \centering
    \begin{subfigure}[b]{0.4\textwidth}
         \centering
         \includegraphics[width=\textwidth]{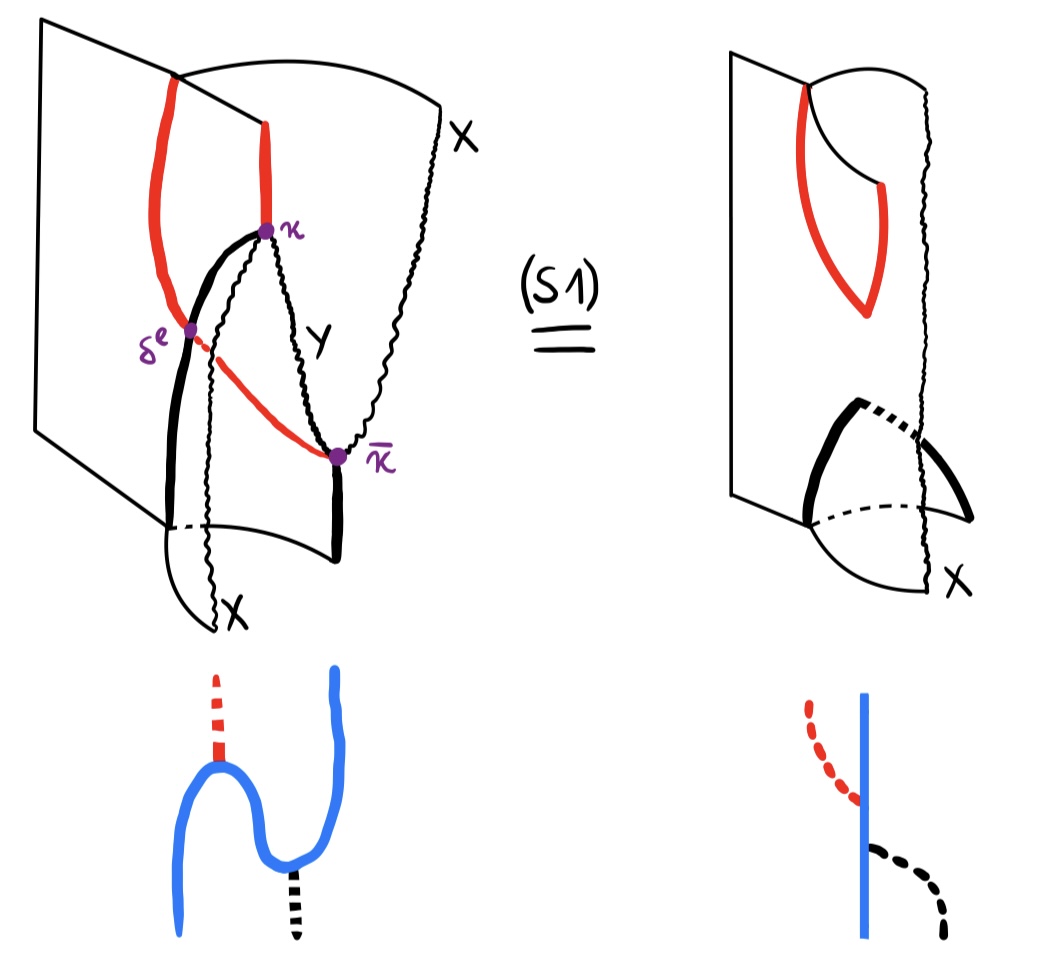}
    \end{subfigure}
     \hspace{4em}
     \begin{subfigure}[b]{0.4\textwidth}
         \centering
         \includegraphics[width=\textwidth]{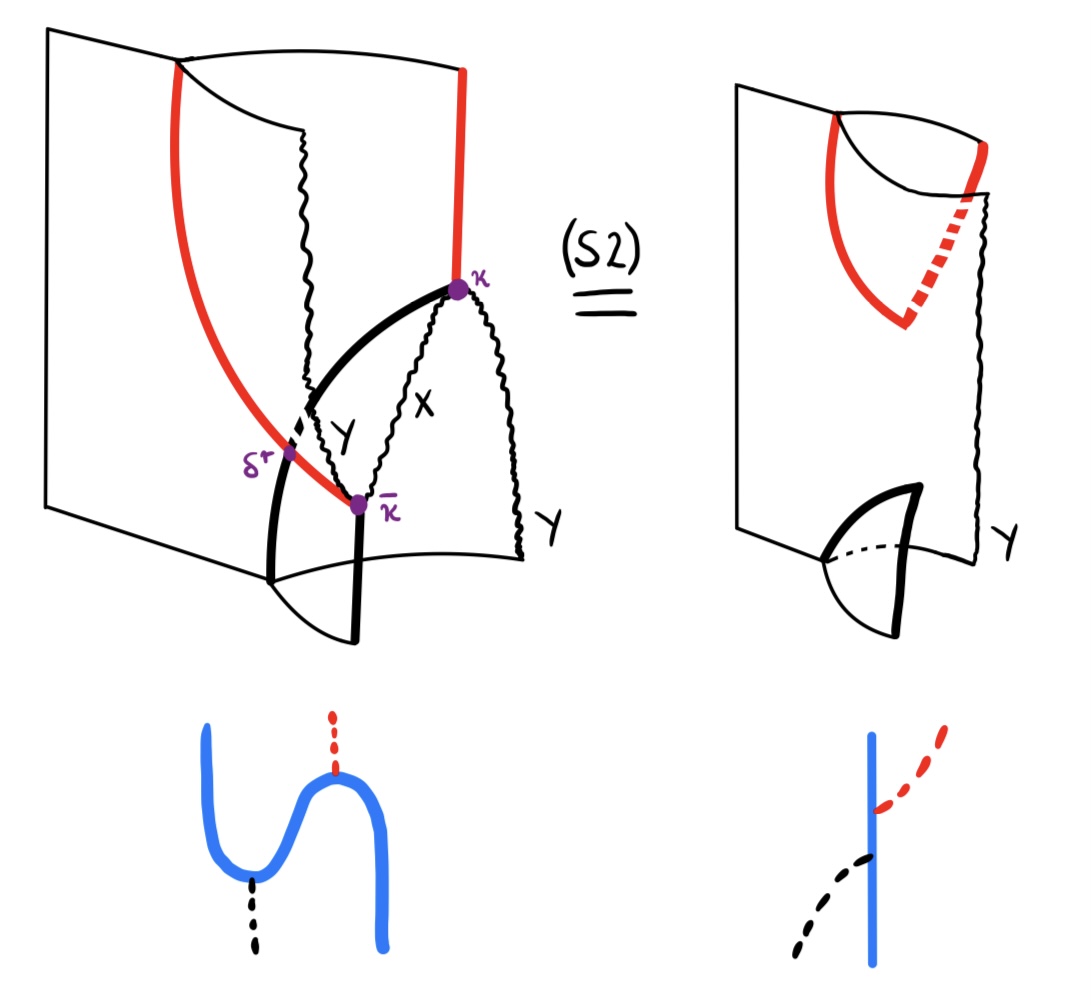}
     \end{subfigure}
\caption{The snake equations. Dotted lines in string diagrams represent monoidal units.}\label{fig:LD-snake equations (S1) and (S2)}
\end{figure}

A HOM file for a side-inverse LD-pairing and LD-copairing is available \href{https://maxdemirdilek.github.io/Research/SurfaceDiagrams}{here}. The website also supplies an STL file for the surface diagram obtained from the leftmost surface diagram in the snake equation (\ref{secondzigzag}) by pre- and post-composing with appropriate unitors.

\begin{remark}
    The source of an LD-pairing and an LD-copairing uses the $\otimes$-monoidal structure, while the target employs the $\parLL$-monoidal structure.
\end{remark}

\begin{example}\label{unitors are side inverse}
    Using the triangle coherence axioms (\ref{eq:A1}) to (\ref{eq:A4}), one finds that the unitors ${l^{\otimes}_K\colon\: 1 \otimes K\ra K}$ and $(l^{\parLL}_1)^{-1}\colon\: 1\ra K\parLL 1$ (similarly, ${r^{\otimes}_K\colon\: K \otimes 1\ra K}$ and ${(r^{\parLL}_1)^{-1}\colon\:1\ra 1\parLL K)}$ are side-inverse to each other.
\end{example}

In the literature, the notions from Definition \ref{pairing/copairing,side-inverse} also appear under a different name:

\begin{remark}\label{nuclearity}
In \cite[Def. A.1.]{LD-functors} two objects $X,Y\in \cC$ together with an LD-pairing ${\kappa\colon X\otimes Y\ra K}$ and an LD-copairing $\overline{\kappa}\colon 1\ra Y\parLL X$ satisfying the snake equation (\ref{firstzigzag}) are called \emph{witness of the right nuclearity} of $\operatorname{id}_X$. Similarly, Cockett and Seely define the notion of \emph{witness of the left nuclearity} by employing the snake equation (\ref{secondzigzag}).
\end{remark}

When viewing monoidal categories as LD-categories with coinciding monoidal products, the following definition specializes to the notion of rigid dual in a monoidal category.

\begin{definition}\label{def:LDduals}(\cite[Def. 4.1.]{WDC}). Let $X$ be an object in $\cC$. An object $X^{\ast}\in \cC$ is said to be a \emph{left LD-dual} of $X$ if there exist an LD-pairing $\operatorname{ev}_X\colon X^{\ast}\otimes X \ra K$ and an LD-copairing $\operatorname{coev}_X\colon 1\ra X\parLL X^{\ast}$, called \emph{evaluation} and \emph{coevaluation}, that are side-inverse to each other. An object is called \emph{left LD-dualizable} if it possesses a left LD-dual. \emph{Right LD-duals} are left duals in $\cC^{\text{rev}}$ and are denoted by $^{\ast}X$, with $\operatorname{ev}'_X\colon X\otimes {^{\ast}X} \ra K$ and $\operatorname{coev}'_X\colon 1\ra {^{\ast}X}\parLL X$.
\end{definition}

\begin{remark}\label{LD-duals are unique} 
    LD-duals are unique up to unique isomorphism \cite[Lemma A.6.]{LD-functors}. This fact can also be verified using surface diagrams. Being LD-dualizable is thus a property of an object.
\end{remark}

\begin{remark}
    Upon inspecting the surface diagrams for the snake equations of Definition \ref{def:LDduals} from the front face, and suppressing the monoidal units, we obtain the familiar string-diagrammatic representations of the snake equations for rigid duals in a monoidal category.
\end{remark}

\begin{example}\label{ex:fd bimodule is LD-dualizable}
Let $M$ be a finite-dimensional bimodule over a finite-dimensional $k$-algebra $A$. In the LD-category $A\operatorname{-bimod}$, the $k$-linear dual $A$-bimodule $DM$ (as defined in Example \ref{bimodule example}) serves as both the left and right LD-dual of $M$. Let $\{m_i\}_i$ be a basis of $M$ with dual basis $\{m^i\}_i$ of $DM$. Define the left evaluation map $\operatorname{ev}_M\colon DM \otimes_A M \to DA$ on representatives of elementary tensors by $\operatorname{ev}_M([f \otimes x])(a) := f(x.a)$ for $f \in DM$, $x \in M$, $a \in A$; this is a well-defined $A$-bimodule morphism. The left coevaluation map $\operatorname{coev}_M\colon A \to M \otimes^A DM$ is given by $\operatorname{coev}_M(a) := \sum_i a.m_i \otimes m^i$. These satisfy the snake equations \eqref{firstzigzag} and \eqref{secondzigzag}. The right evaluation and coevaluation maps are defined analogously.
\end{example}

Let us investigate how Frobenius LD-functors and their morphisms interact with LD-duals. When specializing to Frobenius monoidal functors, the following results are well known \cite[Thm. 2, Prop. 7]{NoteOnFrobMonFunctors}; they carry over to Frobenius LD-functors nearly verbatim. Our proofs are obtained using the three-dimensional graphical calculus.

\begin{prop} \label{FrobPresDuals}
Frobenius LD-functors between LD-categories preserve LD-duals.    
\end{prop}

\begin{proof}
 Let $F\colon \cC\ra \cD$ be a Frobenius LD-functor. Let $X^{\ast}\in \cC$ be a left LD-dual of an object $X\in \cC$ with evaluation $\operatorname{ev}_X$ and coevaluation $\operatorname{coev}_X$. Then the morphisms $F(X^{\ast})\otimes F(X)\ra K$ and $1\ra F(X)\parLL F(X^{\ast})$, given respectively by
\vspace{-0.2cm}
\begin{figure}[H]
     \centering
     \begin{subfigure}[b]{0.14\textwidth}
         \centering
    \includegraphics[width=\textwidth]{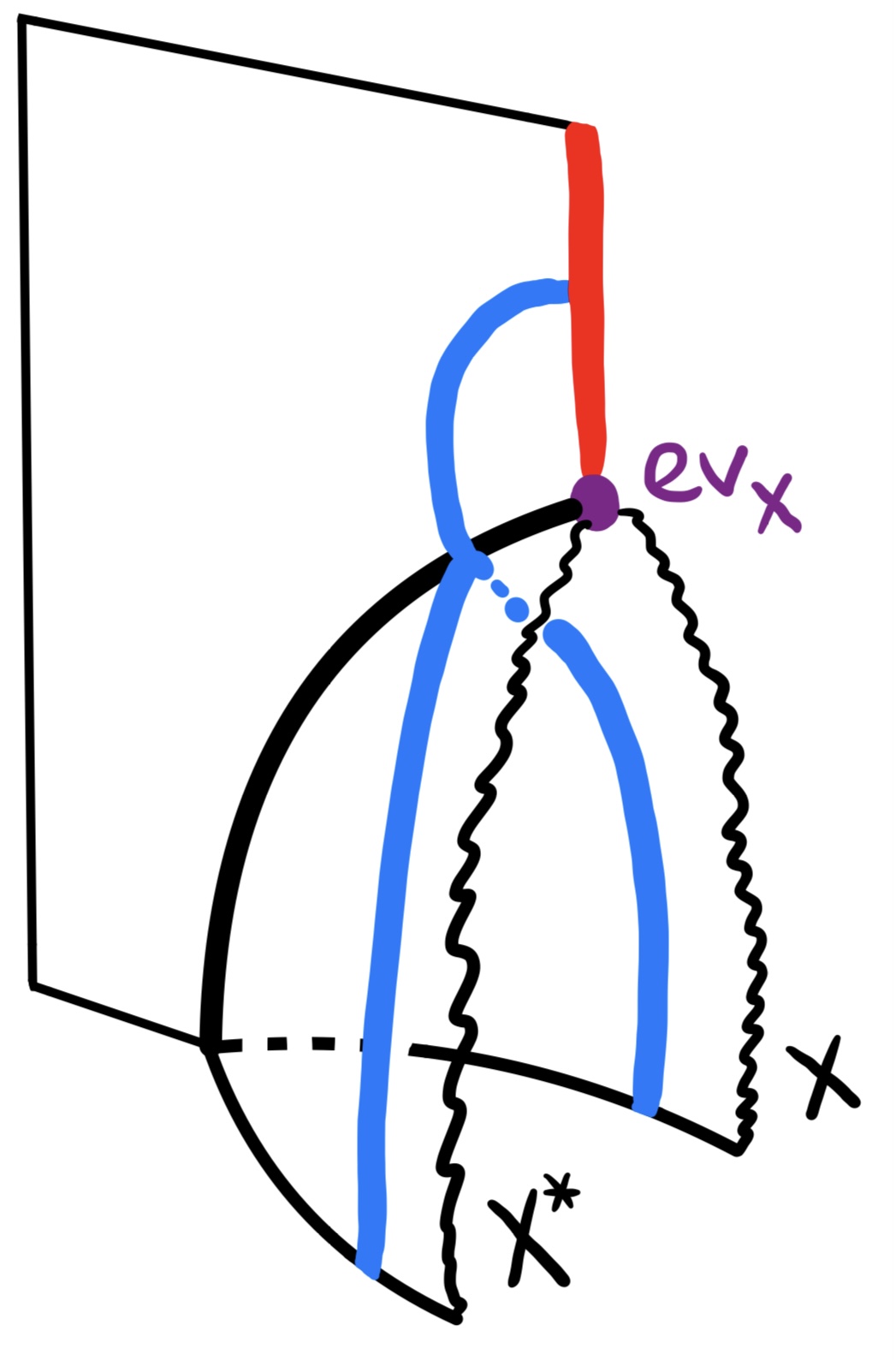}
 \end{subfigure}
     \hspace{4em}
     \begin{subfigure}[b]{0.14\textwidth}
         \centering
         \includegraphics[width=\textwidth]{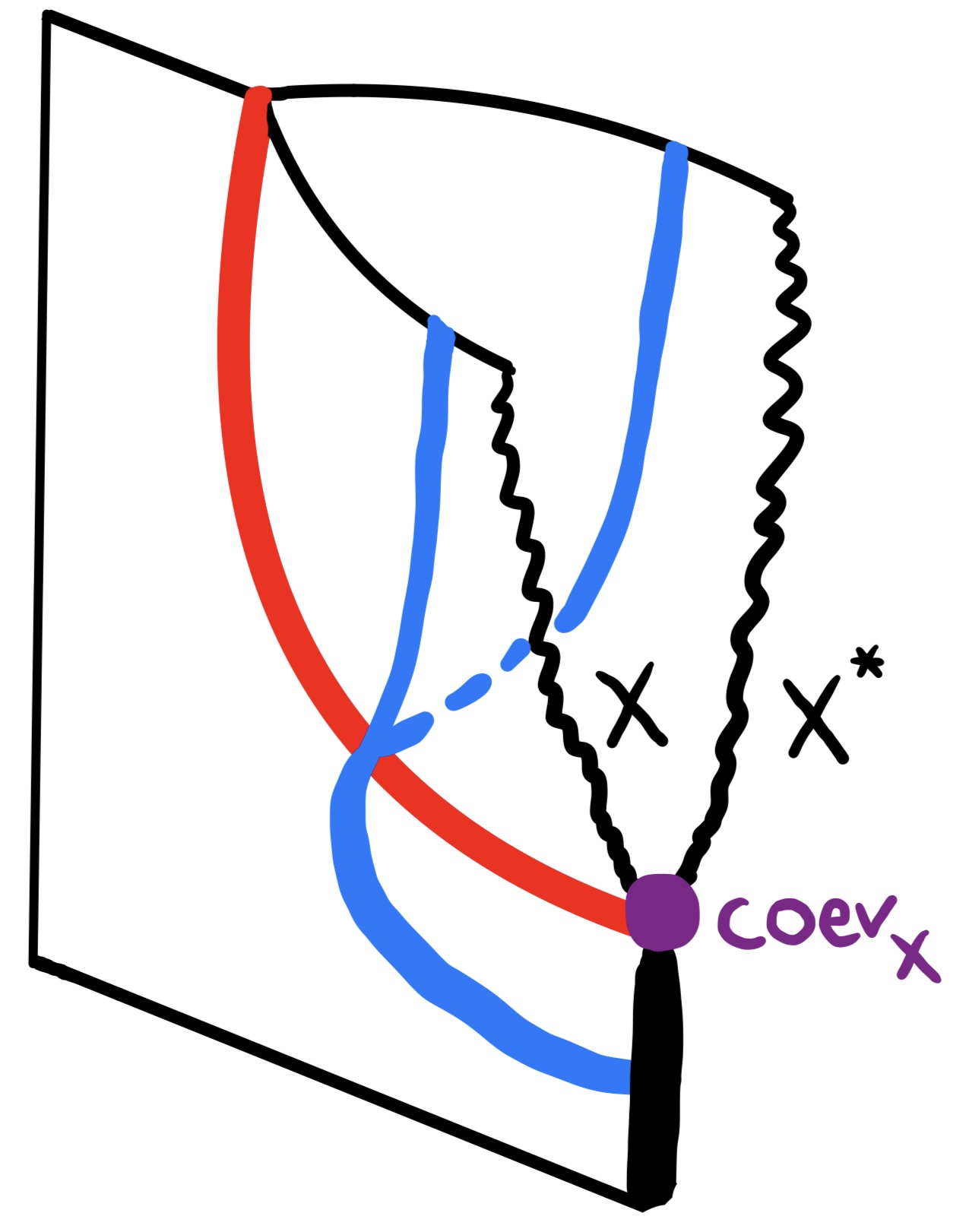}
     \end{subfigure},
\end{figure} 
\vspace{-0.2cm}
\nid witness $F(X^{\ast})$ as a left dual of $F(X)$. To see this, we verify snake equation (\ref{firstzigzag}) as follows:
\vspace{-0.2cm}
\begin{figure}[H]
    \centering
    \includegraphics[width=0.78\textwidth]{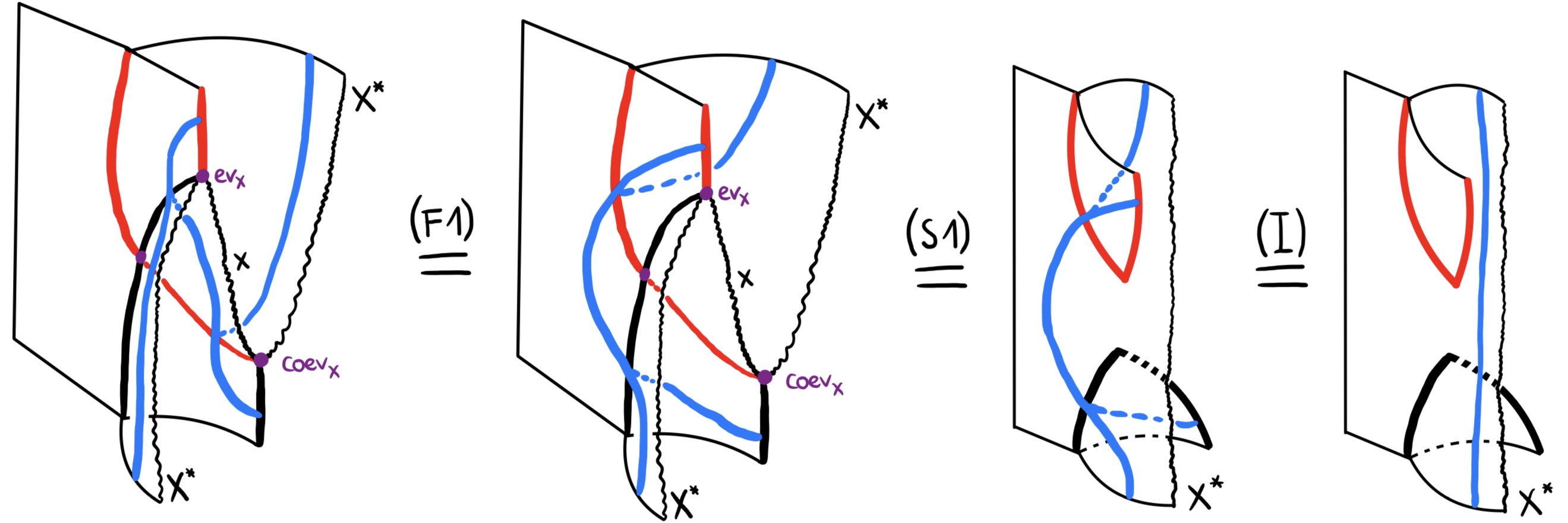}.
\end{figure}
\vspace{-0.2cm}
The first step applies the Frobenius relation \eqref{eq:F1 Frob LD} in Figure \ref{fig: Relations FrobLD-fun} to push the blue functor line of $F$ to the back. In the second step, we then use the snake equation (\ref{firstzigzag}) for $X\in \cC$. Equation (I) follows from the unitality (resp. counitality) axiom of the lax $\otimes$-monoidal structure (resp. oplax $\parll$-monoidal structure) on $F$.

\smallskip

The proof of snake equation (\ref{secondzigzag}) is obtained by applying the $2$-functor $(-)^{\operatorname{rev}}$ to the surface diagrams appearing in the surface diagrammatic proof of snake equation (\ref{firstzigzag}). Since right LD-duals in $\cC$ are left LD-duals in the LD-category $\cC^{\operatorname{rev}}$, this completes the proof.
\end{proof}

\begin{corollary}
Let $A$ be a finite-dimensional $k$-algebra. An $A$-bimodule is LD-dualizable in the LD-category $A\operatorname{-bimod}$ of Example \ref{bimodule example} if and only if it is finite-dimensional.
\end{corollary}
\begin{proof}
Consider the unique morphism of $k$-algebras $\varphi\colon k\rightarrow A$. View the monoidal category $\operatorname{vect}_k$ of $k$-vector spaces as an LD-category where the two monoidal structures coincide, and identify it with the LD-category $k\operatorname{-bimod}$ of $k$-bimodules. By Example \ref{ex:restriction of scalars}.(ii) and Proposition \ref{FrobPresDuals}, the forgetful functor $\operatorname{Res}_{\varphi}\colon A{\operatorname{-bimod}}\rightarrow \operatorname{vect}_k$ preserves LD-duals. Since a $k$-vector space is rigid dualizable (or LD-dualizable) if and only if it is finite-dimensional, any LD-dualizable $A$-bimodule must be finite-dimensional. For the converse, see Example \ref{ex:fd bimodule is LD-dualizable}.
\end{proof}

\begin{remark}
The rigid dualizable objects in the monoidal category of bimodules over a finite-dimensional algebra are the finite-dimensional bimodules that are projective as left and right modules, see \cite[Lemma 5.5]{fuchs2024grothendieckverdierdualitycategoriesbimodules}.
Yet, {\em all} finite-dimensional bimodules are
LD-dualizable.
\end{remark}

It is classical that morphisms between Frobenius algebras in monoidal categories are invertible. As we shall see in Section \ref{sec:LD-FrobAlg}, this statement is a corollary of the following result. This result concerns the interaction of LD-duals with morphisms of Frobenius LD-functors:

\begin{prop}\label{PropMorphLDfuncIso}
Let $\alpha\colon F\ra G$ be a morphism of Frobenius LD-functors $F,G\colon \cC\ra \cD$. 
\mbox{If $X\in \cC$} is left or right LD-dualizable, then the component $\alpha_X\colon\: F(X)\,\ra\, G(X)$ is invertible.
\end{prop}

\begin{proof} 
 Let $^{\ast}X\in \cC$ be a right LD-dual of $X\in \cC$. We claim that the following morphism $\beta_X \colon\: G(X)\,\ra\,F(X)$ is the inverse of $\alpha_X\colon\: F(X)\,\ra\,G(X)$:
 
    \begin{figure}[H]
        \centering
        \includegraphics[scale=0.088]{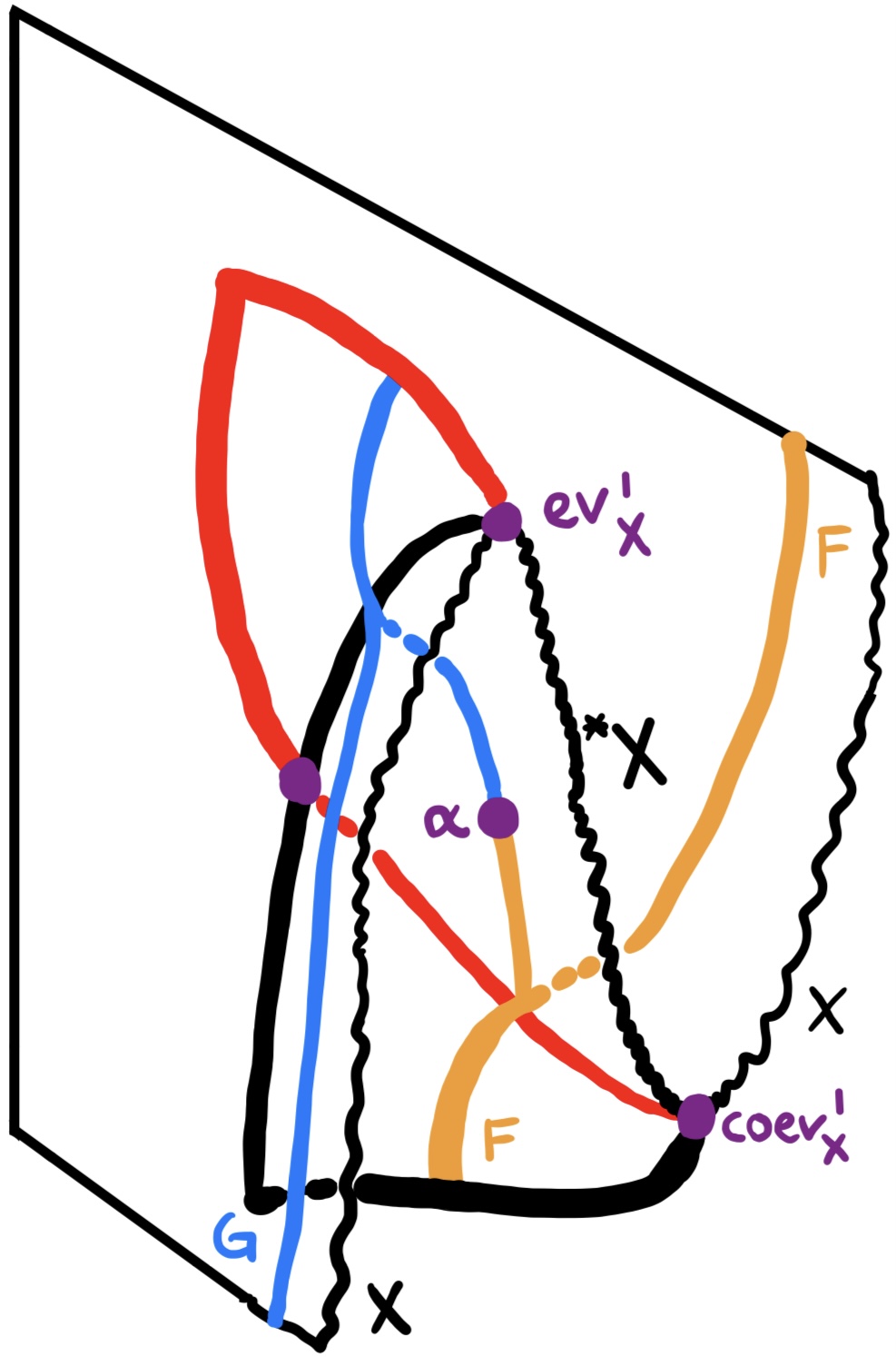}.
    \end{figure}
    
We show that equation $\alpha_X\circ\beta_X\,=\,\operatorname{id}_{G(X)}$ holds: 
    \begin{figure}[H]
        \centering
        \includegraphics[width=0.99\textwidth]{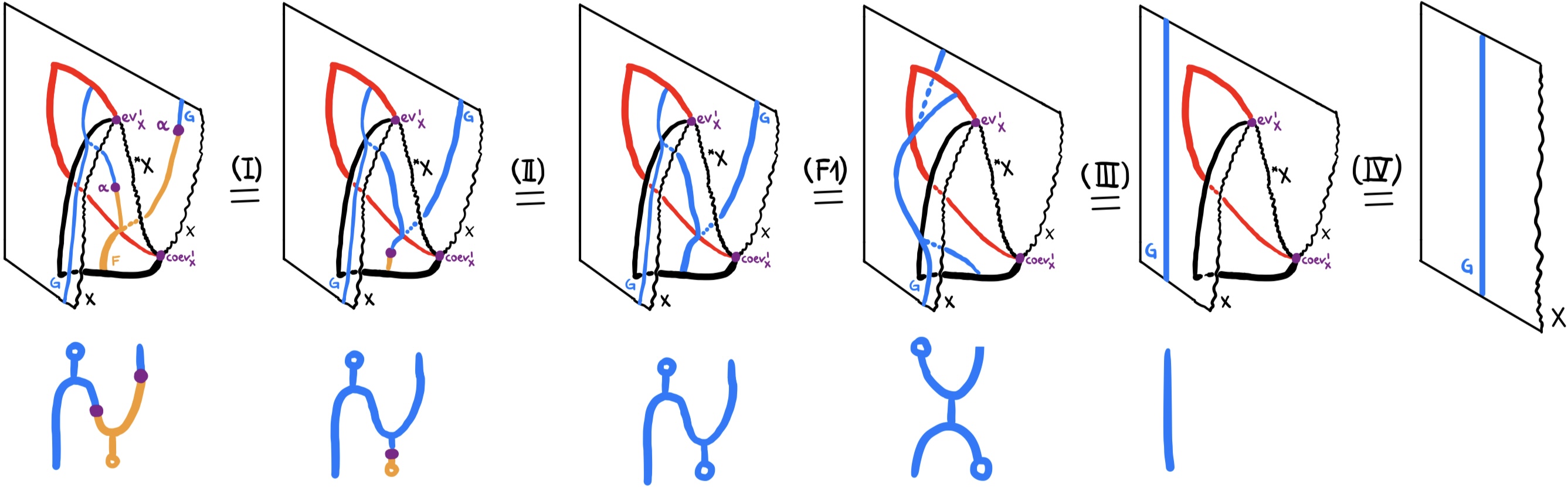}.
    \end{figure}
\vspace{-0.4cm}
Equation (I) is true, since $\alpha$ is, as an $\parLL$-opmonoidal natural transformation, compatible with the $\parLL$-monoidal product. Equation (II) holds since $\alpha$ is, for the same reason, compatible with the $\parLL$-monoidal unit. Equation (III) follows from the unitality (resp. counitality) axiom of the lax $\otimes$-monoidal structure (resp. oplax $\parll$-monoidal structure) on $G$. Equation (IV) holds by the snake equation (\ref{firstzigzag}) and the invertibility of the unitors.

\smallskip

Analogously, one proves $\beta_X\circ \alpha_X=\operatorname{id}_{F(X)}$. Finally, the statement for left LD-duals follows by applying the $2$-functor $(-)^{\operatorname{rev}}$ to the surface diagrams occurring in the proof above.
\end{proof}

\subsection[Grothendieck-Verdier categories (GV-categories)]{Grothendieck-Verdier categories}\label{sec:GV-cats}
Dualizability is a property of an object in a monoidal category. A monoidal category is called rigid if all of its objects are dualizable. The following definition naturally generalizes this notion of rigidity from monoidal categories to LD-categories. Our terminology follows that of Cockett and Seely \cite{WDC} and is motivated by the relationship between multiplicative linear logic and LD-categories; see \cite{WDC, BluteScott}.

\begin{definition}\label{def:LDnegation} A \emph{linearly distributive category with negation}, or \emph{LD-category with negation} for short,  is an LD-category with the property that every object admits a left and right LD-dual. 
\end{definition}

\begin{remark}
 By Proposition \ref{PropMorphLDfuncIso}, every $2$-morphism in the $2$-category of LD-categories with negation, Frobenius LD-functors, and morphisms between them is invertible.   
\end{remark}

As we will see in Theorem \ref{thm:GV correspond LD neg}, LD-categories with negation can also be characterized as follows:

\begin{definition}\label{def:GV-category}
Let $\cC=(\cC,\otimes,1)$ be a monoidal category.

\begin{itemize}
        \item A \emph{dualizing object} of $\cC$ is an object $K\in \cC$ such that for every object $Y\in \cC$, the functor $\operatorname{Hom}_{\cC}(-\otimes Y,K)$ is representable by some object $DY\in \cC$ and the induced contravariant functor $D$ on $\cC$ is an antiequivalence. The functor $D$ is called the \emph{duality functor} induced by the dualizing object $K$.
        \item A \emph{Grothendieck-Verdier structure} on $\cC$ is a choice of a dualizing object. A \emph{Grothendieck-Verdier category}, or \emph{GV-category} for short, is a monoidal category together with a Grothendieck-Verdier structure.
    \item If the monoidal unit $1$ is a dualizing object of $\cC$, we call $\cC$ an \emph{r-category}.
\end{itemize}
\end{definition}

\begin{remark}\label{star autonomous vs GV}
We follow the terminology of \cite{BoDrinfeld}; GV-categories are also known as (non-symmetric) $\ast$-autonomous categories, e.g. \cite{BarrStarAut, BarrNon}. Note, however, that the definition of a $\ast$-autonomous category in \cite[\S 4]{Pastro}, see also \cite[Def. 2.3]{HaLe}, does not completely coincide with our Definition \ref{def:GV-category}.\footnote{We thank Mike Shulman for correspondence on this point.}
\end{remark}

\begin{example}\label{GV-category of bimodules}
    Let $A$ be a finite-dimensional $k$-algebra. The finite-dimensional $A$-bimodule $DA$ from Example \ref{bimodule example} is a dualizing object in the monoidal category $(A\operatorname{-bimod}^{\operatorname{f.d.}},\otimes_A,A)$ of finite-dimensional $A$-bimodules; see \cite{fuchs2024grothendieckverdierdualitycategoriesbimodules}.
\end{example}

For later reference, we recall how a dualizing object $K$ in a monoidal category $(\cC,\otimes,1)$ induces a duality functor:

\begin{remark}\label{rep determines functoriality}
For every object $Y\in \cC$, choose a representing object $DY\in \cC$ and an isomorphism
    \begin{equation}\label{representing transformation}
        \Phi^{Y}_X\colon\operatorname{Hom}_{\cC}(X\otimes Y,K)\,\xrightarrow{\simeq}\, \operatorname{Hom}_{\cC}(X,DY),
    \end{equation}
natural in $X\in \cC$. The duality functor $D$ assigns to an object $Y\in \cC$ its representing object $DY$. On morphisms, it is defined as the unique functor such that the family of maps $(\Phi^{Y}_X)_{X,Y\in \cC}$ is natural in $Y$. In other words, given a morphism $f\in \operatorname{Hom}_{\cC}(X,Y)$, we set 
\begin{equation*}
    D(f)\,:=\,\big(\Phi^{X}_{DY}\circ\operatorname{Hom}_{\cC}(DY\otimes f,K)\circ(\Phi^{Y}_{DY})^{-1}\big)\big(\operatorname{id_{DY}}\big).
\end{equation*} 
The assignment $D$ defines a contravariant functor on $\cC$. Different choices of representing object $DY$ and natural isomorphism $\Phi^Y_{-}$ yield canonically isomorphic duality functors. We, therefore, speak of \emph{the} duality functor.
\end{remark}

Next, we clarify the relationship between dualizing objects and LD-categories with negation. For this purpose, the graphical calculus developed in Section \ref{sec:surfdiagforGray} proves to be useful. This relationship is not surprising, considering \cite[Thm. 4.5]{WDC},
\cite[\S 4]{Pastro}, and \cite[Thm. 2.4]{HaLe}. (The reader should have Remark \ref{star autonomous vs GV} in mind.)

\begin{theorem}\label{thm:GV correspond LD neg}
The following data on a given monoidal category $(\cC,\otimes,1)$ are equivalent:
\vspace{-0.05cm}
\begin{enumerate}[label=(\roman*)]
    \item A choice of a dualizing object $K\in\cC$.
    \vspace{-0.05cm}
    \item A choice of a second monoidal structure $(\parLL,K)$ on $\cC$ together with a choice of natural transformations ${\distl \colon \otimes\circ\, (\operatorname{id_\cC} \times \parLL) \ra {\parLL}\, {\circ}\, (\otimes \times \operatorname{id_\cC})}$ and
    ${\distr \colon {\otimes}\,{\circ}\,(\parLL \times\, {\operatorname{id_\cC}}) \ra {\parLL}\, {\circ}\,(\operatorname{id_\cC}\times \,{\otimes)}}$ endowing the category $\cC$ with the structure of an LD-category with negation.
\end{enumerate}
\end{theorem}
As a consequence of Theorem \ref{thm:GV correspond LD neg}, we will use the terms GV-category and LD-category with negation interchangeably.
 
The proof of Theorem \ref{thm:GV correspond LD neg} will take up a large part of this section. We first show that the $\parLL$-monoidal unit $K$ of an LD-category with negation $\cC$ is a dualizing object of the monoidal category $(\cC,\otimes,1)$. To do so, we require the following technical lemma, which will also be needed later in several places. In particular, it prepares the notion of an LD-dual morphism, which we will discuss in Remark \ref{surface dgms for duality functor D}.
\begin{lemma}\label{computation-heavy lemma LDGV}
   Let $\cC$ be an LD-category. Let $\kappa\colon Z\otimes Y\ra K$ be an LD-pairing that admits a side-inverse LD-copairing $\overline{\kappa}\colon 1\ra Y\parLL Z$. Then, for any LD-pairing $\gamma\colon X\otimes W\ra K$ and any morphism $f\in \operatorname{Hom}_{\cC}(Y,W)$, the morphism
\begin{equation*}
\kappa\circ \Bigg(\Big({\lambda^{\parLL}_{Z}}\circ (\gamma\parLL Z)\circ((X\otimes f)\parLL Z)\circ\distl_{X,Y,Z}\circ(X\otimes\overline{\kappa})\circ(\rho^{\otimes}_X)^{-1}\Big)\otimes Y\Bigg)
\end{equation*}
   is equal to 
   \begin{equation*}
\gamma\circ(X \otimes f).
\end{equation*}
\end{lemma}
\begin{proof}
The first expression corresponds to the first surface diagram in the following calculation. Such surface diagrams are best read from the front face. From this perspective, the following graphical proof can informally be summarized as ‘pull $\gamma$ up, push $\kappa$ down, use a snake equation for $\kappa$ and $\overline{\kappa}$.’ The conceptual interest of this proof lies in the fact that the coherence axioms for LD-categories enter. In fact, they enter for the first time in this paper:
    \begin{figure}[H]
        \centering
        \includegraphics[width=0.99\textwidth]{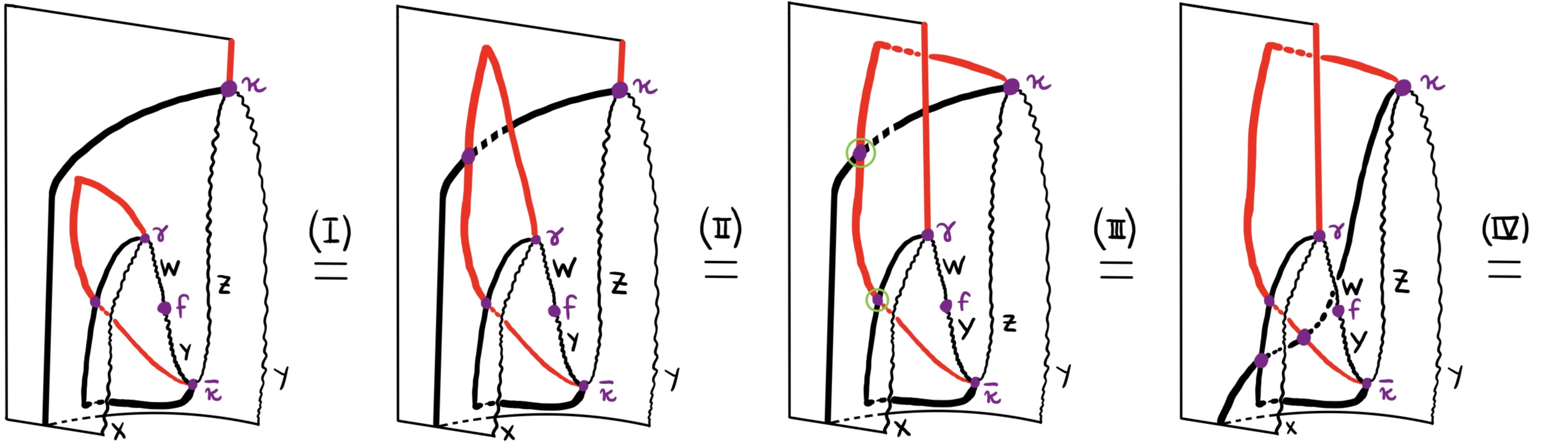}
    \end{figure}
    \begin{figure}[H]
        \centering
        \includegraphics[width=0.99\textwidth]{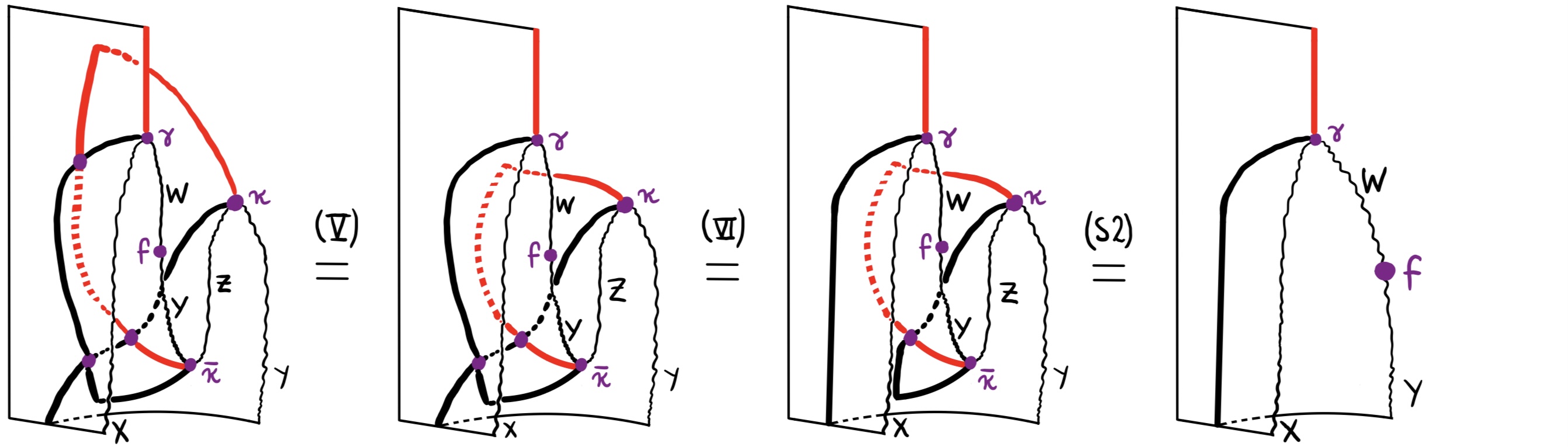}.
    \end{figure}
Equation (I) holds by the triangle coherence axiom (\ref{eq:A3}). Equation (II) follows from Mac Lane's coherence theorem for the $\parLL$-monoidal structure. This amounts to a change of sheets in the top part of the surface diagram. Equation (III) follows from the pentagon coherence axiom (\ref{eq:A9}). Intuitively, Equation (III) amounts to pulling the rightmost sheet towards the front face. We have marked the vertices to which the coherence axiom is applied with light-green circles. Equation (V) follows from the triangle axiom (\ref{eq:A4}). Equation (VI) uses the triangle diagram for the $\otimes$-monoidal product, cf. Figure \ref{fig:MacTriangle}. Equation (IV) holds by the axioms of a strict monoidal $2$-category.
\end{proof}

For the remainder of this section, we consider an LD-category with negation $\cC$. For every object $X$ in $\cC$, fix a left LD-dual $X^{\ast}$ (resp. right LD-dual $\prescript{\ast}{}{X}$) together with a corresponding evaluation $\operatorname{ev}_X$ and coevaluation $\operatorname{coev}_X$ (resp. $\operatorname{ev}'_X$ and $\operatorname{coev}'_X$). As before, let $K\in \cC$ denote the $\parLL$-monoidal unit.

\begin{prop}\label{LDisGV1}
For objects $X,Y\in \cC$, the map 
    \begin{equation*}
        \Phi^{Y}_X\colon\; \operatorname{Hom}_{\cC}(X\otimes Y,K)\,\ra\, \operatorname{Hom}_{\cC}(X,Y^{\ast})
    \end{equation*}
given by $\Phi^{Y}_X(\gamma)\;:=\;l^{\parLL}_{{Y}^{\ast}}\circ(\gamma\parLL Y^{\ast})\circ\distl_{X,Y,Y^{\ast}}\circ(X\otimes\operatorname{coev}_Y)\circ(r^{\otimes}_{X})^{-1}$ and the map 
    \begin{equation*}
        \Psi^{Y}_X\colon\; \operatorname{Hom}_{\cC}(X,Y^{\ast})\,\ra \,\operatorname{Hom}_{\cC}(X\otimes Y,K)
    \end{equation*}
given by $\Psi^{Y}_X(g)\;:=\;\operatorname{ev}_Y\circ \,(g\otimes Y)$ are inverses. Both maps are natural in the component $X\in \cC$. Using the right duality instead, one finds an analogous statement.
\end{prop}

\begin{proof}
The composite $\Phi_X^Y\circ\Psi_X^Y$ is equal to the identity map. This follows from the functoriality of the $\otimes$-tensor product, the naturality of the distributors and unitors, and finally from snake equation (\ref{firstzigzag}). Conversely, for any $\gamma\in \operatorname{Hom}_{\cC}(X\otimes Y,K)$, the morphism $(\Psi_X^Y\circ\Phi_X^Y)(\gamma)$ is equal to $\gamma$. This follows immediately from Lemma \ref{computation-heavy lemma LDGV} by setting $\kappa:=\operatorname{ev}_Y$, $\overline{\kappa}:=\operatorname{coev}_Y$ and $f:=\operatorname{id}_{Y}$.
\end{proof}

By Proposition \ref{LDisGV1}, we can use the arguments from Remark \ref{rep determines functoriality} to define a contravariant functor $D$ on $\cC$. On objects, this functor is given by $D(X):=X^{\ast}$ for $X\in \cC$. Although we will eventually show that $D$ is the duality functor, we have not yet established that it is an antiequivalence at this stage.

\begin{remark}\label{surface dgms for duality functor D}
  Let us describe the surface diagrams of $D$. Using the distributor, the functor $D$ is defined on morphisms as follows:
\begin{figure}[H]
    \centering
    \includegraphics[width=0.95\textwidth]{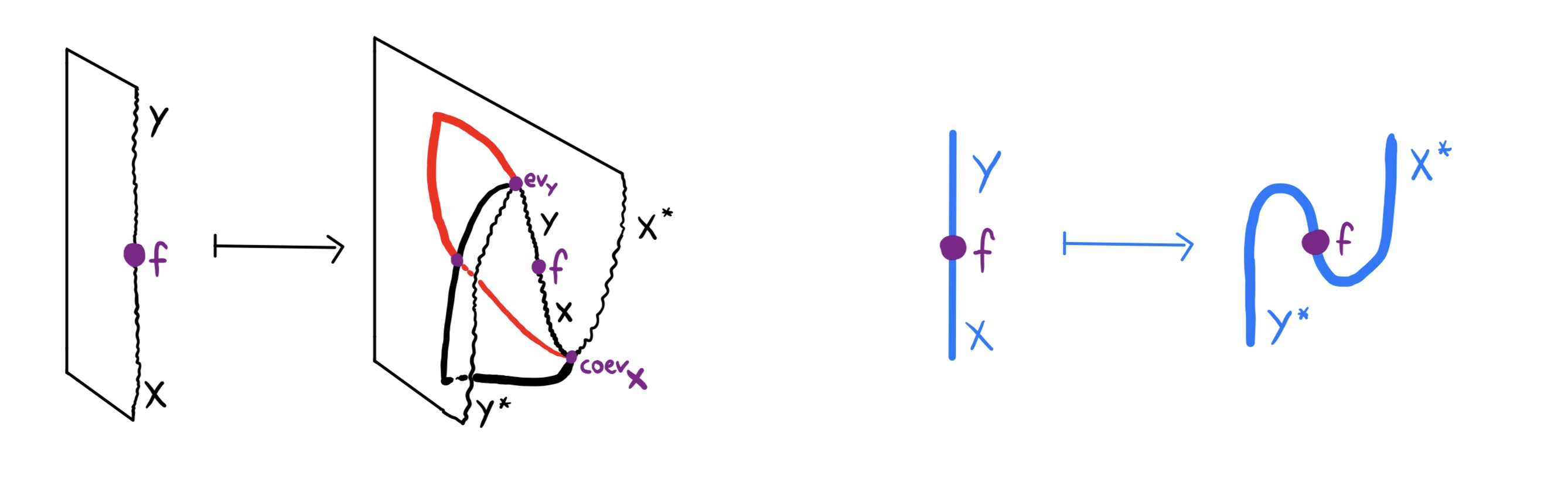}
    \caption{The duality functor $D$ and the extracted string diagrams.}
    \end{figure}
Analogously, the assignment $D^\prime$, given by $D^\prime(X):=\prescript{\ast}{}{X}$ for an object $X\in \cC$ and by 
    \begin{equation*}
        D^\prime(f)\,:=\,r^{\parLL}_{\prescript{\ast}{}{X}}\circ\big({\prescript{\ast}{}{X}}\parLL\operatorname{ev}'_Y\big)\circ\big({\prescript{\ast}{}{X}}\otimes(f\otimes {\prescript{\ast}{}{Y}})\big)\circ\distr_{\prescript{\ast}{}{X},X,\prescript{\ast}{}{Y}}\circ\big(\operatorname{coev}'_X\otimes {^{\ast}Y}\big)\circ(l^{\otimes}_{^{\ast}Y})^{-1}
    \end{equation*}
for a morphism $f\in \operatorname{Hom}_{\cC}(X,Y)$, defines another contravariant functor $D^{\prime}$ on $\cC$. Up to relabelling, the surface diagram of the morphism $D^\prime(f)$ is obtained by applying the $2$-functor $(-)^{\operatorname{rev}}$ to the surface diagram of the morphism $D(f)$.  
\end{remark}

From now on, we will occasionally depict the functors $D\colon \cC^{\operatorname{op}}\ra \cC$ and $D^\prime\colon \cC\ra \cC^{\operatorname{op}}$ without remembering their explicit definition on morphisms:

\begin{remark}\label{opposite cat and functor D}
 The opposite category will be colored light blue. The $1$-cell depicting the duality functor $D$ will also be drawn in blue. As an example, the evaluation
    \begin{equation*}
        \operatorname{ev}_X\colon D(X)\otimes X\, \longrightarrow \, K
    \end{equation*} 
 will often be represented as in Figure \ref{evX}. \nid For $f\in \operatorname{Hom}_{\cC^{\operatorname{op}}}(Y,X)\eqdef \operatorname{Hom}_{\cC}(X,Y)$, the morphism $D(f)\in \operatorname{Hom}_{\cC}(D(Y),D(X))$ will be drawn as in Figure \ref{Df}. Note that natural transformations in $\cC$ are depicted in the reversed direction if located inside the light-blue $2$-cell, i.e. if one considers them in $\cC^{\operatorname{op}}$.
    \begin{figure}[H]
        \centering
    \begin{subfigure}[b]{0.2\textwidth}
        \centering
        \includegraphics[width=0.8\textwidth]{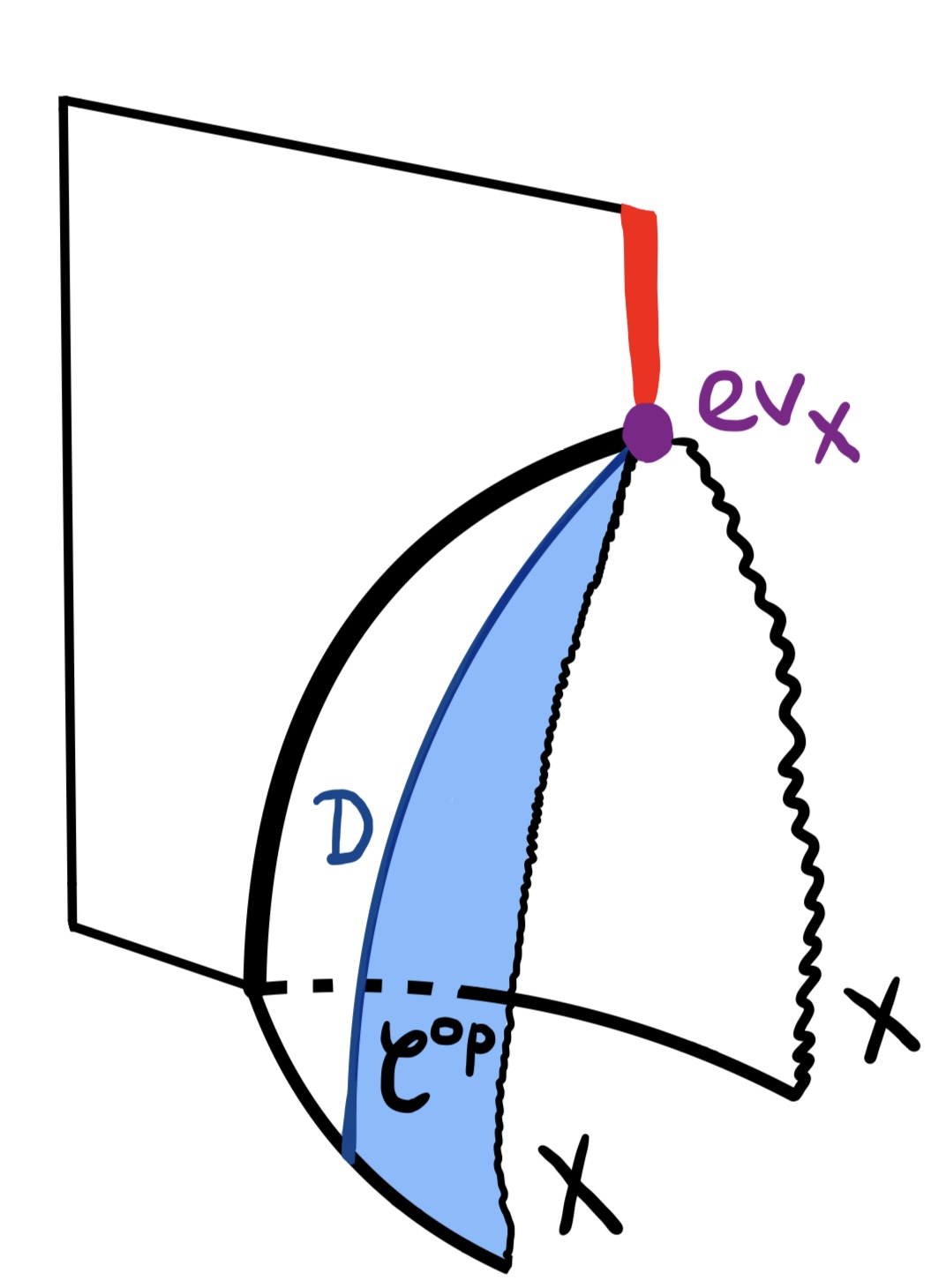}
        \caption{A cap $\operatorname{ev}_X$.}
        \label{evX}
        \end{subfigure}
    \hspace{8em}
    \begin{subfigure}[b]{0.2\textwidth}
         \centering
         \includegraphics[width=0.8\textwidth]{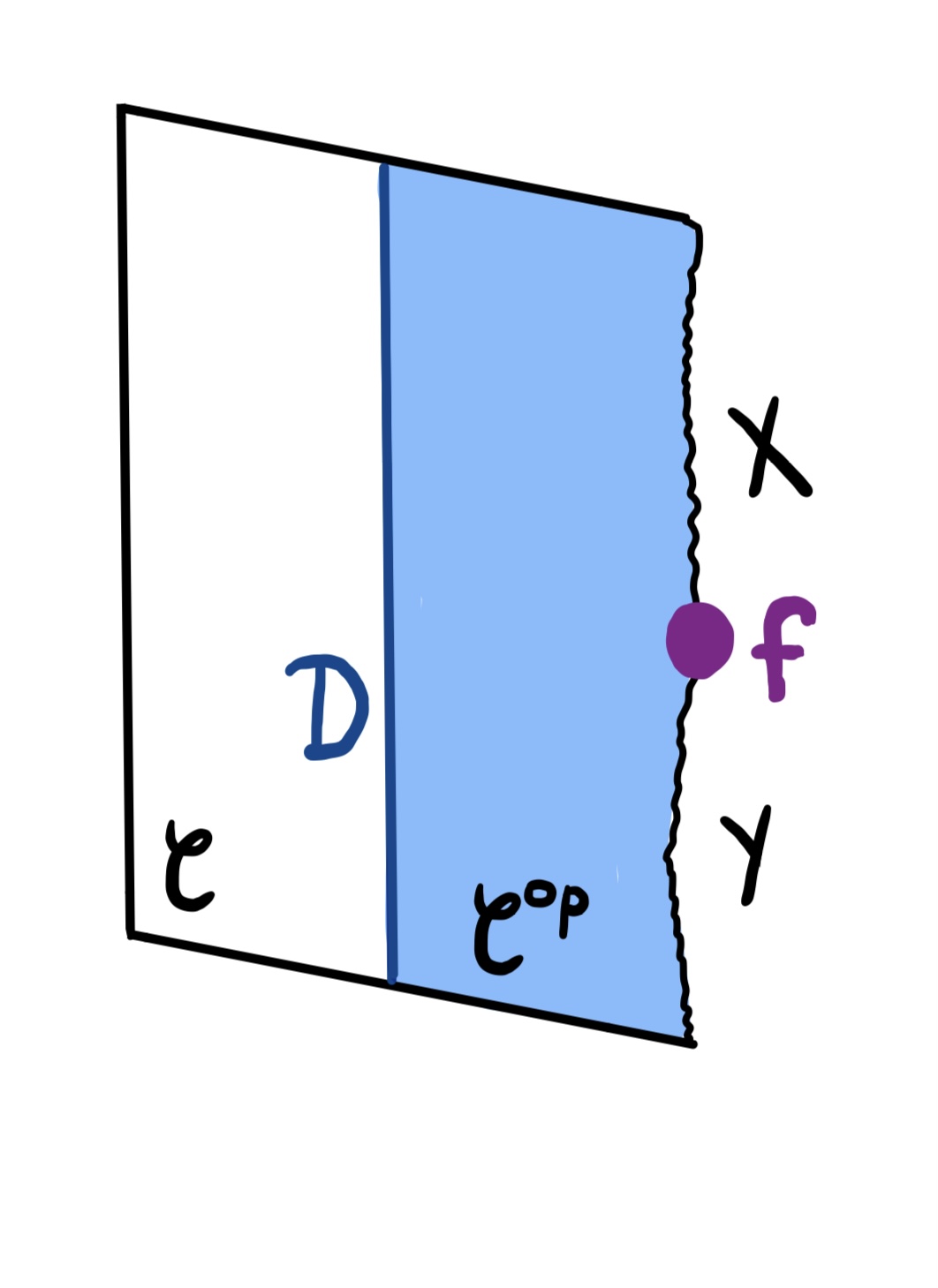}
         \caption{$D(f)$.}
         \label{Df}
    \end{subfigure}
    \caption{The duality functor $D$, continued.}
    \label{fig:Drawing D}
    \end{figure}
\end{remark}

The next proposition is a direct consequence of Lemma \ref{computation-heavy lemma LDGV}. Graphically, it states that we can slide morphisms along cups and caps.

\begin{prop}\label{evalDuality functor}
For any $f\in \operatorname{Hom}_{\cC}(X,Y)$, we have the equation
\begin{equation*}
    {\operatorname{ev}_X} \circ (D(f)\otimes X) \, = \, {\operatorname{ev}_Y} \circ (D(Y)\otimes f).
\end{equation*} Analogous statements hold for the coevaluation morphisms and the contravariant functor $D^\prime$.
\end{prop}

We can now complete the proof that any LD-category with negation $\cC$ carries a GV-structure: 

\begin{prop}\label{quasi-inverses}
The functors $D$ and $D^\prime$ are quasi-inverse to each other. In particular, $D$ is an antiequivalence.
\end{prop}

\begin{proof}
We only show that the composite $DD^\prime$ is isomorphic to the identity functor. Applying the $2$-functor $(-)^{\operatorname{rev}}$ to the surface diagrams occurring in the proof for the composite $DD^\prime$, we then find that the other composite $D^\prime D$ is isomorphic to the identity functor. 

For any object $X$ in $\cC$, consider the morphisms 
\begin{gather*}
\eta_X\,:=\,r^{\parLL}_{{^{\ast}(X^{\ast})}}\circ\big({^{\ast}(X^{\ast})}\parLL\operatorname{ev}_X\big)\circ\distr_{{^{\ast}(X^{\ast})},X^{\ast},X}\circ\big(\operatorname{coev}'_{X^{\ast}}\otimes X\big)\circ(l^{\otimes}_{X})^{-1}
\\
\epsilon_X\,:=\,r^{\parLL}_{X}\circ\big(X\parLL\operatorname{ev}'_{X^{\ast}}\big)\circ\distr_{X,X^{\ast},{^{\ast}(X^{\ast})}}\circ\big(\operatorname{coev}_{X}\otimes {^{\ast}(X^{\ast})}\big)\circ(l^{\otimes}_{{^{\ast}(X^{\ast})}})^{-1}.
\end{gather*}
The equation $\epsilon_X\circ\eta_X=\operatorname{id}_X$ is proved by first using Lemma \ref{computation-heavy lemma LDGV} (applied to the LD-category $\cC^{\operatorname{rev}}$ with $\kappa:=\operatorname{ev}'_{X^{\ast}}$, $\overline{\kappa}:=\operatorname{coev}'_{X^{\ast}}$ and $\gamma:=\operatorname{ev}_X$) and then applying snake equation (\ref{secondzigzag}). It is instructive to perform the proof using surface diagrams. Up to relabelling, the same graphical proof then establishes the other equation $\eta_X \circ \epsilon_X=\operatorname{id}_{^{\ast}(X^{\ast})}$.
By Proposition \ref{evalDuality functor}, both families of morphisms $(\epsilon_X)_{X}$ and $(\eta_X)_{X}$ are natural in $X$.
\end{proof} 

To show the converse, namely that any GV-category carries the structure of an LD-category with negation, we first recall the following lemma.

\begin{lemma}\label{lemma equivalence induces monoidal structure}
Let $\cA$ and $\cB$ be categories, and let $F\colon \cA \ra \cB$ be an equivalence of categories with quasi-inverse $G$. Then, any monoidal structure $(\otimes_{\cA},1_{\cA})$ on $\cA$ induces a monoidal structure $(\otimes_{\cB},1_{\cB})$ on $\cB$, where the monoidal product is given by
    \begin{equation*}
        X\otimes_{\cB}Y\,:=\, F(GX\otimes_{\cA}GY), \qquad \text{for } X,Y\in \cB,
    \end{equation*}
and the monoidal unit is $1_{\cB}\,:=\,F(1_{\cA}).$
The monoidal categories $(\cA,\otimes_{\cA},1_{\cA})$ and $(\cB,\otimes_{\cB},1_{\cB})$ are monoidally equivalent by strong monoidal functors, with underlying functors $F$ and $G$.
\end{lemma}

Let us collect our results:

\begin{proof}[Proof of Theorem \ref{thm:GV correspond LD neg}]
    Propositions \ref{LDisGV1} and \ref{quasi-inverses} have shown that the $\parLL$-monoidal unit of any LD-category with negation is a dualizing object in the underlying $\otimes$-monoidal category.

    For the other direction, let $(\cC,\otimes,1,K)$ be a GV-category. By Lemma \ref{lemma equivalence induces monoidal structure}, the monoidal structure $(\otimes^{\operatorname{rev}},1)$ on the opposite category $\cC^{\operatorname{op}}$ induces on $\cC$ a second monoidal product 
    \begin{equation*}
        \parLL\;:=\;D^{\prime}\circ\otimes^{\operatorname{rev}}\circ (D\times D)
    \end{equation*}
    with monoidal unit $D^{\prime}(1)$. %Note that we have not distinguished between the functor $\otimes^{\operatorname{rev}}$ and its opposite $(\otimes^{\operatorname{rev}})^{\operatorname{op}}$. We will continue not to do so. 
    Using Yoneda's lemma, we identify $D^{\prime}(1)\cong K$. For the definition of the distributors ${\distl \colon \otimes\circ\, (\operatorname{id_\cC} \times \parLL) \ra {\parLL}\, {\circ}\, (\otimes \times \operatorname{id_\cC})}$ and
    ${\distr \colon {\otimes}\,{\circ}\,(\parLL \times\, {\operatorname{id_\cC}}) \ra {\parLL}\, {\circ}\,(\operatorname{id_\cC}\times \,{\otimes)}}$, as well as for the proof that they make $(\cC,\otimes,1,\parLL, K)$ into an LD-category, we refer to \cite[Prop. 4.11]{fuchs2024grothendieckverdierdualitycategoriesbimodules}. In \cite[Prop. 3.35]{fuchs2024grothendieckverdier} it is shown that the linearly distributive structure on a GV-category indeed defines an LD-category with negation. It is immediate that the assignments which assign an LD-category with negation to a GV-category and vice versa are inverses.
\end{proof}

Next, we recall the following observation from \cite[\S 4]{BoDrinfeld}:

\begin{remark}\label{double dual mon structure}
    For any GV-category $(\cC,\otimes,1,K)$, there is a canonical isomorphism
    \begin{equation*}
        \operatorname{Hom}_{\cC}(X\otimes (Y\otimes Z),K)\,\xlongrightarrow{\simeq}\,\operatorname{Hom}_{\cC}(D^2(Y\otimes Z)\otimes X,K)
    \end{equation*}
    natural in $X,Y,Z\in \cC$. It is obtained by pre- and postcomposing the natural isomorphism \begin{equation*}
        \operatorname{Hom}_{\cC}(X, D(Y\otimes Z))\,\xlongrightarrow{\simeq}\,\operatorname{Hom}_{\cC}(D^2(Y\otimes Z),D(X))
        \end{equation*}
        with the natural isomorphism (\ref{representing transformation}).
Similarly, we find a canonical isomorphism 
    \begin{equation*}
        \operatorname{Hom}_{\cC}(X\otimes (Y\otimes Z),K)\,\xlongrightarrow{\simeq}\,\operatorname{Hom}_{\cC}(D^2(Y)\otimes D^2(Z)\otimes X,K)
    \end{equation*}
natural in $X,Y,Z\in \cC.$ With Yoneda's lemma, this yields a natural isomorphism 
    \begin{equation*}\label{D2 strong monoidal}
        \varphi^{2,D^2}\colon\, D^2\circ \otimes\,\xlongrightarrow{\simeq}\, \otimes\circ(D^2\times D^2).
    \end{equation*}
It endows the covariant equivalence $D^2$ with a strong monoidal structure $(\varphi^{2,D^2},\varphi^{0,D^2})$ by \cite[Prop. 4.2]{BoDrinfeld}.
\end{remark}

We use this monoidal structure on the double dual in the proof of the following result:

\begin{prop}\label{duality functor is FrobLD}
The duality functor $D$ of a GV-category $\cC=(\cC,\otimes,1,K)$ carries the structure of a strong Frobenius LD-functor $\cC\ra \cC^{\operatorname{lop}}$. Here, as in Definition \ref{three2-functors}.(iii), $\cC^{\operatorname{lop}}$ denotes the linearly distributive opposite of $\cC$.
\end{prop}

\begin{proof}
Applying Lemma \ref{lemma equivalence induces monoidal structure} to the monoidal category $(\cC^{\operatorname{op}},\otimes^{\text{rev}},1)$ and the equivalence $(D^{\prime},D)$, we obtain the monoidal structure $(\parLL,D^{\prime}(1)\cong K)$ on $\cC$, where $X\parLL Y:=D^{\prime}(DY\otimes DX)$. Lemma \ref{lemma equivalence induces monoidal structure} also yields a strong monoidal structure on the duality functor 
    \begin{equation}\label{first strong mon structure on D}
        D\colon\,(\cC,\parLL,D^{\prime}(1)\cong K) \,\longrightarrow\, (\cC^{\operatorname{op}},\otimes^{\operatorname{rev}},1).
    \end{equation}

Similarly, the monoidal structure $(\otimes,1)$ on $\cC$ together with the equivalence $(D,D^{\prime})$ induces another monoidal structure $(\parLL',D(1))$ on the category $\cC^{\operatorname{op}}$, where $X\parLL^{\prime}Y:=D(D^{\prime}X\otimes D^{\prime}Y)$. Lemma \ref{lemma equivalence induces monoidal structure} also yields a strong monoidal structure on the duality functor 
    \begin{equation}\label{second strong mon structure on D}
        D\colon \,(\cC,\otimes,1)\,\longrightarrow \, (\cC^{\operatorname{op}}, \parLL',D(1)).
    \end{equation}

Recall from Remark \ref{double dual mon structure} the strong $\otimes$-monoidal structure $(\varphi^{2,D^2},\varphi^{0,D^2})$ on the double dual $D^2$. For $X,Y\in \cC$, let $\widetilde{m}_{X,Y}$ denote the isomorphism $D^{\prime}(\varphi^{2,D^2}_{D^{\prime}X, D^{\prime}Y})^{-1}.$ The natural isomorphism
    \begin{align*}
        D^{\prime}(DX\otimes DY)\,\simeq\,D^{\prime}(D^2D^{\prime}X\otimes D^2D^{\prime}Y)\,\xrightarrow{\widetilde{m}_{X,Y}}\,D^{\prime}D^2(D^{\prime}X\otimes D^{\prime}Y)\\\simeq D(D^{\prime}X\otimes D^{\prime}Y),
    \end{align*}
together with the isomorphism $D^{\prime}(1)\,\xrightarrow{D^{\prime}(\varphi^{0,D^2})^{-1}}\,D^{\prime}D^2(1)\cong D(1),$ equips the identity functor on $\cC^{\operatorname{op}}$ with a strong monoidal structure
    \begin{equation}\label{id functor on cop}
        \operatorname{id}_{\cC^{\operatorname{op}}}\colon\,(\cC^{\operatorname{op}},\parLL',D(1))\,\longrightarrow\,(\cC^{\operatorname{op}},\parLL^{\operatorname{rev}},D^{\prime}(1)).
    \end{equation}
By composing the strong monoidal functor (\ref{second strong mon structure on D}) with the strong monoidal functor (\ref{id functor on cop}), we obtain a strong monoidal structure on the duality functor
    \begin{equation}\label{2nd strong mon structure on D}
        D\colon \,(\cC,\otimes,1)\,\longrightarrow \,(\cC^{\operatorname{op}},\parLL^{\operatorname{rev}},D^{\prime}(1)\cong K).
    \end{equation}
The identities proven in \cite[Prop. 3.30]{fuchs2024grothendieckverdier} amount to the statement that the monoidal structures (\ref{first strong mon structure on D}) and (\ref{2nd strong mon structure on D}) endow the duality functor $D$ with the structure of a Frobenius LD-functor $\cC\ra \cC^{\operatorname{lop}}$.
\end{proof} 

Propositions \ref{duality functor is FrobLD} and \ref{quasi-inverse gets Frobenius LD-structure} immediately imply the following result.

\begin{corollary}\label{dual dprime is frobenius ld}
Let $(\cC,\otimes,1,K)$ be a GV-category. The quasi-inverse $D^{\prime}$ of the duality functor $D$ carries the structure of a strong Frobenius LD-functor such that the unit $\eta\colon \operatorname{id}_{\cC}\xrightarrow{\simeq}D^{\prime}D$ and counit $\epsilon\colon DD^{\prime}\xrightarrow{\simeq}\operatorname{id}_{\cC}$ are isomorphisms of Frobenius LD-functors.
\end{corollary}

We will need the following lemma in several places. It is easily verified using string diagrams for the $2$-category of categories, functors, and natural transformations.

\begin{lemma}\label{unit counit under duality}
Let $(\cC,\otimes,1,K)$ be a GV-category. Consider the unit $\eta\colon \operatorname{id}_{\cC}\xrightarrow{\simeq}D^{\prime}D$ and the counit $\epsilon\colon DD^{\prime}\xrightarrow{\simeq}\operatorname{id}_{\cC}$ of the adjoint equivalence $D\dashv D^{\prime}.$
For all $X\in \cC$, we then have 
\begin{equation*}
    D(\eta_X)\;=\;\epsilon_{DX}.
\end{equation*}
\end{lemma}

\begin{remark}\label{DoubleDual FrobLD}
As a composite of Frobenius LD-functors, the double dual $D^2$ carries the structure of a strong Frobenius LD-endofunctor by Proposition \ref{LDC2Cat}. In particular, the double dual $D^2$ is a strong monoidal endofunctor on the monoidal category $(\cC,\otimes,1)$. Using Lemma \ref{unit counit under duality}, one can show that this strong $\otimes$-monoidal structure is precisely the one described by Boyarchenko-Drinfeld \cite[\S 4]{BoDrinfeld}.
\end{remark}

\section{LD-Frobenius algebras}\label{sec:LD-FrobAlg}
In this section, we fix an LD-category $\cC$ and view the unit category $\cI$ as an LD-category with coinciding monoidal products. 

\subsection{Frobenius relations}
Algebras and coalgebras in the LD-category $\cC$ are defined with respect to the two different monoidal structures of $\cC$ as follows:
\begin{definition}\label{definition algebra-coalgebra}
    An \emph{algebra} in $\cC$ is a lax $\otimes$-monoidal functor $\cI\ra \cC$. A \emph{coalgebra} in $\cC$ is an algebra in the LD-category $\cC^{\text{cop}}$, i.e. an oplax $\parLL$-monoidal functor $\cI\ra \cC$.
\end{definition}

Since the unit category is represented as transparent, the surface diagrams for an algebra in $\cC$ are obtained as follows: Cut the surface diagrams for a lax $\otimes$-monoidal functor from Section \ref{application:moncats and functors} along the blue functor lines. These cuts yield two surfaces from each surface diagram. For each pair, we retain only the surface that is furthest from the front face.

\begin{remark}\label{Cutting-off surface dgms AlgCoalg}
    Specifically, by cutting the surface diagrams in Figures \ref{fig:varphi2} and \ref{fig:varphi0}, we obtain the surface diagrams for the multiplication $\mu$ and unit $\eta$ of an algebra $(A,\mu,\eta)$ in $\cC$, as shown in Figures \ref{multi mu} and \ref{unit eta}. Dually, Figures \ref{comulti delta} and \ref{counit epsilon} depict the surface diagrams for the comultiplication $\Delta$ and counit $\epsilon$ of a coalgebra $(C,\Delta,\epsilon)$ in $\cC$.
\begin{figure}[H]
    \centering
    \begin{subfigure}[b]{0.2\textwidth}
    \centering
    \includegraphics[width=0.62\textwidth]{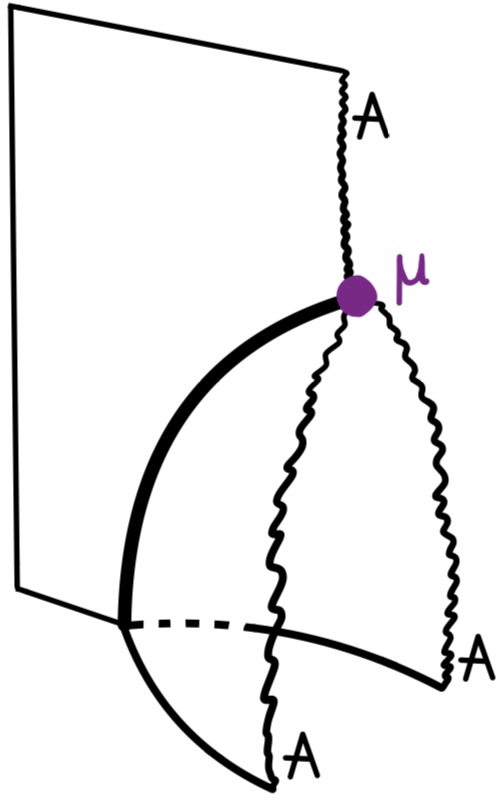}
    \caption{Multiplication $\mu$.}
    \label{multi mu}
    \end{subfigure}
\hspace{1.5em}
    \begin{subfigure}[b]{0.2\textwidth}
    \centering
    \includegraphics[width=0.72\textwidth]{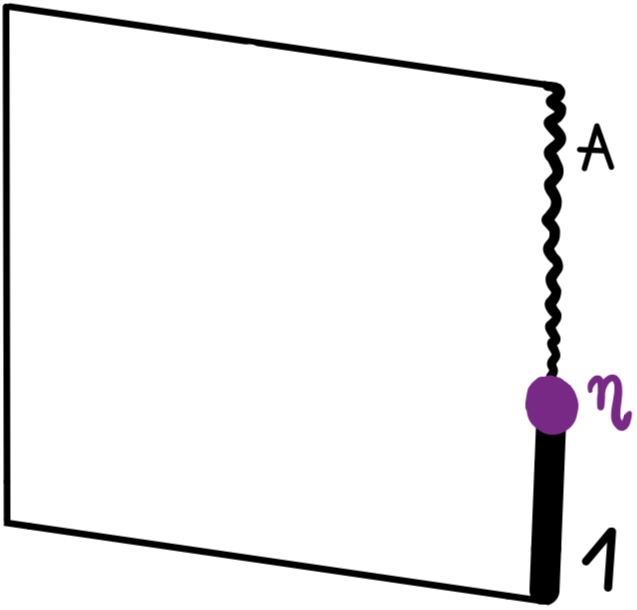}
    \caption{Unit $\eta$.}
    \label{unit eta}
    \end{subfigure}
\hspace{1.5em}
    \begin{subfigure}[b]{0.23\textwidth}
    \centering
    \includegraphics[width=0.57\textwidth]{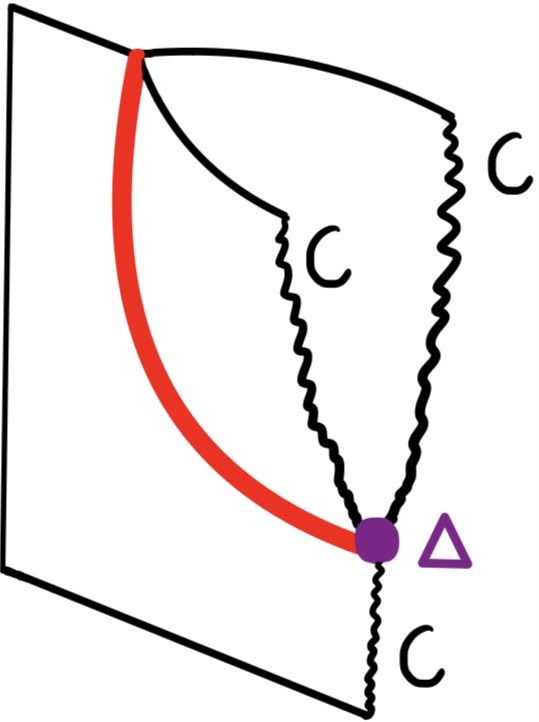}
    \caption{Comultiplication $\Delta$.}
    \label{comulti delta}
    \end{subfigure}
\hspace{1.5em}
    \begin{subfigure}[b]{0.2\textwidth}
    \centering
    \includegraphics[width=0.72\textwidth]{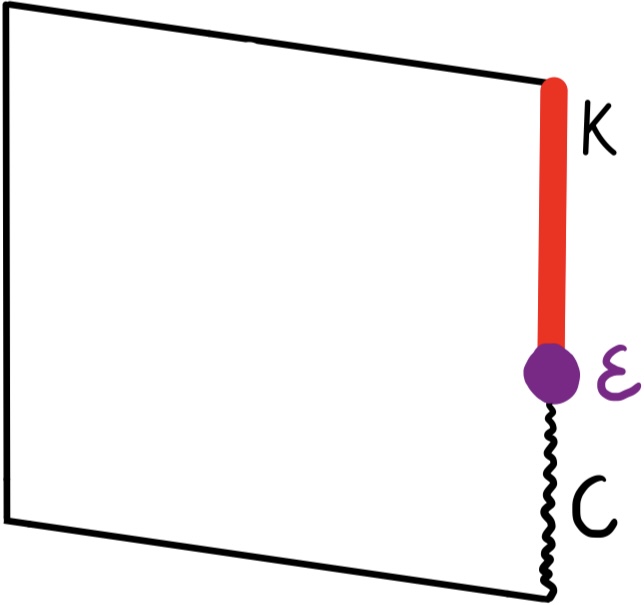}
    \caption{Counit $\epsilon$.}
    \label{counit epsilon}
    \end{subfigure}
\caption{Surface diagrams for an algebra $(A,\mu,\eta)$ and a coalgebra $(C,\Delta,\epsilon)$.}
\label{fig:surface dgms for alg and coalg}
\end{figure}
\end{remark}

\begin{figure}[H]
    \centering
    \begin{subfigure}[b]{0.48\textwidth}
    \centering
    \includegraphics[width=0.75\textwidth]{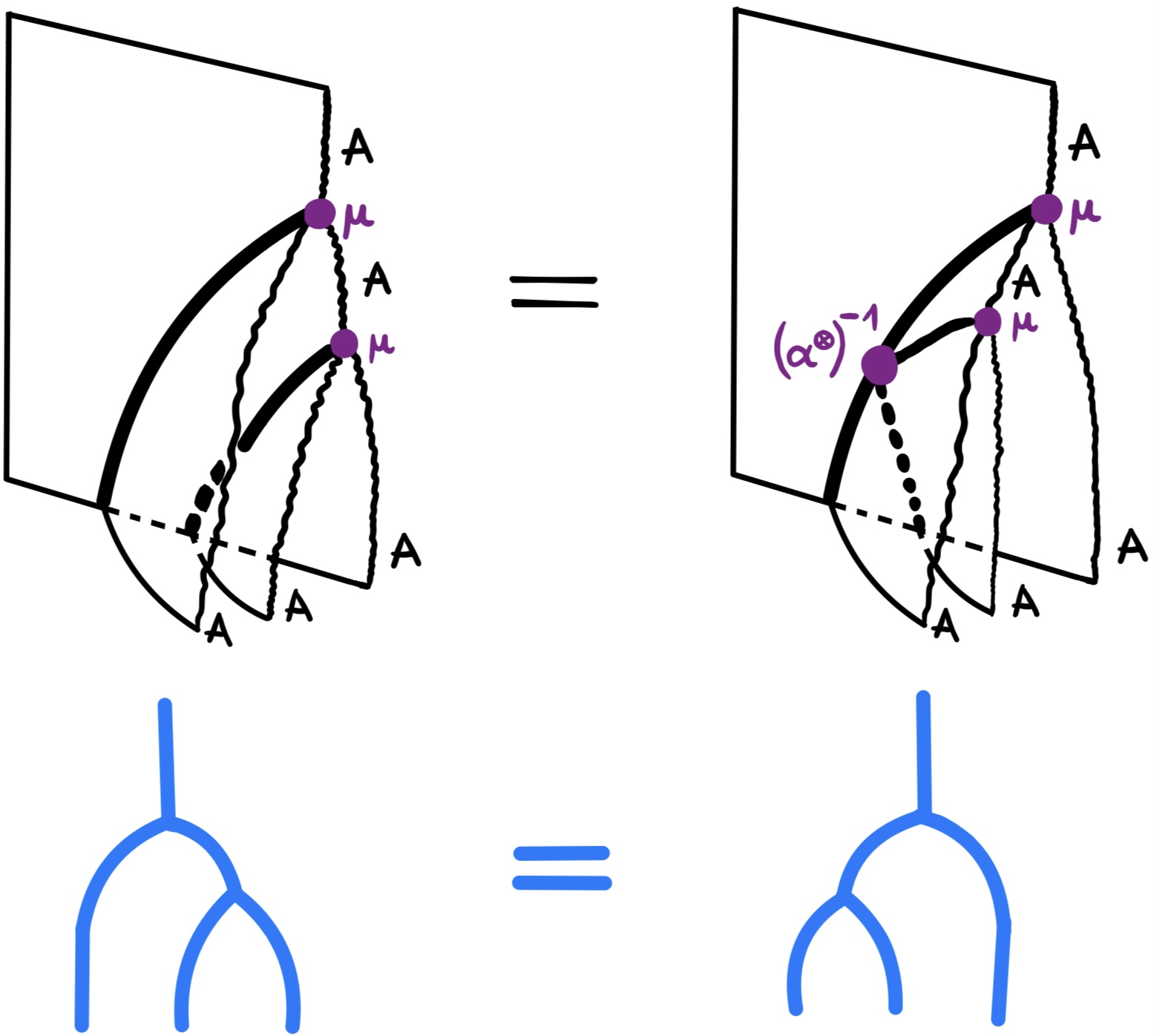}
    \caption{Associativity.}
    \label{Associativity}
    \end{subfigure}
\hfill
    \begin{subfigure}[b]{0.48\textwidth}
    \centering
    \includegraphics[width=0.88\textwidth]{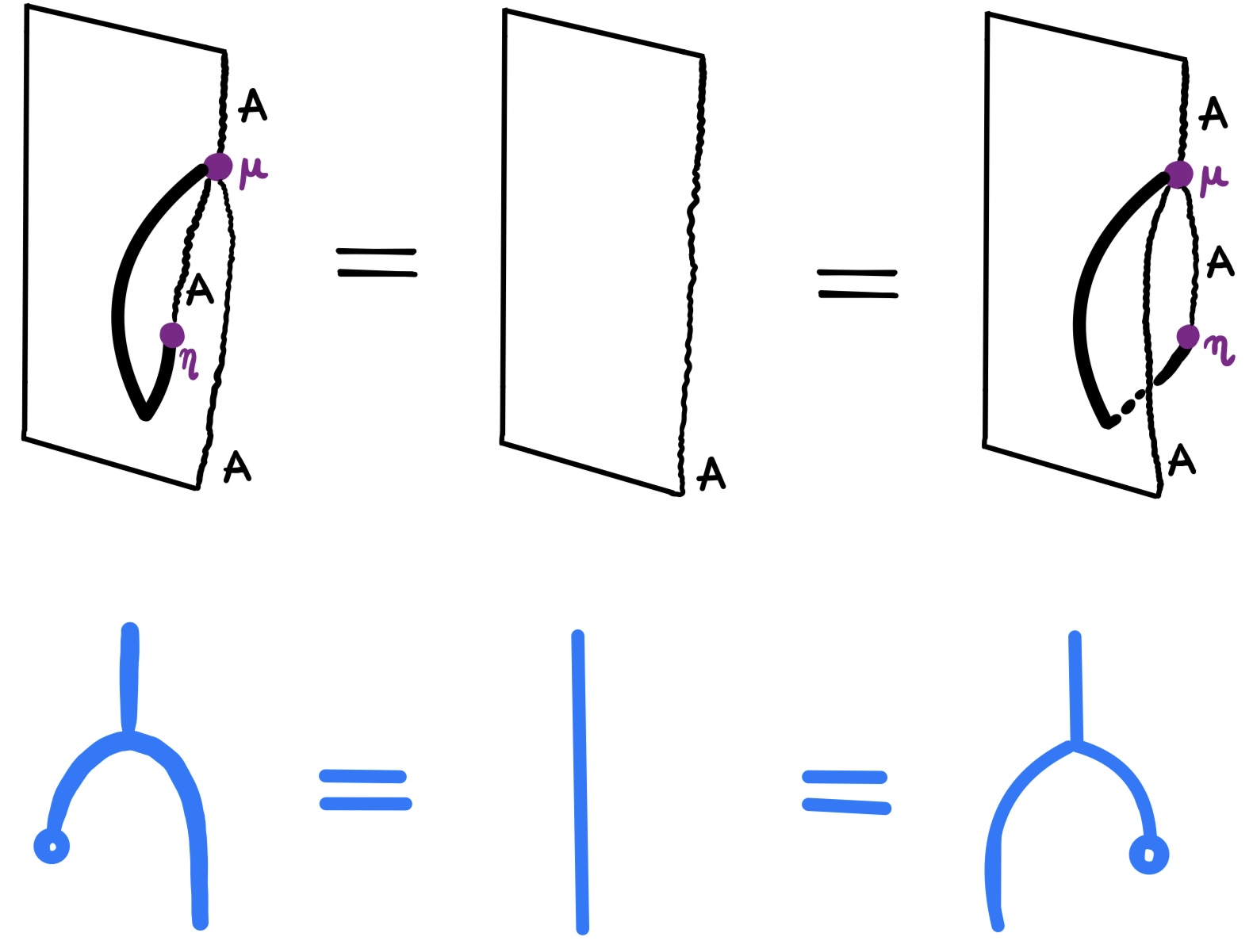}
    \caption{Unitality.}
    \label{unitality}
    \end{subfigure}
\caption{Surface diagrams for an algebra $(A,\mu,\eta)$, continued.}
\label{fig:Associativity and unitality of an algebra}
\end{figure}

\begin{remark}
  Analogously, we cut off the surface diagram for the coherence axioms of a lax $\otimes$-monoidal functor. We then obtain the surface diagrams that represent the associativity and unitality of an algebra; see Figure \ref{fig:Associativity and unitality of an algebra}. When viewing the resulting surfaces from the front face, we detect the well-known string diagrams encoding associativity and unitality for an algebra in a strict monoidal category. Associators and unitors occur as $0$-cells in the negative X-direction behind them; compare Remark \ref{Mac Lane's strictification theorem}. An application of the $2$-functor $(-)^{\operatorname{cop}}$ from Definition \ref{three2-functors}.(ii) to the surface diagrams from Figure \ref{fig:Associativity and unitality of an algebra} yields the surface diagrams encoding the coassociativity and the counitality relations of a coalgebra in $\cC$.
\end{remark}

We describe the surfaces from Figure \ref{fig:Associativity and unitality of an algebra} informally to provide three-dimensional intuition:

\begin{remark}
Moving from the left-hand to the right-hand surface, the associativity relation in Figure \ref{Associativity} amounts to picking the middle sheet at the lower $0$-cell $\mu$, pulling it first towards the upper $0$-cell $\mu$ and the pulling it down onto the sheet closest to the side face. In this process, the functor labelling the bottom face of our canvas is unchanged. As a result, the middle sheet must cross the $\otimes$-functor line closest to the back face to move onto the sheet closest to the front face. This crossing is labelled by an inverse associator. The unitality relations in Figure \ref{unitality} intuitively state that a cupped shell can be attached to (or detached from) a sheet.
\end{remark}

\begin{remark}
The STL files for the surfaces in Figures \ref{fig:surface dgms for alg and coalg} and \ref{fig:Associativity and unitality of an algebra} can be found \href{https://maxdemirdilek.github.io/Research/SurfaceDiagrams}{here}.
\end{remark}

Next, we define bi-(co)modules over (co)algebras.

\begin{definition}\label{def module comodule}
    Given an algebra $A$ in $\cC$, a \emph{(left) $A$-module} is a left $A$-module in the monoidal category $(\cC,\otimes,1).$ Given two algebras $A,B$ in $\cC$, an \emph{$(A,B)$-bimodule} is an $(A,B)$-bimodule in the monoidal category $(\cC,\otimes,1).$ \emph{(Left) $C$-comodules} and \emph{$(C,D)$-bi-comodules} over coalgebras $C,D$ are defined dually by passing to the LD-category $\cC^{\operatorname{cop}}$, i.e. by using the $\parLL$-monoidal structure.
\end{definition}

\begin{remark}
    Similarly to Remark \ref{Cutting-off surface dgms AlgCoalg}, the surface diagrams for (co)modules and bi-(co)modules can be seen as cut-off surface diagrams. For instance, by cutting off the surface diagrams for modules and bi-modules over a lax $\otimes$-monoidal functor, one obtains the diagrams for modules and comodules in the sense of Definition \ref{def module comodule}. See \cite[Def. 39]{YetterModuleLaxMon} for a definition of the former notions.
\end{remark}

\begin{definition}
    Let $A$ be an algebra and $C$ be a coalgebra in $\cC$. Given two modules $(M,\beta)$ and $(N,\gamma)$ over $A$, a morphism $f\colon M \ra N$ is a \emph{morphism of $A$-modules} if $\gamma\circ(A\otimes f)=f\circ\beta$. A \emph{morphism of $C$-comodules} is defined dually. We denote the category of $A$-modules by $_A\cC$ and the category of $C$-comodules by $^C\cC$.
\end{definition}

We now come to Frobenius algebras; here, we have an algebra and a coalgebra with compatibility conditions. Specifically, by cutting off the surface diagrams for a Frobenius LD-functor (Definition \ref{def:FrobLin}), one generalizes Frobenius algebras from monoidal categories to LD-categories:

\begin{definition}\label{def:FrobAlg}
    An \emph{LD-Frobenius algebra} in $\cC$, or \emph{Frobenius algebra} for short, is a Frobenius LD-functor $\cI \ra \cC$. Here, $\cI$ denotes the unit category. A \emph{morphism of LD-Frobenius algebras} is a morphism between Frobenius LD-functors in the sense of Definition \ref{def:MorphFrob}.
\end{definition}

Explicitly, a Frobenius algebra in $\cC$ consists of an object $A\in \cC$ endowed with an algebra structure $(\mu,\eta)$ with respect to $(\otimes,1)$ and a coalgebra structure $(\Delta,\epsilon)$ with respect to $(\parLL,K)$ such that the following \emph{LD-Frobenius relations} hold:
\bigskip
\begin{equation}\label{LD-FrobeniusRelations}
\begin{matrix}
\includegraphics[scale=0.12]{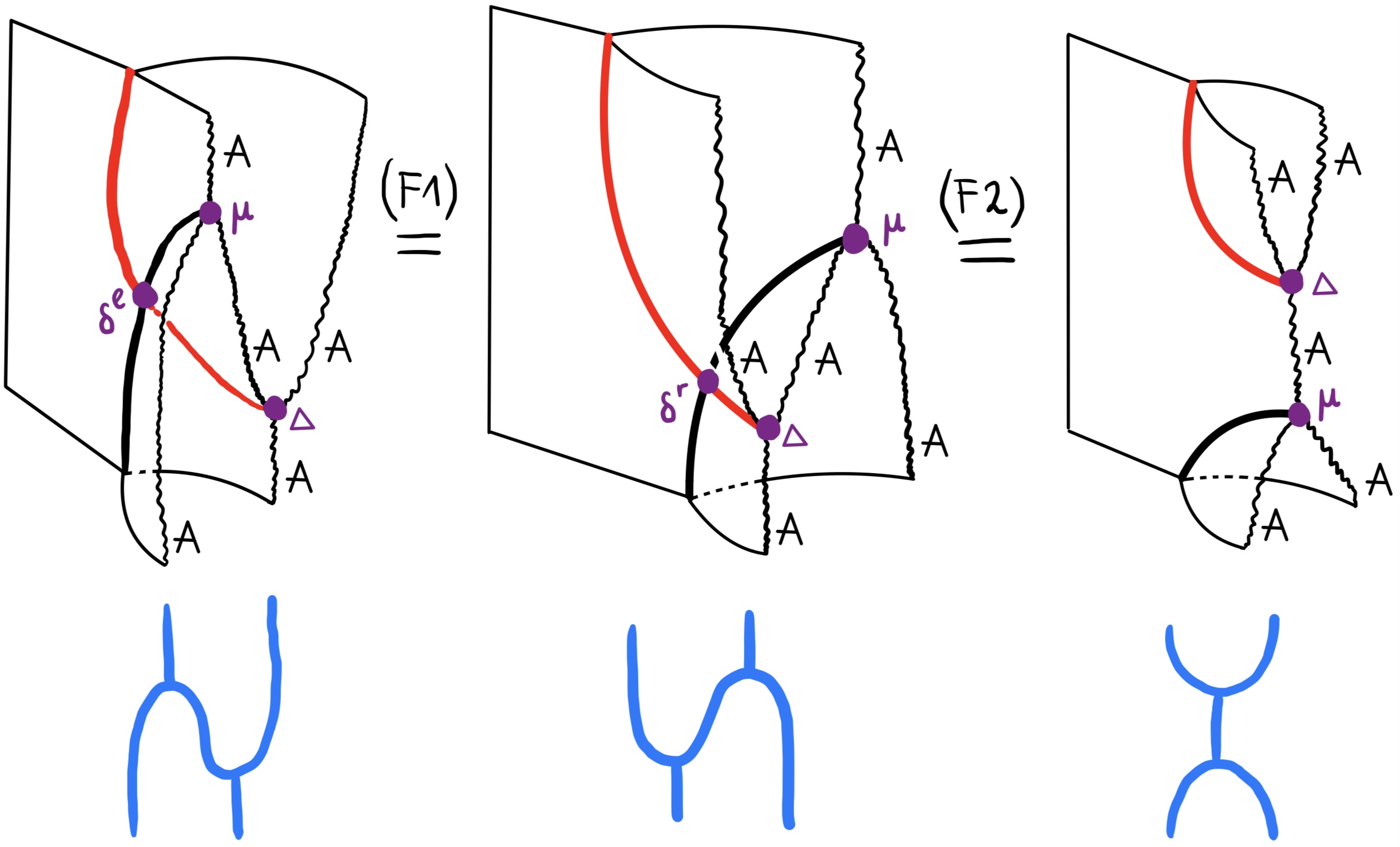}.
\end{matrix}
\end{equation}

From now on, the LD-Frobenius relations will simply be called \emph{Frobenius relations}. 

\begin{remark}
    The STL and HOM files for Frobenius algebras are presented \href{https://maxdemirdilek.github.io/Research/SurfaceDiagrams}{here}.
\end{remark}

\begin{remark}
    The Frobenius relations (\ref{LD-FrobeniusRelations}) generally involve non-invertible distributors. We thus need to include the third dimension and use surface diagrams.
\end{remark}

Cutting off two different surface diagrams can yield the same diagram:

\begin{remark}
    The right-hand surface diagram appearing in the Frobenius relations (\ref{LD-FrobeniusRelations}) is obtained by cutting off, in the manner described in Remark \ref{Cutting-off surface dgms AlgCoalg}, either the left-hand surface diagram of Frobenius relation \eqref{eq:F1 Frob LD} in Figure \ref{rel F1} or the lef-hand one in Figure \ref{rel F2}.
\end{remark}

\begin{remark}
    The Frobenius relations (\ref{LD-FrobeniusRelations}) of a Frobenius algebra $(A,\mu,\eta,\Delta,\epsilon)$ express that the comultiplication $\Delta\colon A\ra A\parLL A$ (equivalently, the multiplication $\mu\colon A\otimes A\ra A$) is a morphism of $A$-bimodules (equivalently, a morphism of $A$-bi-comodules). Note that the $A$-bimodule structure (resp. $A$-bicomodule structure) on the $\parLL$-monoidal product $A\parLL A$ (resp. $\otimes$-monoidal product $A\otimes A$) involves both distributors.
\end{remark}

\begin{remark}\label{category FrobC}
    Frobenius algebras in $\cC$ and their morphisms form a category $\mathsf{Frob}_{\cC}$. By definition, the category $\mathsf{Frob}_{\cC}$ is the hom-category $\mathsf{Frob}(\cI,\cC)$, where $\mathsf{Frob}$ denotes the $2$-category from Proposition \ref{LDC2Cat}. The three $2$-functors from Definition \ref{three2-functors} yield three 1-involutions $\mathsf{Frob}_{\cC}\cong\mathsf{Frob}_{\cC^{\operatorname{rev}}}$, $\mathsf{Frob}_{\cC}\cong\mathsf{Frob}_{\cC^{\operatorname{cop}}}$ and $\mathsf{Frob}_{\cC}\cong\mathsf{Frob}_{\cC^{\operatorname{lop}}}$. 
\end{remark}

By definition, the results on Frobenius LD-functors from Sections \ref{LD-cats and functors} and \ref{duality in LD-categories} specialize to Frobenius algebras:
\begin{remark}
    By Proposition \ref{LDC2Cat}, Frobenius LD-functors preserve Frobenius algebras. By Proposition \ref{FrobPresDuals}, since the single object in the unit category $\cI$ is LD-dualizable, Frobenius algebras are self-dual in the sense of Definition \ref{def:LDduals}. By Proposition \ref{PropMorphLDfuncIso}, the category $\mathsf{Frob}_{\cC}$ is a groupoid.
\end{remark}

Definition \ref{def:FrobAlg} generalizes the following two definitions found in the literature:

\begin{remark}\label{GV-Frobenius algebras}
    An LD-Frobenius algebra in an LD-category with negation is referred to as a \emph{GV-Frobenius algebra} in \cite{fuchs2024grothendieckverdier}. Further specializing, LD-Frobenius algebras in a symmetric GV-category, where all objects are rigid dualizable but the dualizing object is not necessarily chosen as the monoidal unit, are called \emph{(degree-$d$) Frobenius algebra objects} by Fu and Vial \cite[Def. 2.1]{MotivicGlobalTor}.
\end{remark}

Frobenius algebras in LD-categories naturally arise in applications:

\begin{examples}\label{example ld-frobenius algebras}
    \begin{enumerate}[label=(\roman*)]
        \item Trivializations of the partially defined relative Serre functor of a GV-module category over a GV-category $\cC$ yield LD-Frobenius algebras in $\cC$ \cite[Thm. 5.20]{fuchs2024grothendieckverdier}.
        \item The cohomology ring $H^{\ast}(M,\mathbb{Q})$ of a closed, connected, oriented real $d$-manifold $M$ is an LD-Frobenius algebra in the r-category of finite-dimensional graded $\mathbb{Q}$-vector spaces, with dualizing object the $1$-dimensional $\mathbb{Q}$-vector space sitting in degree $d$. The Frobenius structure arises from the fundamental class of $M$ via Poincaré duality \cite[Ex. 2.5]{MotivicGlobalTor}. 
        \item Similarly, if $M$ is a compact Kähler manifold, the cohomology ring $H^{\ast}(M,\mathbb{Q})$ can be endowed with a Frobenius algebra structure in the r-category of pure rational Hodge structures \cite[Ex. 2.6.]{MotivicGlobalTor}.
        \item The Chow motive of any $d$-dimensional smooth projective variety over a field $k$ is naturally a (commutative) LD-Frobenius algebra in the r-category of Chow motives \cite[Lemma 2.7]{MotivicGlobalTor}. In this example, the invertible dualizing object is the motive $\mathds{1}(1)^{-d}$, where $\mathds{1}(1)$ denotes the Tate motive.
    \end{enumerate}
\end{examples}

Similar to the usual Frobenius algebras in monoidal categories, the Frobenius relations for LD-Frobenius algebras are not independent. Put differently, cutting off surface diagrams can render certain relations redundant:

\begin{theorem}\label{FrobRelationsNotIndep}
    Let $A\in \cC$ be an object. Assume that $A$ is endowed with an algebra structure $(\mu,\eta)$ and a coalgebra structure $(\Delta,\epsilon)$ in the sense of Definition \ref{definition algebra-coalgebra}. Assume that the Frobenius relation (F1) in Equation (\ref{LD-FrobeniusRelations}) holds. Then the Frobenius relation (F2) from Equation (\ref{LD-FrobeniusRelations}) follows.
\end{theorem}

Before giving the proof of Theorem \ref{FrobRelationsNotIndep}, we explain our proof strategy, since it will be used frequently in the rest of this paper:

\begin{remark}\label{proof strategy}
    First, we recall the string-diagrammatic proof of the statement specialized to strict monoidal categories. In the LD-setting, we then extend the two-dimensional string diagrams into a third dimension to keep track of the coherence data. When viewed from the front face, the surface diagrams look like the ordinary string diagrams of the proof in the strict monoidal setting. We draw these ordinary string diagrams in blue below their surface diagrams. Finally, we fill in the missing steps in the proof in the LD-setting using the coherence axioms from Appendix \ref{coherenceLD}. The blue horizontal ‘pedal line’ below a surface diagram indicates that the string diagram extracted from the front face of the surface diagram is equal to the string diagram on the left of the pedal line (up to planar isotopy).
\end{remark}

\begin{proof}
By following the proof strategy from Remark \ref{proof strategy} we find the following graphical proof. It is a higher-dimensional version of the proof given in \cite[Lemma A1]{Street-Pastro-WeakHopf}:
    \begin{figure}[H]
    \centering
    \includegraphics[width=0.925\textwidth]{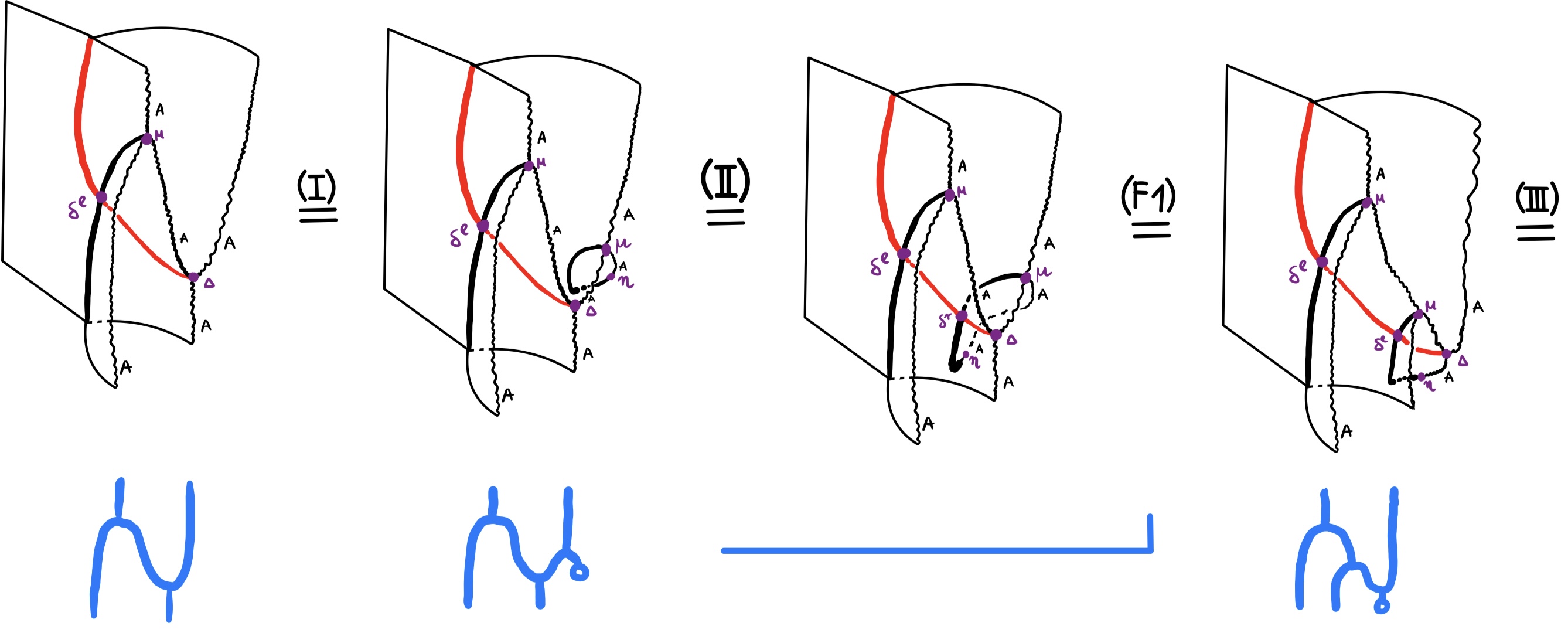}
\end{figure}
\begin{figure}[H]
    \centering
    \includegraphics[width=0.925\textwidth]{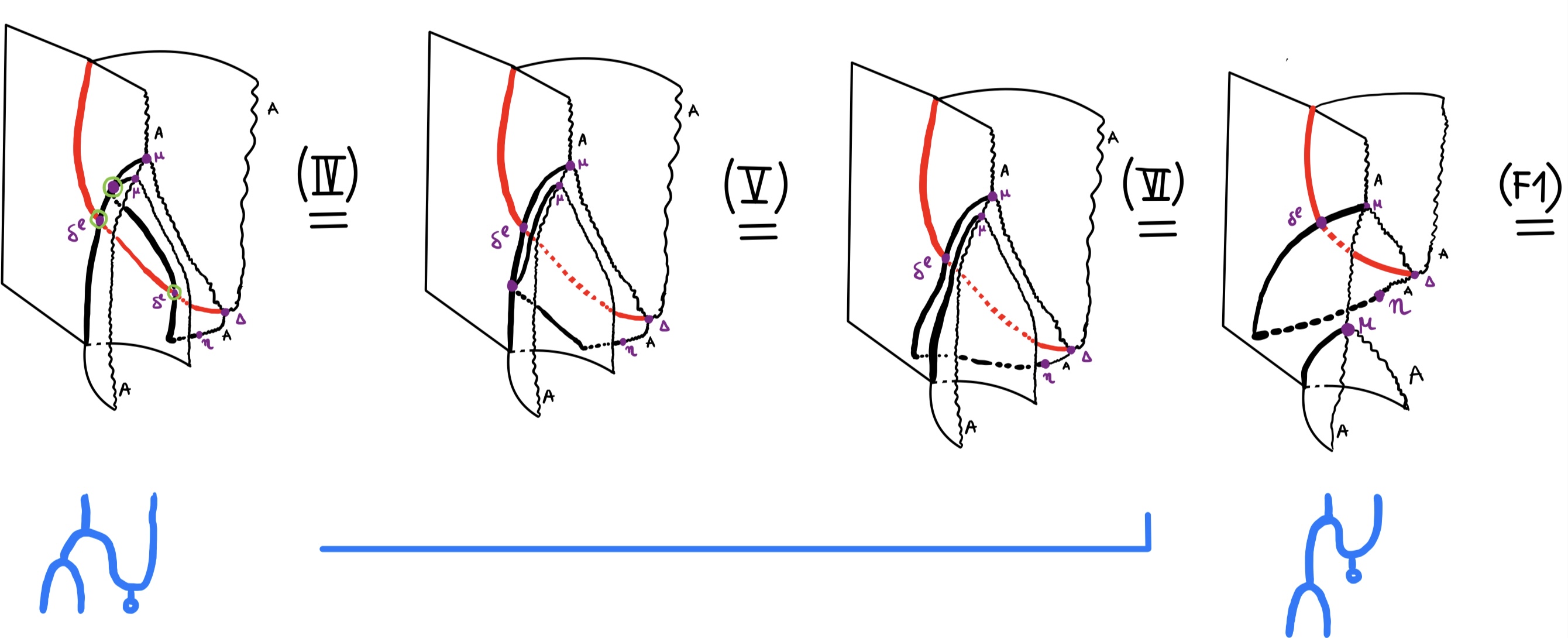}
\end{figure}
\begin{figure}[H]
    \centering
    \includegraphics[width=0.925\textwidth]{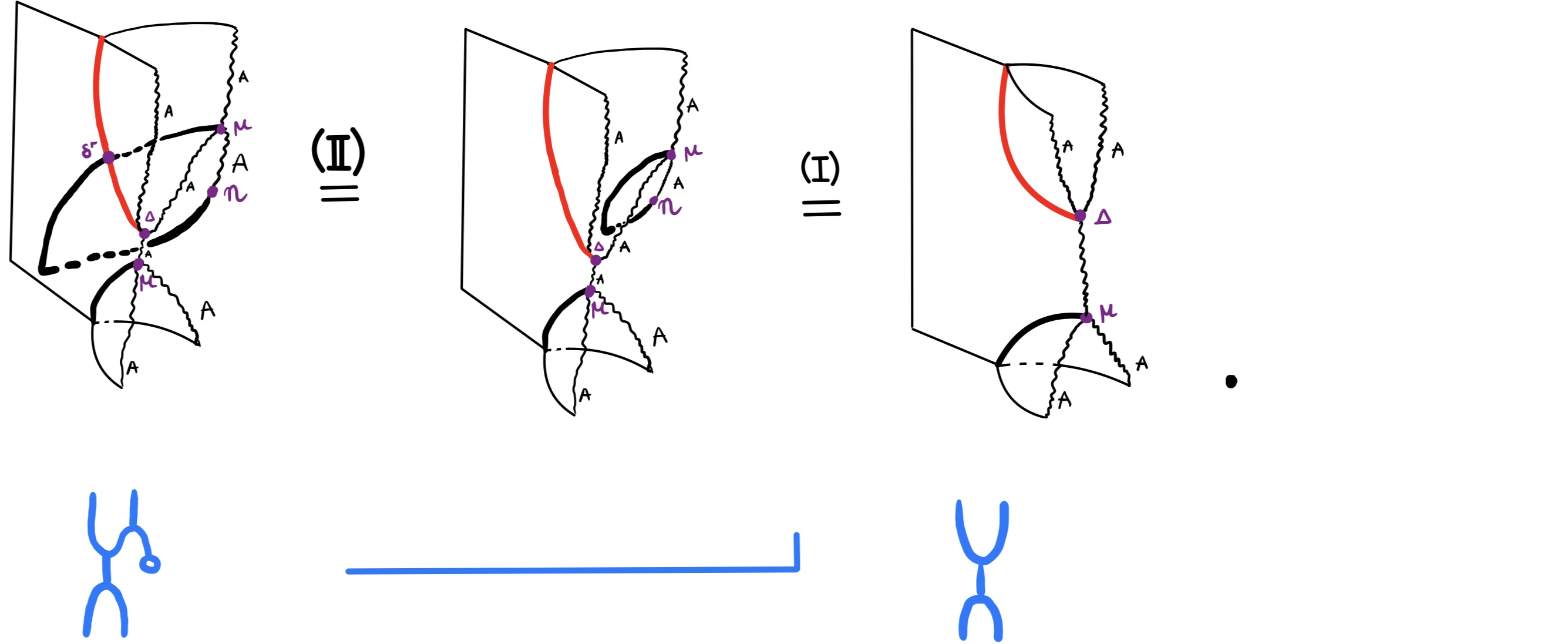}
\end{figure}
Equation (I) follows from the unitality of the multiplication as depicted in Figure \ref{unitality}, while Equation (II) holds by the triangle coherence axiom (\ref{eq:A2}). Equation (III) follows from the associativity of the multiplication $\mu$ as depicted in Figure \ref{Associativity}, while Equation (IV) holds by the pentagon coherence axiom (\ref{eq:A5}). Equation (V) uses the triangle diagram for the $\otimes$-monoidal product. Equation (VI) holds by the axioms of a strict monoidal $2$-category.
\end{proof}

\begin{remark}
Let $(A,\mu,\eta)$ be a unital, but not-necessarily associative, algebra, and $(A,\Delta,\epsilon)$ a counital, but not-necessarily coassociative, coalgebra in an LD-category $\cC$. Suppose that the quintuple $(A,\mu,\eta,\Delta,\epsilon)$ satisfies the Frobenius relations (\ref{LD-FrobeniusRelations}). Using the proof strategy outlined in Remark \ref{proof strategy}, one can then show that $(A,\mu,\eta)$ is necessarily associative and $(A,\Delta,\epsilon)$ is necessarily coassociative.
\end{remark}

\subsection{Frobenius forms}
In this section, we present an equivalent characterization of Frobenius algebras using the following notion:

\begin{definition}
 A \emph{form} on an object $X\in \cC$ is a morphism $\lambda\in \operatorname{Hom}_{\cC}(X,K)$ with target the $\parLL$-monoidal unit $K\in \cC$.
\end{definition}

\begin{definition}\label{def:frobenius morphism}
    Let $(A,\mu,\eta)$ be an algebra in $\cC$. Assume that $A$ is right LD-dualizable. We call a form $\lambda\colon A\ra K$ \emph{(right) Frobenius} if the morphism
    \begin{equation*}
        \psi^r_\lambda\colon\; A\,\longrightarrow \,{^{\ast}A}
    \end{equation*} defined by
    \begin{equation}\label{psir}
    \psi^r_\lambda\;:=\;r_{{^{\ast}A}}^{\parLL}\circ({^{\ast}A}\parLL \lambda)\circ({^{\ast}A}\parLL \mu)\circ\distr_{{^{\ast}A},A,A}\circ(\operatorname{coev}'_A\otimes A)\circ (l^{\otimes}_A)^{-1}
    \end{equation}
is invertible. Assuming that $A$ is left LD-dualizable, the form $\lambda\colon A \ra K$ is called \emph{(left) Frobenius} if it is right Frobenius in $\cC^{\text{rev}}.$
\end{definition}

\begin{figure}[H]
    \centering
    \begin{subfigure}[b]{0.48\textwidth}
    \centering
    \includegraphics[width=0.65\textwidth]{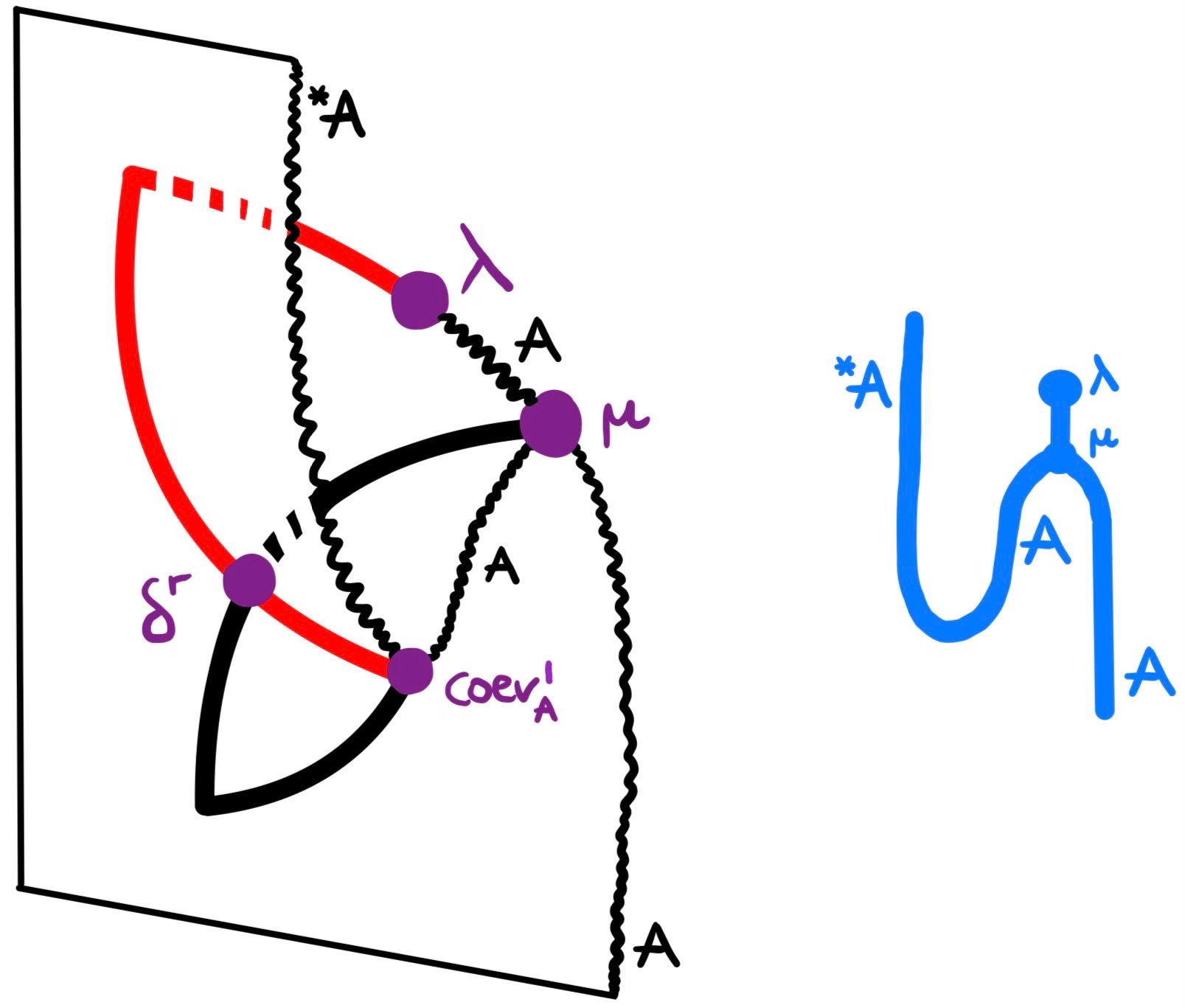}
    \caption{The morphism $\psi^r_{\lambda}$.}
    \label{fig:Frobenius morphism}
    \end{subfigure}
\hfill
    \begin{subfigure}[b]{0.48\textwidth}
    \centering
    \includegraphics[width=0.8\textwidth]{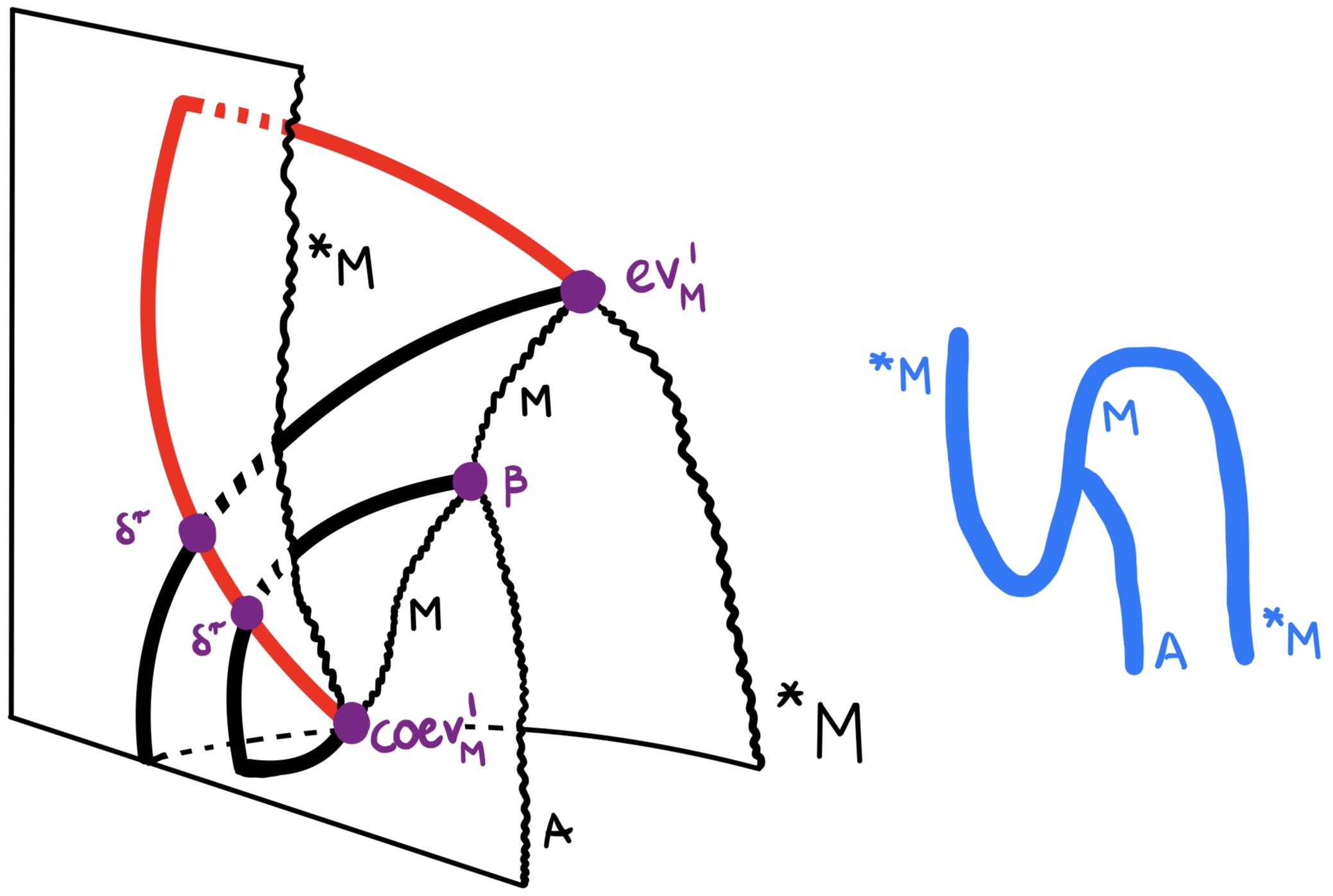}
    \caption{The morphism ${\beta^{\prime}}$.}
    \label{fig:module on astA}
    \end{subfigure}
\caption{Surface diagrams for Definition \ref{def:frobenius morphism} and Lemma \ref{dual of algebra is module}.}
\end{figure}

\begin{remark}
    Figure \ref{fig:Frobenius morphism} shows the surface diagram of the morphism $\psi^r_\lambda$ from Equation (\ref{psir}). 
\end{remark}

\begin{remark}
    Although the morphism $\psi^r_\lambda$ depends on the choice of coevaluation, its invertibility is independent of this choice.
\end{remark}

For later use, we record two lemmas. These lemmas are well known (and easy to prove) when specializing to rigid duals in monoidal categories. The proof in the LD-setting highlights the extent to which the coherence data and their relations are at work in the background.

\begin{lemma}\label{dual of algebra is module}
    Let $(A,\mu,\eta)$ be an algebra in $\cC$. Let $(M,\beta)$ be a right $A$-module. Assume that the underlying object $M$ is right LD-dualizable. Then the morphism ${\beta^{\prime}}\colon A\otimes \prescript{\ast}{}{M}\ra \prescript{\ast}{}{M}$ from Figure \ref{fig:module on astA} endows the right LD-dual $\prescript{\ast}{}{M}$ with a left $A$-module structure.
\end{lemma}

\begin{proof}
Following the proof strategy from Remark \ref{proof strategy}, we graphically prove the compatibility of the action ${\beta^{\prime}}$ with the multiplication as follows:

\smallskip

\bigskip

\begin{figure}[H]
    \centering
    \includegraphics[width=\textwidth]{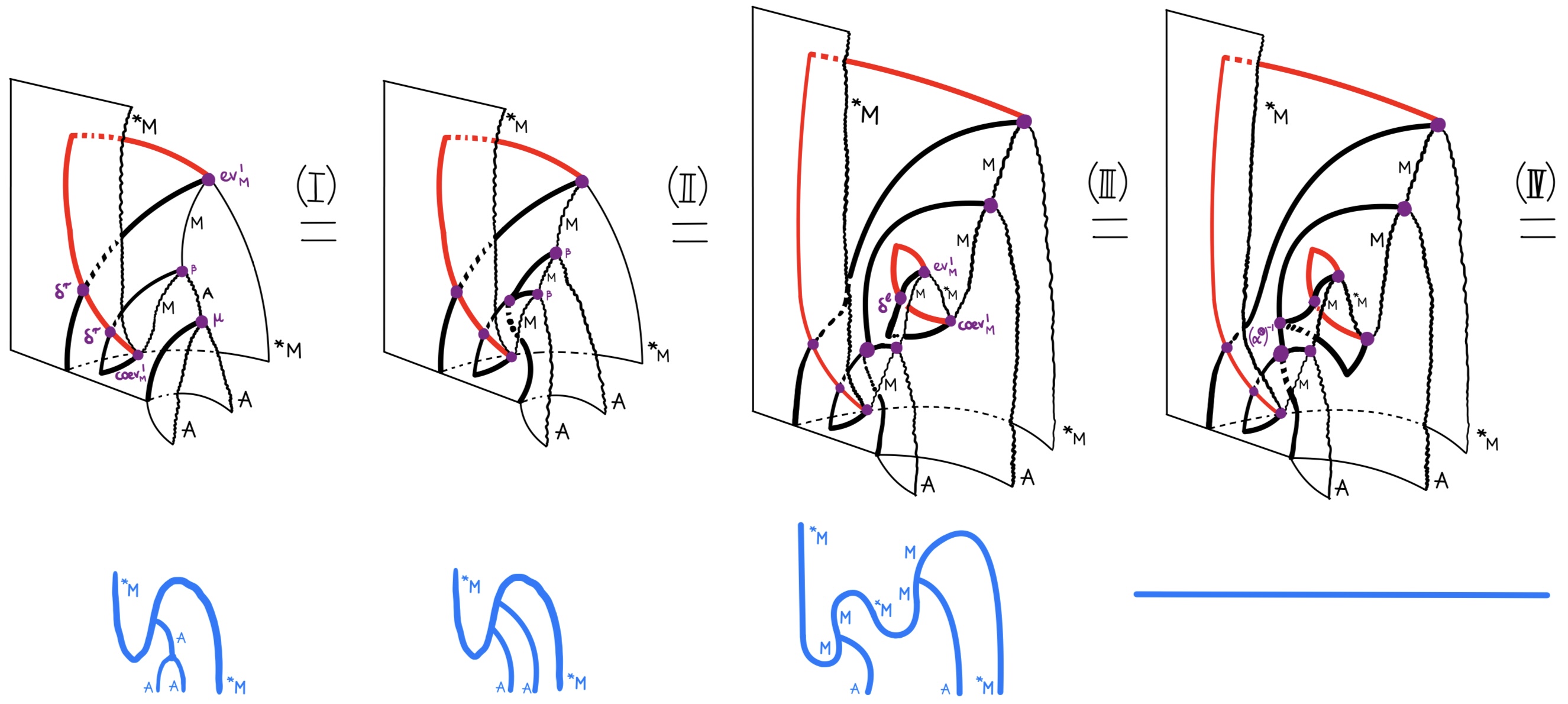}
\end{figure}

\bigskip

\smallskip

\begin{figure}[H]
    \centering
    \includegraphics[width=\textwidth]{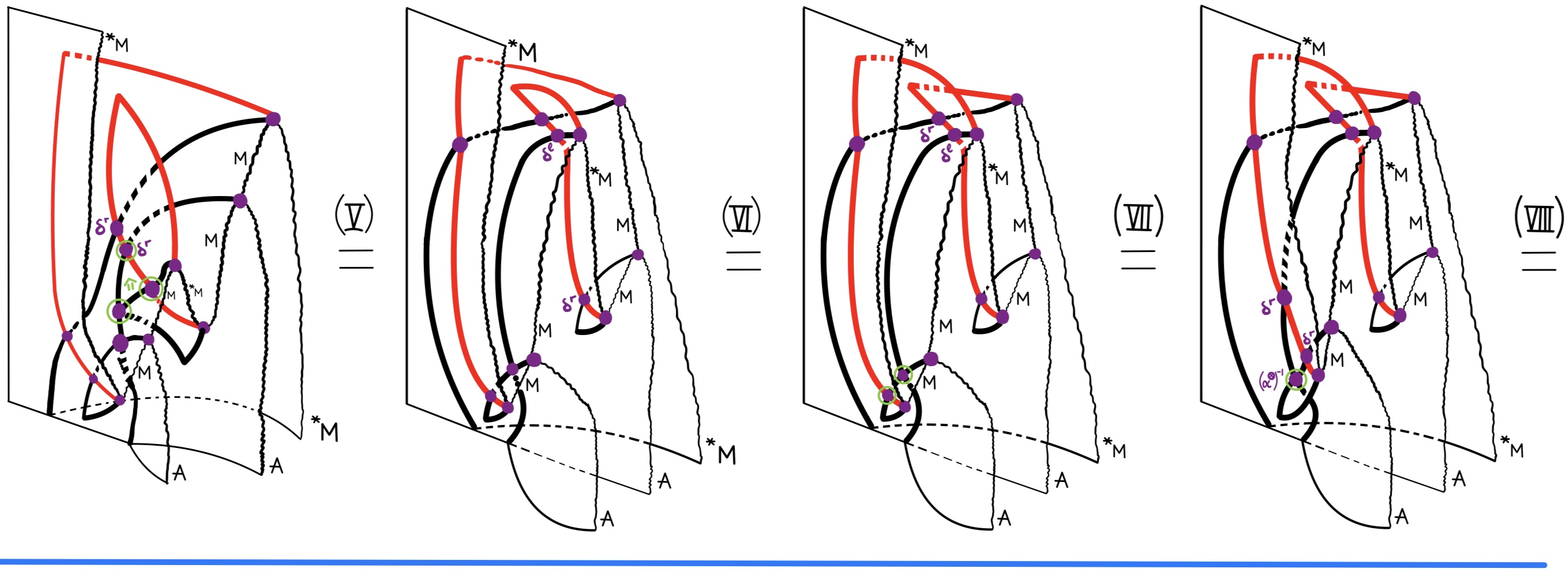}
\end{figure}

\begin{figure}[H]
    \centering
    \includegraphics[width=0.97\textwidth]{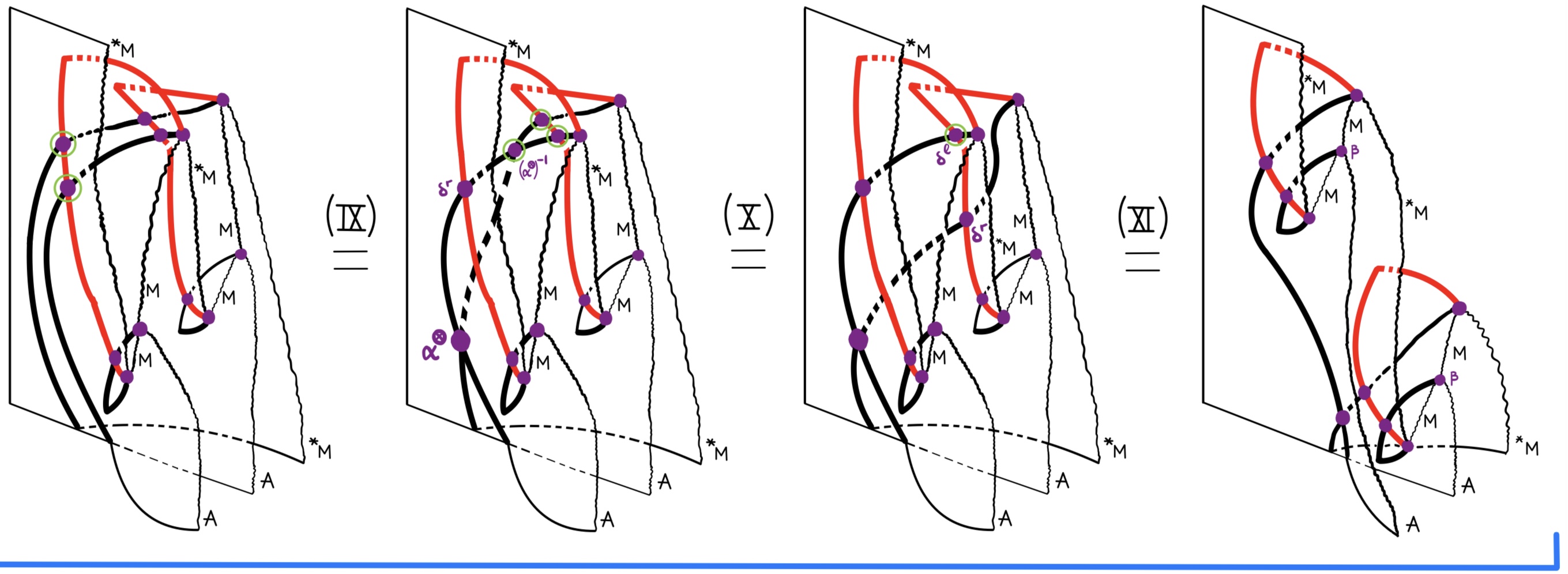}.
\end{figure}

Equation (I) holds since the action $\beta$ is compatible with the multiplication $\mu$; compare Figure \ref{Associativity}. Equation (II) follows from the snake equation (\ref{firstzigzag}); see Figure \ref{fig:LD-snake equations (S1) and (S2)}. Mac Lane's coherence theorem for the $\otimes$-monoidal structure implies Equation (III); compare Figure \ref{fig:MacTriangle}. Equation (IV) follows from the triangle coherence axiom (\ref{eq:A3}). Equation (V) holds by the pentagon coherence axiom (\ref{eq:A9}). As before, we have colored those vertices light green, which are involved in the coherence axiom. We also indicated the direction in which we move a sheet. Equation (VI) follows from an application of Mac Lane's coherence theorem for the $\parLL$-monoidal structure to the top of our surface diagram. Equation (VII) holds by the pentagon coherence axiom (\ref{eq:A6}), while Equation (VIII) is implied by Mac Lane's coherence theorem for the $\otimes$-monoidal structure. Equations (IX) and (X) are both applications of the pentagon coherence axiom (\ref{eq:A6}). Finally, Equation (XI) follows from the triagle coherence axiom (\ref{eq:A4}) and the axioms of a strict monoidal $2$-category.

\smallskip

The action is also unital:
\begin{figure}[H]
    \centering
    \includegraphics[width=0.94\textwidth]{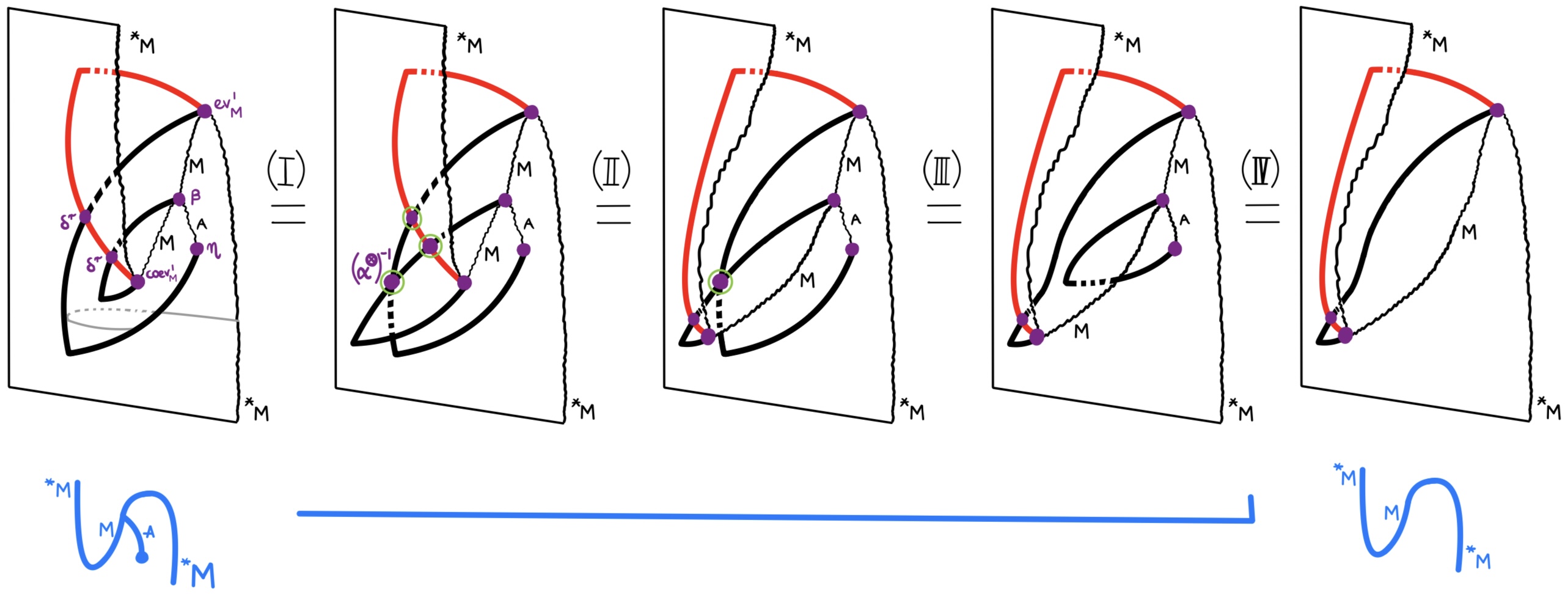}.
\end{figure}
\vspace{-0.2cm}
Equations (I) and (III) follow from Mac Lane's coherence theorem for the $\otimes$-monoidal structure. Equation (II) holds by the pentagon coherence axiom (\ref{eq:A6}), while Equation (IV) holds by the unitality of multiplication; see Figure \ref{unitality}. Finally, we apply snake equation (\ref{secondzigzag}) from Figure \ref{fig:LD-snake equations (S1) and (S2)} to find that the right-hand surface diagram is equal to the identity on the LD-dual $\prescript{\ast}{}{M}$.
\end{proof}

\begin{remark}
    By applying Lemma \ref{dual of algebra is module} to the LD-categories $\cC^{\text{rev}}$, and to $\cC^{\text{cop}}$ and $\cC^{\text{lop}}$, one obtains analogous statements for the left dual of $M$, and for left and right comodules, respectively. 
\end{remark}

\begin{lemma}\label{psi morphism of modules}
    Let $(A,\mu,\eta)$ be an algebra in $\cC$. Assume that $A$ is right LD-dualizable. Let $\lambda\colon A\ra K$ be a form on $A$. Then the morphism $\psi^r_\lambda\colon \, A\,\longrightarrow \,{^{\ast}A}$ from Definition \ref{def:frobenius morphism} is a morphism of left $A$-modules. Here, we consider the left $A$-module structure on $\prescript{\ast}{}{A}$ from Lemma \ref{dual of algebra is module}.
\end{lemma}

\begin{proof}
We again use the proof strategy outlined in Remark \ref{proof strategy}:
\begin{figure}[H]
    \centering
    \includegraphics[width=\textwidth]{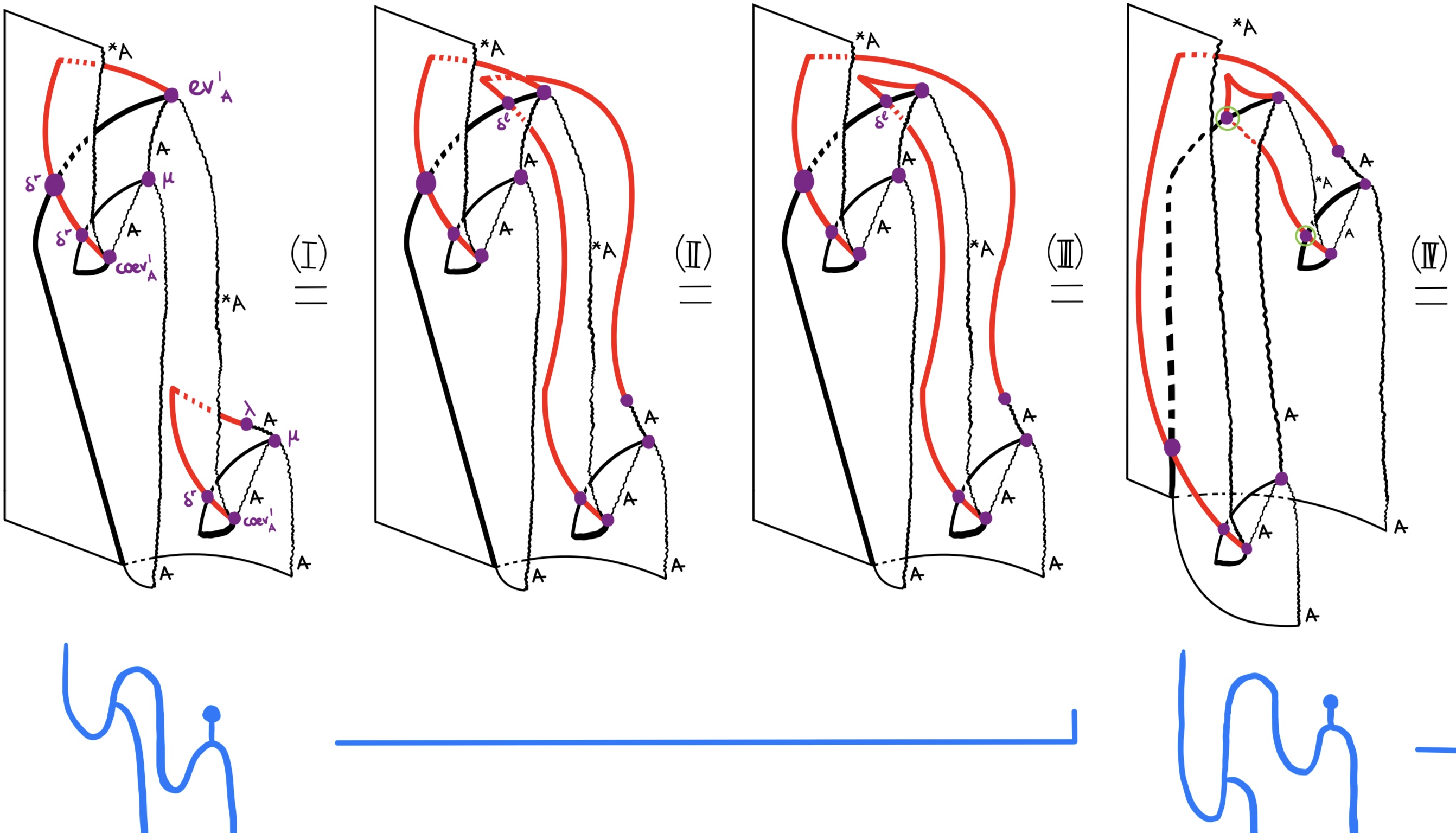}
\end{figure}

\begin{figure}[H]
    \centering
    \includegraphics[width=\textwidth]{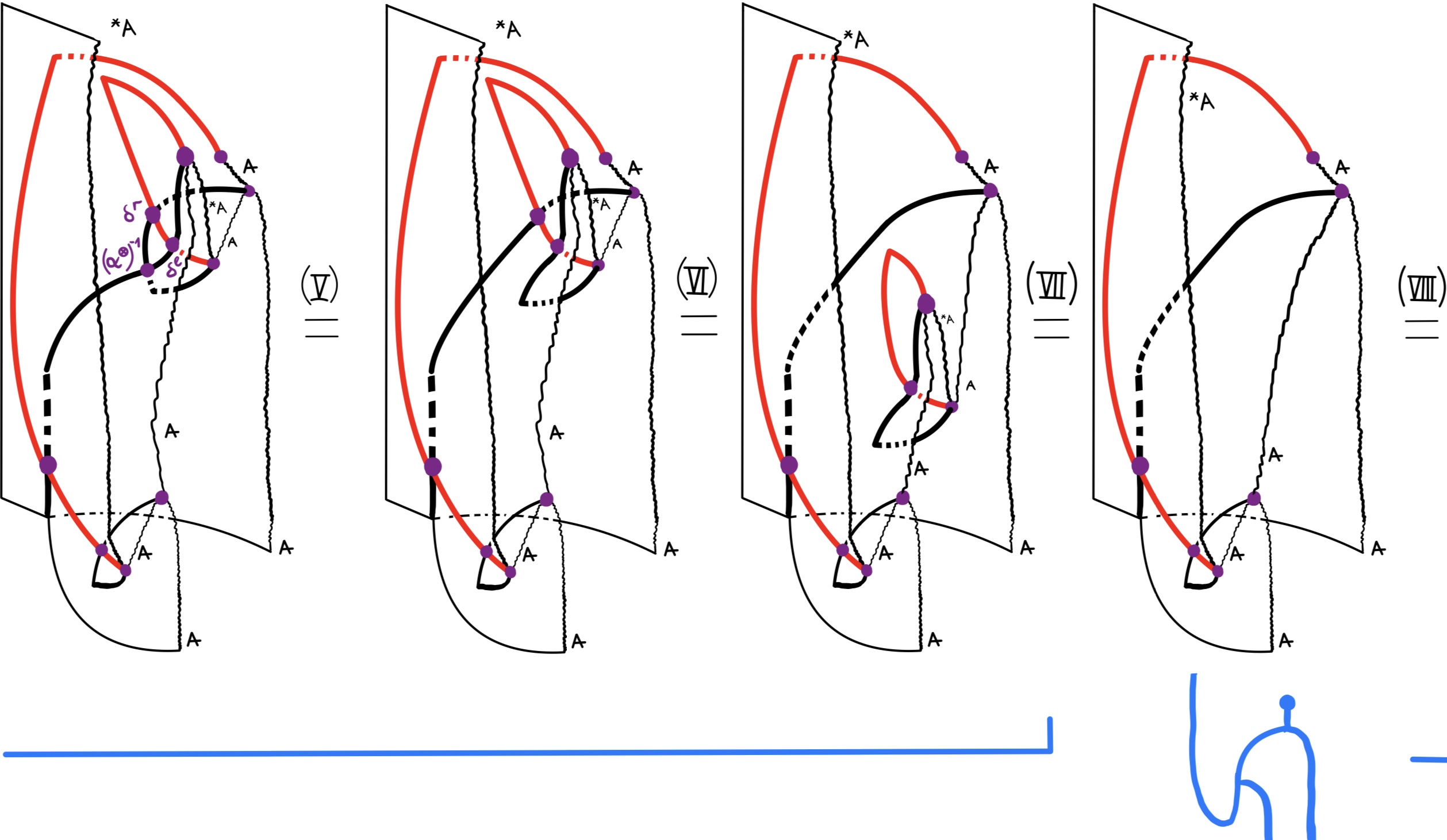}
\end{figure}
\begin{figure}[H]
    \centering
    \includegraphics[width=0.99\textwidth]{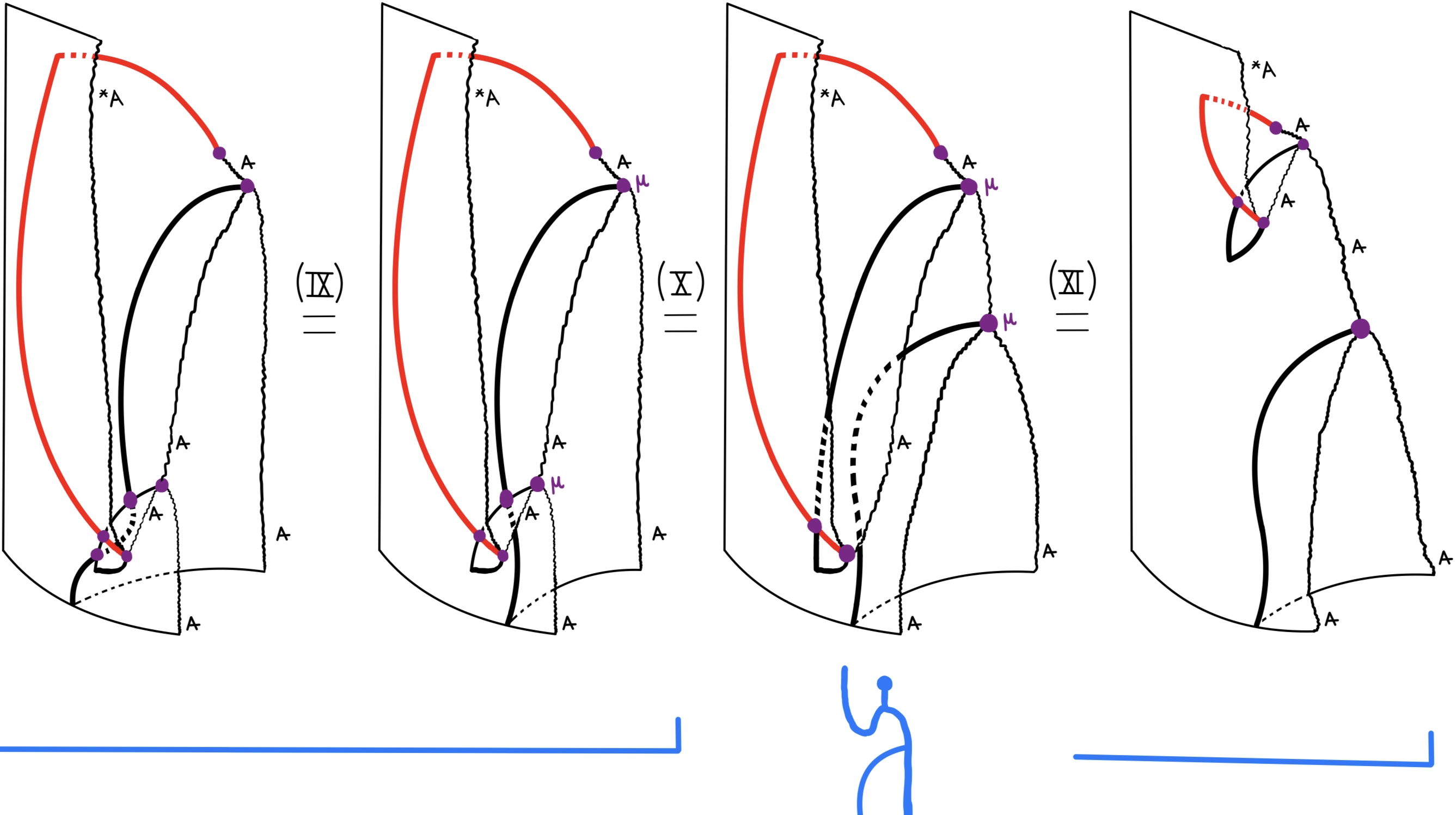}.
\end{figure}
%Kürzungspotential
Equation (I) follows from the triangle coherence axiom (\ref{eq:A4}), while Equation (II) is an application of Mac Lane's coherence theorem for the $\parLL$-monoidal structure to the top of our surface diagram. Equation (III) follows from the axioms of a strict monoidal $2$-category, while equation (IV) holds by the pentagon coherence axiom (\ref{eq:A9}). As in the proof of Lemma \ref{dual of algebra is module}, we have colored those vertices light green, which are relevant for the coherence axiom. Equation (V) follows from Mac Lane's triangle axiom for the $\otimes$-monoidal structure; see Figure \ref{fig:MacTriangle}. Equation (VI) is an application of the triangle coherence axiom (\ref{eq:A3}). Equation (VII) follows from the snake equation (\ref{firstzigzag}); see Figure \ref{fig:LD-snake equations (S1) and (S2)}. Equation (VIII) holds by the pentagon coherence axiom (\ref{eq:A6}), while Equation (IX) follows from Mac Lane's coherence theorem for the $\otimes$-monoidal structure. The associativity of the multiplication $\mu$ from Figure \ref{Associativity} implies Equation (X). Finally, Equation (XI) follows from the axioms of a strict monoidal $2$-category.
\end{proof}

For the main result of this subsection, we also need the following two definitions:
 
\begin{definition}\label{def:ideal}(\cite[Def. 2.7.]{walton2022filtered}).
    Let $A=(A,\mu,\eta)$ be an algebra in a monoidal category $\cC$. A \emph{weak left ideal} of $A$ consists of a left $A$-module $(I,\omega\colon A \otimes I\ra I)$ together with a morphism of left $A$-modules $\varphi\colon I\ra A$. A weak left ideal such that the morphism $\varphi$ is monic is called a \emph{left ideal}. \emph{(Weak) right ideals} are (weak) left ideals in the monoidal category $\cC^{\text{rev}}.$
\end{definition}

In other words, a left ideal of $A$ is a left $A$-submodule of $A$. In the case of algebras over commutative rings, Definition \ref{def:ideal} thus reduces to the standard notion.

\begin{definition}\label{def:abelian LD-cat}
    An \emph{abelian LD-category} is an LD-category whose underlying category is abelian and whose two monoidal products are additive bifunctors.    
\end{definition}

We come to the main theorem in this subsection. It generalizes a result of Walton and Yadav \cite[Thm. 5.3]{walton2022filtered} from algebras in abelian rigid monoidal categories to LD-dualizable algebras in abelian LD-categories. As we shall see in Corollary \ref{frobalgebra via absence of ideals}, the following theorem together with results in \cite{fuchs2024grothendieckverdier} characterizes Frobenius algebras by the absence of non-trivial one-sided ideals.

\begin{theorem}\label{Frobenius via ideals}
    Let $(A,\mu,\eta)$ be an algebra in an abelian LD-category $\cC$. Assume that $A$ is right LD-dualizable. For a form $\lambda\colon A \ra K$ on $A$, the following are equivalent:
    \begin{enumerate}[label=(\roman*)]
        \item The form $\lambda$ is right Frobenius.
        \item For every left \emph{weak} ideal $(I,\omega,\varphi)$ of $A$ that factors through the kernel $\operatorname{ker}(\lambda)$, the morphism $\varphi\colon I \ra A$ is the zero morphism.
        \item For every left ideal $(I,\omega,\varphi)$ of $A$ that factors through the kernel $\operatorname{ker}(\lambda)$, the morphism $\varphi\colon I \ra A$ is the zero morphism.
    \end{enumerate}
\end{theorem}

For the proof of Theorem \ref{Frobenius via ideals}, we will use surface diagrams. We let $(\cV,\otimes,J)$ be the symmetric monoidal category of abelian groups with standard tensor product, and work internally to the strict monoidal $2$-category $\cV{\textnormal-}\mathsf{Cat}$ of preadditive categories.

\begin{proof}
To prove that $(i)$ implies $(ii)$, let $(I,\omega,\varphi)$ be a weak left ideal of $A$ that factors through $\operatorname{ker}(\lambda)$. By the definition of the morphism $\psi_{\lambda}^r\colon A\ra {^{\ast}A}$, as given in Definition \ref{def:frobenius morphism}, we have:
    \begin{equation}
        {\psi_\lambda^r}\circ\varphi \;=\;    {r^{\parLL}_{^{\ast}A}}\circ\big({^{\ast}A}\parLL\, { (\lambda\circ\mu)}\big)\circ\distr_{{^{\ast}A},A,A}\circ({\operatorname{coev}_A^{\prime}}\otimes  A)\circ(l^{\otimes}_A)^{-1}\circ\varphi.
    \end{equation}

Next, we apply the naturality of the left unitor and the right distributor, along with the assumption that $\varphi$ is a morphism of left $A$-modules, which gives:
    \begin{equation}\label{eq: naturality dist unitor varphi}
        {\psi_\lambda^r}\circ\varphi \, = \,{r^{\parLL}_{^{\ast}A}}\circ\big({^{\ast}A}\parLL\, (\lambda\circ \varphi \circ\omega)\big)\circ\distr_{{^{\ast}A},A,I}\circ({\operatorname{coev}_A^{\prime}}\otimes  I)\circ(l^{\otimes}_I)^{-1}.
    \end{equation}

Denote by $\iota\colon \operatorname{ker}(\lambda) \hookrightarrow A$ the canonical monomorphism. By assumption, there exists a morphism $\phi\colon I\ra \operatorname{ker}(\lambda)$ such that $\varphi=\iota\circ \phi.$ Substituting $\varphi=\iota\circ \phi$ into Equation (\ref{eq: naturality dist unitor varphi}), we obtain:
    \begin{equation}
        {\psi_\lambda^r}\circ\varphi \, = \,{r^{\parLL}_{^{\ast}A}}\circ\big({^{\ast}A}\parLL\, (\lambda\circ \iota\circ\phi \circ\omega)\big)\circ\distr_{{^{\ast}A},A,I}\circ({\operatorname{coev}_A^{\prime}}\otimes  I)\circ(l^{\otimes}_I)^{-1}.
    \end{equation}

Finally, since the composite $\lambda\circ \iota$ is by definition zero, it follows from the additivity of the monoidal product $\parLL$ that the composite ${\psi_\lambda^r}\circ\varphi$ is zero. Since $\psi_{\lambda}^r$ is assumed to be invertible, this implies that $\varphi$ must be the zero morphism. Therefore, (ii) follows from (i).

\medskip

\nid Statement (ii) implies (iii) by definition. 

\medskip

\nid For the last implication, assume that (iii) holds. Recall from Lemma \ref{psi morphism of modules} that the morphism $\psi_{\lambda}^r\colon A\ra {^{\ast}A}$ is a morphism of $A$-modules. By standard results, the kernel $I:=\operatorname{ker}(\psi_\lambda^r)$ can thus be given a left $A$-module structure such that the canonical monomorphism $\varphi\colon I \to A$ is a morphism of left $A$-modules. 

To demonstrate that this ideal $I$ factors through $\operatorname{ker}(\lambda)$, we compute the following:
\begin{figure}[H]
    \centering
    \includegraphics[width=0.65\textwidth]{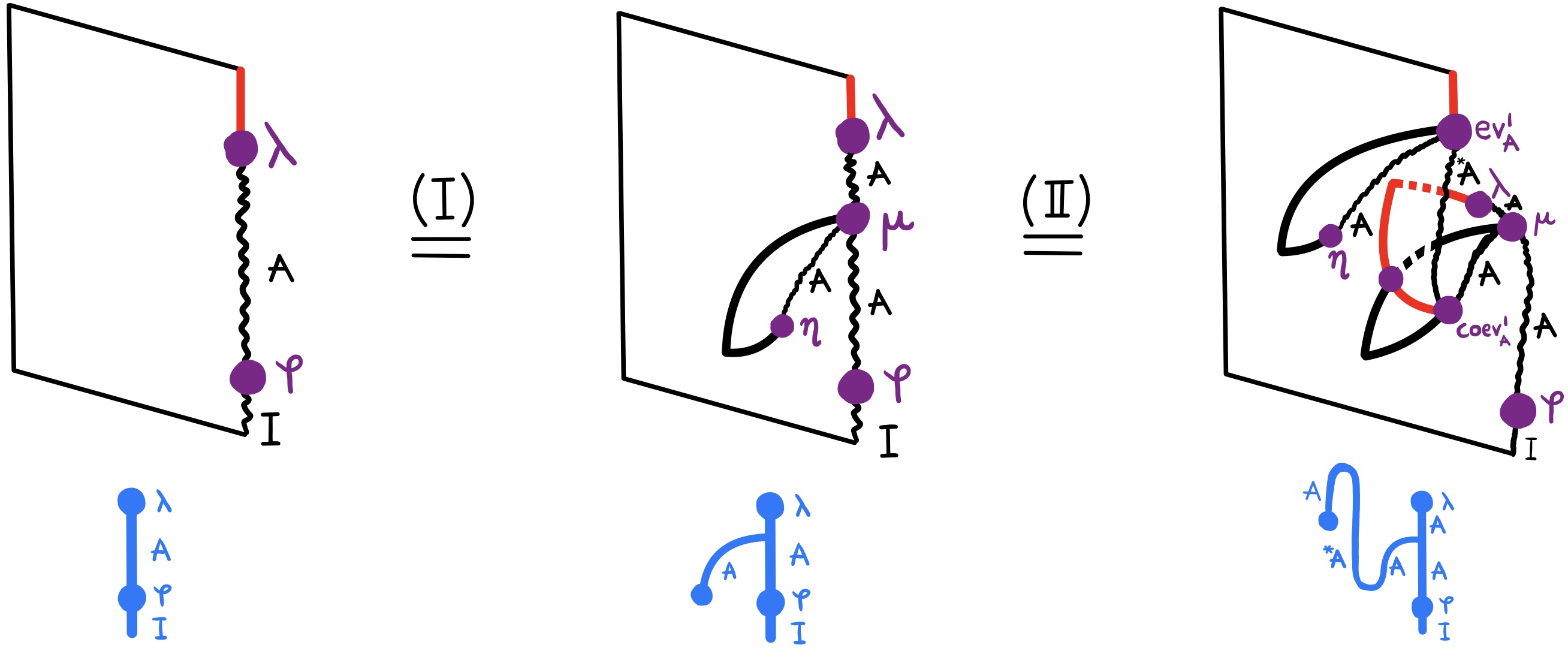}.
\end{figure}
Equation (I) follows from the unitality of the multiplication morphism $\mu$ of $A$, as shown in Figure \ref{unitality}. Equation (II) is an application of Lemma \ref{computation-heavy lemma LDGV} to the LD-category $\cC^{\text{rev}}$, where we set $\gamma:=\lambda\circ \mu$, $\kappa:=\operatorname{ev}_A^{\prime}$, $\overline{\kappa}:=\operatorname{coev}_A^{\prime}$ and $f:=\operatorname{id}_A$. 

Since $\psi^r_\lambda\circ\varphi\;=\;0$ by the definition of $\varphi$, the composite morphism \begin{equation*}
{\operatorname{ev}_A^{\prime}\circ \;{(\eta\otimes {^{\ast}A})}\circ l^{-1}_{{^{\ast}A}}\circ\psi^r_\lambda\circ\varphi}
\end{equation*}
represented by the rightmost surface diagram is the zero morphism. This implies that $\lambda\circ\varphi$ is also zero, which means that the ideal $I$ factors through $\operatorname{ker}(\lambda)$. By assumption (iii), it follows that $\varphi$ must be the zero morphism. Since $\varphi$ is monic, we deduce that the kernel $I$ is zero.

Analogously, it can be proved that the cokernel $\operatorname{coker}(\psi^r_\lambda)$ is zero. Since $\cC$ is an abelian category, this shows that the form $\lambda$ is right Frobenius. 
\end{proof}

\begin{remark}
    By applying Theorem \ref{Frobenius via ideals} to the LD-category $\cC^{\text{rev}}$, we find an analogous characterization of \emph{left} Frobenius forms via \emph{right} (weak) ideals.
\end{remark}

Let us discuss the relationship between Frobenius forms and Frobenius algebras. To do so, recall the notions of LD-pairing and LD-copairing from Definition \ref{pairing/copairing,side-inverse}.

\begin{definition}
An LD-pairing $\kappa$ on an object $X\in \cC$ is called \emph{non-degenerate} if there exists an LD-copairing $\overline{\kappa}$ on $X$ that is side-inverse to $\kappa$. An LD-pairing $\kappa$ on an algebra $(A,\mu,\eta)$ in $\cC$ is called \emph{invariant} if we have
    \begin{equation*}
        \kappa\circ(\mu \otimes A)\;=\;\kappa\circ(A\otimes \mu)\circ \alpha^{\otimes}_{A,A,A}.
    \end{equation*}
\end{definition}

\begin{remark}
An invariant LD-pairing $\kappa$ on an algebra $(A,\mu,\eta)$ is automatically compatible with the unit $\eta$, that is, we have 
    \begin{equation*}
        \kappa\circ(\eta\otimes A)\circ l^{\otimes}_A\;=\;\kappa\circ(A\otimes \eta)\circ r^{\otimes}_A.
    \end{equation*}
This is a purely $\otimes$-monoidal statement and can, therefore, be proved using string diagrams.
\end{remark}

We define two categories $\mathsf{Alg}^{\kappa}_{\cC}$ and $\mathsf{Alg}^{\lambda}_{\cC}$:

\begin{definition}\label{category Alg}
     The objects of $\mathsf{Alg}^{\kappa}_{\cC}$ are pairs $(A,\kappa)$ consisting of an algebra $A$ in $\cC$ and an invariant non-degenerate LD-pairing $\kappa$ on $A$. A morphism $f\colon (A,\kappa_A)\ra (B,\kappa_B)$ in $\mathsf{Alg}^{\kappa}_{\cC}$ is a unital algebra morphism $f\colon A\ra B$ such that $\kappa_B\circ(f\otimes f)=\kappa_A.$
\end{definition}

\begin{definition}\label{category AlgLambda}
    The objects of $\mathsf{Alg}^{\lambda}_{\cC}$ are pairs $(A,\lambda)$ consisting of a right LD-dualizable algebra $A$ in $\cC$ and a Frobenius form $\lambda\colon A \ra K.$ A morphism $f\colon (A,\lambda_A)\ra (B,\lambda_B)$ in $\mathsf{Alg}^{\lambda}_{\cC}$ is an algebra morphism $f\colon A \ra B$ such that $\lambda_B\circ f=\lambda_A.$
\end{definition}

We recall the three equivalent characterizations of Frobenius algebras from \cite{fuchs2024grothendieckverdier}:

\begin{theorem}\label{equivalent characterizations LD-Frobenius algebras}{\ \\}\vspace{-0.7cm}
    \begin{enumerate}
        \item For the LD-category $\cC$, the following two groupoids are equivalent:
        \begin{enumerate}[label=(\roman*)]
            \item The category $\mathsf{Frob}_{\cC}$ as described in Remark \ref{category FrobC}.
            \item The category $\mathsf{Alg}^{\kappa}_{\cC}$ as described in Definition \ref{category Alg}.
        \end{enumerate}
        \item If $\cC$ is an LD-category with negation, the groupoids from Part 1. are both equivalent to the category $\mathsf{Alg}^{\lambda}_{\cC}$ as described in Definition \ref{category AlgLambda}.
    \end{enumerate}
\end{theorem}

\begin{proof}
Part 1. of Theorem \ref{equivalent characterizations LD-Frobenius algebras} is proved in \cite[Prop. 4.11]{fuchs2024grothendieckverdier}. Although \cite[Prop. 4.11]{fuchs2024grothendieckverdier} is formulated in the setting of GV-categories, the same proof goes through for LD-categories. Part 2. of Theorem \ref{equivalent characterizations LD-Frobenius algebras} is shown in \cite[Prop. 4.13]{fuchs2024grothendieckverdier}.
\end{proof}

Theorems \ref{equivalent characterizations LD-Frobenius algebras} and \ref{Frobenius via ideals} together immediately imply the following result:

\begin{corollary}\label{frobalgebra via absence of ideals}
If $\cC$ is an LD-category with negation, the following data on an algebra $(A,\mu,\eta)$ in $\cC$ are equivalent:
 \begin{enumerate}[label=(\roman*)]
    \item A coalgebra structure $(\Delta,\epsilon)$ on $A$ such that $(A,\mu,\eta,\Delta,\epsilon)$ is a Frobenius algebra in $\cC$.
    \item A form $\lambda\colon A\ra K$ such that, for every left ideal $(I,\omega,\varphi)$ of $A$ that factors through the kernel $\operatorname{ker}(\lambda)$, the morphism $\varphi$ is the zero morphism.
    \item A form $\lambda\colon A\ra K$ such that, for every right ideal $(I,\omega,\varphi)$ of $A$ that factors through the kernel $\operatorname{ker}(\lambda)$, the morphism $\varphi$ is the zero morphism.
\end{enumerate}
\end{corollary}

%%%%%%%%%%%%%%%%%%%%%%%%%%%%%%%%%%

\subsection{The category of (co)modules}
In this subsection, we show that given an LD-Frobenius algebra $A\in \cC$, the category ${_A\cC}$ of left $A$-modules is isomorphic to the category ${^A\cC}$ of left $A$-comodules. This result is well known for monoidal categories; e.g. see \cite{Abrams_ModulesComodules} for the monoidal category of vector spaces. 

The statement for a general LD-category $\cC$ is non-obvious since the module and comodule structures of an LD-Frobenius algebra are defined with respect to different monoidal products. Consequently, distributors appear in our proof, which makes it necessary to deal with the coherence axioms from Appendix \ref{coherenceLD}.

\begin{figure}[H]
     \centering
     \begin{subfigure}[b]{0.48\textwidth}
         \centering
         \includegraphics[width=0.7\textwidth]{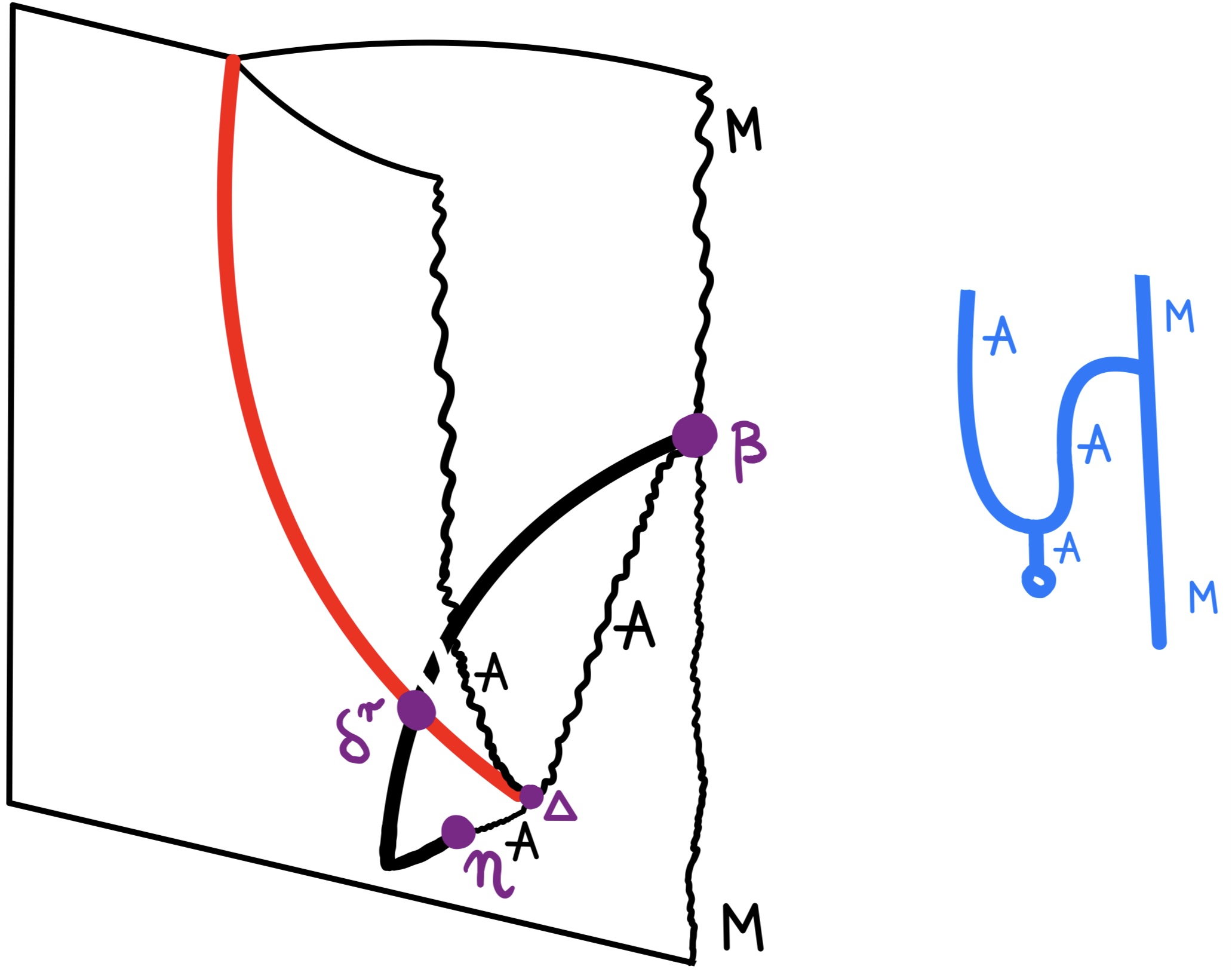}
         \caption{The morphism $\overline{\beta}$ from Lemma \ref{module-induces-comodule}.}
         \label{fig:overlineBeta}
     \end{subfigure}
\hfill
     \begin{subfigure}[b]{0.48\textwidth}
         \centering
         \includegraphics[width=0.7\textwidth]{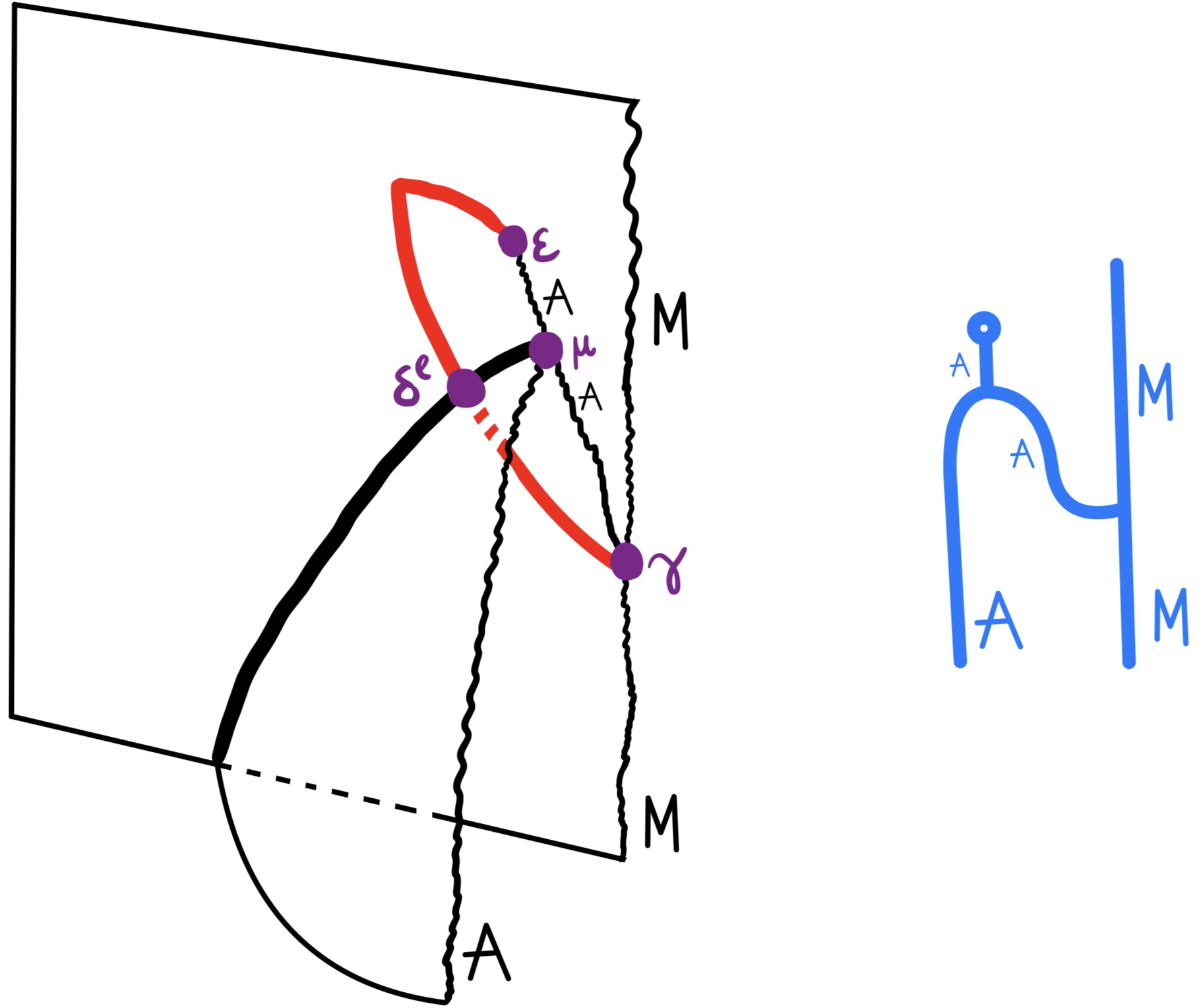}
         \caption{The morphism $\overline{\gamma}$ from Lemma \ref{comodule-induces-module}.}
         \label{fig:overlineGamma}
     \end{subfigure}
        \caption{Induced module and comodule structures.}
\end{figure}

We split our proof into various lemmas.

\begin{lemma}\label{module-induces-comodule}
    Let $A$ be a Frobenius algebra in the LD-category $\cC$. Given an $A$-module ${(M,\beta)\in{_A\cC}}$, the morphism $\overline{\beta}\colon M \ra A \parLL M$ from Figure \ref{fig:overlineBeta} endows $M$ with the structure of an $A$-comodule. This assignment yields a functor $F\colon _A \cC\ra {^A \cC}$. The functor $F$ commutes with the forgetful functors to $\cC$.
\end{lemma}

\begin{proof}The morphism $\overline{\beta}$ is compatible with the comultiplication:

\begin{figure}[H]
    \centering
    \includegraphics[width=0.96\textwidth]{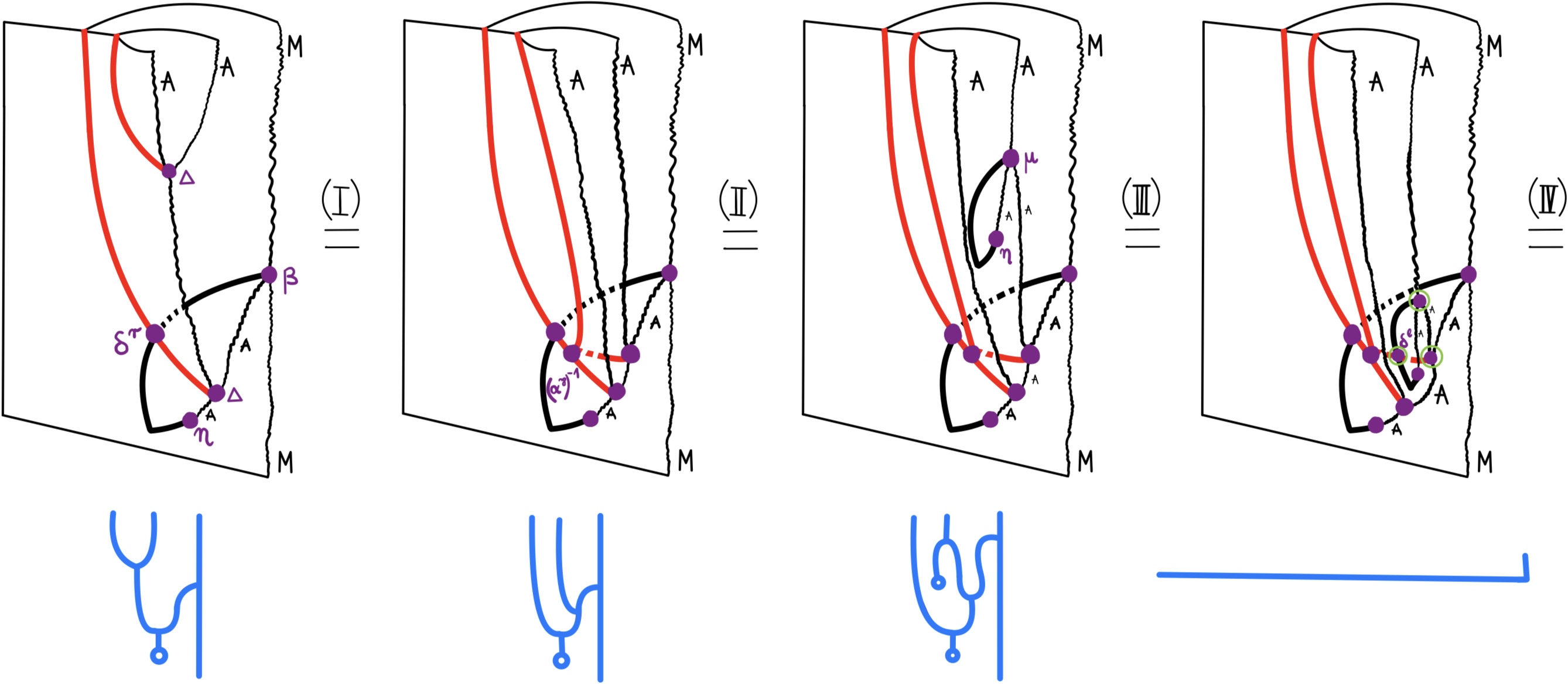}
\end{figure}

\bigskip

\begin{figure}[H]
    \centering
    \includegraphics[width=0.99\textwidth]{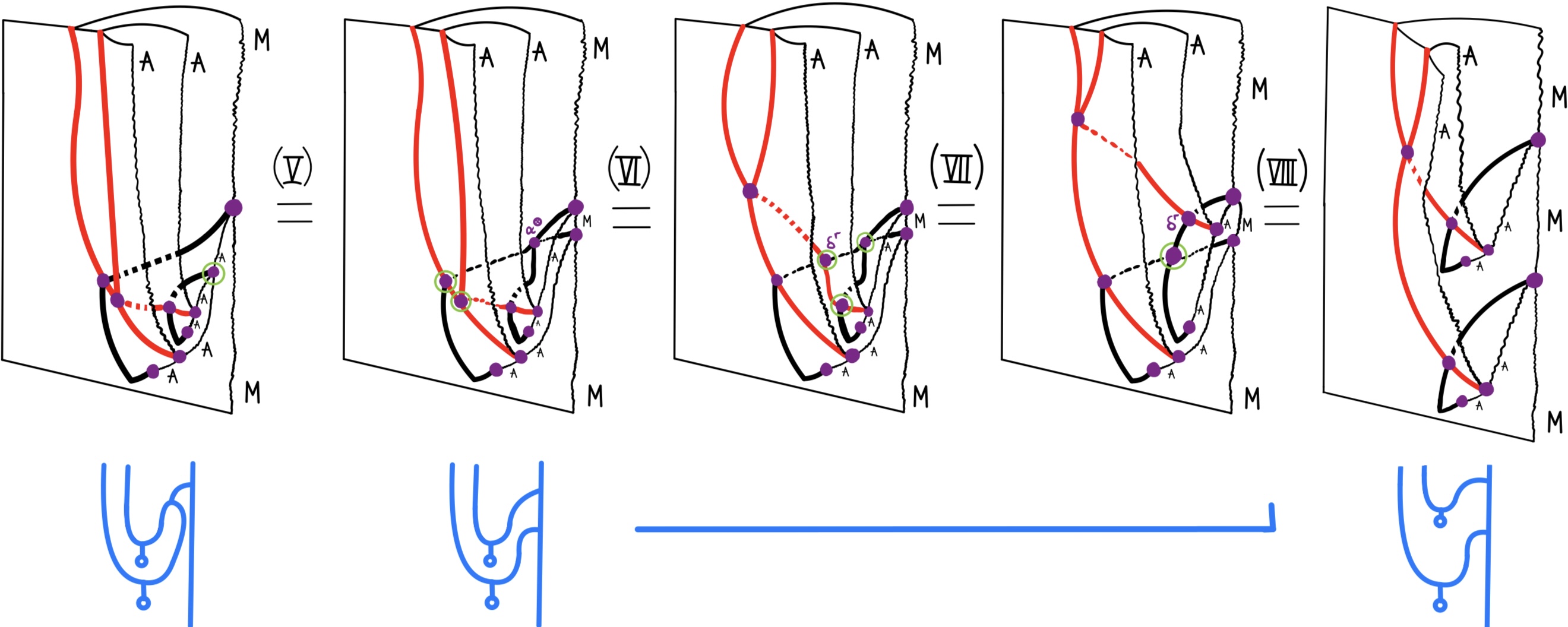}.
\end{figure}

\nid Equation (I) holds by the coassociativity of the comultiplication $\Delta$; c.f. Figure \ref{Associativity}. Equation (II) follows from the unitality of the multiplication $\mu$; see Figure \ref{unitality}. The triangle coherence axiom (\ref{eq:A1}) implies Equation (III). Equation (IV) follows from the Frobenius relation (F1); see Figure \ref{LD-FrobeniusRelations}. Equation (V) holds since the morphism $\beta$ is compatible with the multiplication $\mu$ by definition; cf. Figure \ref{Associativity}. Equations (VI) and (VII) are instances of the pentagon coherence axioms (\ref{eq:A7}) and (\ref{eq:A6}) from Appendix \ref{coherenceLD}, respectively. Finally, Equation (VIII) follows from Mac Lane's coherence theorem for the $\otimes$-monoidal structure.

The proof of counitality involves the triangle axiom (\ref{eq:A3}) and is omitted. The fact that the assignment $(M,\beta)\mapsto (M,\overline{\beta})$ extends to a functor is immediate.
\end{proof}

Applying Lemma \ref{module-induces-comodule} to the LD-category $\cC^{\text{cop}}$ yields the following claim:

\begin{lemma}\label{comodule-induces-module}
    Let $(N,\gamma)\in {^A\cC}$ be an $A$-comodule. Then the morphism $\overline{\gamma}$ from Figure \ref{fig:overlineGamma} endows the object $N$ with the structure of an $A$-module. This assignment yields a functor $G\colon ^A\cC\ra {_A\cC}$. As in Lemma \ref{module-induces-comodule}, the functor $G$ commutes with the forgetful functors ${^A \cC\ra \cC}$ and $_A \cC\ra \cC$.
\end{lemma}

We come to the main result of this subsection:

\begin{theorem}\label{thm:Frob Alg modules isomorphic to comodules}
    Let $A\in \cC$ be a Frobenius algebra in an LD-category $\cC$. The category of $A$-modules $_A\cC$ is isomorphic to the category of $A$-comodules ${^A\cC}$, where the isomorphism can be chosen to commute with the forgetful functors to $\cC$.
\end{theorem}

\begin{proof}
To prove that the functors $F\colon _A\cC \ra {^A\cC}$ from Lemma \ref{module-induces-comodule} and $G\colon {^A\cC} \ra {_A\cC}$ from Lemma \ref{comodule-induces-module} are inverses, it suffices to show that the composites $GF$ and $FG$ act as identities on objects. We only prove that, for any left $A$-module $(M,\beta)\in {_A \cC}$, we have $GF(M,\beta)=(M,\beta)$; the other proof is similar. Our proof involves the pentagon axiom (\ref{eq:A9}) in the background: 

\bigskip
    \begin{figure}[H]
    \centering
    \includegraphics[width=0.98\textwidth]{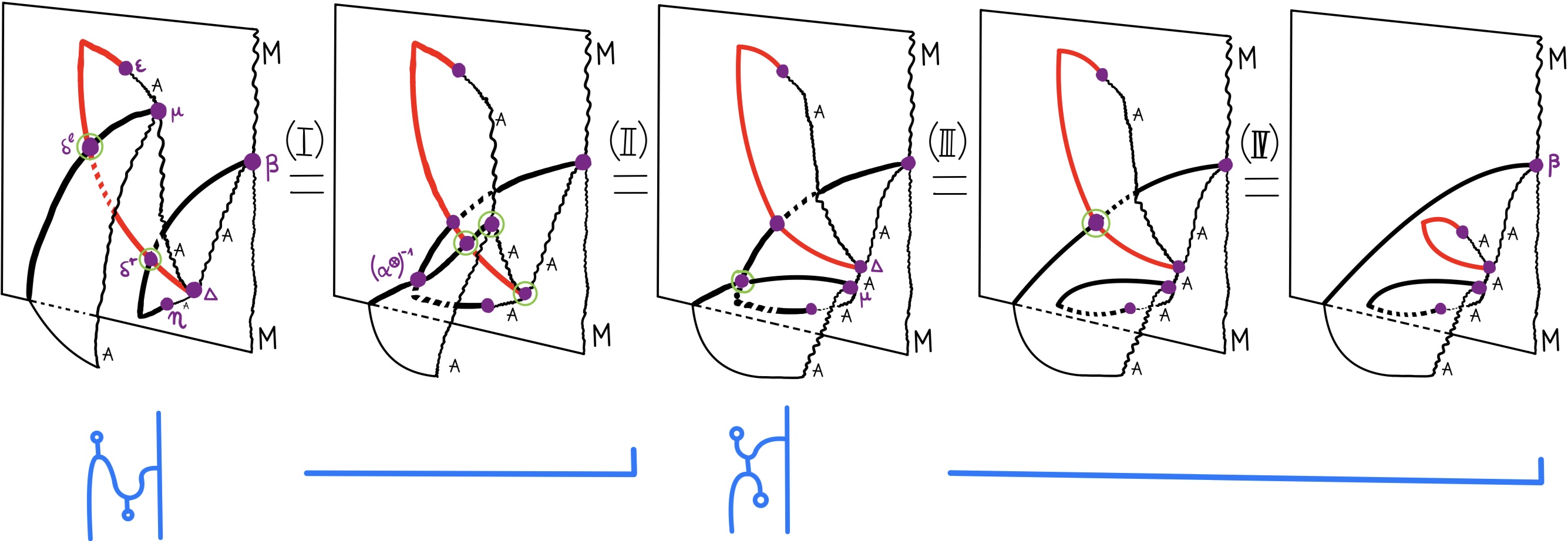}.
\end{figure}

Equation (I) follows from the pentagon coherence axiom (\ref{eq:A9}). Equation (II) holds by the Frobenius relations (F1) and \eqref{eq:F2 Frob LD}; see Equation (\ref{LD-FrobeniusRelations}). Mac Lane's triangle axiom for the $\otimes$-monoidal structure implies Equation (III); see Figure \ref{fig:MacTriangle}. Equation (IV) follows from the triangle coherence axiom (\ref{eq:A3}). Finally, the morphism represented by the right-hand surface diagram equals $\beta$ by the unitality of the multiplication $\mu$ and the counitality of the comultiplication $\Delta$, cf. Figure \ref{unitality}.
\end{proof}

%%%%%%%%%%%%%%%%%%%%%%%%%%

\section{Higher Frobenius-Schur indicators for pivotal GV-categories}\label{sec:FrobSchurIndic}
In this section, we generalize higher Frobenius-Schur indicators from $k$-linear pivotal rigid monoidal categories to $k$-linear pivotal GV-categories. Invariants of monoidal categories are hard to obtain. It is thus remarkable that, in the Grothendieck-Verdier setting, higher Frobenius-Schur indicators can be defined and are categorical invariants. As in previous sections, the three-dimensional graphical calculus will be an important tool.

As categorical invariants, higher Frobenius-Schur indicators can serve to distinguish categories-with-structure in the setting of GV-categories. For rigid monoidal categories, Schauenburg and Ng \cite{CentralInvariantsAndHigher} have used them to tell apart representation categories of semisimple quasi-Hopf algebras, which have the same number of (isomorphism classes) of simple objects and the same fusion rules. Applying Frobenius-Schur indicators to study concrete examples of non-rigid GV-categories goes beyond the scope of this paper, and is left to future investigation. 

\break

\subsection{Pivotal GV-categories}\label{sec:pivotalGV}
We first discuss the notion of pivotality. Pivotal rigid monoidal categories are well known, e.g. \cite{Freyd-Yetter}, \cite[Def. 4.7.7]{EGNO}. Boyarchenko and Drinfeld \cite[Def. 5.1]{BoDrinfeld} give the following generalization to the GV-setting.

\begin{definition}\label{def:pivotalStructure}
A \emph{pivotal structure} on a GV-category $(\cC,\otimes,1,K)$ is an isomorphism
\begin{equation*}
    \psi_{X,Y}\colon\, \operatorname{Hom}_{\cC}(X\otimes Y,K)\,\xlongrightarrow{\simeq}\,\operatorname{Hom}_{\cC}(Y\otimes X,K),
\end{equation*}
that is natural in $X,Y\in \cC$, and satisfies the following two properties:
\begin{enumerate}[label=(\roman*)]
    \item The diagram
        \begin{equation*}
        \adjustbox{max width=\textwidth,center}{
        \begin{tikzcd}
        {\operatorname{Hom}_{\cC}((Z\otimes X)\otimes Y,K)}&&&{\operatorname{Hom}_{\cC}(Z\otimes (X\otimes Y),K)}\\{\operatorname{Hom}_{\cC}((Y\otimes(Z\otimes X),K)}&&&{\operatorname{Hom}_{\cC}((X\otimes Y)\otimes Z,K)}\\
        {\operatorname{Hom}_{\cC}((Y\otimes Z)\otimes X,K)}&&&{\operatorname{Hom}_{\cC}(X\otimes(Y\otimes Z),K)}
        \arrow["\simeq"', "\operatorname{Hom}_{\cC}(\alpha^{-1}_{Z,X,Y}{,}K)" {yshift=2pt}, from=1-1, to=1-4]
        \arrow["\psi_{Z\otimes X,Y}"' {xshift=-2pt}, "\simeq", from=1-1, to=2-1]
        \arrow["\operatorname{Hom}_{\cC}(\alpha_{Y,Z,X}{,}K)"' {xshift=-2pt},"\simeq", from=2-1, to=3-1]
        \arrow["\psi_{Y\otimes Z,X}"' {yshift=-2pt}, "\simeq", from=3-1, to=3-4]
        \arrow["\simeq", "\operatorname{Hom}_{\cC}(\alpha_{X,Y,Z}{,}K)"' {xshift=2pt}, from=3-4, to=2-4]
        \arrow["\psi_{X\otimes Y,Z}"' {xshift=2pt}, "\simeq", from=2-4, to=1-4]
        \end{tikzcd}}
\end{equation*}
commutes for all objects $X,Y,Z\in \cC$.
    \item We have 
\begin{equation*}
\psi_{X,Y}\;=\;\psi^{-1}_{Y,X},
\end{equation*}
for all objects $X,Y\in \cC$.
\end{enumerate}
A \emph{pivotal GV-category} is a GV-category together with a choice of a pivotal structure.
\end{definition}

We now provide an equivalent definition of a pivotal GV-category that is better suited for the surface diagrammatic calculus. To do so, recall from \cite[Lemma 5.6]{BoDrinfeld} that for any GV-category $(\cC,\otimes,1,K)$, there is a one-to-one correspondence between natural isomorphisms $\operatorname{id}_{\cC}\xrightarrow{\simeq}D^2$ and isomorphisms \begin{equation*}
\operatorname{Hom}_{\cC}(X\otimes Y,K)\,\xlongrightarrow{\simeq}\,\operatorname{Hom}_{\cC}(Y\otimes X,K),\end{equation*} natural in $X,Y\in \cC$. Additionally, recall from Remark \ref{DoubleDual FrobLD} that the double dual $D^2$ is equipped with the structure of a strong Frobenius LD-endofunctor.

\begin{theorem}
    \label{equivalent def of pivotal structure}
    Let $(\cC,\otimes,1,K)$ be a GV-category. A natural isomorphism $\rho\colon \operatorname{id}_{\cC}\xrightarrow{\simeq} D^2$ corresponds to a pivotal structure on $\cC$ in the sense of Definition \ref{sec:pivotalGV} if and only if $\rho$ is a morphism of Frobenius LD-functors.
\end{theorem}

Theorem \ref{equivalent def of pivotal structure} and Proposition \ref{LDC2Cat} directly imply the following result:

\begin{corollary}
    Pivotal structures on a GV-category $\cC$ correspond to isomorphisms of Frobenius LD-functors $D^{\prime}\xrightarrow{\simeq} D$. Here, the duality functor $D$ and its quasi-inverse $D^{\prime}$ are endowed with the strong Frobenius LD-structures from Proposition \ref{duality functor is FrobLD} and Corollary \ref{dual dprime is frobenius ld}, respectively.
\end{corollary}

\begin{remark}
\begin{enumerate}[label=(\roman*)]
    \item Theorem \ref{equivalent def of pivotal structure} shows, in particular, that the notion of pivotality in the sense of Definition \ref{sec:pivotalGV} treats the $\otimes$-monoidal structures and the $\parLL$-monoidal structures on an equal footing.
    \item Given a GV-category $\cC$, any morphism of Frobenius LD-functors $\rho\colon \operatorname{id}_{\cC}\ra D^2$ is invertible by Proposition \ref{PropMorphLDfuncIso}.
    \item By Theorem \ref{equivalent def of pivotal structure}, the set of pivotal structures on a GV-category $\cC$ is either empty or a torsor over the group $\operatorname{Aut}^{\otimes,\parLL}(\operatorname{id}_{\cC})$ of automorphisms of the Frobenius LD-functor $\operatorname{id}_\cC$. %In the special case of finite-dimensional modules over a finite-dimensional Hopf algebra, this is the group of central grouplike elements.
\end{enumerate}
\end{remark}

After the proof of Theorem \ref{equivalent def of pivotal structure}, we will use the following terminology and notation:

\begin{remark}
    \begin{enumerate}[label=(\roman*)]
    \item For simplicity and to maintain consistency with the standard terminology for rigid monoidal categories, we will refer to either of the equivalent structures in Theorem \ref{equivalent def of pivotal structure} as a pivotal structure. 
    \item Similarly, we often denote a pivotal GV-category $(\cC,\otimes,1,K,\psi)$ simply by $(\cC,\rho)$, where $\rho\colon \operatorname{id}_{\cC}\xrightarrow{\simeq} D^2$ is the natural isomorphism corresponding to the pivotal structure $\psi$.
    \end{enumerate}
\end{remark}

\begin{examples}\label{pivotality of bimodule example}
\begin{enumerate}[label=(\roman*)]
    \item Let $A$ be a finite-dimensional $k$-algebra. Recall from Example \ref{GV-category of bimodules} the GV-category of finite-dimensional $A$-bimodules $A\operatorname{-bimod}^{\operatorname{f.d.}}$ with dualizing object $DA$. Let $\operatorname{vect}_k$ denote the rigid monoidal category of finite-dimensional $k$-vector spaces. The canonical identification of a finite-dimensional $k$-vector space $V$ with its double dual $V^{\ast\ast}$ yields a monoidal isomorphism $\operatorname{id}_{\operatorname{vect}_k}\xrightarrow{\simeq}(-)^{\ast\ast}$. This monoidal isomorphism lifts to an isomorphism $\operatorname{id}_{A\operatorname{-bimod}^{\operatorname{f.d.}}}\xrightarrow{\simeq}D^2$ of Frobenius LD-functors, thus providing a pivotal structure on the GV-category of finite-dimensional bimodules.
\item If the finite-dimensional $k$-algebra $A$ is commutative, any left $A$-module is canonically an $A$-bimodule. The left $A$-module $DA$ is a dualizing object in the full monoidal subcategory $A\operatorname{-mod}^{\operatorname{f.d.}}\subset A\operatorname{-bimod}^{\operatorname{f.d.}}$ of finite-dimensional left $A$-modules. Thus, the pivotal structure on the GV-category $A\operatorname{-bimod}^{\operatorname{f.d.}}$ restricts to a pivotal structure on the GV-category $A\operatorname{-mod}^{\operatorname{f.d.}}$.
\item More generally, in \cite{HLZ}, a class of vertex operator algebras and their categories of modules has been identified, which naturally carry a pivotal GV-structure; see also \cite{allen2021duality}.
\end{enumerate}
\end{examples}

We now turn to the proof of Theorem \ref{equivalent def of pivotal structure}. To begin, we establish the following notation: For the remainder of this subsection, we fix a GV-category $(\cC,\otimes,1,K)$. Let ${\eta\colon \operatorname{id}_{\cC} \xrightarrow{\simeq}D^{\prime}D}$ and ${\epsilon \colon DD^{\prime}\xrightarrow{\simeq}\operatorname{id}_{\cC}}$ denote the unit and the counit of the adjoint equivalence $D \dashv D^{\prime}$. Let $f\colon K\xrightarrow{\simeq} D1$ and $g\colon D^{\prime}1\xrightarrow{\simeq} K$ be the isomorphisms obtained via Yoneda's lemma from the natural isomorphisms 
\begin{equation*}
    \operatorname{Hom}_{\cC}(X,K)\cong \operatorname{Hom}_{\cC}(X\otimes 1, K)\cong \operatorname{Hom}_{\cC}(X,D1),\text{ and}
    \end{equation*}
\begin{equation*}
    \operatorname{Hom}_{\cC}(X,D^{\prime}1)\cong \operatorname{Hom}_{\cC}(1,DX) \cong \operatorname{Hom}_{\cC}(1\otimes X,K)\cong \operatorname{Hom}_{\cC}(X,K).
\end{equation*}

\bigskip

Theorem \ref{equivalent def of pivotal structure} will be a consequence of the following three results.

\begin{prop}\label{pivotality via pivot}(\cite[Prop. 5.7]{BoDrinfeld}). A natural isomorphism $\rho\colon \operatorname{id}_{\cC}\xrightarrow{\simeq}D^2$ corresponds to a pivotal structure on $\cC$ if and only if it satisfies the following two conditions:
\begin{enumerate}[label=(\roman*)]
    \item $\rho$ is $\otimes$-monoidal, and
    \item $\rho_K\colon K \xrightarrow{\simeq}D^2K$ equals the isomorphism
    \begin{equation*}
        K\xrightarrow{f}D1\xrightarrow{D(\epsilon_1)}D^2D^{\prime}1\xrightarrow{D^2(g)}D^2K.
    \end{equation*}
\end{enumerate}
Moreover, a natural isomorphism $\rho\colon \operatorname{id}_{\cC}\xrightarrow{\simeq}D^2$ satisfying these two conditions also satisfies \begin{equation}\label{pivot under duality}
    \rho_{DX}\;=\;D(\rho_X)^{-1} \quad \forall X\in \cC.
\end{equation}
\end{prop}

Let $(\varphi^{2,D^2},\varphi^{0,D^2},\nu^{2,D^2},\nu^{0,D^2})$ denote the strong Frobenius LD-structure on the double duality functor $D^2\colon \; \cC \,\ra\, \cC$ from Remark \ref{DoubleDual FrobLD}.

\begin{lemma}\label{triangle diagram for rho}
Let $\rho\colon \operatorname{id}_{\cC}\xrightarrow{\simeq}D^2$ be a natural isomorphism. Then $\rho$ satisfies condition (ii) from Proposition \ref{pivotality via pivot} if and only if we have \begin{equation}\label{triange dgm rho par monoidal}
    \operatorname{id}_K\;=\;\nu^{0,D^2}\circ\rho_K.
\end{equation}
\end{lemma}

\begin{proof}
We freely use the notation introduced before Proposition \ref{pivotality via pivot}. By the definition of the $\parLL$-unit $\nu^{0,D^2}$ (see Remark \ref{DoubleDual FrobLD} and the proof of Proposition \ref{duality functor is FrobLD}), we have
\begin{equation}
\nu^{0,D^2} \;=\; g \circ {D^{\prime}(\varphi^{0,D^2})} \circ \eta_{D(1)} \circ {D(\epsilon_1^{-1})} \circ {D^2(g^{-1})}.
\end{equation}
Thus, Equation (\ref{triange dgm rho par monoidal}) amounts to the the commutativity of the following outer diagram:

 \begin{equation*}
  \adjustbox{max width=\textwidth}{
  \begin{tikzcd}
	{D^{\prime}(1)} & {D^{\prime}D^2(1)} & {D(1)} & {D^2 D^{\prime}(1)} \\
	& {} & {K} & {D^2(K).}
 \arrow["D^{\prime}(\varphi^{0,D^2})^{-1}" {yshift=5pt},"\simeq"', from=1-1, to=1-2]
 \arrow["\eta_{D(1)}^{-1}" {yshift=5pt}, "\simeq"', from=1-2, to=1-3] 
 \arrow["D(\epsilon_1)" {yshift=5pt}, "\simeq"', from=1-3, to=1-4]
 \arrow["f" {xshift=-3pt}, "\simeq"', from=2-3, to=1-3]
 \arrow["D^2(g)" {xshift=5pt},"\simeq"', from=1-4, to=2-4]
 \arrow["\rho_K"' {yshift=-5pt},"\simeq", from=2-3, to=2-4]
\arrow["g"'{yshift=-1pt},"\simeq", bend right=15, from=1-1, to=2-3]
\arrow[phantom,"\textup{(I)}"{xshift=-10pt}, bend left=5, from=1-2, to=2-3]
\arrow[phantom,"\textup{(II)}"{yshift=-1pt}, from=1-3, to=2-4]
 \end{tikzcd}
}
\end{equation*}

\nid By Boyarchenko-Drinfeld \cite[Prop. 4.2]{BoDrinfeld}, we have \begin{equation}\label{varphi0 equal to}
   \epsilon_1^{-1} \, = \, D(g)\circ D(f)\circ \varphi^{0,D^2}.\end{equation} 
   A careful inspection shows that Equation (\ref{varphi0 equal to}), together with Lemma \ref{unit counit under duality} and the naturality of the unit $\eta\colon \operatorname{id}_{\cC}\xrightarrow{\simeq}D^{\prime}D$, imply that diagram (I) is commutative. Thus, diagram (II), which is just condition (ii) in Proposition \ref{pivotality via pivot}, commutes if and only if the outer diagram does.
\end{proof}

\begin{lemma}\label{pivot is par monoidal}
Let $\rho\colon\operatorname{id}_{\cC}\xrightarrow{\simeq}D^2$ be a natural isomorphism that satisfies conditions (i) and (ii) of Proposition \ref{pivotality via pivot}. Then $\rho$ is also a $\parLL$-monoidal natural transformation.
\end{lemma}

\begin{proof}
By Lemma \ref{triangle diagram for rho} it suffices to show \begin{equation}\label{par is compatible with comultiplication}
\nu^{2,D^2}_{X,Y}\circ\rho_{X\parLL Y}\;=\;\rho_X\parLL\rho_Y \quad \forall X,Y\in \cC.
\end{equation}
Writing out the definition of the comultiplication morphism $\nu^{2,D^2}$ and the $\parLL$-monoidal product explicitly, we find that Equation (\ref{par is compatible with comultiplication}) is equivalent to the commutativity of the following outer diagram:
\begin{equation*}
  \adjustbox{max width=\textwidth}{
  \begin{tikzcd}
	{D^2D^{\prime}(DY\otimes DX)} & {D(DY\otimes DX)} &{D(D^{\prime}D^2Y\otimes D^{\prime}D^2X)} & {D^{\prime}D^2(D^{\prime}D^2Y\otimes D^{\prime}D^2X)} \\
	& {D^{\prime}D^2(DY\otimes DX)} & & {D^{\prime}(D^2D^{\prime}D^2Y\otimes D^2D^{\prime}D^2X)} \\
    & {D^{\prime}(DY\otimes DX)} & & {D^{\prime}(D^3Y\otimes D^3X)}.
    \arrow["D(\epsilon^{-1})" {yshift=5pt},"\simeq"', from=1-1, to=1-2]
    \arrow["D({\eta^{-1}}\otimes {\eta^{-1}})" {yshift=5pt}, "\simeq"', from=1-2, to=1-3] \arrow["\eta" {yshift=5pt}, "\simeq"', from=1-3, to=1-4]
    \arrow["D^{\prime}\big(\varphi^{2,D^2}\big)" {xshift=5pt},"\simeq"', from=1-4, to=2-4]
    \arrow["D^{\prime}\big(D^2\eta\otimes D^2\eta\big)" {xshift=5pt},"\simeq"', from=2-4, to=3-4]
    \arrow["D^{\prime}(D\rho_Y\otimes D\rho_X)"' {yshift=-5pt},"\simeq", from=3-2, to=3-4]
    \arrow["\rho_{D^{\prime}(DY\otimes DX)}","\simeq"', bend left=35, from=3-2, to=1-1]
    \arrow["\eta"'{xshift=-4pt},"\simeq", from=1-2, to=2-2]
    \arrow["D^{\prime}D^2(\eta\otimes \eta)","\simeq"', bend right=2, from=1-4, to=2-2]
    \arrow["D^{\prime}\big(\varphi^{2,D^2}\big)","\simeq"', from=2-2, to=3-4]
    \arrow["D^{\prime}(\rho^{-1})" {xshift=-4pt},"\simeq"', from=3-2, to=2-2]
    \arrow[phantom,"\textup{(I)}"{xshift=80pt,yshift=5pt}, from=1-2, to=2-2]
    \arrow[phantom,"\textup{(II)}"{xshift=-10pt,yshift=-21pt}, from=1-3, to=2-4]
    \arrow[phantom,"\textup{(III)}"{xshift=-140pt,yshift=-50pt}, from=1-3, to=2-4]
    \arrow[phantom,"\textup{(IV)}"{xshift=-290pt,yshift=-20pt}, from=1-3, to=2-4]
   \end{tikzcd}}
\end{equation*}
By the naturality of the unit $\eta\colon \operatorname{id}_{\cC}\xrightarrow{\simeq}D^{\prime}D$, diagram (I) commutes. Similarly, by the naturality of the multiplication morphism $\varphi^{2,D^2}$, diagram (II) commutes. Diagram (III) commutes by Equation (\ref{pivot under duality}) and since $\rho$ is an $\otimes$-monoidal natural transformation by condition (i). The diagram (IV) can be shown to commute by using the naturality of the unit $\eta\colon \operatorname{id}_{\cC}\xrightarrow{\simeq}D^{\prime}D$ and the counit $\epsilon \colon DD^{\prime}\xrightarrow{\simeq}\operatorname{id}_{\cC}$ together with Lemma \ref{unit counit under duality} and Equation (\ref{pivot under duality}).
\end{proof}

Let us collect our results:

\begin{proof}[Proof of Theorem \ref{equivalent def of pivotal structure}]
By Proposition \ref{pivotality via pivot} and Lemma \ref{pivot is par monoidal}, any natural isomorphism ${\operatorname{id}_{\cC}\xrightarrow{\simeq}D^2}$ that corresponds to a pivotal structure on $\cC$ is a morphism of Frobenius LD-functors.
Conversely, by Proposition \ref{pivotality via pivot} and Lemma \ref{triangle diagram for rho}, any isomorphism of Frobenius LD-functors ${\operatorname{id}_{\cC}\xrightarrow{\simeq}D^2}$ corresponds to a pivotal structure.
\end{proof}

\begin{figure}[H]
    \centering
    \includegraphics[width=0.3\textwidth]{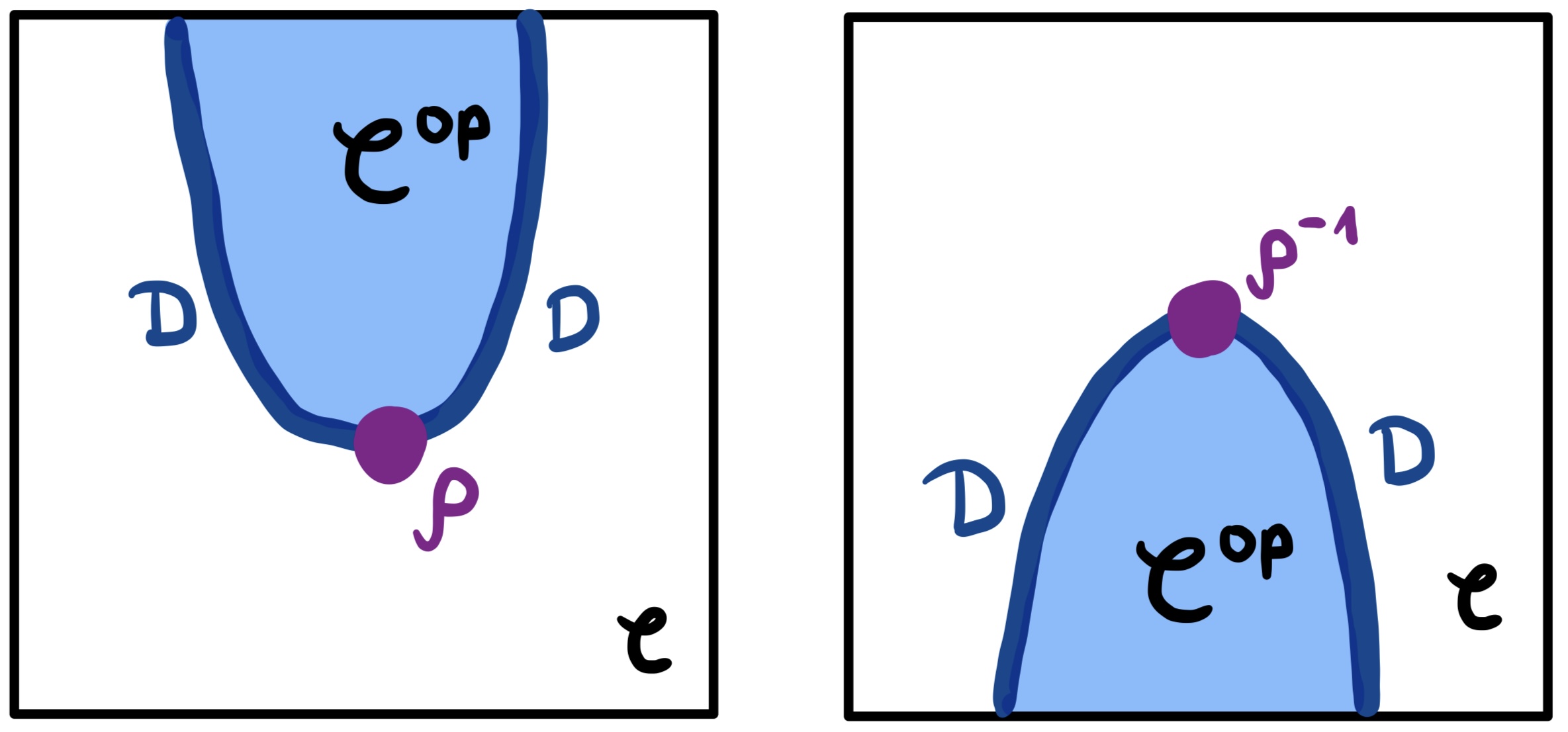}
    \caption{A pivotal structure $\rho\colon \operatorname{id}_{\cC}\xrightarrow{\simeq} D^2$ and its inverse $\rho^{-1}$.}
    \label{pivot}
\end{figure}

\begin{figure}[H]
\centering
\begin{subfigure}[b]{0.7\textwidth}
         \centering
         \captionsetup{justification=centering}
         \includegraphics[width=0.92\textwidth]{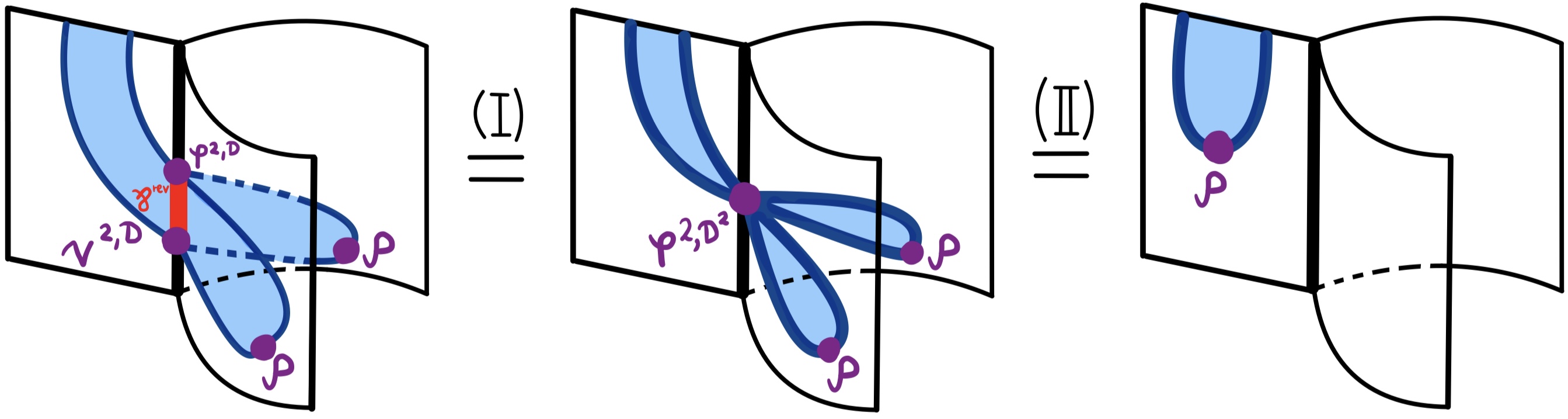}
         \caption{$\rho$ is compatible with the multiplication morphism $\varphi^{2,D^2}$.}
         \label{f monoidal}
\end{subfigure}
\\[2ex]
\begin{subfigure}[b]{0.6\textwidth}
        \centering
        \captionsetup{justification=centering}
         \includegraphics[width=0.92\textwidth]{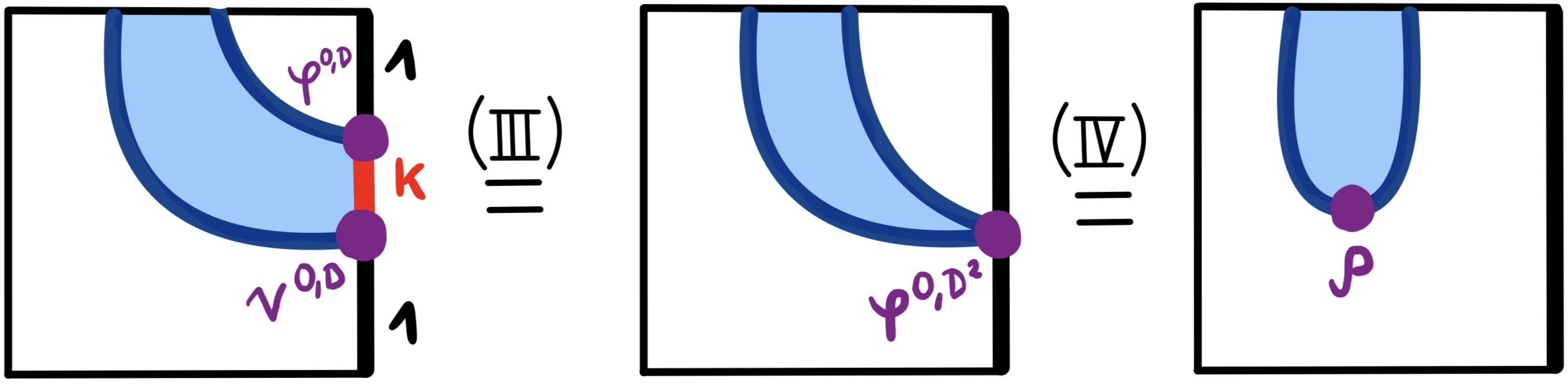}
         \caption{$\rho$ is compatible with the unit morphism $\varphi^{0,D^2}$.}
         \label{f monoidal2}
\end{subfigure}     
\caption{A pivotal structure $\rho\colon \operatorname{id}_{\cC}\xrightarrow{\simeq} D^2$ is an $\otimes$-monoidal natural transformation.}
\label{fig:surface dgms for pivot}
\end{figure}

Pivotal structures $\rho\colon \operatorname{id}_{\cC}\xrightarrow{\simeq} D^2$ can be depicted graphically, as shown in Figures \ref{pivot} and \ref{fig:surface dgms for pivot}. The coherence morphisms of the strong Frobenius LD-functor $D\colon \cC \ra\cC^{\operatorname{op}}$ are denoted by
\begin{equation*}
\varphi^{2,D}_{X,Y}\in\operatorname{Hom}_{\cC^{\operatorname{op}}}\big({D(X)}\parLL^{\operatorname{rev}}{D(Y)},D(X\otimes Y)\big), \qquad \text{and} \qquad {\varphi^{0,D}_{X,Y}\in\operatorname{Hom}_{\cC^{\operatorname{op}}}\big(K,D(1)\big)},
\end{equation*}
\begin{equation*}
{\nu^{2,D}_{X,Y}\in\operatorname{Hom}_{\cC^{\operatorname{op}}}\big(D({X}\parLL{Y}),{D(X)}\otimes^{\operatorname{rev}}{D(Y)}\big)}, \qquad \text{and} \qquad {\nu^{0,D}_{X,Y}\in\operatorname{Hom}_{\cC^{\operatorname{op}}}\big(D(K),1\big)}.
\end{equation*}

Equations (II) and (IV) in Figure \ref{fig:surface dgms for pivot} express that $\rho$ is an $\otimes$-monoidal natural transformation. Equations (I) and (III) in Figure \ref{fig:surface dgms for pivot} hold by the definition of the strong $\otimes$-monoidal structure of the composite $D^2$, see Remark \ref{double dual mon structure}.

\medskip

The surface diagrams that encode $\rho$ as an oplax $\parLL$-monoidal transformation are similar.

%%%%%%%%%%%%%%%%%%%%%%%%

\subsection{Pivotal functors}\label{sec:pivotalFunctor}
We turn to the discussion of the notion of morphism between pivotal GV-categories. By abuse of notation, given two GV-categories $\cC$ and $\cD$, we denote both of their duality functors by $D$. Let $F\colon \cC\ra \cD$ be a Frobenius LD-functor with unit morphism $\varphi^0\colon 1 \ra F(1)$ and comultiplication morphism ${\upsilon^2\colon F\circ \parLL \ra {\parLL} \circ {(F\times F)}}$. LD-duals are unique up to a unique isomorphism, see Remark \ref{LD-duals are unique}. Furthermore, by Proposition \ref{FrobPresDuals}, $FD(X)$ is a left LD-dual of $F(X)$. Hence, for any object $X\in \cC$, there exists a unique isomorphism
\begin{equation}\label{xiFX}
    \xi^F_X\in \operatorname{Hom}_{\cD}\big(FD(X), DF(X)\big)
\end{equation}
such that
\begin{equation}\label{charact property duality transfo}
    {(F(X)\parLL\,\xi^F_X)}\circ \widetilde{\operatorname{coev}}_{F(X)}\;=\;\operatorname{coev}_{F(X)},    
\end{equation}
where we have defined
\begin{equation*}
    \widetilde{\operatorname{coev}}_{F(X)}\colon \: 1\, \xrightarrow{\varphi^{0}}\, F(1)\,\xrightarrow{F(\operatorname{coev_X})}\, F(X\parLL D(X))\,\xrightarrow{\upsilon^2_{X,D(X)}}\,F(X)\parLL FD(X).    
\end{equation*}

\begin{definition}\label{definition of duality transformation}
    For any Frobenius LD-functor $F$ between GV-categories $\cC$ and $\cD$, the morphisms $(\xi^F_X)_{X\in\cC}$ assemble into an isomorphism of functors
        \begin{equation}
            \xi^F\colon F\circ D\,\xlongrightarrow{\simeq}\, D\circ F,
        \end{equation}
    called the \emph{duality transformation} for $F$.
\end{definition}

\begin{figure}[H]
     \centering
    \begin{subfigure}[c]{0.48\textwidth}
         \centering
         \includegraphics[width=0.82\textwidth]{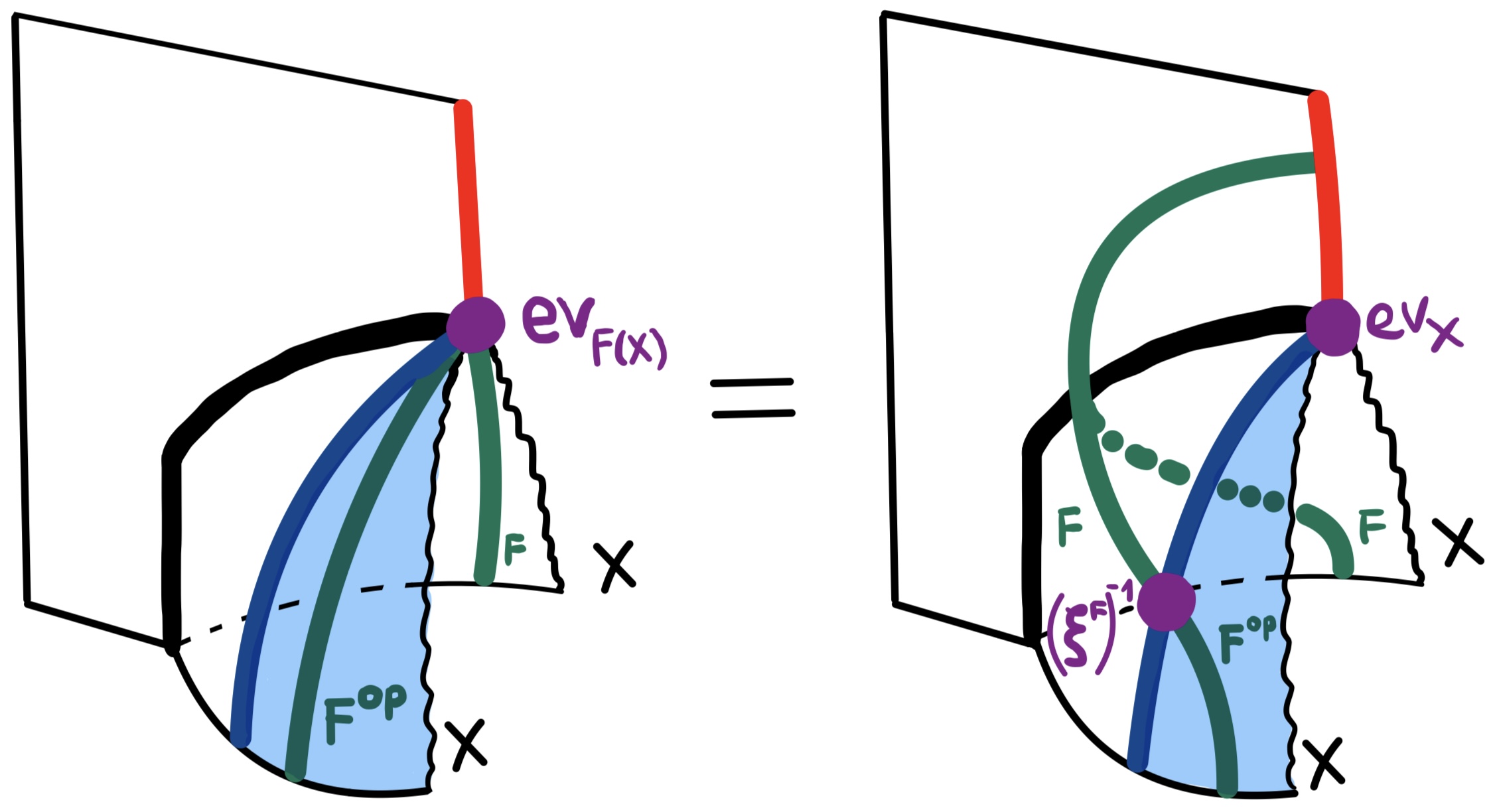}
         \caption{Characterizing $(\xi^F)^{-1}$.}
         \label{charac of duality transfo}
     \end{subfigure}
     \begin{subfigure}[c]{0.48\textwidth}
         \centering
    \includegraphics[width=0.82\textwidth]{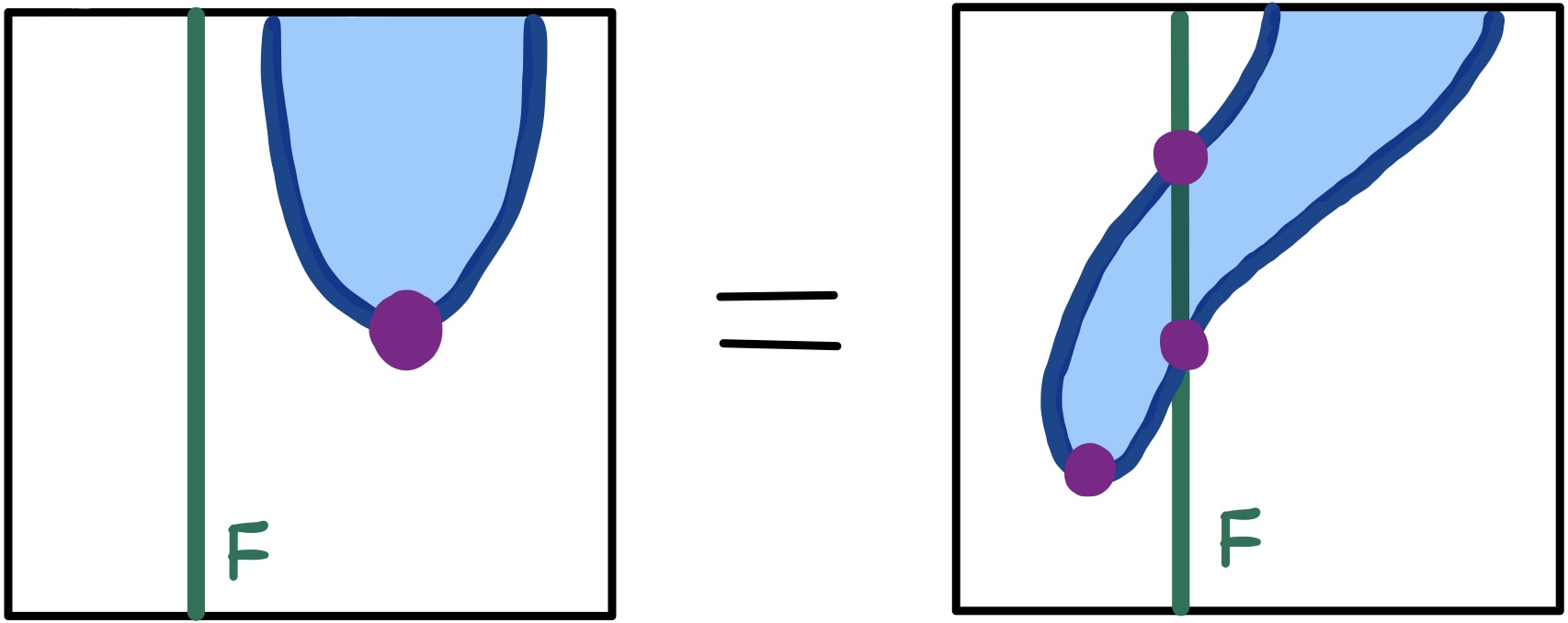}
    \caption{Equation (\ref{preserve pivotality}).}
    \label{fig:pivotality preserving}
 \end{subfigure}
     \caption{Surface diagrams for the inverse duality transformation and the pivotal functor $F$.}
\end{figure}

\begin{remark}\label{duality transfo Remark}
\begin{enumerate}[label=(\roman*)]
    \item There is no asymmetry between evaluations and coevaluations in the definition of the duality transformation. The inverse $(\xi^F)^{-1}$ of the duality transformation can be characterized as the unique natural isomorphism $D \circ F \xrightarrow{\simeq} F \circ D$, such that the equality of surface diagrams in Figure \ref{charac of duality transfo} is satisfied.
    \item As composites of Frobenius LD-functors, both $F \circ D$ and $D \circ F$ carry the structure of Frobenius LD-functors, according to Proposition \ref{LDC2Cat}. By adapting the proof of \cite[Lemma 1.1.]{SchauNgHigherFrob} in the rigid setting to surface diagrams, one can show that the duality transformation $\xi^F$ is a morphism of these Frobenius LD-functors.
\end{enumerate}
\end{remark}

Ng and Schauenburg \cite[\S 1]{SchauNgHigherFrob} define pivotal strong monoidal functors between pivotal rigid monoidal categories. We generalize their definition to Frobenius LD-functors between pivotal GV-categories.

\begin{definition}
    Let $(\cC,\rho^{\cC})$ and $(\cD,\rho^{\cD})$ be pivotal GV-categories, and let $F\colon \cC \ra \cD$ be a Frobenius LD-functor between the underlying GV-categories. Let $\xi^F$ denote the duality transformation of $F$. We call $F$ \emph{pivotal} if the diagram 
    \begin{equation}\label{preserve pivotality}
\begin{tikzcd}
	{F(X)} && {D^2F(X)} \\
	{FD^2(X)} && {DFD(X)}
	\arrow["{\rho^{\cD}_{F(X)}}"{yshift=2pt}, from=1-1, to=1-3]
	\arrow["{F(\rho^{\cC}_{X})}"'{xshift=-2pt}, from=1-1, to=2-1]
	\arrow["{D(\xi^F_X)}"{xshift=2pt}, from=1-3, to=2-3]
	\arrow["{\big(\xi^F_{D(X)}\big)^{-1}}"{yshift=-2pt}, from=2-3, to=2-1]
\end{tikzcd}
\end{equation}
commutes for any object $X\in \cC.$
\end{definition}

\begin{remark}
As shown in Figure \ref{fig:pivotality preserving}, Equation (\ref{preserve pivotality}) expresses graphically that we can pull pivotal structures (and their inverses) over the functor line of $F$.
\end{remark}

%Kürzungspotential
We state the following lemma for future reference. It follows directly from the snake equations in Figure \ref{fig:LD-snake equations (S1) and (S2)} and the invertibility of the duality transformation $\xi^F$.
\begin{lemma}\label{evFDV and coevFDV}
Let $F\colon \cC\ra \cD$ be a Frobenius LD-functor between GV-categories. For any object $X\in \cC$, the LD-pairing ${DFD(X)}\otimes {FD(X)}\,\ra\, K$, defined by \begin{equation}\label{FrobSchur LD-pairing}
{\operatorname{ev}_{DF(X)}}\circ{({D((\xi_X^F)^{-1})}\otimes{\xi_X^F})},\end{equation}
and the LD-copairing
$1\ra {FD(X)}\parLL {DFD(X)}$, defined by
    \begin{equation}\label{FrobSchur LD-copairing}
        {((\xi_X^F)^{-1}}\parLL {D(\xi_X^F))}\circ{\operatorname{coev}_{DF(X)}},
    \end{equation} are side-inverse to each other in the sense of Definition \ref{pairing/copairing,side-inverse}.
\end{lemma}

%%%%%%%%%%%%%%%%%%%%%%%%%%%%%%%

\subsection{Higher Frobenius-Schur indicators}\label{sec:FrobSchur}
Let us establish some conventions for $k$-linear GV-categories over a field $k$. We assume that both monoidal products of a $k$-linear GV-category $\cC$ are $k$-bilinear functors. Additionally, we assume that each morphism space in $\cC$ is finite-dimensional. In other words, we take $(\cV,\otimes,J)$ to be the symmetric monoidal category of finite-dimensional $k$-vector spaces, and we work internally to the strict monoidal $2$-category $\cV{\textnormal-}\mathsf{Cat}$.

Next, we turn to the definition of our generalized Frobenius-Schur indicators. They are defined for a $k$-linear pivotal GV-category $(\cC,\rho)$. We fix an object $V\in \cC,$ and let $V^{\parLL n}$ denote the right-bracketed $n$-fold $\parLL$-monoidal product of $V$. The $k$-vector space $\operatorname{Hom}_{\cC}(1,V^{\parLL n})$ is called the \emph{space of invariant $n$-tensors}, or \emph{space of invariant tensors} for short. 

\begin{figure}[H]
\centering
    \begin{subfigure}[b]{0.27\textwidth}
         \centering
         \includegraphics[width=0.96\textwidth]{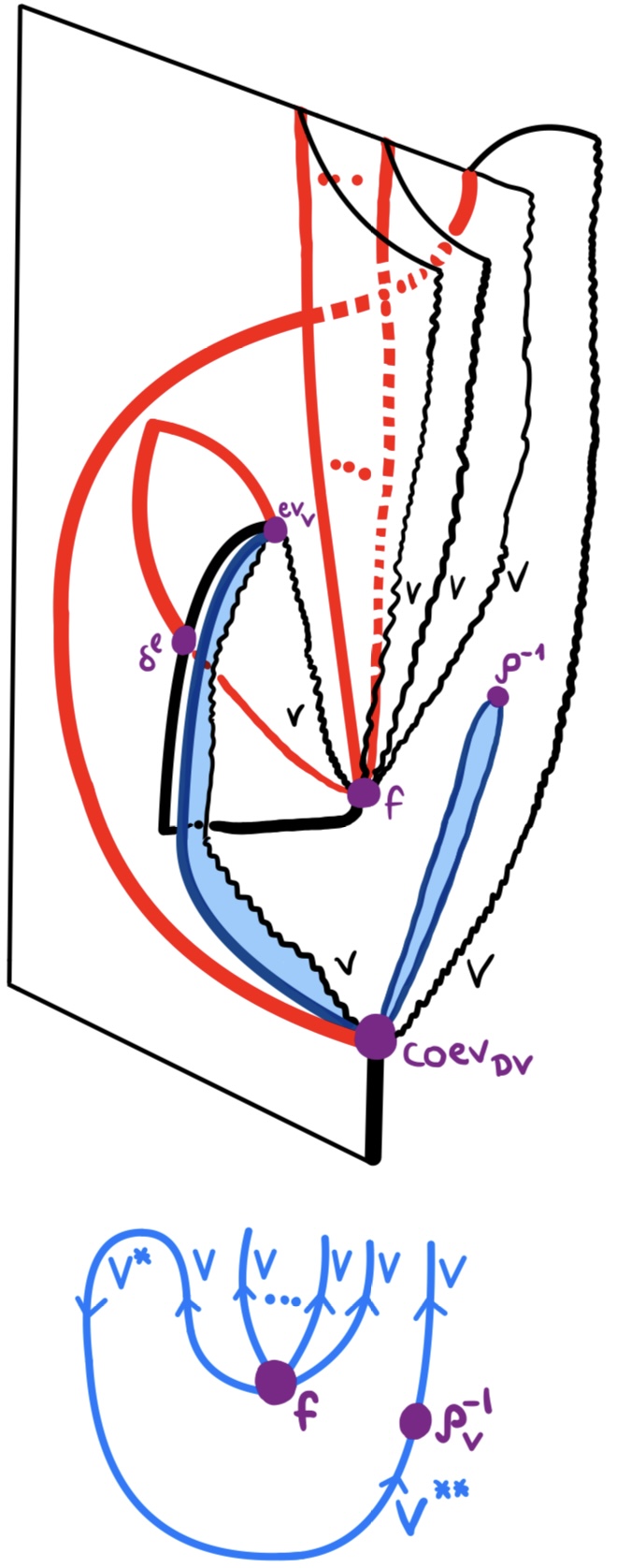}
         \caption{$E_V^{(n)}(f)$.}
         \label{fig:EnV}
    \end{subfigure}
\hspace{6em}
    \begin{subfigure}[b]{0.23\textwidth}
    \centering
    \includegraphics[width=0.96\textwidth]{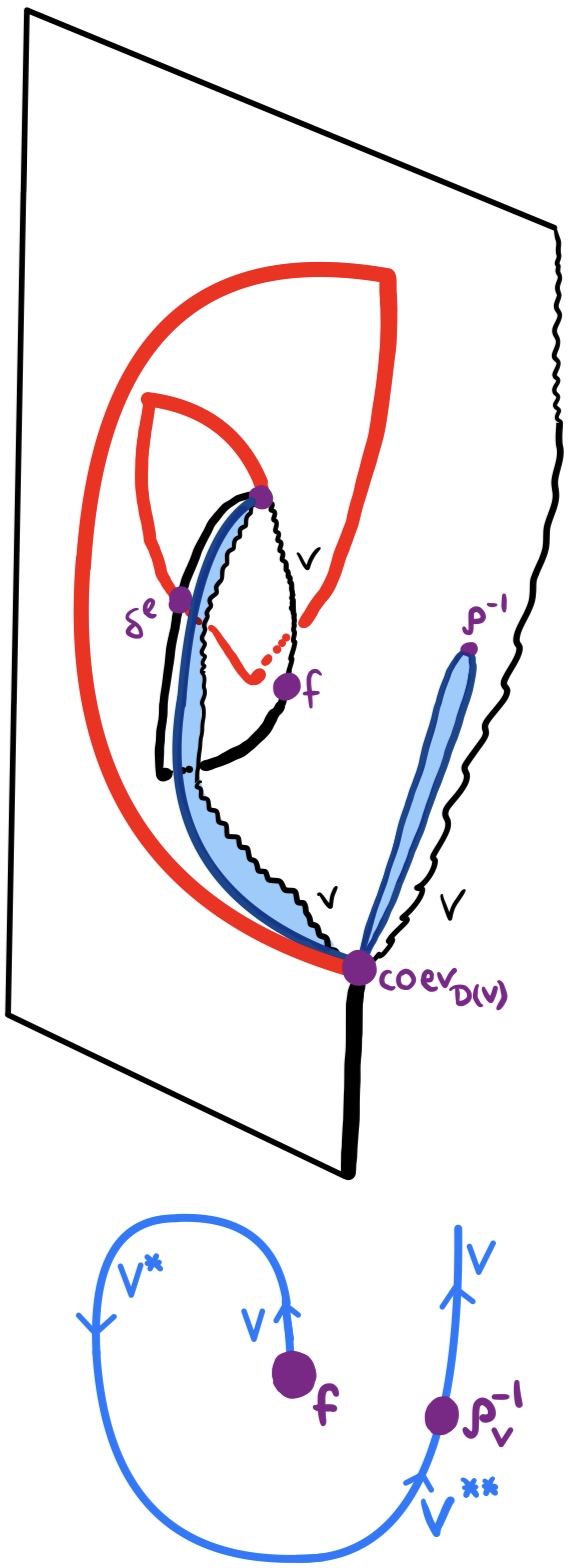}
    \caption{$E_V^{(1)}(f)$.}
    \label{fig:E1V}
    \end{subfigure}
\caption{Definition of the endomorphisms $E_V^{(n)}$ and $E_V^{(1)}$.}
\label{fig:frobschurendo}
\end{figure}

Figure \ref{fig:frobschurendo} shows the surface diagrammatic definition of the $k$-linear endomorphism\, \begin{equation*}
    E_V^{(n)}\colon \operatorname{Hom}_{\cC}(1,V^{\parLL n}) \,\longrightarrow\, \operatorname{Hom}_{\cC}(1,V^{\parLL n}); \quad f\mapsto E_V^{(n)}(f)
\end{equation*}
on the space of $n$-tensors. When viewed from the front face, the surface diagrams in Figure \ref{fig:frobschurendo} yield string diagrams which are familiar from the definition of higher Frobenius-Schur indicators for pivotal rigid monoidal categories \cite[Def. 3.1]{SchauNgHigherFrob}.

\smallskip

Following \cite[Def. 3.1.]{SchauNgHigherFrob} and \cite{KaSoZh} in the rigid monoidal case, we define the following:

\begin{definition}
Let $\cC$ be a $k$-linear pivotal GV-category. For integers $r,n\geq 1$, we call the $k$-linear trace \begin{equation*}
    \nu_{n,r}(V)\,:=\, \operatorname{Tr}((E^{(n)}_V)^r)
\end{equation*} 
the \emph{$(n,r)$-th Frobenius-Schur indicator} of $V\in \cC$. Here, $E^{(n)}_V$ is the endomorphism introduced in Figure \ref{fig:frobschurendo}.
\end{definition}

We need the following technical lemma to show the invariance of Frobenius-Schur indicators under equivalences of $k$-linear pivotal GV-categories. This lemma generalizes Ng and Schauenburg's \cite[Lemma 4.1.]{SchauNgHigherFrob} from strong monoidal functors between \emph{rigid} monoidal categories to Frobenius LD-functors between \emph{GV-categories}. From this lemma, it follows that Ng and Schauenburg's result (which applies specifically to the rigid pivotal case) holds not only for strong monoidal functors but also for Frobenius monoidal functors.

\begin{lemma}\label{lemma:technical lemma}
    Let $F\colon \cC\ra \cD$ be a Frobenius LD-functor between GV-categories. As in Definition \ref{LD-cats and functors}, denote its comultiplication and unit morphism by $\upsilon^2$ and $\varphi^0$, respectively. For all $V,W\in \cC$, the following diagram commutes:
\begin{equation}\label{hexagon diagram technical lemma}
    \begin{tikzcd}
	   {\operatorname{Hom}_{\cC}(1,V\parLL W)} & {\operatorname{Hom}_{\cD}(F(1),F(V\parLL W))} & {\operatorname{Hom}_{\cD}(1,F(V)\parLL F(W))}\\
	   {\operatorname{Hom}_{\cC}(D(V),W)} & {\operatorname{Hom}_{\cD}(FD(V),F(W))} & {\operatorname{Hom}_{\cD}(DF(V),F(W)),}
	   \arrow["F", from=1-1, to=1-2]
	   \arrow["A", from=1-2, to=1-3]
	   \arrow["F", from=2-1, to=2-2]
	   \arrow["C","\simeq"', from=2-2, to=2-3]
	   \arrow["B^{V,W}", from=1-1, to=2-1]
        \arrow["B^{F(V),F(W)}", from=1-3, to=2-3]
    \end{tikzcd}
\end{equation}
where, for $f\in\operatorname{Hom}_{\cC}(1,V\parLL W)$, we used the relevant structure morphisms. Concretely, we set
\begin{equation*}
    A\,:=\, \operatorname{Hom}_{\cD}(\varphi^0,\upsilon^2_{V,W}),
\end{equation*}
\begin{equation*}
    B^{V,W}(f)\,:=\,l^{\parLL}_W\circ(\operatorname{ev}_V\parLL W)\circ\distl_{D(V),V,W}\circ(D(V)\otimes f)\circ(r^{\otimes}_{D(V)})^{-1},
\end{equation*}
\begin{equation*}
    C\,:=\,\operatorname{Hom}_{\cD}((\xi^F_V)^{-1},F(W)),    
\end{equation*}
where the duality transformation $\xi^F$ was introduced in Definition \ref{definition of duality transformation}.
\end{lemma}

\begin{proof}
The first surface diagram in the following figure shows the composite of the upper horizontal and the right-hand vertical map in the diagram (\ref{hexagon diagram technical lemma}), applied to a morphism ${f\in \operatorname{Hom}_{\cC}(1,V\parLL W)}$:
\begin{figure}[H]
    \centering
    \includegraphics[width=0.98\textwidth]{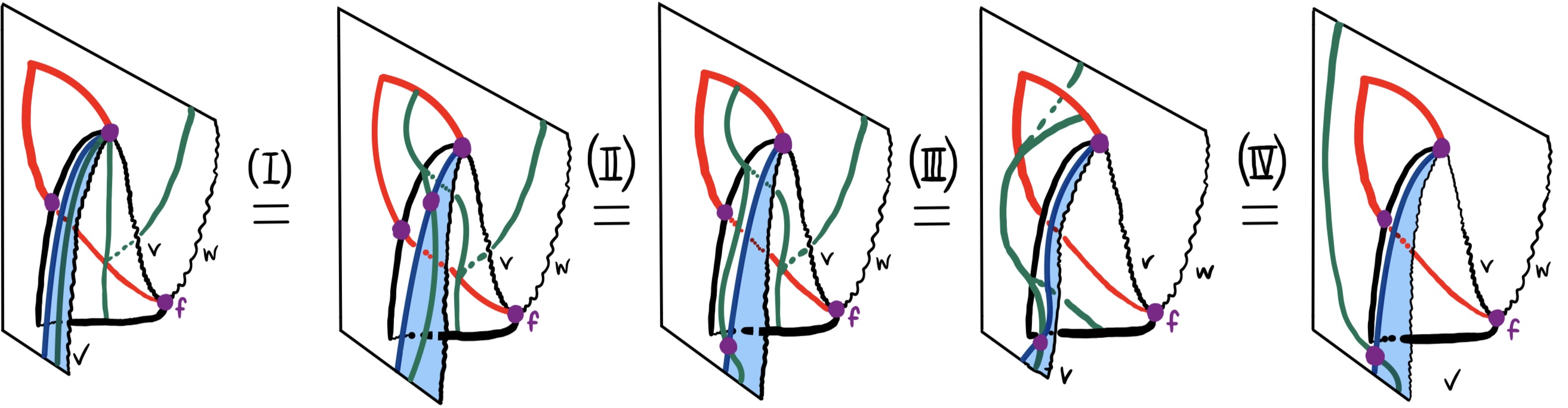}.
\end{figure}
Our strategy is to move the upper part of the green line depicting the functor $F$ towards the back. Explicitly, Equation (I) is satisfied by Remark \ref{duality transfo Remark}.(i). Equation (II) amounts to sliding the inverse of the duality transformation down. This is permissible by the axioms of a strict monoidal $2$-category. Equation (III) holds by the Frobenius relation \eqref{eq:F1 Frob LD} from Definition \ref{def:FrobLin}. To see this, the reader might wish to focus on the green functor lines. Equation (IV) follows from the unitality (resp. counitality) axiom of the lax $\otimes$-monoidal structure (resp. oplax $\parll$-monoidal structure) on $F$. We thus end up with the composite of the left-hand map with the bottom horizontal map, applied to $f.$
\end{proof}

The next theorem generalizes a result of Shimizu \cite[Thm. 5.2]{ShiPivotalCover} from strong monoidal functors between pivotal rigid monoidal categories to pivotal Frobenius LD-functors between pivotal GV-categories:

\begin{theorem}\label{main thm FrobSchur}
    Let $F\colon(\cC,\rho^{\cC})\ra (\cD,\rho^{\cD})$ be a $k$-linear pivotal Frobenius LD-functor between $k$-linear pivotal GV-categories. Let $V\in \cC$. For an integer $n\geq 1$, consider the $k$-linear map
    \begin{equation*}
    \setlength{\arraycolsep}{0pt}
    \renewcommand{\arraystretch}{1.2}
    \tilde{F}_n\colon 
    \left\{ \quad
    \begin{array}{ c c l }
    \operatorname{Hom}_{\cC}(1,V^{\parLL n}) & {} \longrightarrow {} & \operatorname{Hom}_{\cD}(1,F(V)^{\parLL n}) \\
    f & {} \longmapsto {} & \big(1\xrightarrow{\varphi^0}F(1)\xrightarrow{F(f)}F(V^{\parLL n})\xrightarrow{\upsilon^n_{V,\ldots,V}}F(V)^{\parLL n}\big).
    \end{array}
    \right.
    \end{equation*}
Here, $\varphi^0$ is the unit morphism of the Frobenius LD-functor $F$. The symbol $\upsilon^n_{V,\ldots, V}$ denotes the unique morphism $F(V^{\parLL n})\ra F(V)^{\parLL n}$ composed of instances of the opmonoidal structure of $F$.

Then the following diagram commutes:
\begin{equation}\label{ShimizuThm}
\begin{tikzcd}
	{\operatorname{Hom}_{\cC}(1,V^{\parLL n})} && {\operatorname{Hom}_{\cD}(1,F(V)^{\parLL n})} \\
	{\operatorname{Hom}_{\cC}(1,V^{\parLL n})} && {\operatorname{Hom}_{\cD}(1,F(V)^{\parLL n}).}
	\arrow["\tilde{F}_n", from=1-1, to=1-3]
	\arrow["E^{(n)}_V"', from=1-1, to=2-1]
	\arrow["E^{(n)}_{F(V)}", from=1-3, to=2-3]
	\arrow["\tilde{F}_n", from=2-1, to=2-3]
\end{tikzcd}
\end{equation} 
\end{theorem}

\begin{remark}
    Informally, Theorem \ref{main thm FrobSchur} states that $k$-linear pivotal Frobenius LD-functors preserve the cyclic structure of the invariant tensors.
\end{remark}

\begin{proof}
We choose the evaluation $\operatorname{ev}_{FD(V)}$ and coevaluation $\operatorname{coev}_{FD(V)}$ as in Lemma \ref{evFDV and coevFDV}.

For $n=1$, the first surface diagram in the following figure is the composition of the upper horizontal and the right-hand vertical map in the diagram (\ref{ShimizuThm}), applied to a morphism $f\in \operatorname{Hom}_{\cC}(1,V)$:

\begin{figure}[H]
    \centering
    \includegraphics[width=0.99\textwidth]{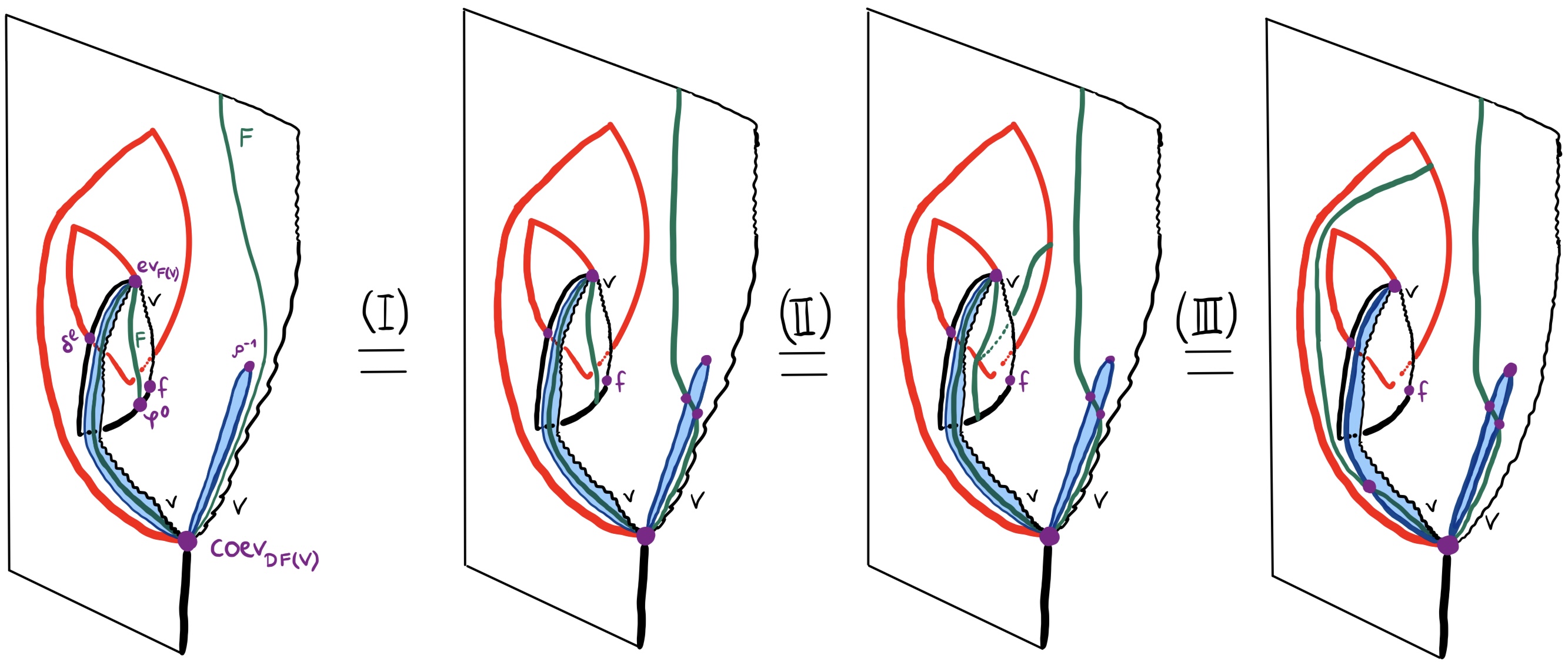}
\end{figure}
\begin{figure}[H]
    \centering
    \includegraphics[width=0.93\textwidth]{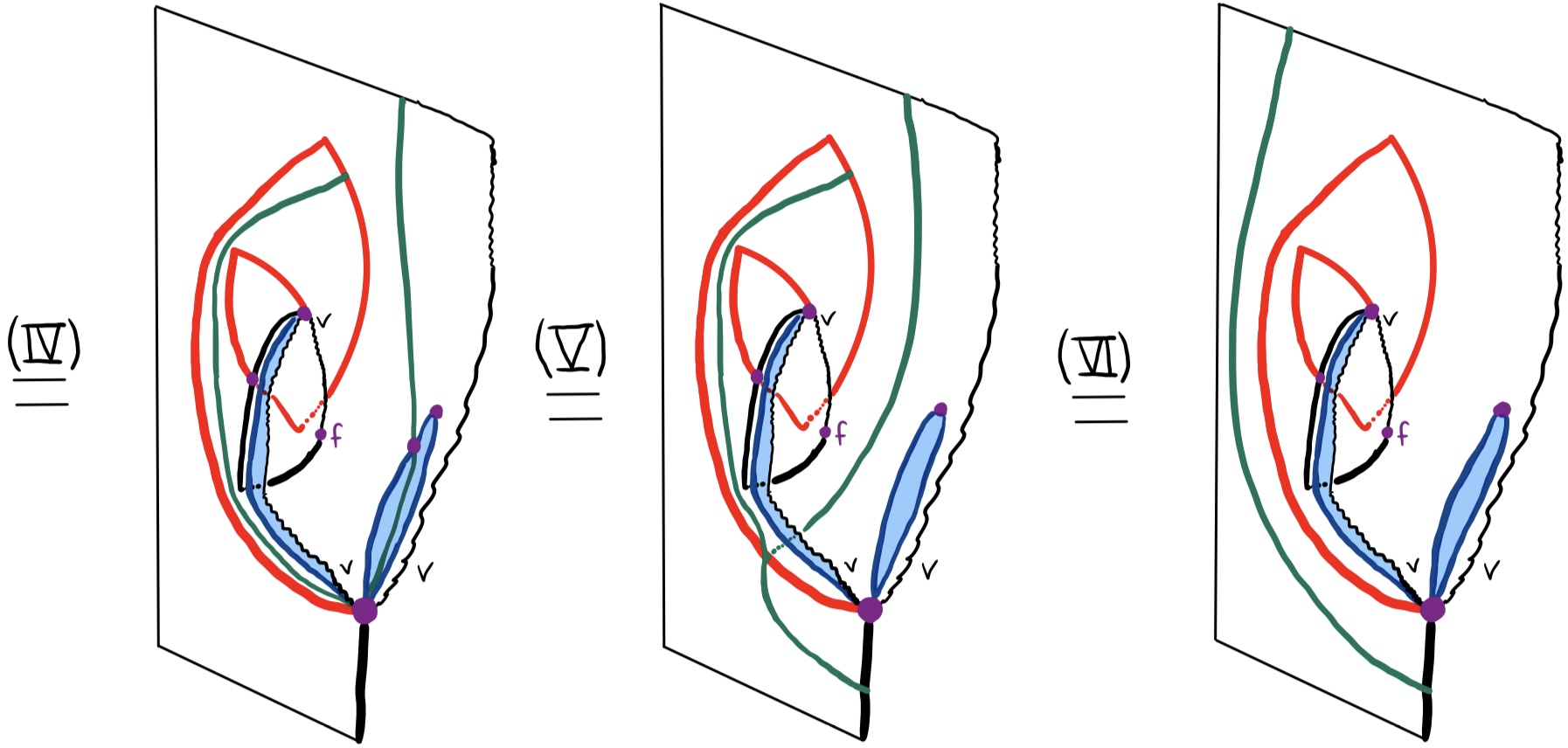}.
\end{figure}
Our proof strategy is to move the green functor line of $F$ behind the red functor lines of $\parLL$. We begin by first using the pivotality of $F$ in Equation (I). Equation (II) follows from the counitality axiom of the oplax $\parll$-monoidal structure on $F$. Equation (III) is an immediate consequence of Lemma \ref{lemma:technical lemma}, while Equation (IV) holds by our choice of coevaluation $\operatorname{coev}_{FD(V)}$ as in Lemma \ref{evFDV and coevFDV}. Equation (V) is an application of Equation (\ref{charact property duality transfo}), which characterizes the duality transformation for $F$. Equation (VI) follows from the counitality axiom of the oplax $\parll$-monoidal structure on $F$.

\bigskip

Similarly, for an integer $n\geq 2$ and a morphism $f\in \operatorname{Hom}_{\cC}(1,V^{\parLL n})$, we have:
\begin{figure}[H]
    \centering
    \includegraphics[width=0.99\textwidth]{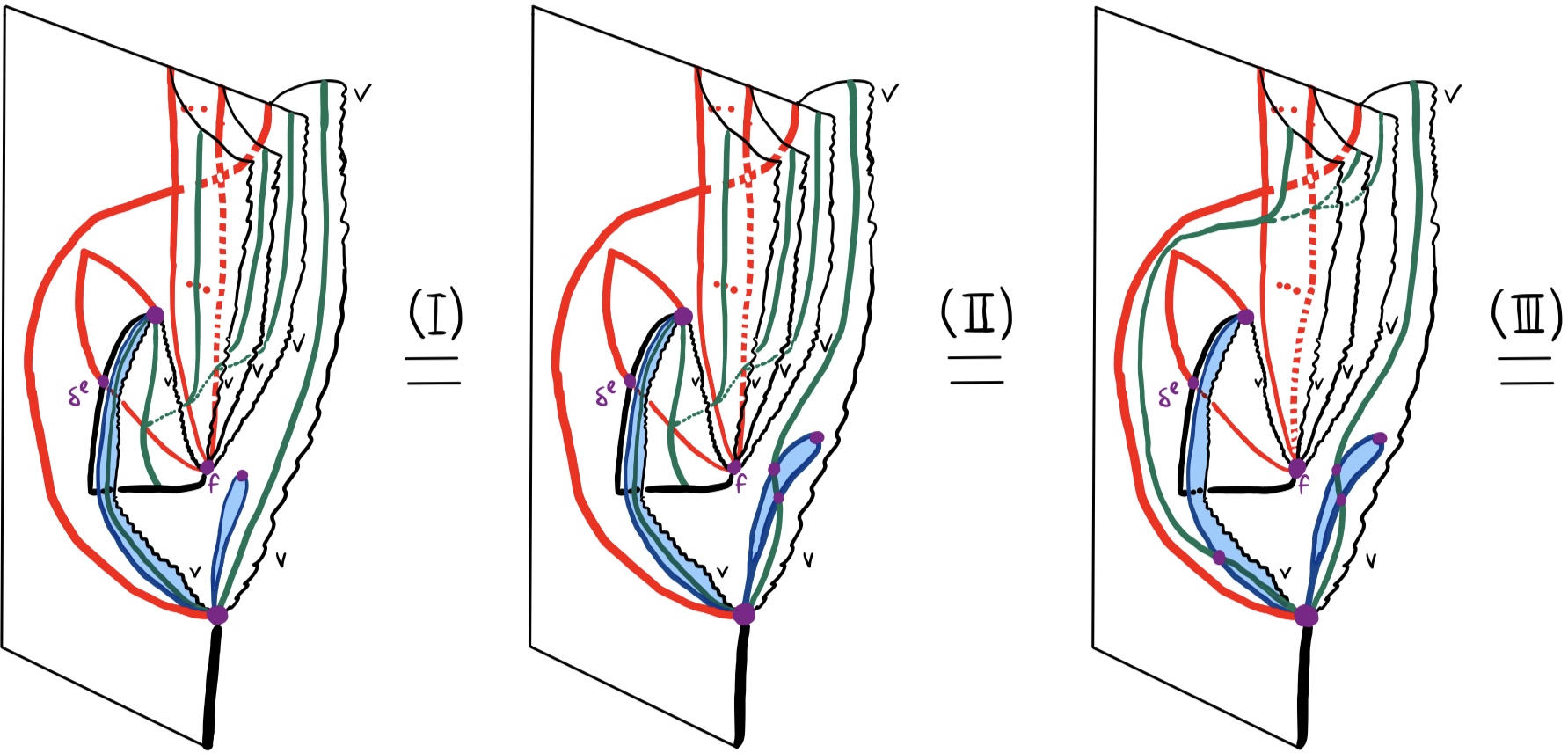}
\end{figure}
    \medskip
\begin{figure}[H]
    \centering
    \includegraphics[width=0.9\textwidth]{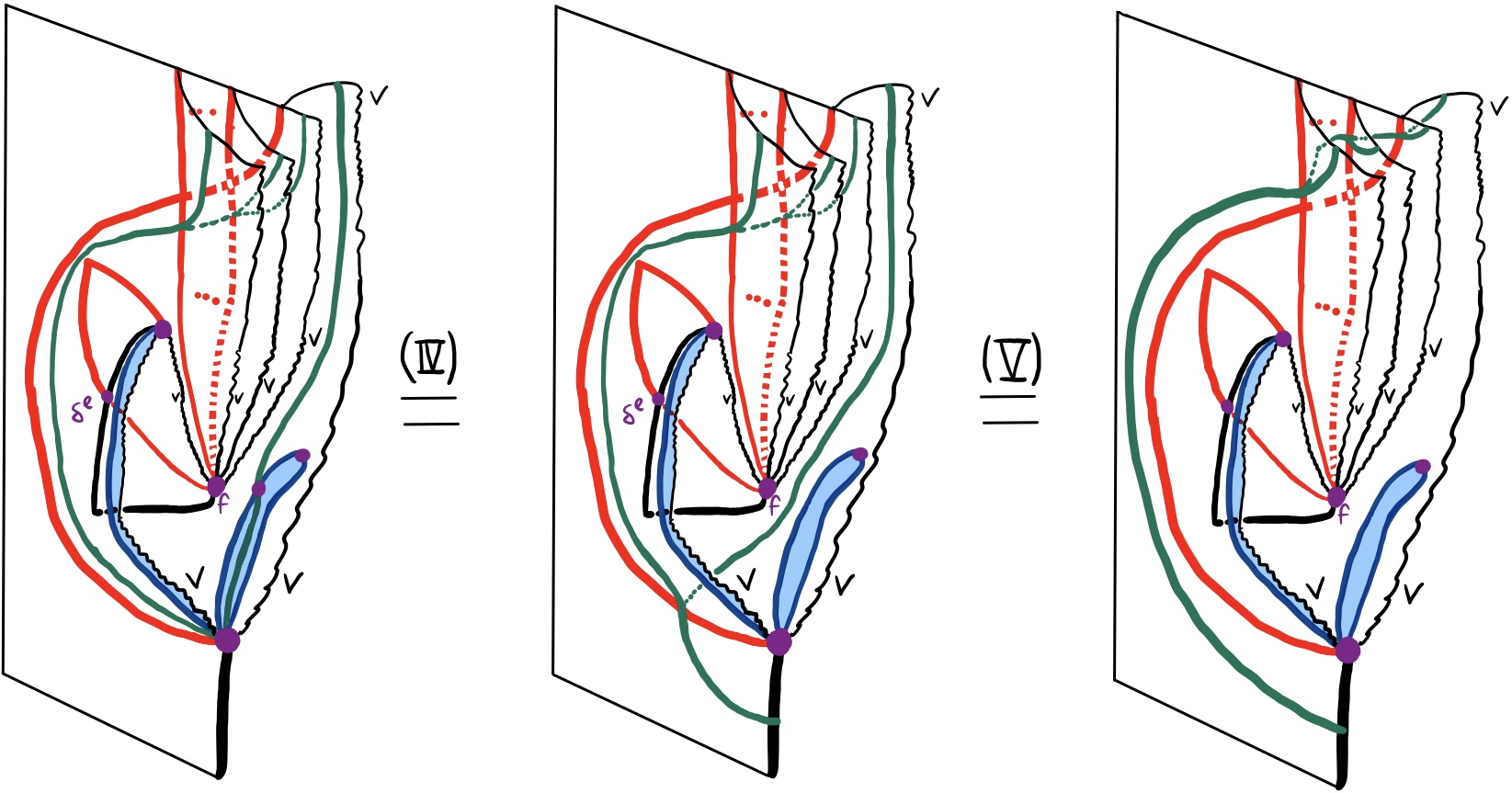}.
\end{figure}
    
Equation (I) holds since $F$ is pivotal, while Equation (II) is an immediate consequence of Lemma \ref{lemma:technical lemma}. Equation (III) holds by our choice of coevaluation $\operatorname{coev}_{FD(V)}$. Equation (IV) is an application of Equation (\ref{charact property duality transfo}), which characterizes the duality transformation. Equation (V) follows from the coassociativity axiom of the oplax $\parll$-monoidal structure on $F$.
\end{proof}

Recall from Definition \ref{def:strong Frobenius LD-functor} that a Frobenius LD-functor is called strong if all four of its coherence morphisms are isomorphisms. For the rest of this section, the following weakening of the notion of a strong Frobenius LD-functor is relevant:

\begin{definition}
    A Frobenius LD-functor $F\colon \cC\ra \cD$ between LD-categories $\cC$ and $\cD$ is \emph{$(\parLL,1)$-semistrong} if the comultiplication ${\upsilon^2\colon F\circ \parLL \ra {\parLL} \circ {(F\times F)}}$ and the unit $\varphi^0\colon 1 \ra F(1)$ are both invertible.
\end{definition}

The next claim is immediate:

\begin{lemma}\label{prop half-strong}
    Let $F\colon\cC\ra\cD$ be a $(\parLL,1)$-semistrong Frobenius LD-functor between LD-categories. Assume that the underlying functor of $F$ is fully faithful. For objects $X,Y\in \cC$, denote by $\hat{F}_{X,Y}$ the inverse of the map $F\colon {\operatorname{Hom}_{\cC}(X,Y)\ra\operatorname{Hom}_{\cD}(F(X),F(Y))}$. For any integer $n\geq1$ and any object $V\in \cC$, the $k$-linear map
    \begin{equation*}
    \setlength{\arraycolsep}{0pt}
    \renewcommand{\arraystretch}{1.2}
    \tilde{F}_n\colon 
    \left\{ \quad
    \begin{array}{ c c l }
    \operatorname{Hom}_{\cC}(1,V^{\parLL n}) & {} \longrightarrow {} & \operatorname{Hom}_{\cD}(1,F(V)^{\parLL n}) \\
    f & {} \longmapsto {} & \big(1\xrightarrow{\varphi^0}F(1)\xrightarrow{F(f)}F(V^{\parLL n})\xrightarrow{\upsilon^n_{V,\ldots,V}}F(V)^{\parLL n}\big)
    \end{array}
    \right.
    \end{equation*}
possesses a two-sided inverse
\begin{equation*}
    \setlength{\arraycolsep}{0pt}
    \renewcommand{\arraystretch}{1.2}
    (\tilde{F}_n)^{-1}\colon 
    \left\{ \quad
    \begin{array}{ c c l }
    \operatorname{Hom}_{\cD}(1,F(V)^{\parLL n}) & {} \longrightarrow {} & \operatorname{Hom}_{\cC}(1,V^{\parLL n}) \\
    f & {} \longmapsto {} & \hat{F}\big(F(1)\xrightarrow{(\varphi^0)^{-1}}1\xrightarrow{f}F(V)^{\parLL n}\xrightarrow{(\upsilon^n_{V,\ldots,V})^{-1}}F(V^{\parLL n})\big).
    \end{array}
    \right.
    \end{equation*}
\end{lemma}

Since the trace is cyclically invariant, Theorem \ref{main thm FrobSchur} and Lemma \ref{prop half-strong} immediately imply the following:

\begin{theorem}\label{thm: FrobSchur categorical invariants}
    Let $F\colon\cC\ra\cD$ be a $k$-linear pivotal $(\parLL,1)$-semistrong Frobenius LD-functor between $k$-linear pivotal GV-categories. Assume that the underlying functor of $F$ is fully faithful. For all integers $n,r\geq 1$ and for every object $V\in \cC$, we then have 
    \begin{equation*}
        \nu_{n,r}(V)\;=\;\nu_{n,r}(F(V)).
    \end{equation*}
    In particular, the $(n,r)$-th Frobenius-Schur indicator is invariant under $k$-linear pivotal Frobenius LD-equivalence in the sense of Definition \ref{def:strong Frobenius LD-functor}.
\end{theorem}

%Kürzungspotenzial
Lastly, we briefly discuss the $n$-th power of the endomorphism $E_V^{(n)}$.

\begin{prop}\label{E1V is id}
    Let $(\cC,\rho)$ be a pivotal GV-category. We have $E^{(1)}_V=\operatorname{id}.$ Thus, $\nu_{1,1}(V)=1.$
\end{prop}

\begin{proof}
By Example \ref{unitors are side inverse}, we can assume that we have chosen LD-duals, their evaluations and coevaluations such that $D(1)=K$, $D(K)=1$, $\operatorname{coev}_{D(1)}=(l^{\parLL}_1)^{-1}$ and $\operatorname{ev}_{1}=r^{\otimes}_K$.

The first surface diagram depicts the composite $E^{(1)}_V(f)$ for a morphism $f\in \operatorname{Hom}_{\cC}(1,V)$:
\begin{figure}[H]
    \centering
    \includegraphics[width=0.77\textwidth]{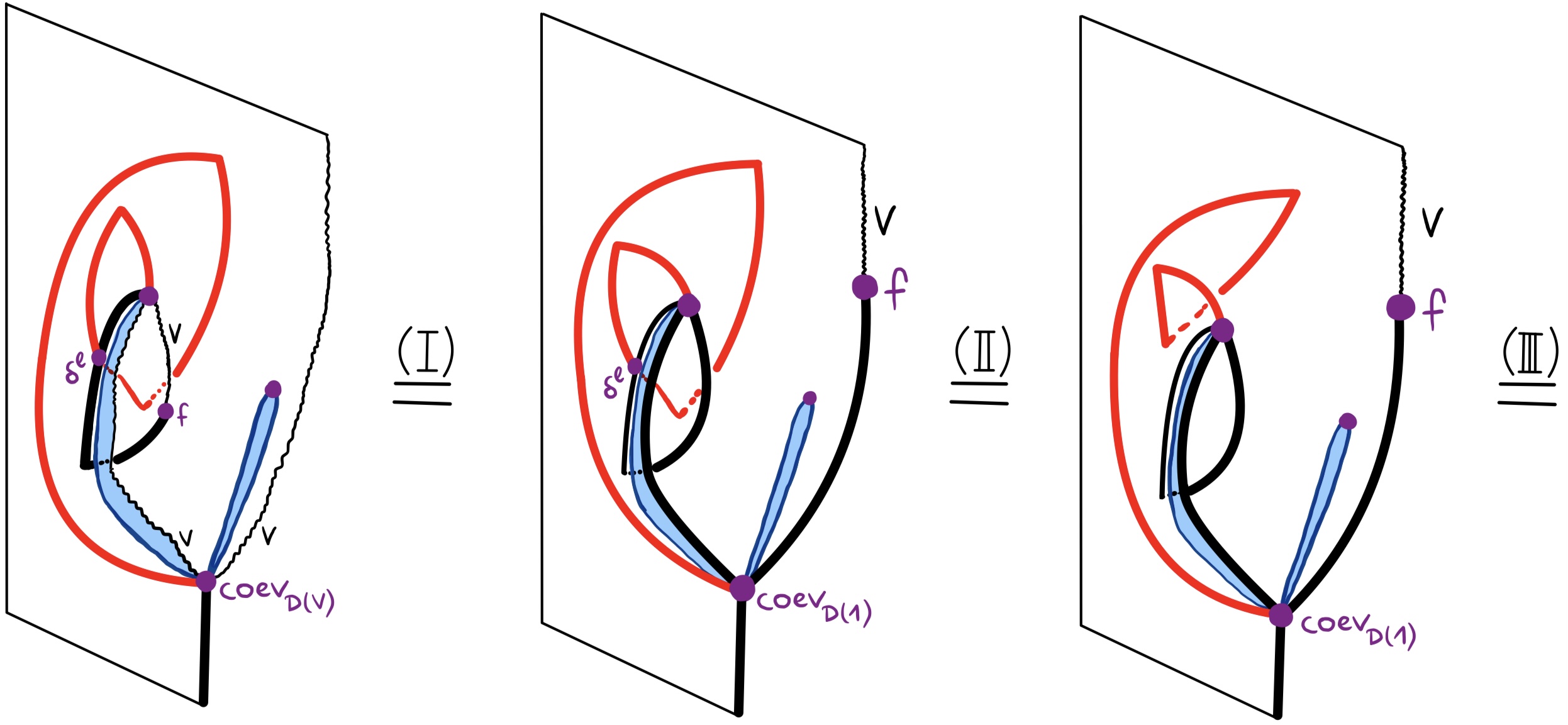}
\end{figure}
\begin{figure}[H]
    \centering
    \includegraphics[width=0.77\textwidth]{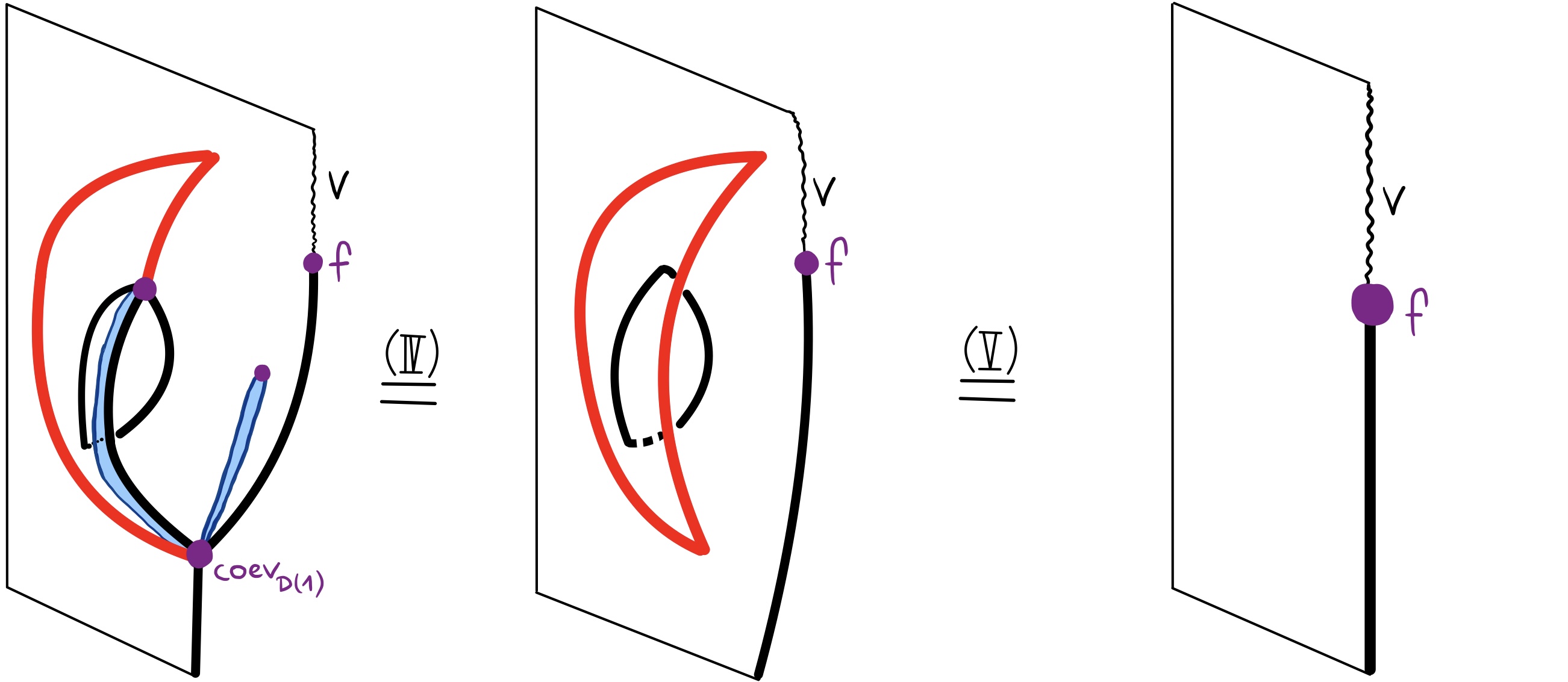}.
\end{figure}
    Equation (I) follows from Proposition \ref{evalDuality functor}, while Equation (II) holds by the triangle coherence axiom (\ref{eq:A4}). Equation (III) follows from Mac Lane's coherence theorem for the $\parLL$-monoidal structure, while Equation (IV) holds by our choice of $\operatorname{coev}_{D(1)}$ and $\operatorname{ev}_{1}$. Finally, Equation (V) holds since the unitors $l^{\parLL}$ and $(l^{\parLL})^{-1}$, as well as the unitors $r^{\otimes}$ and $(r^{\otimes})^{-1}$, are inverses.
\end{proof}

\begin{remark}
Proposition \ref{E1V is id} immediately implies Theorem \ref{main thm FrobSchur} for the case $n=1.$     
\end{remark}

\begin{remark}
    The arguments in the proof of Proposition \ref{E1V is id} can be generalized to show that $(E^{(n)})^n=\operatorname{id},$ for all integers $n\geq 1.$ Thus, for an algebraically closed field $k$ of characteristic zero, the $(n,r)$-th Frobenius-Schur indicator $\nu_{n,r}(V)$ is a cyclotomic integer in $\mathbb{Q}_n\subset k$. The proofs are similar to those in \cite[Thm. 5.1]{SchauNgHigherFrob}.
\end{remark}

%%%%%%%%%%%%%%%%%%%%%%%%%%%%%

\appendix
\section{Coherence axioms for LD-categories}\label{coherenceLD}
This section presents the ten coherence axioms for the distributors in an LD-category, as described in Definition \ref{def:LDcat}, using surface diagrams. The coherence axioms are taken from the (corrected version) of \cite{WDC}.

\medskip

Before presenting the coherence axioms, let us briefly discuss them informally: Graphically, they amount to moving a gluing boundary of a sheet along a surface. That is, the coherence axioms introduce moves that change the configuration of sheets within the canvas. In this section, we color a sheet that is being moved light grey. We do this to improve readability only, and not to indicate a different category.

\bigskip

\begin{definition}[Continuation of Definition \ref{def:LDcat}]\vphantom\newline

\begin{enumerate}[label=(\roman*)]
    \item There are four compatibility conditions involving unitors and distributors:
        \begin{equation}\label{eq:A1}
        \begin{aligned}
        \centering
        \includegraphics[width=0.5\textwidth]{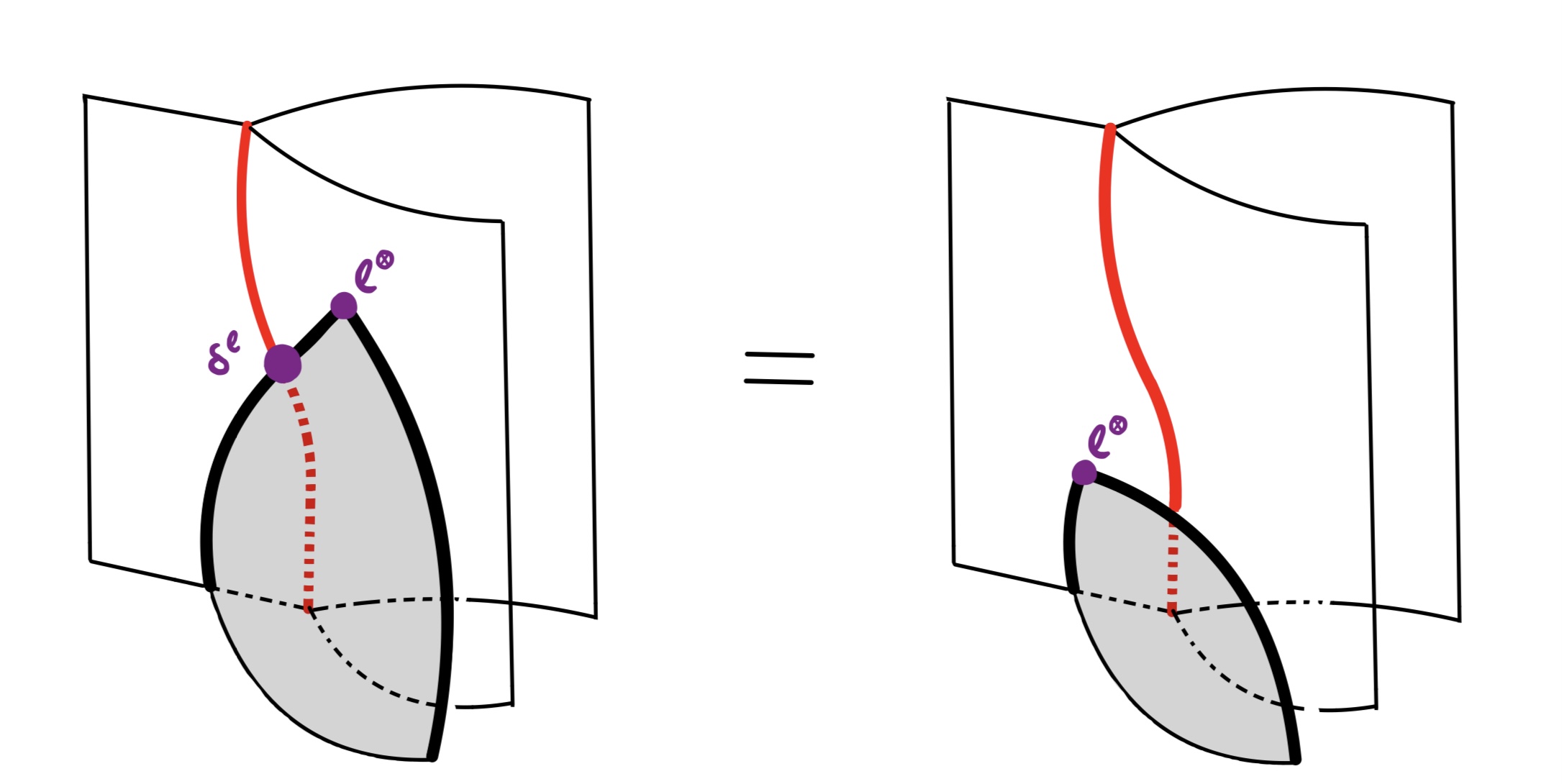}
        \end{aligned}
    \end{equation}

    \begin{equation}\label{eq:A2}
        \begin{aligned}
        \centering
        \includegraphics[width=0.5\textwidth]{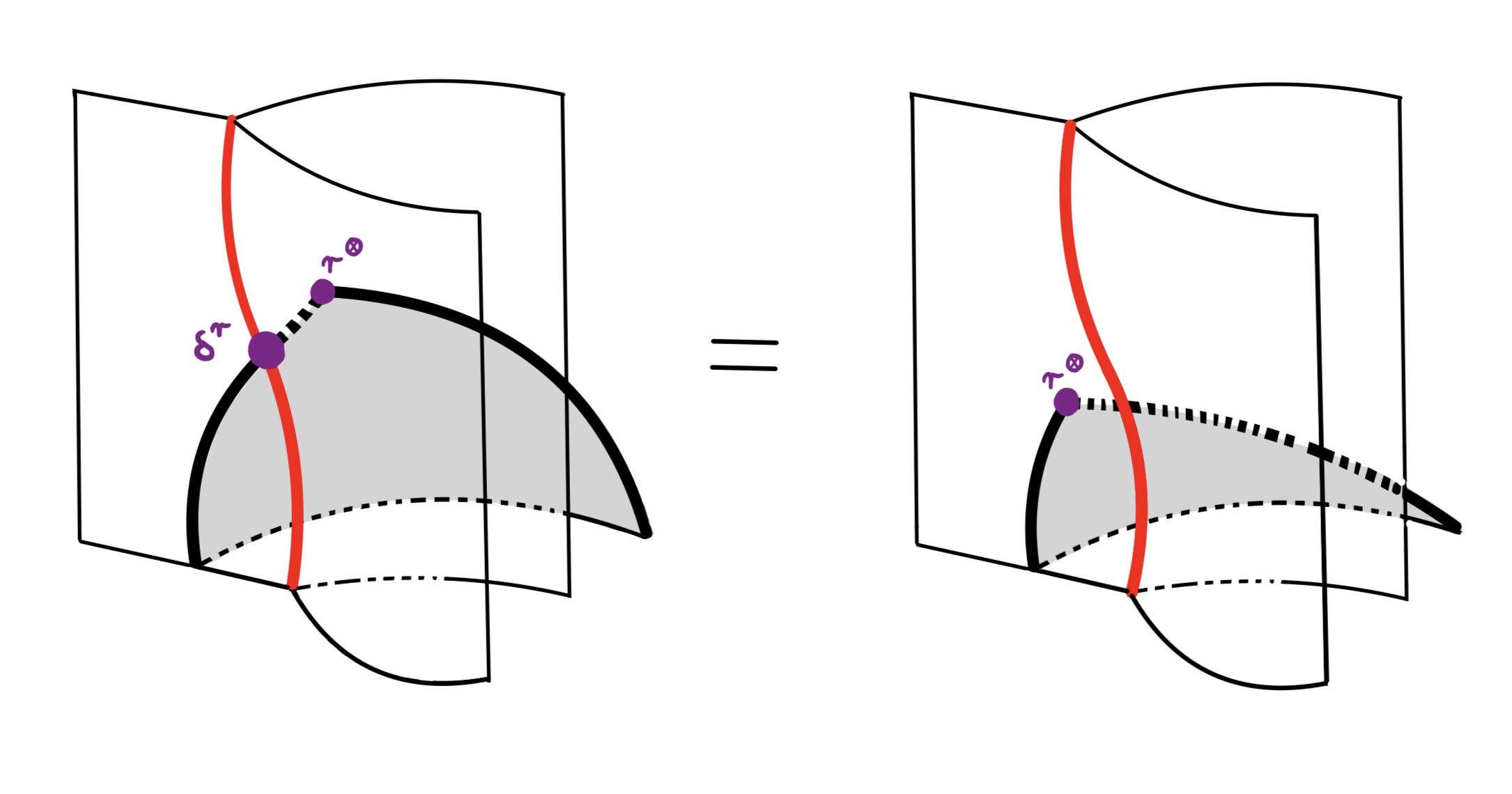}
        \end{aligned}
    \end{equation}

    \begin{equation}\label{eq:A3}
        \begin{aligned}
        \centering
        \includegraphics[width=0.5\textwidth]{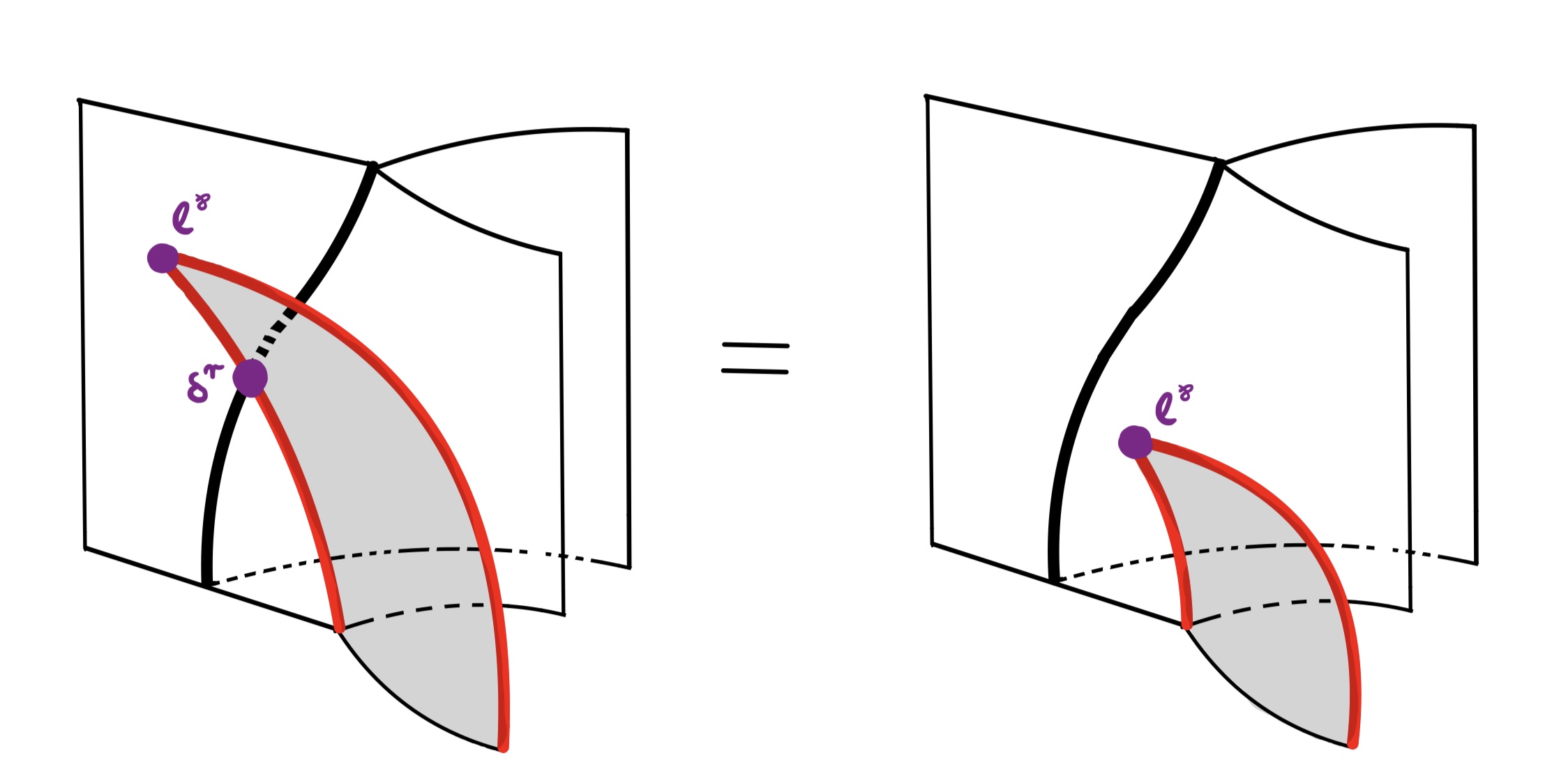}
        \end{aligned}
    \end{equation}

    \begin{equation}\label{eq:A4}
        \begin{aligned}
        \centering
        \includegraphics[width=0.5\textwidth]{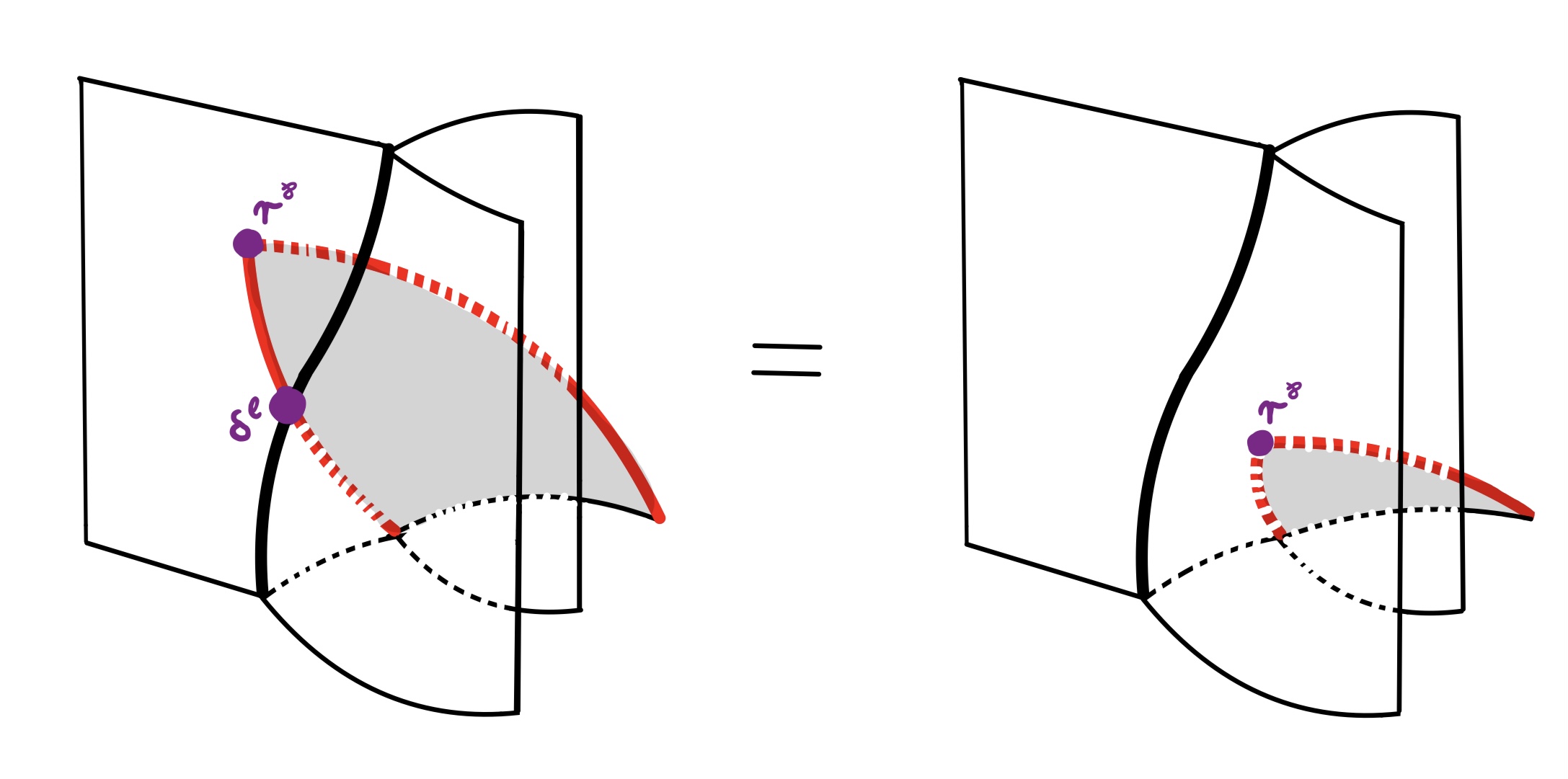}.
        \end{aligned}
    \end{equation}
Informally, the coherence conditions (\ref{eq:A1}) to (\ref{eq:A4}) amount to pulling a (light grey) triangular sheet, bounded by lines of one color, above a functor line of a different color. The surface diagrams for each of these coherence conditions contain a total of three purple vertices. In other words, the coherence conditions (\ref{eq:A1}) to (\ref{eq:A4}) correspond to triangular commutative diagrams.

\item There are four compatibility conditions involving associators and distributors:

\bigskip

\begin{equation}\label{eq:A5}
    \begin{aligned}
   \centering
    \includegraphics[width=0.545\textwidth]{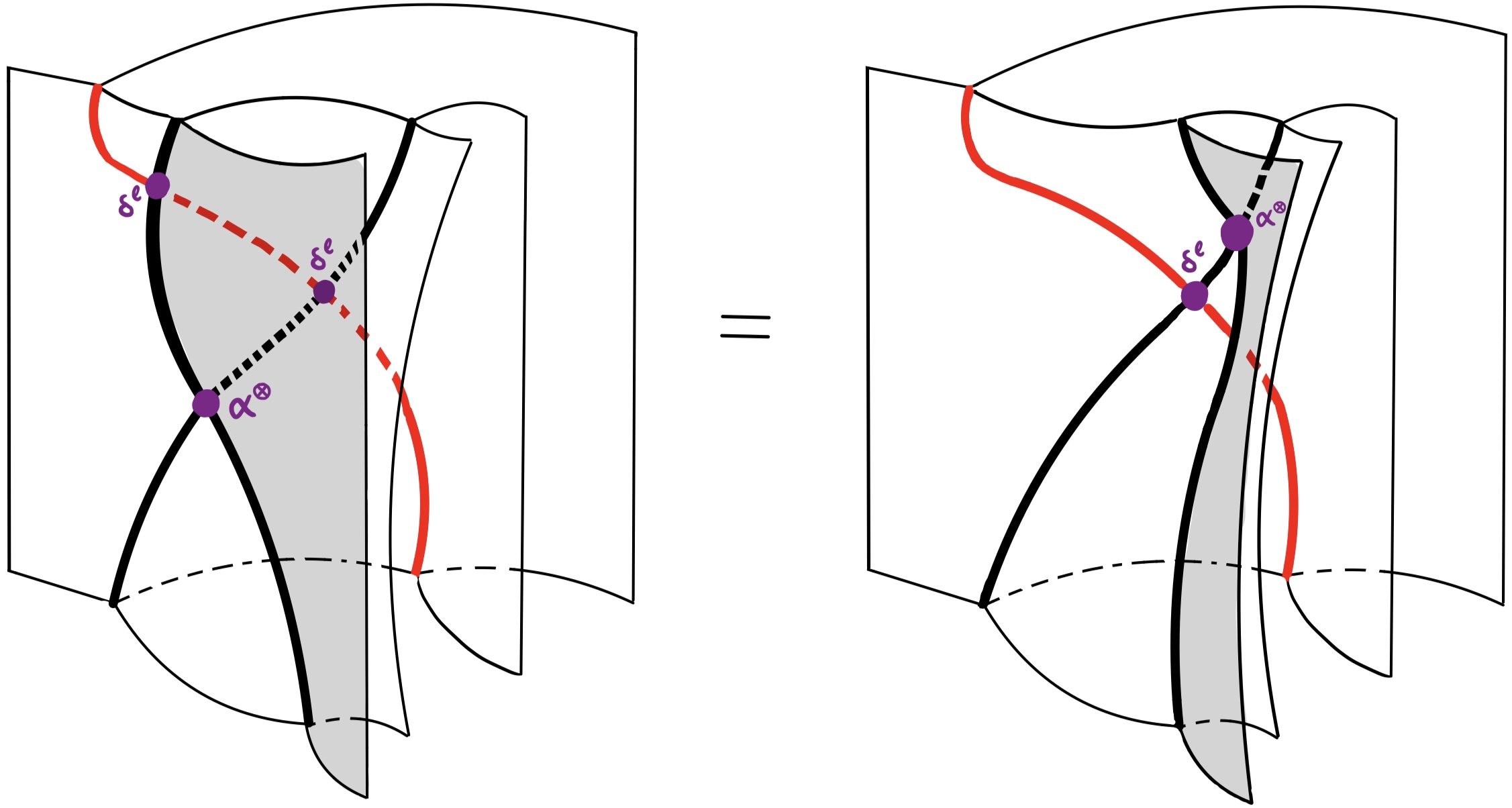}
    \end{aligned}
\end{equation}

\begin{equation}\label{eq:A6}
    \begin{aligned}
   \centering
    \includegraphics[width=0.545\textwidth]{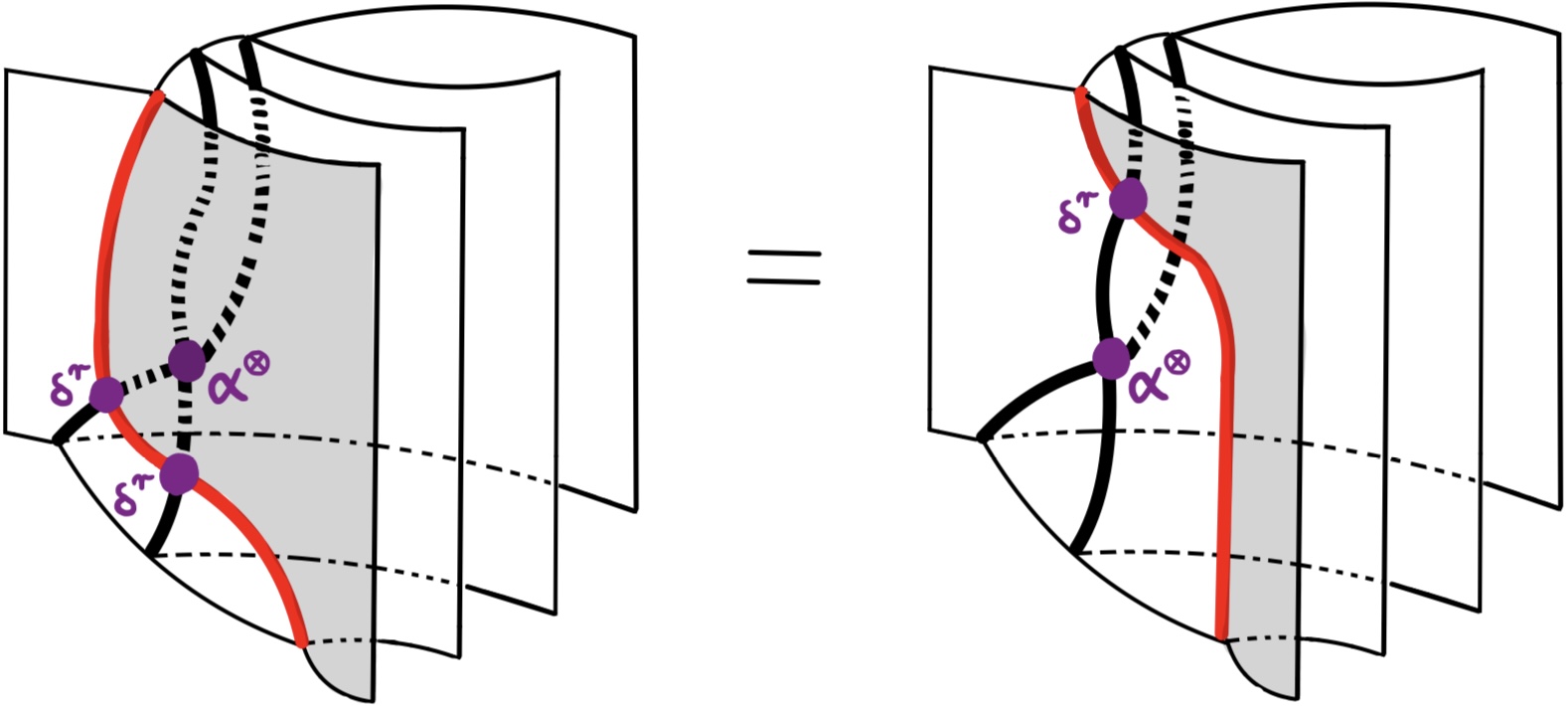}
    \end{aligned}
\end{equation}

\begin{equation}\label{eq:A7}
    \begin{aligned}
    \centering
    \includegraphics[width=0.545\textwidth]{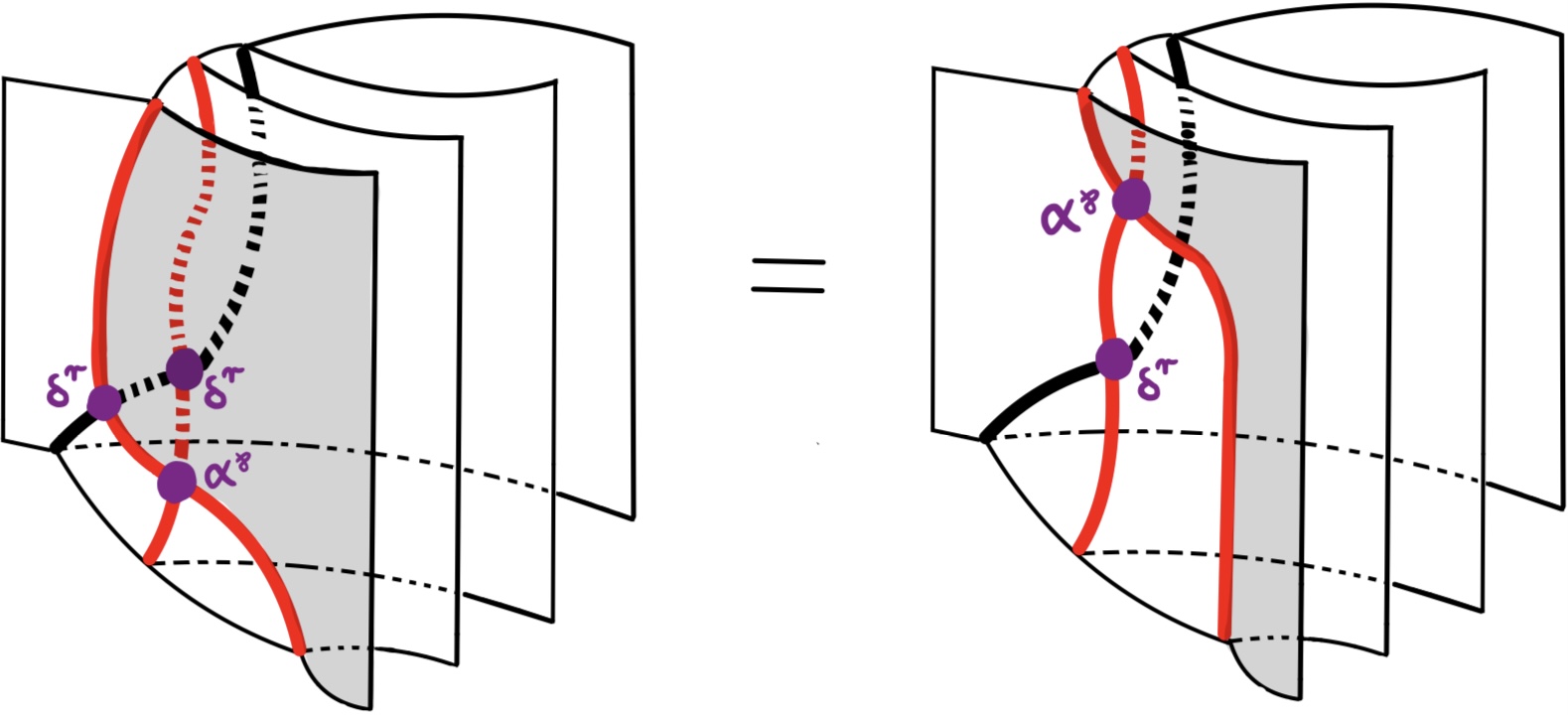}
    \end{aligned}
\end{equation}

\begin{equation}\label{eq:A8}
    \begin{aligned}
    \centering
    \includegraphics[width=0.545\textwidth]{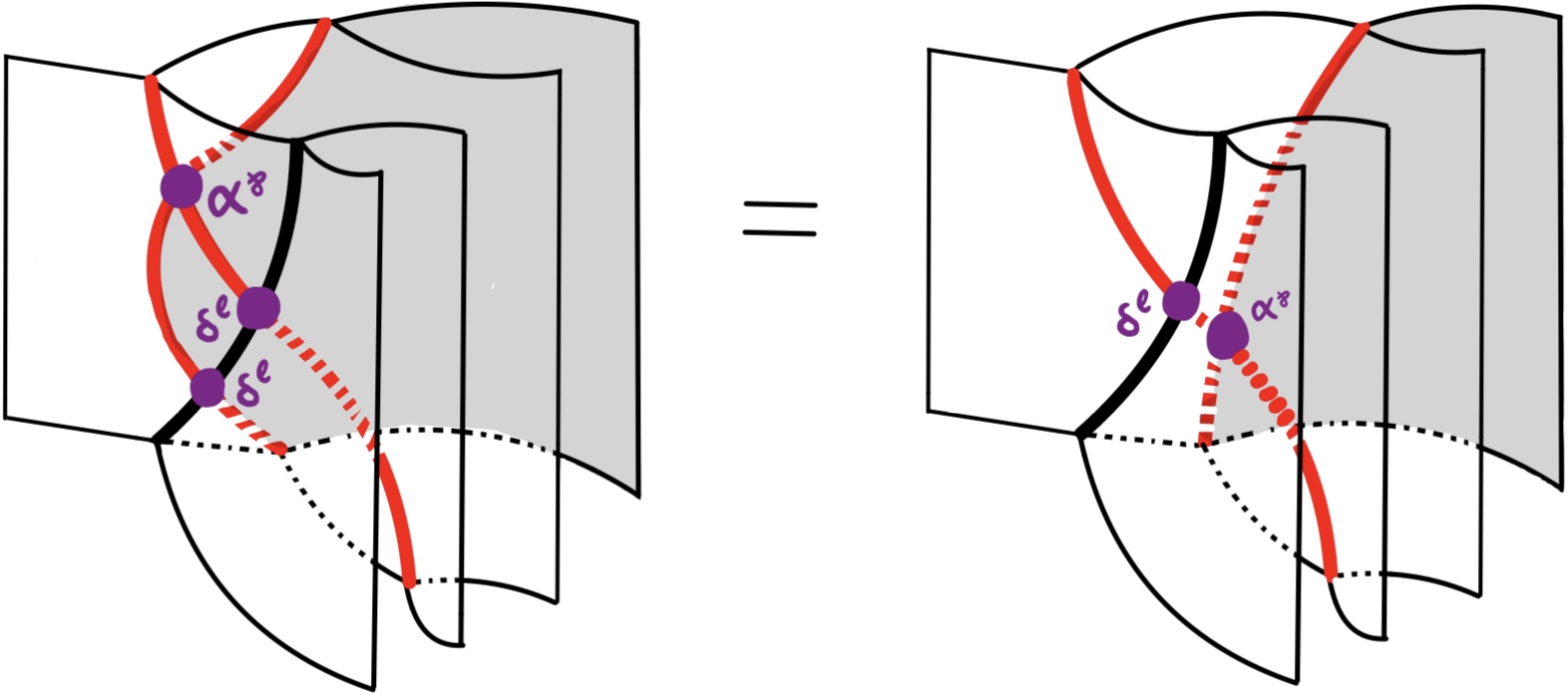}.
    \end{aligned}
\end{equation}

\item There are two compatibility conditions involving only the distributors:

\begin{equation}\label{eq:A9}
    \begin{aligned}
    \centering
    \includegraphics[width=0.545\textwidth]{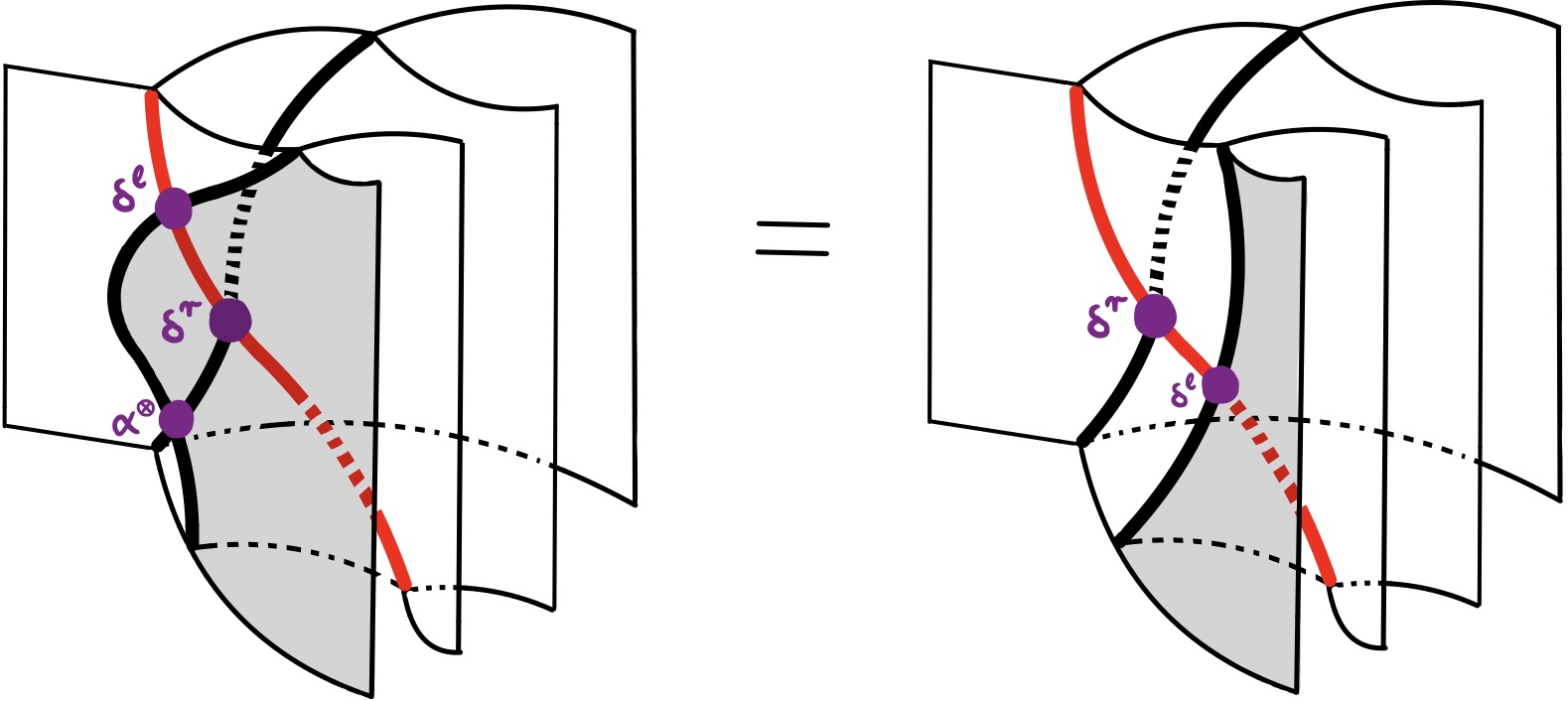}
    \end{aligned}
\end{equation}

\begin{equation}\label{eq:A10}
    \begin{aligned}
    \centering
    \includegraphics[width=0.549\textwidth]{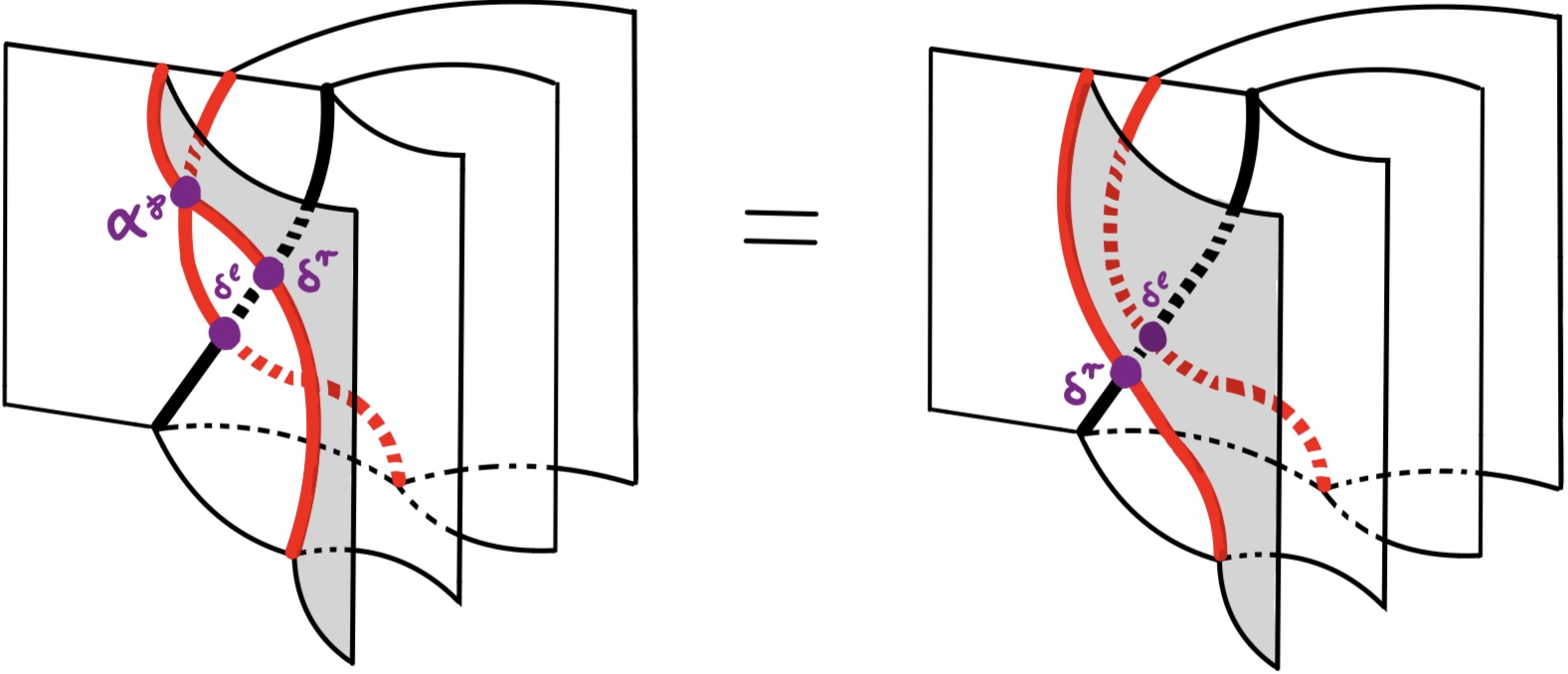}.
    \end{aligned}
\end{equation}
Informally, the coherence conditions (\ref{eq:A5}) to (\ref{eq:A10}) amount to pulling a (light grey) sheet along a surface above a singular vertex. The surface diagrams for each of these coherence conditions contain a total of five purple vertices. In other words, the coherence conditions (\ref{eq:A5}) to (\ref{eq:A10}) correspond to pentagonal commutative diagrams. 
\end{enumerate}
\end{definition}

\section{Auxiliary proofs}\label{sec:appendixproof}

\begin{proof}[Proof sketch of Proposition \ref{quasi-inverse gets Frobenius LD-structure}.]
Without loss of generality, we assume that the underlying equivalence $(F,G)$ is an adjoint equivalence with unit $\eta\colon \operatorname{id}_{\cC}\xrightarrow{\simeq}FG$ and counit $\epsilon\colon GF\xrightarrow{\simeq}\operatorname{id}_{\cC}$. 

To show uniqueness, assume that $G$ carries a strong Frobenius LD-structure 
\begin{equation*}
(\varphi^{2,G},\varphi^{0,G},\nu^{2,G},\nu^{0,G})
\end{equation*}
such that the unit $\eta\colon \operatorname{id}_{\cC}\xrightarrow{\simeq}GF$ and counit $\epsilon\colon FG\xrightarrow{\simeq}\operatorname{id}_{\cD}$ are isomorphisms of Frobenius LD-functors. Denote the strong Frobenius LD-structure on $F$ by $(\varphi^{2,F},\varphi^{0,F},\nu^{2,F},\nu^{0,F})$.

For $X,Y\in\cC$, we then have:
\begin{align}\label{varphi2G}
    \begin{split}
        \varphi^{2,G}_{X,Y}& \;\eqabove{(I)} \; \varphi^{2,G}_{X,Y}\circ (G(\epsilon_X)\otimes G(\epsilon_Y))\circ(\eta_{G(X)}\otimes\eta_{G(Y)})\\ &\;\eqabove{(II)} \; G(\epsilon_X\otimes \epsilon_Y)\circ\varphi^{2,G}_{FG(X),FG(Y)}\circ(\eta_{G(X)}\otimes\eta_{G(Y)})\\ &\;\eqabove{(III)} \; G(\epsilon_X\otimes \epsilon_Y)\circ G(\varphi^{2,F}_{G(X),G(Y)})^{-1}\circ G(\varphi^{2,F}_{G(X),G(Y)})\circ\varphi^{2,G}_{FG(X),FG(Y)}\circ(\eta_{G(X)}\otimes\eta_{G(Y)})\\&\;\eqabove{(IV)} \; G(\epsilon_X\otimes\epsilon_Y)\circ G(\varphi^{2,F}_{G(X),G(Y)})^{-1}\circ\eta_{G(X)\otimes G(Y)}.
    \end{split}
\end{align}
Equation (I) follows from one snake equation for the unit $\eta\colon \operatorname{id}_{\cC}\xrightarrow{\simeq}GF$ and the counit ${\epsilon\colon FG\xrightarrow{\simeq}\operatorname{id}_{\cD}}$, while Equation (II) holds by the naturality of the multiplication morphism $\varphi^{2,G}$. Equation (III) inserts an identity morphism by using that the $\otimes$-monoidal structure on $F$ is strong, while Equation (IV) follows from the unit $\eta\colon \operatorname{id}_{\cC}\xrightarrow{\simeq}GF$ being an $\otimes$-monoidal natural isomorphism. 

\smallskip

For $X,Y\in \cC$, we similarly find:
\begin{align}
\varphi^{0,G}\,& =\,G\big(\varphi^{0,F}\big)^{-1}\circ\eta_1,\label{varphi0G}\\
\nu^{2,G}_{X,Y}\,&=\,\eta^{-1}_{G(X)\parLL G(Y)}\circ G(\nu^{2,F}_{G(X),G(Y)})^{-1}\circ G(\epsilon^{-1}_X\parLL\epsilon_Y^{-1}),\label{nu2G}\\
\nu^{0,G}\,&=\,\eta^{-1}_K\circ G(\nu^{0,F})^{-1}. \label{nu0G}
\end{align}
This shows that, if it exists, a Frobenius LD-structure on $G$ such that the unit $\eta\colon \operatorname{id}_{\cC}\xrightarrow{\simeq}GF$ and counit $\epsilon\colon FG\xrightarrow{\simeq}\operatorname{id}_{\cD}$ are isomorphisms of Frobenius LD-functors is unique. 

\smallskip

Next, one proves that the coherence morphisms $\varphi^{2,G}$ and $\varphi^{0,G}$, as defined in Equations (\ref{varphi2G}) and (\ref{varphi0G}), yield a strong $\otimes$-monoidal structure on $G$. This is well known; see, e.g. \cite[\S 1.4.4]{CatTann}. Analogously, one shows that the coherence morphisms $\nu^{2,G}$ and $\nu^{0,G}$, as defined in Equations (\ref{nu2G}) and (\ref{nu0G}), yield a strong $\parLL$-monoidal structure on $G$. It is an easy verification that the unit and counit are isomorphisms of Frobenius LD-functors if the Frobenius LD-structure on $G$ is defined as in Equations (\ref{varphi2G}), (\ref{varphi0G}), (\ref{nu2G}) and (\ref{nu0G}).

\smallskip

Finally, we prove the Frobenius relation \eqref{eq:F1 Frob LD} for $G$ surface-diagrammatically. The functors $F$ and $G$ are colored in blue and orange, respectively:
    \begin{figure}[H]
        \centering
        \includegraphics[width=0.99\textwidth]{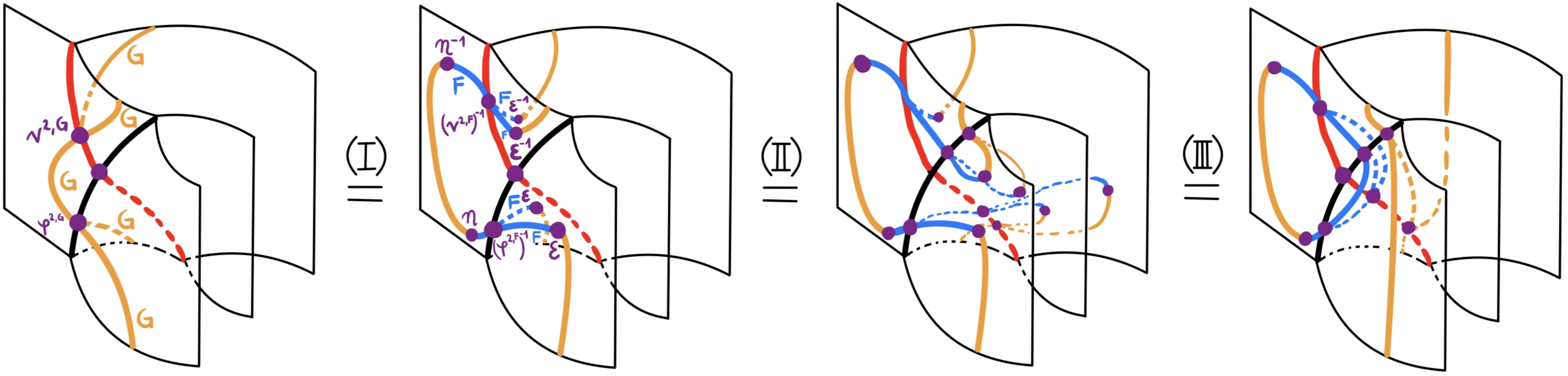}
    \end{figure}
    \begin{figure}[H]
        \centering
        \includegraphics[width=0.99\textwidth]{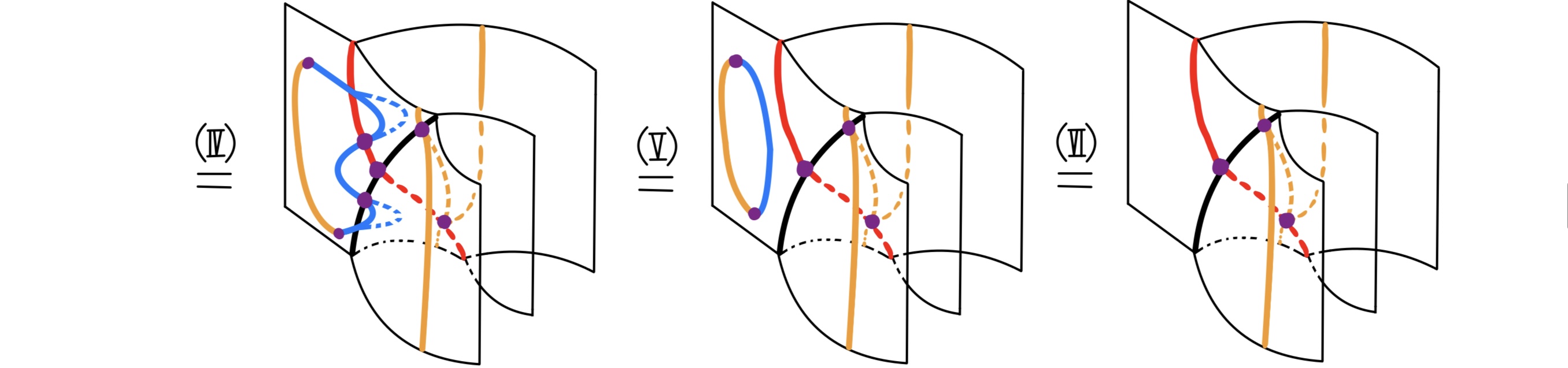}.
    \end{figure}
Equation (I) follows from the definition (\ref{varphi2G}) of the multiplication morphism $\varphi^{2,G}$ and the definition (\ref{nu2G}) of the comultiplication morphism $\nu^{2,G}$ of $G$, while Equation (II) holds since the counit $\epsilon\colon FG\xrightarrow{\simeq}\operatorname{id}_{\cD}$ is an isomorphism of Frobenius LD-functors. Equations (III) and (VI) use the invertibility of the counit $\epsilon\colon FG\xrightarrow{\simeq}\operatorname{id}_{\cD}$ and unit $\eta\colon \operatorname{id}_{\cC}\ra GF$, respectively. Since the Frobenius LD-structure on $F$ is strong, Equation (V) holds. Equation (IV) amounts to an application of Frobenius relation \eqref{eq:F1 Frob LD} for $F$.

The proof of Frobenius relation \eqref{eq:F2 Frob LD} is obtained from the above surface-diagrammatic proof by applying the $2$-functor $(-)^{\operatorname{rev}}$. 
\end{proof}

%%%%%%%%%%%%%%%%%%%%%%%%%%%%%%%%%%%%%%%%%%%%%%%%%%%%%%%%%%%%%%%%%%%%%%%%

\vfill

\nid
\nid\textbf{Acknowledgments:}\\[.3em]
We thank Sam Bauer, Aaron David Fairbanks, Jürgen Fuchs, JS Lemay, Christian Reiher, David Reutter, Gregor Schaumann, Mike Shulman, Lukas Woike, and Simon Wood for correspondence and valuable discussions.

The authors acknowledge support from the DFG through the CRC 1624 \emph{Higher Structures, Moduli Spaces and Integrability}, project number 506632645. C.S. is partially funded by the Excellence Cluster EXC 2121 \emph{Quantum Universe}, project number 390833306.

\bibliographystyle{alphaurl}
\bibliography{references}
\end{document}